\documentclass[12pt]{article}

\usepackage{amsmath, amssymb, eucal, amscd, color}
\begin{document}

\newtheorem{theorem}{Theorem}
\renewcommand{\thetheorem}{}
\newtheorem{proposition}{Proposition}
\renewcommand{\theproposition}{}
\newtheorem{lemma}{Lemma}
\renewcommand{\thelemma}{}
\newtheorem{corollary}{Corollary}
\newtheorem{remark}{Remark}
\newcommand{\Def}{{\bf Definition.}\quad}
\newcommand{\Thm}{{\bf Theorem.}\quad}
\newcommand{\Prop}{{\bf Proposition.}\quad}
\newcommand{\Lem}{{\bf Lemma.}\quad}
\newcommand{\Ex}{{\bf Example.}\quad}
\newcommand{\Rem}{{\bf Remark.}\quad}
\newtheorem{example}{Example}
\newcommand{\ca }{\mathcal}
\newcommand{\CC}{{\mathbb{C}}}
\newcommand{\QQ}{{\mathbb{Q}}}
\newcommand{\RR}{{\mathbb{R}}}
\newcommand{\ZZ}{{\mathbb{Z}}}
\newcommand\Hom{\operatorname{Hom}}
\newcommand\End{\operatorname{End}}
\newcommand\Spec{\operatorname{Spec}}
\newcommand\CH{\operatorname{CH}}
\newcommand\Cone{\operatorname{Cone}}
\newcommand{\cA}{{\cal A}} \newcommand{\cB}{{\cal B}}\newcommand{\cC}{{\cal C}}
\newcommand{\cD}{{\cal D}}\newcommand{\cE}{{\cal E}}\newcommand{\cF}{{\cal F}}
\newcommand{\cG}{{\cal G}}\newcommand{\cH}{{\cal H}}\newcommand{\cI}{{\cal I}}
\newcommand{\cJ}{{\cal J}} \newcommand{\cK}{{\cal K}} \newcommand{\cL}{{\cal L}}\newcommand{\cM}{{\cal M}} \newcommand{\cN}{{\cal N}} \newcommand{\cO}{{\cal O}}\newcommand{\cP}{{\cal P}} \newcommand{\cQ}{{\cal Q}} \newcommand{\cR}{{\cal R}}\newcommand{\cS}{{\cal S}} \newcommand{\cT}{{\cal T}} 
\newcommand{\cZ}{{\cal Z}} \newcommand{\Z}{{\cal Z}}
\newcommand{\cU}{{\cal U}} \newcommand{\cV}{{\cal V}}
\newcommand{\bF}{{\mathbf{F}}}
\newcommand{\bU}{ {\mathbf{U}} }
\newcommand{\bK}{ {\mathbf{K}} }
\newcommand{\bL}{ {\mathbf{L}} }

\newcommand{\chZ}{ {\cal Z}} 

\newcommand{\al}{\alpha} \newcommand{\be}{\beta} \newcommand{\ga}{\gamma} 
\newcommand{\de}{\delta} \newcommand{\eps}{\epsilon} \newcommand{\la}{\lambda}
\newcommand{\De}{\Delta}

\newcommand{\BA}{{\mathbb{A}}} \newcommand{\BB}{{\mathbb{B}}}
\newcommand{\BC}{{\mathbb{C}}} \newcommand{\BD}{{\mathbb{D}}}
\newcommand{\BE}{{\mathbb{E}}} \newcommand{\BF}{{\mathbb{F}}}
\newcommand{\BG}{{\mathbb{G}}} \newcommand{\BH}{{\mathbb{H}}} 
\newcommand{\BI}{{\mathbb{I}}} 
\newcommand{\BIc}{{{I}}}
\newcommand{\bbI}{{\mathbb{I}}}
\newcommand{\BJ}{{\mathbb{J}}} \newcommand{\BK}{{\mathbb{K}}} 
\newcommand{\BL}{{\mathbb{L}}} \newcommand{\BM}{{\mathbb{M}}} 
\newcommand{\BN}{{\mathbb{N}}} \newcommand{\BO}{{\mathbb{O}}} 
\newcommand{\BP}{{\mathbb{P}}} \newcommand{\BQ}{{\mathbb{Q}}}
\newcommand{\BR}{{\mathbb{R}}} \newcommand{\BS}{{\mathbb{S}}}
\newcommand{\BT}{{\mathbb{T}}} \newcommand{\BU}{{\mathbb{U}}} 
\newcommand{\BV}{{\mathbb{V}}} 
\newcommand{\BW}{{\mathbb{W}}} 
\newcommand{\Bf}{{\mathbf{f}}}  \newcommand{\Bg}{{\bold{g}}}
\newcommand{\Bu}{{\bold{u}}}\newcommand{\Bv}{{\bold{v}}}
\newcommand{\Bw}{{\bold{w}}}

\newcommand{\sq}{\square} 
\newcommand{\ul}{\underline}
\newcommand{\ulm}{ {\ul{m}\,} } \newcommand{\uln}{ {\ul{n}} }  
\newcommand{\ot}{\otimes}
\newcommand{\scirc}{\circ} 
\newcommand{\hts}{\hat{\otimes}}
\newcommand{\htimes}{\hat{\times}}
\newcommand{\ts}{\otimes}
\newcommand{\vphi}{{\varphi}}  
\newcommand{\ddd}{\cdots}
\newcommand{\hlim}{\operatornamewithlimits{holim}} %
\newcommand{\bd}{\partial}
\newcommand{\sAb}{s Ab}       %
\newcommand{\stm}{\times}
\newcommand{\Tot}{\operatorname{Tot}}
\newcommand{\Ker}{\operatorname{Ker}}
\newcommand{\Cok}{\operatorname{Cok}}
\newcommand{\Gr}{\operatorname{Gr}}
\newcommand{\bop}{\oplus}
\newcommand{\moplus}{\mathop{\textstyle\bigoplus}}
\newcommand{\bts}{\mathop{\textstyle\bigotimes}}
\newcommand{\bhts}{\mathop{\textstyle \widehat{\bigotimes}}}
\newcommand{\Cyc}{\Z} \newcommand{\C}{\Z}
\newcommand{\injto}{\hookrightarrow}
\newcommand{\isoto}{\overset{\sim}{\to}}
\newcommand{\J}{\cJ}
\newcommand{\ti}{\tilde}
\newcommand{\CHM}{{C\!H\!\cal M}}
\newcommand{\dotsigma}{\dot{\sigma} }
\newcommand{\dotdelta}{\dot{\delta} }
\newcommand{\dotvphi}{\dot{\vphi} }
\newcommand{\maruphi}{\varphi^{in}}

\newcommand{\term}{\operatorname{tm}}
\newcommand{\tm}{\operatorname{tm}}
\newcommand{\ctop}{\overset{\,\,\circ}}
\newcommand{\init}{\operatorname{in}}

\newcommand{\Symb}{\mathop{Symb}}
\newcommand{\bound}{\partial}
\newcommand{\proj}{\pi}
\newcommand{\surjto}{\twoheadrightarrow}
\newcommand{\sgn}{\operatorname{sgn}}
\newcommand{\tbar}{\,\rceil \!\! \lceil}
\newcommand{\dbar}{|} 
\newcommand{\tp}{\mbox{top}}
\newcommand{\bsigma}{\mbox{\boldmath $\sigma$}}
\newcommand{\smallbsigma}{\mbox{\boldmath $\scriptstyle{\sigma}$}}
\newcommand{\bbrace}[1]{\mbox{\boldmath \{ }#1
\mbox{\boldmath \} }}

\newcommand{\bvphi}{\mbox{\boldmath $\vphi$}}
\newcommand{\bpartial}{\mbox{\boldmath $\partial$}}
\newcommand{\incl}{\operatorname{incl}}
\newcommand{\bDelta}{\mbox{\boldmath $\Delta$}}
\newcommand{\diag}{{\sl diag}}

\newcommand{\bskip}{\bigskip}
\newcommand{\sskip}{\smallskip}

\newcommand{\Seq}{{\BS eq}}

\newcommand{\mapright}[1]{%
  \smash{\mathop{%
    \hbox to
     1cm{\rightarrowfill}}\limits^{#1} } } 
\newcommand{\mapr}{\mapright}
\newcommand{\smapr}[1]{%
  \smash{\mathop{%
    \hbox to 0.5cm{\rightarrowfill}}\limits^{#1} } } 
\newcommand{\maprb}[1]{%
  \smash{\mathop{%
    \hbox to 1cm{\rightarrowfill}}\limits_{#1} } } 
\newcommand{\mapleft}[1]{%
  \smash{\mathop{%
    \hbox to 1cm{\leftarrowfill}}\limits^{#1} } }
\newcommand{\mapl}{\mapleft}
\newcommand{\maplb}[1]{%
  \smash{\mathop{%
    \hbox to 1cm{\leftarrowfill}}\limits_{#1} } }
\newcommand{\mapdown}[1]{\Big\downarrow
  \llap{$\vcenter{\hbox{$\scriptstyle#1\, $}}$ }}
\newcommand{\mapd}{\mapdown}
\newcommand{\mapdownr}[1]{\Big\downarrow
  \rlap{$\vcenter{\hbox{$\scriptstyle#1\, $}}$ }}
\newcommand{\mapdr}{\mapdownr}
\newcommand{\mapup}[1]{\Big\uparrow
  \llap{$\vcenter{\hbox{$\scriptstyle#1\, $}}$ }}
 \newcommand{\mapu}{\mapup}
\newcommand{\mapupr}[1]{\Big\uparrow
  \rlap{$\vcenter{\hbox{$\scriptstyle#1\, $}}$ }}
 \newcommand{\mapur}{\mapupr}

\newcommand{\sss}[1]{\noindent(#1)}
\newcommand{\mts}{\overset{\bullet}{\ts}}
\newcommand{\tildets}{{\tilde\ts}}
\newcommand{\shift}{\sl{shift}}
\newcommand{\tildetimes}{{\tilde\times}}
\newcommand{\tbullet}{^{\bullet\,\bullet\,\bullet} }
\newcommand{\dbullet}{^{\bullet\,\bullet} }

\newcommand{\partI}{\RefHaRelcorresp}

\newcommand{\RefBlone}{{[Blo-1]}}
\newcommand{\RefBltwo}{{[Blo-2]}}
\newcommand{\RefBlthree}{{[Blo-3]}}
\newcommand{\RefCH}{{[CH]}}
\newcommand{\RefFu}{{[Fu]}}
\newcommand{\RefHaoneone}{{[Ha-1]}}
\newcommand{\RefHaonetwo}{{[Ha-2]}}
\newcommand{\RefHaonethree}{{[Ha-3]}}
\newcommand{\RefHaDescent}{{[Ha-4]}}
\newcommand{\RefHaRelcorresp}{{[Ha-5]}}
\newcommand{\RefKa}{{[Ka]}}
\newcommand{\RefKaSh}{{[KS]}}
\newcommand{\RefLev}{{[Le]}}
\newcommand{\RefMacL}{{[Ma-1]}}
\newcommand{\RefMacLtwo}{{[Ma-2]}}
\newcommand{\RefTera}{{[Te]}}
\newcommand{\RefVe}{{[Ve]}}
\newcommand{\RefVo}{{[Vo]}}
\newcommand{\RefWe}{{[We]}}

\title{\bf Quasi DG categories and mixed motivic sheaves}
\author{Masaki Hanamura}
\date{}
\maketitle

\begin{abstract} 
We introduce the notion of a quasi DG category, 
generalizing that of a DG category. 
To a quasi DG category satisfying certain additional conditions, we associate
another quasi DG category, the quasi DG category of $C$-diagrams. 
We then show the homotopy category of the quasi DG category of $C$-diagrams has 
the structure of a triangulated category.  
This procedure is then applied to  produce a triangulated category 
of mixed motives over a base variety.
\end{abstract}

\renewcommand{\thefootnote}{\fnsymbol{footnote}}
\footnote[0]{ 2010 {\it Mathematics Subject Classification.} Primary 14C25; Secondary 14C15, 14C35.

 Key words: algebraic cycles, Chow group, motives. } 

\newcommand{\secBN}{0}
\newcommand{\secDG}{1}
\newcommand{\secfcn}{2}
\newcommand{\secCd}{3}
\newcommand{\secDS}{4}
\newcommand{\Multicpx}{0.1}
\newcommand{\Finiteorderedset}{0.2}
\newcommand{\Tensorproductcpx}{0.3}
\newcommand{\Phicpx}{0.4}
\newcommand{\DefDG}{\secDG.1}
\newcommand{\wqDG}{\secDG.2}
\newcommand{\DGwqDG}{\secDG.3}
\newcommand{\FtbarS}{\secDG.4}
\newcommand{\HowqDG}{\secDG.5}
\newcommand{\DefqDG}{\secDG.6} 
\newcommand{\RmkDefqDG}{\secDG.7}
\newcommand{\ExqDG}{\secDG.8}
\newcommand{\SymbqDG}{\secDG.9}
\newcommand{\Cdiag}{\secDG.10}

\newcommand{\lsum}{s}
\newcommand{\rsum}{t}
\newcommand{\Acyclicitysigma}{(\DefqDG), (6)}
\newcommand{\Compatibilitydiagsigma}{(\DefqDG), (7)}
\newcommand{\Propintersection}{(\DefqDG), (8)}
\newcommand{\propPIone}{1}
\newcommand{\propPItwo}{2}
\newcommand{\propPIthree}{3}
\newcommand{\propPIfour}{4}
\newcommand{\propPIfive}{5}
\newcommand{\AdditivityGISigma}{3.b}
\newcommand{\FMM}{\secfcn.3}
\newcommand{\DefCdiagram}{\secfcn.4}
\newcommand{\Defbsigma}{\secfcn.5}
\newcommand{\GKL}{\secfcn.6}
\newcommand{\HIGI}{\secfcn.7}
\newcommand{\acyclicGI}{\secfcn.8}
\newcommand{\htsPhicpx}{\secfcn.9}
\newcommand{\GISigma}{\secfcn.10} 
\newcommand{\rhoPi}{\secfcn.11}
\newcommand{\CompDoublecpx}{\secfcn.12}
\newcommand{\GIT}{\secfcn.13}
\newcommand{\GITSigma}{\secfcn.14}
\newcommand{\BFIS}{\secfcn.15}
\newcommand{\DefCdzero}{\secfcn.16}
\newcommand{\propFI}{\secfcn.17}
\newcommand{\Proofacyclicity}{\secfcn.18}
\newcommand{\PlanCd}{\secCd.1}
\newcommand{\ZFKL}{\secCd.2}
\newcommand{\GKLM}{\secCd.3}
\newcommand{\Compos}{\secCd.4}
\newcommand{\Identity}{\secCd.5}
\newcommand{\IdentityProp}{\secCd.6}
\newcommand{\Sigmaprolong}{\secCd.7}
\newcommand{\Propintersectiontheta}{\secCd.8}
\newcommand{\AdditivityBGKL}{\secCd.9}
\newcommand{\AdditivityBG}{\secCd.10}
\newcommand{\AdditivityBF}{\secCd.11}
\newcommand{\ThmCdqDG}{\secCd.12}
\newcommand{\ThmCdqDGtriang}{\secCd.13}
\newcommand{\Shiftfunctor}{\secCd.14}
\newcommand{\Conemorphism}{\secCd.15}
\newcommand{\Axiomproof}{\secCd.16}
\newcommand{\Axiomprooftwo}{\secCd.17}
\newcommand{\DefDS}{\secDS.1}
\newcommand{\ThmDS}{\secDS.2}

\newcommand{\PH}{P\!H}
\newcommand{\PBG}{P\!\BG}
\newcommand{\ph}{\phantom}
\newcommand{\vecta}{\mbox{\boldmath $a$}}
\newcommand{\vectb}{\mbox{\boldmath $b$}}

In this paper we introduce the notion of a {\it quasi DG category\/},  and give
a procedure to construct a triangulated category associated to it. 
This method is then applied to the construction of the triangulated category 
of mixed motivic sheaves over a base variety.

If the base variety is the Spec of the ground field, this coincides with the 
triangulated category of motives as in \RefHaoneone\, and \RefHaonetwo.
As for the construction of the triangulated category of mixed motives, 
there have been approaches by M. Levine and by V. Voevodsky as well
(see \RefLev\, and \RefVo).  

The main idea in \RefHaoneone\,, \RefHaonetwo\,,   is as follows. 
 We start with the class of smosoth projective
varieties $X$ over the base field, 
the cycles complexes of the product of those varieties
$\Z(X\times Y, \bullet)$, and the composition maps 
$\Z(X\times Y, \bullet)\ts \Z(Y\times Z, \bullet)\to \Z(X\times Z, \bullet)$.
We refer to \RefBlone, \RefBltwo\, and  \RefBlthree\  for the cycle complex and higher 
Chow groups.
Although the composition is only partially defined, it is defined on 
a quasi-isomorphic subcomplex.
In other words, we have almost a DG category where $\Z(X\times Y, \bullet)$ is the group of homomorphisms from $X$ to $Y$. 
Out of these we construct the triangualted
category in which the objects (which are called {\it $C$-diagrams})
and morphisms can be explicitly described in terms of varieties and cycles.

In order to extend this idea and construct the category of mixed motives over
any base variety $S$, we encounter the following problem. 
With the construction of the category of pure motives over $S$ 
 in mind (see \RefCH\,), 
one takes the class of smooth quasi-projective varieties $X$, equipped with 
projective maps $X\to S$, and considers the cycle complex $Z(X\times_S Y, \bullet)$, 
where $X$ and $Y$ are both smooth varieties with projective maps to $S$. 
There is, however,  no partially defined composition map 
$\Z(X\times_S Y, \bullet)\ts \Z(Y\times_S Z, \bullet)\to \Z(X\times_S Z, \bullet)$.
So the construction in \RefHaoneone\,, \RefHaonetwo\,
cannot be immediately generalized. 
The present paper provides a solution to this problem on the categorial aspect:
 If we are given a quasi DG category, 
we provide a method to produce an associated triagulated category.

The notion of a quasi DG category is a generalization of that of 
a DG category.  Recall that a DG category is an additive category $\cC$, such that
for a pair of objects $X, Y$ the group of homomorphisms $F(X, Y)$ has the structure 
of a complex of abelian groups, 
and the composition $F(X, Y)\ts F(Y, Z)\to F(X, Z)$ is a map of 
complexes.  
(Our sign convention for the tensor product of complexes
differs from the usual one, see \S 0. 
Accordingly, the domain of the composition map is $F(X, Y)\ts F(Y, Z)$
instead of $F(Y, Z)\ts F(X, Y)$.)

A quasi DG category (which is not really a category)
 consisits of an underlying category $\cC_0$, 
 and to a pair of objects there corresponds a complex $F(X, Y)$. 
One would like to think of $F(X, Y)$ as the group of homomorphisms.
But there is no 
composition map $F(X, Y)\ts F(Y, Z)\to F(X, Z)$. 
Instead, we have the following structure:
\smallskip 

(1) There is 
a quasi-isomorphic double subcomplex of 
$F(X, Y)\ts F(Y, Z)$, which we denote by 
$ F(X, Y)\hts F(Y, Z)$; so there is an injective quasi-isomorphism
$$\iota: F(X, Y)\hts F(Y, Z)\injto F(X, Y)\ts F(Y, Z)\,\,.$$

(2) There is another complex $F(X, Y, Z)$ and a surjective quasi-isomorphism 
$$\sigma: F(X, Y, Z)\to F(X, Y)\hts F(Y, Z)\,\,.$$

(3) There is a map of complexes 
$\vphi: F(X, Y, Z)\to F(X, Z)\,\,.$
\smallskip 

\noindent In the derived category at least, one has an induced map 
$$\psi: F(X, Y)\ts F(Y, Z)\to F(X, Z)$$
obtained by inverting the quasi-isomorphism $\iota\sigma: F(X, Y, Z)\to 
F(X, Y)\ts F(Y, Z)$, and 
composing it with $\vphi$. 
In particular one has  the induced map on cohomology 
$\psi: 
H^0F(X, Y)\ts H^0F(Y, Z)\to 
H^0F(X, Z)$. 

We naturally 
require that the above pattern persists for more than three objects as follows. 
\smallskip 

(1) For each sequence of objects 
 $X_1, \ddd , X_n$ ($n\ge 2$), 
we are given  a complex  of  abelian groups 
$F(X_1, \cdots , X_n)$. 

Let $(1, n)$ denote the subset $\{2, \cdots, n-1\}$. 
 For a subset of integers 
 $S=\{i_1, \cdots, i_{a-1}\}\subset (1, n)$, let 
$i_0=1$, $i_a=n$ and define an $a$-tuple complex by 
$$F(X_1, \cdots , X_n\tbar S):=F(X_{i_0}, \cdots, X_{i_1})\ts
F(X_{i_1}, \cdots, X_{i_2})\ts\cdots\ts
F(X_{i_{a-1}}, \cdots, X_{i_a})\,\,.
$$
We are also given  an $a$-tuple  complex  $F(X_1, \cdots , X_n| S)$ 
and an injective 
quasi-isomorphism 
$$\iota_S:F(X_1, \cdots , X_n| S)\injto F(X_1, \cdots , X_n\tbar  S)\,\,.$$
We assume  $F(X_1, \cdots , X_n\dbar \emptyset)
=F(X_1,\cdots, X_n)$.

(2) To an inclusion $S\subset S'$ there corresponds a surjective quasi-isomorphism 
$$\sigma_{S\, S'}: F(X_1, \cdots , X_n\dbar S)\to 
F(X_1, \cdots , X_n\dbar S')\,\,.$$
One has $\sigma_{S\, S}=id$, and for 
 $S\subset S'\subset S''$, $\sigma_{S\,S''}=
\sigma_{S'\,S''}\sigma_{S\,S'}$.  In particular we have 
a surjective quasi-isomorphism
$$\sigma_S:=\sigma_{\emptyset\, S} : F(X_1, \cdots , X_n)\to 
F(X_1, \cdots , X_n\dbar S)\,.$$

(3)  For a subset 
$K=\{k_1, \ddd , k_b\}\subset (1, n)$ disjoint from $S$,  there is a
 map of complexes
$$\vphi_K: F(X_1, \ddd, X_n | S)\to F(X_1, \ddd,\widehat{X_{k_1}}, 
\ddd , \widehat{X_{k_b}}, \ddd ,  X_n | S)\,\,.$$
If $K$ is the disjoint union of $K'$ and $K''$, one has 
$\vphi_K=\vphi_{K'}\vphi_{K''}$. 

 If $K$ and $S'$ are disjoint $\sigma_{S\, S'}$ and $\vphi_K$ commute.  
The injection $\iota_S$ and the maps $\sigma$, $\vphi$ are compatible 
(see (\DefqDG) for what this means).  
From these requirements it follows that the composition 
map $\psi: 
H^0F(X, Y)\ts H^0F(Y, Z)\to 
H^0F(X, Z)$ 
is associative. 
\bigskip 

Thus a quasi DG category consists of the category $\cC_0$, 
the complexes $F(X_1, \cdots, X_n|S)$, and the maps 
$\iota$, $\sigma$, and $\vphi$ satisfying the above conditions.
A few other conditions are required, as we briefly explain
the idea below, and we refer to (\DefqDG) for details. 
\smallskip 

(4) The category $\cC_0$ is a groupoid, and it is also a 
symmmetric monoidal category with respect to a functor 
$(X,  Y)\mapsto X\oplus Y$ called the direct sum, and 
an object $O$ called the zero object. 
 The multiple complexes $F(X_1, \ddd, X_n | S)$
are covariantly  functorial in each variable $X_i$. 
The maps $\iota$, $\sigma$, and $\vphi$ are compatible with 
 the covariant functoriality.
 In addition, $F(X_1, \ddd, X_n)$ is assumed additive in each variable 
in an appropriate sense, elaborated in (\DefqDG). 
 
(5) For each object $X$ there is the identity element $1_X$ in $H^0F(X, X)$. 
\smallskip 

The first example of a quasi DG category 
comes from a DG category.
 Indeed given a DG category, and a sequence of objects $X_1, \cdots, X_n$, let 
$$F(X_1, \cdots, X_n)=F(X_1, X_2)\ts\cdots\ts F(X_{n-1}, X_n)\,,$$
and $F(X_1, \cdots, X_n|S)=F(X_1, \cdots, X_n)$ for any $S$;
take $\sigma_{S\, S'}$ to be identities, and 
$\vphi_K$ to be the composition at $X_k$, $k\in K$.

To a quasi DG category $\cC$ one can associate an additive category 
 $Ho(\cC)$, called the 
{\it homotopy category} of $\cC$; it is  the  category in which 
the objects are the same as for $\cC$, and the homomorphism 
group is $H^0 F(X, Y)$, and the composition is the map $\psi$ above. 

Recall that given an additive category $\cA$, one has the DG category 
of complexes with values in $\cA$, and its homotopy category is a 
triangulated category. 
We will give an analogous construction for a quasi DG category.
In \S\S 2 and 3, we take a quasi DG category $\cC$ and construct 
a related quasi DG category denoted $\cC^\Delta$. 
(For this we require on $\cC$ two extra structure (\DefqDG), (iv), (v). 
These are the existence of diagonal elements and diagonal extension 
 (\DefqDG)(iv), and 
the existence of a generating set, notion of proper intersection, and distinguished subcomplexes, (\DefqDG)(v).)
 For such $\cC$, we produce a related 
 quasi DG category $\cC^\Delta$, where the objects are what 
 we call {\it $C$-diagrams\/} in $\cC$, see (1.9). 
 A $C$-diagram is of the form 
$K= (K^m; f(m_1, \cdots  , m_\mu)\, )$, 
where $(K^m)$ is a sequence of objects of $\cC$
indexed by $m\in \ZZ$, almost all
of which are zero, and 
$$f(m_1, \cdots , m_\mu)\in F(K^{m_1}, \cdots , K^{m_\mu})^{-(m_\mu-m_1-\mu+1)}$$
is a collection of elements indexed by sequences $(m_1<m_2<\cdots<m_\mu)$
with $\mu\ge 2$, satisfying the following conditions:
\smallskip 

(i) For each $j$ with $1<j<\mu$, 
$$\sigma_{K^{m_j}}(f(m_1, \cdots , m_\mu)\, )
=f(m_1, \cdots , m_j)\ts f(m_j, \cdots , m_\mu)\,\, $$
in $F(K^{m_1}, \cdots , K^{m_j})\ts F(K^{m_j}, \cdots , K^{m_\mu})$. 

(ii) For each $(m_1, \cdots , m_\mu)$, one has 
$$\bound f(m_1, \cdots , m_\mu)+\sum_{1\le t<\mu}\,
\sum_{m_t<k<m_{t+1}}
(-1)^{\,m_\mu+\mu+k+t}\vphi_{K^{m_k}}(f(m_1, \cdots ,m_t, k, m_{t+1}, \cdots,  m_\mu)\, )=0\,\,.$$
(Here $\partial$ is the differential of the complex $F(K^{m_1}, \cdots, K^{m_\mu})$.)
\smallskip 

\noindent One observes that $K$ is like a complex with terms $K^m$ in degree $m$, 
and with differentials $f(m, n)$ from $K^m$ to $K^n$ for $m<n$. 
Note also that given
an object $X$ of $\cC$ and $n\in \ZZ$, one can define a $C$-diagram $K$, by setting 
$K^n=X$, $K^m=O$ for $m\neq n$, and $f(m_1, \cdots , m_\mu)=0$
for all sequences $(m_1<m_2<\cdots<m_\mu)$.  We write $X[-n]$ for $K$.

For the class of $C$-diagrams to form a quasi DG category, we
must define the complexes $\BF(K_1, \cdots , K_n)$ for  a sequence of $C$-diagrams $K_1, \cdots, K_n$, 
 together with the maps $\sigma$ and $\vphi$. 
 We carry out the construction of these  in \S 2,  and verify the
 axioms of a quasi DG category in \S\S 2 and 3. 
 If $X$, $Y$ are objects of $\cC$ that can be viewed as $C$-diagrams, 
 one has the identity $\BF(X, Y)=F(X, Y)$. 
 With the quasi DG category $\cC^\Delta$ thus obtained, 
one can  consider its homotopy category. 
The main result of this paper is stated as follows 
(see (\ThmCdqDG) and (\ThmCdqDGtriang)\,).
 \bigskip 
 
 {\bf Main Theorem.}\quad {\it Let $\cC$ be a quasi DG 
category satisfying the conditions (iv),(v) of (\DefqDG). 
 Let $\cC^\Delta$ be the quasi DG category of $C$-diagrams in $\cC$, and $Ho(\cC^\Delta)$ its homotopy
 category. Then we have 

(I) For an object $X$ of $\cC$ and $n\in \ZZ$, there corresponds an object
$X[n]$ in $\cC^\Delta$.
For two objects $X, Y $ of $\cC$ and $m, n\in \ZZ$ we have 
$$\Hom_{Ho(\cC^\Delta)}(X[m], Y[n])=H^{n-m}F(X, Y)\,\,.$$
The composition
of morphisms
between three objects $X[m]$, $Y[n]$, $Z[\ell]$ (where 
$X, Y, Z$ are objects of $\cC$ and $m, n, \ell \in \ZZ$), 
$$H^0\BF(X[m], Y[n])\ts H^0\BF(Y[n], Z[\ell])\to H^0\BF(X[m], Z[\ell])$$
is identified via the above isomorphisms with 
the map 
$$H^{n-m}F(X, Y)\ts H^{\ell-n}F(Y, Z)\to H^{\ell -m}F(X, Z)\,.$$
It coincides with the map induced on cohomology from the map in the derived category 
$$F(X, Y)\ts F(Y, Z) \to F(X, Z)\,, $$ 
obtained by composing the inverse of $\sigma$ and $\vphi$. 

(II) The additive category $Ho(\cC^\Delta)$ has the structure of a triangulated category. 
 }\bigskip

So far there is no geometry involved. 
For us the main example of a quasi DG category is that of symbols over a quasi-projective variety $S$, 
 denoted $Symb(S)$, see (\SymbqDG) for details. 
A typical object of $Symb(S)$ is of the form $(X, r)$, where $r$ is an integer and $X$ a smooth variety with 
a projective map to $S$. 
For two such objects $(X, r)$ and $(Y, s)$, the corresponding complex $F((X, r), (Y, s)\, )$ is quasi-isomorphic to 
$\Z_{\dim Y -s+r}(X\times_S Y, \bullet)$, the cycle complex of the fiber product $X\times_S Y$. 
We refer to \partI \, for the  construction of the complexes 
$$F((X_1, r_1), \cdots, (X_n, r_n) |S)$$
and the maps $\iota_S$, $\sigma_{S\, S'}$, and $\vphi_K$. 
The additional conditions (\DefqDG), (iv) and (v) are satisfied for $Symb(S)$. 

In \S 4 we apply the construction of \S \S 2 and 3 to $Symb(S)$. 
The resulting triangulated category $\cD(S):=Ho(Symb(S)^\Delta)$ 
is by definition the triangulated category of mixed motives over 
$S$.  See (\ThmDS) for the properties of $\cD(S)$. 
For  $S=\Spec k$ the construction of the triangulated category 
in \RefHaoneone \,, \RefHaoneone \, is similar
 but simpler since then $\Symb(S)$ is (almost) a
DG category, and the notion of $C$-diagrams is simpler. 
(Essentially the same idea appeared in \RefKa,\, preceding \RefHaoneone.)
It is useful to have a construction of $\cD(k)$ via $C$-diagrams, 
because one can construct objects concretely using cycles; 
see for example \RefTera. 

We collected basic notions in \S 0 regarding multiple complexes and finite
totally ordered sets. 
 
Regarding the technical aspects, we point out two problems for the reader. 
The first is the delicate question of signs; we give the basic argument in \S 0, and 
wherever there is additional issue we elaborated on it. 
The signs for the complex $\BH(K, L)$ in (\GKL) require the most care. 
The second is the question of general positions (choice of distinguished subcomplexes); 
the notion of proper intersection and distinguished subcomplexes in 
(\DefqDG) are designed to resolve such problems, allowing us to take distinguished 
subcomplexes as suited for our purposes. The places where we use this are
(\GKL), (\HIGI),  (\GISigma) and (\GITSigma).

In case $S=\Spec k$, the work \RefHaoneone, \RefHaonetwo, \RefHaonethree, 
\RefHaDescent\, deals with not only the 
construction of the triagulated category $\cD(k)$, but also the cohomology 
realization functor and the functor of cohomological motives (which associates to
each quasi-projective variety $X$ its motive $h(X)$ in $\cD(k)$). 
We will discuss these problems in a separate paper. 
\bigskip 

{\bf Acknowledgements.}
We would like to thank S. Bloch,  B. Kahn  and P. May for helpful discussions. 
We also gratefully acknowledge the helpful comments of the referee, which led
to a significant improvement of the paper.

\setcounter{section}{-1}

\section{Basic notions.}

Subsections (\Multicpx) and (\Finiteorderedset) are used throughout this paper, 
(\Tensorproductcpx) and (\Phicpx) are needed in \S\S 2 and 3. 
\bigskip

\sss{\Multicpx} {\it Multiple complexes.}\quad
By a complex of abelian groups we mean a graded abelian group $A^\bullet$
with a map $d$ of degree one satisfying $dd=0$. 
If $u: A\to B$ and $v: B\to C$ are maps of complexes, 
we define $u\cdot v: A\to C$ by $(u\cdot v)(x)=v(u(x))$. 
So $u\cdot v$ is $v\scirc u$ in the usual notation. 
As usual we also write $vu$ for $v\scirc u$ (but not for $v\cdot u$).

A double complex $A=(A^{i, j}; d', d'')$ is a bi-graded abelian group with 
differentials $d'$ of degree $(1, 0)$, $d''$ of degree $(0, 1)$, satisfying
$d'd''+ d''d'=0$. Its total complex $\Tot(A)$ is the complex with 
$\Tot(A)^k=\bop_{i+j=k}A^{i, j}$ and the differential $d=d'+d''$. 
In contrast a ``double" complex 
$A=(A^{i, j}; d', d'')$ is a bi-graded abelian group with 
differentials $d'$ of degree $(1, 0)$, $d''$ of degree $(0, 1)$, satisfying
$d'd''=d''d'$. Its total complex $\Tot(A)$ is given by 
$\Tot(A)^k=\bop_{i+j=k}A^{i, j}$ and the differential $d=d'+(-1)^i d''$ on 
$A^{i, j}$.  
(Note that the totalization depends on the ordering of the 
two gradings; if we reverse the order, the corresponding totalization has 
differential $(-1)^j d'+ d''$ on $A^{i, j}$.)
A ``double" complex can be viewed as a double complex by taking the differentials 
to be $(d', (-1)^i d'')$ (or, when we reverse the order, $(\,(-1)^j d', d'')$).

Let $(A, d_A)$ and $(B, d_B)$ be complexes. Then 
$(A^{i, j}=A^{j}\ts B^i; d_A\ts 1, 1\ts d_B)$ is a ``double" complex.
Its total complex has 
differential $d$ given by 
$$d(x\ts y)= (-1)^{\deg y} dx\ts y + x\ts dy\,\,$$
(for the tensor product complex, 
we always take the reverse order of the gradings for the totalization). 
Note this differs from the usual convention. 

More generally for
$n\ge 2$ one has the notion of $n$-tuple complex and ``$n$-tuple" 
complex. An $n$-tuple (resp. ``$n$-tuple") complex is 
a $\ZZ^n$-graded abelian group $A^{i_1, \cdots, i_n}$ with differentials
$d_1, \cdots, d_n$, $d_k$ raising $i_k$ by 1, such that 
for $k\neq \ell$, $d_kd_\ell+ d_\ell d_k=0$ (resp. $d_kd_\ell=d_\ell d_k$).
An ``$n$-tuple" complex $A^{i_1, \cdots, i_n}$ is an $n$-tuple complex
with respect to the differentials $(d_1, (-1)^{\al_1} d_2, \cdots, (-1)^{\al_n} d_n)$
where $\al_k=\sum_{j<i} \deg_j$.  
In this way we turn an ``$n$-tuple" into an $n$-tuple complex. 
(As for double complexes, one may reverse the order of the gradings; in that case 
we will explicitly mention it.)
A single
complex $\Tot(A)$, called the total complex, is defined in either case. 

As a variant one can define partial totalization. 
To explain it, 
let $S_1, \cdots, S_m$ be an ordered set of non-empty
subsets of $[1, n]:=\{1, \cdots, n\}$
 such that $S_i\cap S_j=\emptyset$ for $i\neq j$ and $\cup S_i=[1, n]$.
Such data corresponds to a surjective map $f: [1, n]\to [1, m]$. 
Then, given an $n$-tuple complex (resp. an ``$n$-tuple" complex)
$A^{i_1, \cdots, i_n}$
 one can ``totalize" in degrees in $S_i$, and form an 
 $m$-tuple (resp. ``$m$-tuple") complex denoted 
$\Tot^{S_1, \cdots, S_m}(A)$ or $\Tot^f(A)$. 
Given surjective maps $f: [1, n]\to [1, m]$ and $g: [1, m]\to [1, \ell]$, one has 
$\Tot^g\Tot^f(A)=\Tot^{gf}(A)$.  
For example, if a subset $S=[k, \ell]\subset [1, n]$ is specified, 
one can ``totalize"
in degrees in $S$, so the result $\Tot^S(A)$ is an $m$-tuple 
(resp. ``$m$-tuple") complex, where $m=n-|S|+1$. 

For $n$ complexes $A_1^\bullet, \cdots, A_n^\bullet$, the tensor
product $A_1^\bullet\ts \cdots \ts A_n^\bullet$ is an ``$n$-tuple" 
complex.  For the tensor product, we view it as an $n$-tuple complex 
with respect to the 
reverse order of the gradings.  This $n$-tuple complex and its 
total complex will be still denoted $A_1^\bullet\ts\cdots\ts A_n^\bullet$. 

The only difference between $n$-tuple and ``$n$-tuple" complexes is 
that of signs. 
For the rest of this section, what we say for 
$n$-tuple complexes equally applies to ``$n$-tuple" complexes. 

If $A$ is an $n$-tuple complex and $B$ an $m$-tuple complex, and when 
$S=[k, \ell]\subset [1, n]$ with $m=n-|S|+1$ is specified, one can talk of 
maps of $m$-tuple complexes $\Tot^S(A)\to B$.  When the choice of 
$S$ is obvious from the context, we just say maps of multiple complexes
$A\to B$. For example if $A$ is an $n$-tuple complex and $B$ an 
$(n-1)$-tuple complex, for each set $S=[k, k+1]$ in $[1, n]$ one can 
speak of maps of $(n-1)$-tuple complexes $\Tot^S(A)\to B$; if $n=2$
there is no ambiguity. 
\bigskip 

(\Multicpx.1) {\it Multiple subcomplexes of a tensor product complex.}\quad
Let $A$ and $B$ be complexes. A  double subcomplex $C^{i, j}\subset
A^i\ts B^j$ is a submodule closed under the two differentials. 
If $\Tot(C)\injto \Tot(A\ts B)$ is a quasi-isomorphism, we say 
$C^{\bullet\, \bullet}$ is a quasi-isomorphic subcomplex. 
It is convenient to let $A^\bullet\hts B^\bullet$ denote
such a subcomplex. 
(Note it does not mean the tensor product of subcomplexes of $A$ and $B$.)
Likewise a quasi-isomorphic 
multiple subcomplex of $A_1^\bullet\ts\cdots \ts A_n^\bullet$
is denoted $A_1^\bullet\hts\cdots \hts A_n^\bullet$. 
\bigskip 

\sss{\Finiteorderedset}
{\it Finite ordered sets, partitions and segmentations.}\quad
Let $I$ be a non-empty finite totally ordered 
set (we will simply  say a finite ordered set), so
$I=\{i_1, \cdots, i_n\}, i_1<\cdots <i_n$, where $n=|I|$. 
 The {\it initial} (resp. {\it terminal})
element of $I$ is $i_1$ (resp. $i_n$); let $\init(I)=i_1$, $\term(I)=i_n$.
If $n\ge 2$, let $\ctop{I}=I-\{\init(I), \term(I)\}$. 

If $I=\{i_1, \cdots, i_n\}$, a subset $I'$ of the form 
$[i_a, i_b]=\{i_a, \cdots, i_b\}$ is called a {\it sub-interval}. 

In the main body of the paper, for the sake of concreteness 
 we often assume $I=[1, n]=\{1, \cdots, n\}$, a subset of 
$\ZZ$. More generally a finite subset of $\ZZ$ is an example of a finite 
ordered set. 

A {\it partition} of $I$ is a disjoint decomposition into 
sub-intervals $I_1, \cdots, I_a$ such that there is a sequence of 
elememtns $\init(I)=i_0< i_1< \cdots < i_{a-1}< i_a=\term(I)$ so that 
$I_k=[i_{k-1}, i_k -1]$. 

 So far we have assumed $I$ and $I_i$ to be of cardinality $\ge 1$. 
In some contexts we allow only finite ordered sets 
with at least two elements. 
There instead of partition the following notion plays a role.
Given a subset of $\ctop{I}$, $\Sigma= \{i_1, \cdots, i_{a-1}\}$,
where $i_1<i_2<\cdots< i_{a-1}$, one has 
a decomposition of $I$ into the sub-intervals $I_1, \cdots, I_a$, where
$I_k=[i_{k-1}, i_k]$, with $i_0=\init(I)$, $i_a=\term(I)$. 
Thus the sub-intervals satisfy 
$I_k\cap I_{k+1}=\{i_k\}$ for $k=1, \cdots, a-1$. 
 The sequence of sub-intervals $I_1, \cdots, I_a$ is
  called the {\it segmentation } of $I$ corresponding to $\Sigma$. 
(The terminology is adopted  to distinguish it from 
 the  partition). 
\bigskip

\sss{\Tensorproductcpx} {\it  Tensor product of ``double'' complexes.}
\quad Let $A\dbullet=(A^{a,p}; d'_A, d''_A)$ be a ``double'' complex
 (so  $d'$ has degree $(1, 0)$, $d''$ has degree $(0,1)$, and $d'd'=0$, 
$d''d''=0$ and $d'd''=d''d'$). The associated total complex $\Tot(A)$
has differential $d_A$ given by $d_A= d'+(-1)^a d''$ on $A^{a,p}$. 
The association $A\mapsto \Tot(A)$ forms a functor.
Let $(B^{b,q}; d'_B, d''_B)$ be another ``double'' complex. 
Then the tensor product of $A$ and $B$  as ``double'' complexes, 
denoted $A\dbullet\times B\dbullet$, 
is by definition the ``double'' complex $(E^{c,r}; d'_E, d''_E)$, where
$$E^{c,r}=\bop_{a+b=c\,,p+q=r}A^{a,p}\ts B^{b,q}$$
and $d'_E=(-1)^b d'_A\ts 1+  1\ts d'_B$, $d''_E=(-1)^q 
d''_A\ts 1+ 1\ts 
d''_B$. 

The tensor product complex $\Tot(A)\ts \Tot(B)$ and 
the total complex of $A\dbullet\times B\dbullet$ are 
related as follows. There is an isomorphism 
of complexes 
$$u: \Tot(A)\ts \Tot (B)\to \Tot(A\dbullet\times B\dbullet)$$
given by 
$u=(-1)^{aq}\cdot id$ on the summand $A^{a,p}\ts B^{b,q}$.

Let $A\dbullet$, $B\dbullet$, $C\dbullet$ be ``double'' complexes. One has
an obvious isomorphism of ``double'' 
complexes
$(A\dbullet\times B\dbullet)\times C\dbullet=A\dbullet\times 
(B\dbullet\times C\dbullet)$; it is denoted 
$A\times B\times C$. 
We will often suppress the double dots for simplicity. 
The following diagram commutes:
$$
\begin{array}{ccc}
\Tot(A)\ts\Tot(B)\ts\Tot(C)&\mapr{u\ts 1} &\Tot(A\times B)\ts\Tot(C) \\
\mapd{1\ts u}& &\mapdr{u} \\
\Tot(A)\ts\Tot(B\times C)&\mapr{u} &\Tot(A\times B\times C)\,\,.
\end{array}
$$
The composition defines an isomorphism 
 $u: \Tot(A)\ts \Tot(B)\ts \Tot(C)\isoto \Tot(A\times B\times C)$.
 
One can generalize this to the case of tensor product of more
than two ``double'' complexes. 
If $A_1, \cdots, A_n$ are ``double'' complexes, there is an isomorphism
of complexes 
$$u_n:\Tot(A_1)\ts\cdots\ts\Tot(A_n)\to \Tot(A_1\times\cdots\times A_n)$$
which coincides with the above $u$ if $n=2$, and is in general a composition
of $u$'s in any order.
As in case $n=3$, one has commutative diagrams involving $u$'s;
we leave the details to the reader. 

Let $A$, $B$, $C$ be ``double'' complexes and 
 $\rho: A\dbullet\times B\dbullet\to 
C\dbullet$
be a map of ``double'' complexes, 
namely it is bilinear and 
for $\al\in A^{a,p}$ and $\be\in B^{b,q}$,
$$d' \rho (\al\ts\be)= \rho((-1)^b d'\al\ts\be+  \al\ts d'\be)$$
and 
$$d'' \rho (\al\ts\be)= \rho((-1)^q d''\al\ts\be+  \al\ts d''\be)\,\,.$$
Composing $\Tot(\rho): \Tot(A\times B)\to \Tot(C)$ with $u:
\Tot(A)\ts\Tot(B)\isoto \Tot(A\times B)$, one obtains the map 
$$\hat{\rho}: \Tot(A)\ts \Tot(B)\to \Tot(C)\,\,;$$
it is given given by $(-1)^{aq}\cdot 
\rho$ on the summand $A^{a,p}\ts B^{b,q }$. 

The same holds for a map of ``double'' complexes 
$\rho: A_1\times \cdots\times A_n\to C$.
\smallskip 

{\it Remark.}\quad One could discuss more general sign rules for 
the change of ordering of the set of gradings of multiple complexes. 
We have restricted our discussions to the case we will need in \S 2. 
\bigskip 

\sss{\Phicpx}{\it The subcomplex $\Phi A$.}\quad
Let $(A\dbullet; d_1, d_2)$ be a ``double" complex satisfying

(i) $A^{a, p}=0$ if $p<0$, 

(ii) The sequence of complexes 
$$A^{\bullet\, 0}\mapr{d_2}A^{\bullet\, 1}\mapr{d_2}\cdots $$
is exact. 
In other words, for each $a$, the complex $A^{a, \bullet}$ sastisfies:
$H^i(A^{a, \bullet})=0$ for $i\neq 0$. 
\smallskip 

We then let $\Phi( A)^\bullet $ be the kernel of $d_2:A^{\bullet\, 0}\to A^{\bullet\, 1}$. 
Then one has an exact sequence of complexes
$$0\to \Phi( A)^\bullet\to A^{\bullet, 0}\mapr{d_2} A^{\bullet, 1}\mapr{d_2}\cdots\,\,.$$
Thus $\Phi( A)^\bullet$ is a complex with differential $d_1$, and the inclusion 
$\Phi( A)^\bullet\injto \Tot(A^{\bullet\,\bullet})$ is a quasi-isomorphism. 
The association $A\mapsto \Phi A$ forms an exact functor from the 
category of ``double" complexes satisfying the conditions (i), (ii)
to the category of complexes.  
$\Phi A$ is so to speak the peripheral complex of $A^{\bullet\,\bullet}$ in the 
second direction. 

This can be generalized to the case of ``$n$-tuple" complexes $(A^{\bullet\,\cdots\, \bullet}; d_1, \cdots, 
d_n)$ satisfying the conditions similar to (i), (ii) with respect to the last degree and differential. 
Then $\Phi (A) =\Ker (d_n)$ is an ``$(n-1)$-tuple" complex, and the inclusion 
$\Phi(A)\injto A^{\bullet\,\cdots\, \bullet}$ is a quasi-isomorphism on total complexes. 
If $A^{\bullet\,\bullet\, \bullet}$ is a ``triple" complex, for example, we put double 
dots, as in $\Phi(A)^{\dbullet}$, to indicate it is a ``double" complex. 
\bigskip

The operation $\Phi$ is compatible with tensor product as follows. 
If $A\dbullet$ and $B\dbullet$ are ``double" complexes satisfying (i), (ii), then $A\dbullet\ts B\dbullet$ is a ``quadruple" complex.  Taking totalization with respect to the second and fourth 
degree, one obtains a ``triple" complex, with three differentials 
$d_1\ts 1$, $1\ts d_1$ and $d_2\ts 1\pm 1\ts d_2$. 
The ``triple" complex 
 $\Tot_{24}(A\dbullet\ts B\dbullet)$ thus obtained satisfies the condition (i), (ii) with respect
to the third degree, and one has 
$$\Phi(A)^\bullet\ts \Phi (B)^\bullet= \Phi (\Tot_{24}(A\dbullet\ts B\dbullet)\,)
\eqno{(\Phicpx.a)}$$
as a ``double" complex. 
Indeed $H^i(A^{a, \bullet}\ts B^{b, \bullet})=0$ for $i\neq 0$ and 
$H^0(A^{a, \bullet}\ts B^{b, \bullet})=H^0(A^{a, \bullet}\,)\ts H^0(B^{b, \bullet})$
by  the K\"unneth formula. 

The natural map 
$$\Phi(A)^{\bullet}\ts \Phi(B)^{\bullet} \to A^{\dbullet}\ts B^{\dbullet}\,,$$
obtained as the tensor product of the maps $\Phi(A)^\bullet \to A^{\dbullet}$ and 
$\Phi(B)^\bullet \to B^{\dbullet}$, 
is a quasi-isomorphism on total complexes. 
Indeed, using (\Phicpx.a) it is identified with the inclusion 
$$\Phi(\Tot_{24}(A\dbullet\ts B\dbullet)\,)\to \Tot_{24}(A\dbullet\ts B\dbullet)\,,\eqno{(\Phicpx.b)}$$
which is a quasi-isomorphism.

This construction generalizes to  the case of a finite sequence of 
``$n$-tuple" complexes 
$$A_1^{\bullet\cdots \bullet}, A_2^{\bullet\,\cdots \bullet}, 
\cdots, A_c^{\bullet\,\cdots \bullet}\,,$$
 the above being the case $c=n=2$. 
 For example if $A\tbullet$ and $B\tbullet$ are ``triple" complexes
satisfying (i), (ii), then  $A\tbullet\ts B\tbullet$ is a ``6-tuple" complex.
$\Tot_{36}(A\tbullet\ts B\tbullet)$ is a ``5-tuple" complex satisfying (i), (ii) with respect
to the last degree. One  has 
$$\Phi(A)^{\dbullet}\ts \Phi (B)^{\dbullet}= \Phi (\Tot_{36}(A\tbullet\ts B\tbullet)\,)$$
as a ``quadruple" complex, and the 
natural map 
$$\Phi(A_1)^{\dbullet}\ts \Phi(B)^{\dbullet} \to A^{\tbullet}\ts B^{\tbullet}\,,$$
is a quasi-isomorphism on total complexes. 
For general $c$ and $n$, we get a ``$c(n-1)$-tuple" complex
$\Phi(A_1)\ts\cdots\ts\Phi(A_c)$. 
\bigskip

\section{Quasi DG categories.}

We refer to (\Multicpx) for multiple complexes and tensor product of complexes. 
In this section we will consider sequences of objects indexed by 
$[1, n]=\{1, \cdots, n\}$, or more generally by a finite (totally) ordered set $I$. 
For notions related to finite ordered sets see (\Finiteorderedset). 
For $n\ge 2$, let $(1, n)$ be the set $\{2, \cdots, n-1\}$;
if $I=[1, n]$, then one has $\ctop{I}=(1, n)$. 
\bigskip 

\sss{\DefDG}\quad 
A {\it DG category} $\cC$ is an additive category such that for a pair of objects 
$X$, $Y$  the group of homomorphisms 
$\Hom_\cC(X, Y)$ has the structure of a complex of abelian groups, 
written $F(X, Y)^\bullet$, and the composition of arrows
$$F(X, Y)^\bullet\ts F(Y, Z)^\bullet\to 
F(X, Z)^\bullet$$
which sends $u\ts v$ to $u\cdot  v$, is a map of complexes. 
Here to  $u:X\to Y$ and $v: Y\to Z$ there corresponds the product 
$u\cdot v: X\to Z$, which is the composition 
$v\scirc u$ in the usual notation.

A complex of abelian groups $A^\bullet$ is said to be (degree-wise) 
$\ZZ$-free if for each $p$ there is a set $\cS_A^p$ such that 
$A^p=\ZZ\cS_A^p$.  
In this section all complexes will be assumed $\ZZ$-free and bounded above. 

The following facts will be often used in this section. 
Let $A^\bullet$, $B^\bullet$ and $C^\bullet$ be $\ZZ$-free bounded-above 
complexes, $u: A^\bullet\to B^\bullet$ be  a map of complexes and 
$u\ts 1: A^\bullet\ts C^\bullet\to B^\bullet\ts C^\bullet$ be the induced map. 
If $u$ is injective, then $u\ts 1$ is injective. 
If $u$ is a quasi-isomorphism, then $u\ts 1$ is a quasi-isomorphism. 

We now give the definition of a weak quasi DG category. 
To understand the conditions, the reader may look at the example given 
after the definition. 
\bigskip 

\sss{\wqDG} \Def
 A {\it  weak quasi $DG$ category $\cC$}
consists of the following data (i)-(iii), satisfying the conditions 
(1)-(4). 

(i) A symmetric monoidal category $\cC_0$, which is a groupoid
(namely all morphisms are isomorphisms).
Thus we have the functor $(X, Y)\mapsto X\oplus Y$, the zero object 
$O$, and the commutativity constraint $\gamma_{XY}: X\oplus Y\isoto 
Y\oplus X$, the associativity constraint $\al_{XYZ}: X\oplus (Y\oplus Z)
\isoto (X\oplus Y)\oplus Z$, and the unit isomorphisms 
$\lambda_X: O\oplus X\isoto X$, $\rho_X: X\oplus O\isoto X$, 
that are subject to the axioms of symmetric monoidal category.
See \RefMacL, Chap. VII,  for symmetric monoidal categories.
The ``product" operation is written additively, and the unit object
is called the zero object.
In what follows, objects and morphisms will mean objects and 
morphisms in  $\cC_0$. 
\smallskip 

(ii) For each 
sequence of objects of $\cC_0$, 
$X_1, \cdots ,  X_n$ ($n\ge 2$), a ($\ZZ$-free, bounded above)
complex of abelian groups $F(X_1, \cdots , X_n)$. 
\smallskip 

 (iii) Two types of maps as follows. 
For $1<k<n$ a map of complexes
$$\tau_{k}(X_1, \cdots, X_n): F(X_1, \cdots , X_n)\to F(X_1, \cdots , X_k)\ts
 F(X_k, \cdots , X_n)\,\,,
$$
often just written $\tau_{X_k}$ or $\tau_k$;
the map $\tau_{k}$ is assumed to be a quasi-isomorphism. 
For $1<\ell< n$ a map of complexes 
$$\vphi_{\ell} (X_1, \cdots, X_n): F(X_1, \cdots , X_n)\to 
F(X_1, \cdots , \widehat{X_\ell}, \cdots , X_n)\,\,, 
$$
also written $\vphi_{X_\ell}$ or $\vphi_\ell$.
\bigskip 

These complexes and maps satisfy 
the conditions below. 
\smallskip

(1)(Functoriality.)\quad 
 The complex $F(X_1, \cdots, X_n)$ is assumed covariantly 
functorial. 
Namely, given a sequence of morphisms $f=(f_1, {\scriptstyle\cdots}, f_n): (X_1, {\scriptstyle\cdots}, X_n)
\to (Y_1, {\scriptstyle\cdots}, Y_n)$, where each $f_i$ is a morphism
in $\cC_0$,
there corresponds an isomorphism of complexes
$$f_*:  F(X_1, \cdots, X_n)\to F(Y_1, \cdots, Y_n)$$
that is  covariantly functorial in $f$. 
Also, if $X_i=O$ for some $i$, the complex $F(X_1, \cdots, X_n)$
is assumed to be zero.

We require  that $\tau_k$ is covariantly functorial:
For a sequence of morphisms
 $f: (X_1, \cdots, X_n)$
\newline 
$\to (Y_1, \cdots, Y_n)$ the following 
square commutes:
$$\begin{array}{ccc}
 F(X_1, \cdots, X_n)   &\mapr{\tau_k}  &F(X_1, \cdots , X_k)\ts
 F(X_k, \cdots , X_n)  \\
\mapd{f_*} & &\mapdr{f_*\ts f_*} \\
F(Y_1, \cdots, Y_n)   &\mapr{\tau_k}  &F(Y_1, \cdots , Y_k)\ts
 F(Y_k, \cdots , Y_n)  
 \end{array}
$$  
where the right vertical map $f_*\ts f_*$ is the tensor product of the maps 
$f_*: F(X_1, \cdots, X_k)\newline\to F(Y_1, \cdots, Y_k)$ and 
$f_*: F(X_k, \cdots, X_n)\to F(Y_k, \cdots, Y_n)$. 
In short, one has $(f_*\ts f_*)\tau_k= \tau_k f_*$.

Also $\vphi_\ell$ is assumed covariantly functorial in each 
$X_i$:
For a sequence of morphisms $f: (X_1, \cdots, X_n) \to (Y_1, \cdots, Y_n)$,
 the square
$$\begin{array}{ccc}
 F(X_1, \cdots , X_n)   &\mapr{\vphi_\ell}  &F(X_1, \cdots , \widehat{X_\ell}, \cdots , X_n)  \\
\mapd{f_*} & &\mapdr{ f_*} \\
F(Y_1, \cdots , Y_n)   &\mapr{\vphi_\ell}  &F(Y_1, \cdots , \widehat{Y_\ell}, \cdots , Y_n)  
 \end{array}
$$  
commutes, namely $f_*\vphi_\ell =\vphi_\ell f_*$. 
\smallskip 

(2)(Commutation identities.)\quad\,\, 
For two elements $k<\ell$ in $(1, n)$, we have the identity
\newline $(1\otimes \tau_{X_\ell}) \tau_{X_k}
=(\tau_{X_k}\otimes 1)\tau_{X_\ell}$,  
namely the following square commutes:
$$\begin{array}{ccc}
F(X_1, \cdots , X_n)&\mapr{\tau_{X_k}} &F(X_1, \cdots , X_k)\ts
 F(X_k, \cdots , X_n) \\ 
\mapd{\tau_{X_\ell}}& &\mapdr{1\ts\tau_{X_\ell}} \\
F(X_1, \cdots , X_\ell)\ts F(X_\ell, \cdots , X_n) &\mapr{\tau_{X_k}\ts 1}
&\!\!\!F(X_1, \cdots , X_k)\ts F(X_k, \cdots , X_\ell)\ts  F(X_\ell, \cdots , X_n)\,.
\end{array}
$$
Note the maps $1\ts\tau_{X_\ell}$ and $\tau_{X_k}\ts 1$ are quasi-isomorphisms. 
Writing $\tau_k$, $\vphi_k$ for 
$\tau_{X_k}$, $\vphi_{X_k}$, the identity reads
$(1\otimes \tau_{\ell})\tau_{k}=(\tau_{k}\otimes 1)\tau_{\ell}$.
Note $\tau_\ell$ is not necessarily the $\ell$-th $\tau$-map.

 For two elements $k<\ell$ in $(1, n)$,
$\vphi_{X_\ell}\vphi_{X_k}
=\vphi_{X_k}\vphi_{X_\ell}$, namely the following
commutes:
 $$\begin{array}{ccc}
F(X_1, \cdots , X_n)&\mapr{\vphi_{X_k}} &F(X_1, \cdots , 
\widehat{X_k}, \cdots , X_n) \\ 
\mapd{\vphi_{X_\ell}}& &\mapdr{\vphi_{X_\ell}} \\
F(X_1, \cdots , \widehat{X_\ell}, \cdots , X_n) &\mapr{\vphi_{X_k}}
&F(X_1, \cdots , \widehat{X_k}, \cdots  ,\widehat{X_\ell}, \cdots ,  X_n)\,\,.
\end{array}
$$

  For distinct elements $k$ and $\ell$ in $(1, n)$,
$\tau_{X_\ell}\vphi_{X_k}=(\vphi_{X_k}\otimes 1)\tau_{X_\ell}$ if 
$k<\ell$, and 
$\tau_{X_\ell}\vphi_{X_k}=(1\otimes \vphi_{X_k})\tau_{X_\ell}$ if 
$k>\ell$. The following diagram is for $k<\ell$. 
$$\begin{array}{ccc}
F(X_1, \cdots , X_n)&\mapr{\vphi_{X_k}} &F(X_1, \cdots , 
\widehat{X_k}, \cdots , X_n) \\ 
\mapd{\tau_{X_\ell}}& &\mapdr{\tau_{X_\ell}} \\
F(X_1, \cdots ,{X_\ell})\ts F(X_\ell,  \cdots , X_n) &\mapr{\vphi_{X_k}
\otimes 1}
&F(X_1, \cdots , \widehat{X_k}, \cdots  ,X_\ell)\ts F(X_\ell, \cdots , X_n)
\end{array}
$$
\smallskip

(3)(Additivity of $F(X_1, \cdots, X_n)$.)\quad
We will assume  that the complex $F(X_1, \cdots, X_n)$ is {\it additive} in each 
variable, in the sense formulated below. 
 
For each $i$, $1\le i\le n$, and a sequence of objects 
$(X_1, \cdots, X_{i-1}, Y_i, Z_i, X_{i+1}, \cdots, X_n)$, 
 we are given maps of complexes 
$$\lsum_i(Y_i, Z_i): 
F(X_1,\cdots, X_{i-1}, Y_i, X_{i+1}, \cdots, X_n)\to 
F(X_1,\cdots, X_{i-1}, Y_i\bop Z_i, X_{i+1},\cdots, X_n)$$
and 
$$\rsum_i(Y_i, Z_i): 
F(X_1,\cdots, Z_i, \cdots, X_n)\to 
F(X_1,\cdots, Y_i\bop Z_i, \cdots, X_n)\,.$$
For $1<i<n$, also given a map of complexes 
$$\pi_i(Y_i, Z_i): F(X_1, \cdots, X_{i-1}, Y_i)\ts F(Z_i,X_{i+1},  \cdots, X_n) 
\to F(X_1, \cdots, Y_i\bop Z_i, \cdots, X_n)\, .$$
The following conditions should be satisfied for 
the maps $s_i$, $t_i$ and $\pi_i$;  
all the conditions are natural and rather obvious ones, 
except possibly the requirements that 
$s_i$ (resp. $\pi_i$) be compatible with $\tau_i$ and $\vphi_i$, 
and that $\theta$ be a quasi-isomorphism.
\smallskip 

(ADD-1)
 Functoriality and commutativity.\quad 
The maps $\lsum_i$, $\rsum_i$ are functorial: 
A sequence of morphisms 
$$f: (X_1, \cdots,  X_{i-1}, Y_i, Z_i, X_{i+1},\cdots, X_n) \to 
(X'_1, \cdots,  X'_{i-1}, Y'_i, Z'_i, X'_{i+1},\cdots, X'_n)$$
induces a commutative square 
$$\begin{array}{ccc}
F(X_1,  \cdots, Y_i,  \cdots, X_n) 
&\mapr{\lsum_i (Y_i, Z_i)} 
&F(X_1, \cdots, Y_i\bop Z_i, \cdots, X_n)  \\
\mapd{f_*  }  &     &\mapdr{f_*}  \\
F(X'_1, \cdots, Y'_i,  \cdots, X'_n) &\mapr{\lsum_i (Y'_i, Z'_i)} 
&\phantom{\,.}F(X'_1, \cdots, Y'_i\bop Z'_i, \cdots, X'_n)\,.
\end{array}$$ 
The same for the map $\rsum_i(Y_i, Z_i)$. 

Similarly the map $\pi_i(Y_i, Z_i)$ is functorial in 
$(X_1, \cdots, Y_i, Z_i, \cdots, X_n)$. 

The maps $\lsum_i$ commute with each other:
For $i<j$, the following square commutes
(for $k\neq i, j$, the variable in the $k$-th spot is $X_k$):
$$\begin{array}{ccc}
F(X_1,  \cdots, Y_i,  \cdots, Y_j, \cdots,  X_n) 
&\mapr{\lsum_i (Y_i, Z_i)} 
&F(X_1, \cdots, Y_i\bop Z_i, \cdots, Y_j, \cdots,  X_n)  \\
\mapd{\lsum_j (Y_j, Z_j) }  &     &\mapdr{\lsum_j (Y_j, Z_j)}  \\
F(X_1, {\scriptstyle\cdots}, Y_i, {\scriptstyle\cdots}
, Y_j\bop Z_j,{\scriptstyle\cdots},  X_n) &\mapr{\lsum_i (Y_i, Z_i)} 
&\phantom{\,.}F(X_1, {\scriptstyle\cdots}, Y_i\bop Z_i, {\scriptstyle\cdots}
, Y_j\bop Z_j,{\scriptstyle\cdots},  X_n)\,.
\end{array}$$ 
Similarly,  $\rsum_i$ and $\rsum_j$ commute.

The maps $\pi_i$ commute with each other:
For $1<i<j<n$, the following
diagram commutes:
$$\begin{array}{ccc}
F(X_1,{\scriptstyle\cdots}, Y_i)\ts F(Z_i, {\scriptstyle\cdots}
, Y_j)\ts F(Z_j, {\scriptstyle\cdots}, X_n)
 &\mapr{\pi_i \ts 1} 
&F(X_1,{\scriptstyle\cdots}, Y_i\bop Z_i, {\scriptstyle\cdots}
, Y_j)\ts F(Z_j, {\scriptstyle\cdots}, X_n)
 \\
\mapd{ 1\ts \pi_j }   &     &\mapdr{\pi_j } \\
F(X_1,{\scriptstyle\cdots}, Y_i)\ts F(Z_i, {\scriptstyle\cdots}
, Y_j\bop Z_j, {\scriptstyle\cdots}, X_n)
&\mapr{\pi_i }    
&
F(X_1,{\scriptstyle\cdots}, Y_i\bop Z_i, {\scriptstyle\cdots}
, Y_j\bop Z_j, {\scriptstyle\cdots}, X_n)
\,.
\end{array}$$ 

For $i\neq j$, the maps $\lsum_i$ and $\pi_j$ commute, 
namely the following square commutes(if $i <j$):
$$\begin{array}{ccc}
F(X_1,  {\scriptstyle\cdots}, Y_i,  {\scriptstyle\cdots},  Y_j)\ts F(Z_j, {\scriptstyle\cdots}, X_n) 
&\mapr{\lsum_i (Y_i, Z_i)\ts 1} 
&F(X_1, {\scriptstyle\cdots}, Y_i\bop Z_i,  {\scriptstyle\cdots},  Y_j)\ts
 F(Z_j, {\scriptstyle\cdots}, X_n)  \\
\mapd{\pi_j}   &     &\mapdr{\pi_j}  \\
F(X_1,  {\scriptstyle\cdots}, Y_i,  {\scriptstyle\cdots},  Y_j\bop Z_j, {\scriptstyle\cdots}, X_n)
&\mapr{\lsum_i(Y_i, Z_i)}
&F(X_1, {\scriptstyle\cdots}, Y_i\bop Z_i,  {\scriptstyle\cdots},  Y_j\bop Z_j, {\scriptstyle\cdots}, X_n) \,.
\end{array}$$

(ADD-2)
Compatibility with the constraint maps.\quad
The map $\lsum_i$, $\rsum_i$ are compatible with the maps 
induced by the constraint maps, namely the following diagrams
all commute.
$$\begin{array}{ccc}
F(X_1,  \cdots, Y_i,  \cdots,  X_n) 
&\mapr{\lsum_i (Y_i, Z_i)} 
&F(X_1, \cdots, Y_i\bop Z_i,  \cdots,  X_n)  \\
\Vert   &     &\mapdr{\gamma(Y_i, Z_i)}  \\
F(X_1,  \cdots, Y_i,  \cdots,  X_n)
&\mapr{\rsum_i(Y_i, Z_i)}
&F(X_1, \cdots, Z_i\bop Y_i,  \cdots,  X_n)\,.
\end{array}$$ 
Here the right vertical map is the isomorphism $\gamma(Y_i, Z_i)$
induces.  Because of this, the properties (in this and subsequent subsections)
for $s_i$ will imply the analogous properites for $t_i$. 
$$\begin{array}{ccc}
F(X_1,  \cdots, Y_i,  \cdots,  X_n) 
&\mapr{\lsum_i(Y_i, Z_i\oplus W_i)} 
&F(X_1, \cdots, Y_i\oplus( Z_i\oplus W_i),  \cdots,  X_n)  \\
\mapd{\lsum_i(Y_i, Z_i)}   &     &\mapdr{\alpha(Y_i, Z_i, W_i)}  \\
F(X_1,  \cdots, Y_i\oplus Z_i,  \cdots,  X_n)
&\mapr{\lsum_i(Y_i\oplus Z_i, W_i)}
&F(X_1, \cdots, (Y_i\oplus Z_i)\oplus W_i ,  \cdots,  X_n)\,,
\end{array}$$ 
where the right vertical map is the isomorphism 
induced by $\alpha(Y_i, Z_i, W_i)$.
The map 
$$\lsum_i(Y_i,  O):
F(X_1,  \cdots, Y_i,  \cdots,  X_n)\to 
F(X_1,  \cdots, Y_i\oplus O,  \cdots,  X_n)$$
coincides with the isomorphism induced by $\rho: Y_i\to Y_i\oplus O$. 
This says that $\lsum_i(Y_i,  O)$ is the identity map, when 
$Y_i\oplus O$ is identified with $Y_i$. 
Similarly the map 
$$\rsum_i(O, Y_i):
F(X_1,  \cdots, Y_i,  \cdots,  X_n)\to 
F(X_1,  \cdots, O\oplus Y_i,  \cdots,  X_n)$$
coincides with the isomorphism
 induced by $\lambda: Y_i\to O\oplus Y_i$. 
\bigskip

(ADD-3)
 Compatibility with $\tau$.\quad
The maps $\lsum_i$ and $\tau_j$ commute. 
If $i=j$, it means the commutativity of the following 
square:
$$\begin{array}{ccc}
F(X_1,  \cdots, Y_i,  \cdots, X_n) 
&\mapr{\lsum_i (Y_i, Z_i)} 
&F(X_1, \cdots, Y_i\bop Z_i, \cdots, X_n)  \\
\mapd{\tau_i  }  &     &\mapdr{\tau_i}  \\
F(X_1, \cdots, Y_i)\ts F(Y_i,  \cdots, X_n)
 &\mapr{\lsum_i \ts \lsum_i} 
&F(X_1, \cdots, Y_i\bop Z_i)\ts F(Y_i\bop Z_i,  \cdots, X_n)\,, 
\end{array}$$ 
where the lower horizontal map is the tensor product of 
$\lsum_i:F(X_1, \cdots, Y_i)\to F(X_1, \cdots, Y_i\bop Z_i)$
and $\lsum_i: F(Y_i,  \cdots, X_n)\to 
F(Y_i\bop Z_i,  \cdots, X_n)$. 
For $i\neq j$ it means the commutativity of following square 
(if, say, $i<j$):
$$\begin{array}{ccc}
F(X_1,  \cdots, Y_i,  \cdots, X_n) 
&\mapr{\lsum_i (Y_i, Z_i)} 
&F(X_1, \cdots, Y_i\bop Z_i, \cdots, X_n)  \\
\mapd{\tau_j }  &     &\mapdr{\tau_j}  \\
F(X_1,{\scriptstyle\cdots} , Y_i,{\scriptstyle\cdots}, X_j) 
\ts F(X_j,  \cdots, X_n)
 &\mapr{\lsum_i\ts 1} 
&F(X_1,{\scriptstyle\cdots}, Y_i\bop Z_i,{\scriptstyle\cdots},X_j)\ts
F(X_j,  {\scriptstyle\cdots},  X_n)\,.
\end{array}$$ 
Similarly,  $\rsum_i$ and $\tau_j$ commute.

The maps $\pi_i$ and $\tau_i$ are compatible,  namely the following
diagram commutes:
$$\begin{array}{ccc}
F(X_1, \cdots, Y_i)\ts F(Z_i, \cdots, X_n) &\mapr{\pi_i (Y_i, Z_i)} 
&F(X_1, \cdots, Y_i\bop Z_i, \cdots, X_n)  \\
   \phantom{aaaaaaaaaaaaaaaaaaaaa}  {\scriptstyle \lsum_i\ts \rsum_i}\!\! \!\!\!\!\!\!\!\!&\searrow      &\mapdr{\tau_i}   \\
    &&F(X_1, \cdots, Y_i\bop Z_i)\ts F(Y_i\bop Z_i, \cdots, X_n)
\end{array}$$ 
where for example $\lsum_i$ is the map 
$F(X_1, \cdots, Y_i)\to F(X_1, \cdots, Y_i\bop Z_i)$. 
If $1<i<j<n$, the maps $\pi_i$ and $\tau_j$ commute, i.e., the following
diagram commutes:
$$\begin{array}{ccc}
F(X_1, \cdots, Y_i)\ts F(Z_i, \cdots, X_n) &\mapr{\pi_i (Y_i, Z_i)} 
&F(X_1, \cdots, Y_i\oplus Z_i, \cdots, X_n)  \\
\mapd{ 1\ts \tau_j }   &     &\mapdr{\tau_j } \\
F(X_1, {\scriptstyle\cdots}, Y_i)\ts F(Z_i, {\scriptstyle\cdots}
, X_j)\ts F(X_j,{\scriptstyle\cdots}, X_n)   &
\mapr{\pi_i (Y_i, Z_i)\ts 1}    &\!\!\!F(X_1, {\scriptstyle \cdots}, Y_i\oplus Z_i, {\scriptstyle \cdots}, X_j)\!\!\ts 
\!\!F(X_j, {\scriptstyle\cdots},  X_n)\,.
\end{array}$$ 
Similarly, if $1<j<i<n$, the maps $\pi_i$ and $\tau_j$ commute.

(ADD-4)
 Compatibility with $\vphi$.\quad
The maps $\lsum_i$ and $\vphi_j$ commute
(similarly, $\rsum_i$ and $\vphi_j$ commute).
Namely, if $i=j$, the following
diagram commutes:
$$\begin{array}{ccc}
F(X_1,  \cdots, Y_i,  \cdots, X_n)\!\!\!\!\!\!\!
 &\mapr{\lsum_i (Y_i, Z_i)} 
&F(X_1, \cdots, Y_i\bop Z_i, \cdots, X_n)  \\
   \phantom{aaaaaaaaaaaaaaaaaaaaaaaaa}  
{{\scriptstyle\vphi_i}}\!\!\!\!\!\!\!&\searrow     
 &\mapdr{\vphi_i}   \\
    &&F(X_1, \cdots, X_{i-1}, X_{i+1},  \cdots, X_n)\,.
\end{array}$$ 
If $i\neq j$, it means the commutativity of the following 
square (assume $i<j$):
$$\begin{array}{ccc}
F(X_1,  \cdots, Y_i,  \cdots, X_n) 
&\mapr{\lsum_i (Y_i, Z_i)} 
&F(X_1, \cdots, Y_i\bop Z_i, \cdots, X_n)  \\
\mapd{\vphi_j }  &     &\mapdr{\vphi_j}  \\
F(X_1, {\scriptstyle\cdots}, Y_i,{\scriptstyle\cdots},
\widehat{X_j},  {\scriptstyle\cdots},  X_n)
 &\mapr{\lsum_i(Y_i, Z_i)} 
&\phantom{\,.}F(X_1, {\scriptstyle\cdots}, Y_i\bop Z_i,{\scriptstyle\cdots},
\widehat{X_j},   {\scriptstyle\cdots},  X_n)\,.
\end{array}$$ 
The same for $\rsum_i$ and $\vphi_j$. 

The maps $\pi_i(Y_i, Z_i)$ and $\vphi_j$ are compatible.
It means that, if $i=j$, the composition of $\pi_i(Y_i, Z_i)$ and the map 
$$\vphi_i: F(X_1, \cdots, Y_i\bop Z_i, \cdots, X_n)\to 
F(X_1, \cdots, X_{i-1}, X_{i+1}, \cdots, X_n)$$
is zero. 
If $i\neq j$, 
the maps $\pi_i(Y_i, Z_i)$ and $\vphi_j$ commute, meaning 
the commutativity of 
the following digram (assume, say,  $i<j$):
$$\begin{array}{ccc}
F(X_1, \cdots, Y_i)\ts F(Z_i, \cdots, X_n) &\mapr{\pi_i (Y_i, Z_i)} 
&F(X_1, \cdots, Y_i\bop Z_i, \cdots, X_n)  \\
\mapd{ 1\ts \vphi_j  }  &     &\mapdr{\vphi_j}  \\
F(X_1, \cdots, Y_i)\ts F(Z_i, \cdots, \widehat{X_j}, \cdots, X_n)   &
\mapr{\pi_i (Y_i, Z_i) }    &\phantom{\,.}
F(X_1, {\scriptstyle\cdots}, Y_i\bop Z_i,{\scriptstyle\cdots}, 
\widehat{X_j},{\scriptstyle\cdots},  X_n)\,.
\end{array}$$ 

(ADD-5)
  The  map  of additivity $\theta$.\quad 
If $1<i<n$,  we define the {\it map of additivity}  
\begin{eqnarray*}
&\theta_i(Y_i, Z_i):& F(X_1, \cdots, Y_i, \cdots, X_n) \bop
 F(X_1, \cdots, Z_i, \cdots, X_n) \\
&& \bop F(X_1, \cdots, Y_i)\ts F(Z_i, \cdots, X_n) \\
&&\bop
  F(X_1, \cdots, Z_i)\ts F(Y_i, \cdots, X_n) \\
&\mapr{}&   F(X_1, \cdots, Y_i\bop Z_i, \cdots, X_n) \qquad\qquad\qquad\qquad\qquad 
 \qquad \qquad  \mbox{(\wqDG.a)}
\end{eqnarray*}
as the sum of the maps $\lsum_i(Y_i, Z_i)$, $\rsum_i(Y_i, Z_i)$, $\pi_i(Y_i, Z_i)$, and 
$\pi_i(Z_i, Y_i)$. 
(The last map is, to be precise, the composition of $\pi_i(Z_i, Y_i)$ with the isomorphism
induced by $\gamma(Z_i, Y_i)$. )
If $i=1$ or $n$, let 
$$\theta_i(Y_i, Z_i): F(X_1, \cdots, Y_i, \cdots, X_n)\bop F(X_1, \cdots, Z_i, \cdots, X_n)
\to F(X_1, \cdots, Y_i\bop Z_i, \cdots, X_n)$$
be the sum of $\lsum_i$ and $\rsum_i$.   
In either case we {\it require} that $\theta$ be a quasi-isomorphism. 
We often refer to  the last two terms in the source of $\theta$
in (\wqDG.a) as the {\it cross terms}. 
\smallskip

{\it Remarks to (3).}\quad 
$\bullet$\quad 
Note if $F(X, Y)$ is a complex,  additive in each variable, then the tensor product 
$F(X_1, X_2)\ts F(X_2, X_3)\ts \cdots\ts F(X_{n-1}, X_n)$ is additive in each 
variable in the above sense.   
So additivity here means ``quadratic additivity", so to speak. 

$\bullet$\quad 
 It follows that the map $\theta_i(Y_i, Z_i)$ is compatible with 
$\tau_j$ and with $\vphi_j$.  
For example, the compatibility of $\theta_i$ and $\tau_i$ means the 
commutativity of the following diagram 

\vspace*{0.5cm}
\hspace*{0.5cm}
\unitlength 0.1in
\begin{picture}( 41.3700, 15.6000)(  0.0000,-18.8000)
\put(0.2000,-4.5000){\makebox(0,0)[lb]{${\textstyle F(X_1, \cdots, Y_i,\cdots, X_n)}$}}%
\put(0.0000,-6.2000){\makebox(0,0)[lb]{${\textstyle \oplus F(X_1, \cdots, Z_i,\cdots, X_n)}$}}%
\put(0.0000,-7.9000){\makebox(0,0)[lb]{${\textstyle \oplus F(X_1, \cdots, Y_i)\ts F(Z_i,\cdots, X_n)}$}}%
\put(0.0000,-9.7000){\makebox(0,0)[lb]{${\textstyle \oplus F(X_1, \cdots, Z_i)\ts F(Y_i,\cdots, X_n)}$}}%
\put(0.3000,-14.4000){\makebox(0,0)[lb]{${\textstyle F(X_1, \cdots, Y_i)\ts F(Y_i,\cdots, X_n)}$}}%
\put(0.0000,-16.3000){\makebox(0,0)[lb]{${\textstyle \oplus F(X_1, \cdots, Z_i)\ts F(Z_i,\cdots, X_n)}$}}%
\put(0.0000,-18.2000){\makebox(0,0)[lb]{${\textstyle \oplus F(X_1, \cdots, Y_i)\ts F(Z_i,\cdots, X_n)}$}}%
\put(0.0000,-20.1000){\makebox(0,0)[lb]{${\textstyle \oplus F(X_1, \cdots, Z_i)\ts F(Y_i,\cdots, X_n)}$}}%
%
{\color[named]{Black}{%
\special{pn 8}%
\special{pa 790 970}%
\special{pa 790 1290}%
\special{fp}%
\special{sh 1}%
\special{pa 790 1290}%
\special{pa 810 1224}%
\special{pa 790 1238}%
\special{pa 770 1224}%
\special{pa 790 1290}%
\special{fp}%
\special{pa 790 1290}%
\special{pa 790 1290}%
\special{fp}%
}}%
%
{\color[named]{Black}{%
\special{pn 8}%
\special{pa 2480 590}%
\special{pa 3070 590}%
\special{fp}%
\special{sh 1}%
\special{pa 3070 590}%
\special{pa 3004 570}%
\special{pa 3018 590}%
\special{pa 3004 610}%
\special{pa 3070 590}%
\special{fp}%
}}%
%
{\color[named]{Black}{%
\special{pn 8}%
\special{pa 2470 1570}%
\special{pa 3060 1570}%
\special{fp}%
\special{sh 1}%
\special{pa 3060 1570}%
\special{pa 2994 1550}%
\special{pa 3008 1570}%
\special{pa 2994 1590}%
\special{pa 3060 1570}%
\special{fp}%
}}%
\put(33.4000,-6.6000){\makebox(0,0)[lb]{${\textstyle F(X_1, \cdots, Y_i\oplus Z_i,\cdots, X_n)}$}}%
%
{\color[named]{Black}{%
\special{pn 8}%
\special{pa 3998 940}%
\special{pa 3998 1290}%
\special{fp}%
\special{sh 1}%
\special{pa 3998 1290}%
\special{pa 4018 1224}%
\special{pa 3998 1238}%
\special{pa 3978 1224}%
\special{pa 3998 1290}%
\special{fp}%
\special{pa 3998 1290}%
\special{pa 3998 1290}%
\special{fp}%
}}%
\put(41.3700,-11.6000){\makebox(0,0)[lb]{$\tau_i$}}%
\put(31.5000,-16.2000){\makebox(0,0)[lb]{${\textstyle  F(X_1, \cdots, Y_i\oplus Z_i)\ts F(Y_i\oplus Z_i,\cdots, X_n)}$}}%
\put(26.1000,-5.4000){\makebox(0,0)[lb]{$\theta_i$}}%
\end{picture}%
\vspace{0.8cm}

\noindent 
where the left vertical map is the diagonal sum of $\tau_i$, $\tau_i$, $id$, and $id$, 
and the lower horizontal map is the diagonal sum of $s\ts s$, 
$t\ts t$, 
$s\ts t$ and $t\ts s$.  
The compatibility of $\theta_i$ and $\vphi_i$
means that the composition of $\theta_i$ and 
$\vphi_i: F(X_1, \cdots, Y_i\bop Z_i, \cdots, X_n)
\to F(X_1, \cdots, X_{i-1}, X_{i+1},  \cdots, X_n)$
coincides with  the sum of $\vphi_i$ on 
$F(X_1, \cdots, Y_i, \cdots, X_n)$,  $\vphi_i$ on 
$F(X_1, \cdots, Z_i, \cdots, X_n)$, and the zero maps
on the cross terms.
%
\bigskip

(4)(Existence of the identity in the ring $H^0F(X, X)$.)\quad   
This condition will be stated in (\HowqDG). 
\bigskip 

For a finite ordered set $I$ and a collection of objects $(X_i)_{i\in I}$ 
indexed by $I$, one can define the complex $F( (X_i)_{i\in I}\, )$, also denoted 
$F(X; I)$, $F(X^I)$ or $F(I)$ for short. Then the conditions (1)-(4) can be stated more naturally. 
\bigskip

\sss{\DGwqDG}
{\bf Example.}\quad A DG category $\cC$, 
in which the complexes $F(X, Y)^\bullet$ are $\ZZ$-free and 
bounded above, can be viewed as a weak quasi DG category.

We have the category $Z^0\cC$, which has the same objects as $\cC$, and 
where the homomorphism
groups are given by 
$$\Hom(X, Y)=Z^0F(X, Y)^\bullet\,,$$
namely the cocycles of degree 0 in the complex of homomorphisms $F(X, Y)^\bullet$.
It is an additive category, so in particular it has the structure
of a symmetric monoidal category with the direct sum and the zero 
object. 

We take as $\cC_0$ the subcategory of $Z^0\cC$, with the same objects as $\cC$
and with morphisms the invertible ones in  $Z^0\cC$. 

For objects $X, X'$, take $F(X, X')$ to be the complex $F(X, X')^\bullet$, and 
for a sequence of objects $X_1, \cdots, X_n$, let 
$$F(X_1, \cdots, X_n)=F(X_1, X_2)\ts F(X_2, X_3)\ts\cdots \ts F(X_{n-1}, X_n)\,.$$
The complex $F(X, X')$ is contravariantly  functorial in $X$ 
(resp. covariantly functorial in $X'$)
on the category $Z^0\cC$:  if $f: Y\to X$ and $f': X'\to Y'$ are morphisms
in  $Z^0\cC$, 
we have the induced map of complexes 
$F(X, X')\to F(Y, Y')$ given by $u\mapsto f\cdot u\cdot f'$. 
So $F(X, X')$ is covariantly functorial 
in both variables on $\cC_0$ 
by the map $u\mapsto f^{-1}\cdot u\cdot f'$. 
  It thus follows that $F(X_1, \cdots, X_n)$
is  functorial on $\cC_0$. 

Let the map $\tau_k$ be the identity, and $\vphi_\ell$ be the composition 
in $X_\ell$.  
One immediately verifies  that the conditions for a weak quasi DG category are 
satisfied (the map $\theta$ is the identity). 
\bigskip

\sss{\FtbarS} Let $\cC$ be a weak quasi DG category. Assume given a sequence of
objects $X_1, \cdots, X_n$. 
 For a subset $I=\{\ell_1, \cdots , \ell_a\}\subset [1, n]$, write 
$F(I)$ in place of
 $F(X_{\ell_1}, \cdots , X_{\ell_a})$ for abbreviation. 
Set $(1, n)=\{2, \cdots, n-1\}$. 
For a subset $S=\{i_1, \cdots , i_a\}\subset (1, n)$, let 
$$F(X_1, \cdots, X_n\tbar S):= F(X_1, \cdots , X_{i_1})\ts 
F(X_{i_1}, \cdots , X_{i_2})
\ts \cdots \ts F(X_{i_a}, \cdots , X_n)\,\,$$
be the tensor product complex; in other words, if $I_1, \cdots, I_a$
is the segmentation corresponding to $S$, 
$F(X_1, \cdots, X_n\tbar S)=F(I_1)\ts \cdots\ts F(I_a)$. 
It is an $a$-tuple complex where the ordered set of gradings  $[1, a]$
is identified with the set $\{I_1, \cdots, I_a\}$. 
Note $F(X_1, \cdots, X_n\tbar\emptyset)=F(X_1, \cdots, X_n)$. 
We also write $F([1, n]\tbar S)$ for $F(X_1, \cdots, X_n\tbar S)$. 

More generally, for a finite ordered set $I$, a sequence of 
objects indexed by $I$, $(X_i)_{i\in I}$, and a subset $S$ of 
$\ctop{I}=I-\{\init(I), \term(M)\}$, 
one defines the complex
$F(I\tbar S)$ in a similar manner. 

For subsets $S\subset S'$ of $\ctop{I}$ we will 
define a map of $a$-tuple  complexes
$$\tau_{S\, S'}: F(I\tbar S)\to  \Tot^f (F(I\tbar S')\,)\,\,.
$$
Here, if $I'_1, \cdots, I'_{a'}$ is the segmentation of $[1, n]$ by $S'$, there is a 
map $f: [1, a']\to [1, a]$ such that $I'_{f(i)}\subset I_i$ for each $i$, thus 
one  has an $a$-tuple complex $\Tot^f (F(I\tbar S')\,)$. 
With this understood, we will just write $\tau_{S\, S'}: F(I\tbar S)\to F(I\tbar S')$
more often than not.

 Let 
$$\tau_S:
 F(I)\to F(I\tbar S)\quad\mbox{or}\quad F(I)\to \Tot F(I\tbar S) $$
 be 
the composition of $\tau_{X_k}$'s for $k\in S$
(if $S=\emptyset$, $\tau_S=id$). 
Define for $S\subset S'$ the map $\tau_{S\, S'}$ as follows. 
If $S=\emptyset$, let $\tau_{\emptyset\, S'}=\tau_{S'}$.  In general 
let $I_1, \cdots, I_a$ be the segmentation of $I$ corresponding to 
$S$, $S'_i=S'\cap \ctop{I_i}$, and 
$\tau_{S'_i}: F(I_i)\to F(I_i\tbar S'_i)$ be the map just defined. 
Then $$\tau_{SS'}:=\bts_i\tau_{ S'_i}: \bts_i F(I_i)\to \bts_i 
F(I_i\tbar S'_i\,)=F(I\tbar S)\,\,.$$
Note that  $\tau_{S, S}=id$. 

 For $K=\{k_1, \cdots , k_b\}\subset (1, n)$ disjoint from $S$, we define a map
$$\vphi_K: F(X_1, \cdots, X_n\tbar S)\to F(X_1, \cdots,\widehat{X_{k_1}}, \cdots , \widehat{X_{k_b}}, \cdots ,  X_n\tbar S)\,\,. $$
More generally for $S\subset \ctop{I}$ and $K\subset \ctop{I}$ disjoint 
from $S$, we have  $\vphi_K: F(I \tbar S)\to F(I-K\tbar S)$.
 If $S=\emptyset$, $\vphi_K$ is the composition of $\vphi_{k}$
for $k\in K$; in general, let $I_1, \cdots, I_a$ be the segmentation of $[1,n]$ corresponding to 
$S$, and 
$$\vphi_K:= \bts_{i}\vphi_{K\cap I_i}:\bts_i F(I_i) \to \bts_i F(I_i-K )\,\,.$$
Note $\vphi_K=id$ if $K=\emptyset$.

The complexes $F(X_1, \cdots, X_n\tbar S)$ and the above maps 
$\tau$, $\vphi$ satisfy the following properties, extending 
(1)-(3) of the previous subsection.
\smallskip

(1)(Functoriality.)\quad 
A sequence of morphisms $f=(f_1, \cdots, f_n): (X_1, \cdots, X_n)$
\newline
$\to (Y_1, \cdots, Y_n)$ induces an isomorphism of complexes
$$f_*:  F(X_1, \cdots, X_n\tbar S)\to F(Y_1, \cdots, Y_n\tbar S)$$
given by the formula 
$f_*(u_1\ts\cdots\ts u_a)=(f_*u_1)\ts\cdots\ts(f_*u_a)$ 
for $u_i\in F(I_i)$. 
The $f_*$ is  covariantly functorial in $f$. 

The map $\tau_{S, S'}$ is covariantly functorial for morphisms.
Also $\vphi_K$ is  covariantly functorial.  

(2)(Commutation identities.)\quad\,\, 
 $\tau_{S\, S'}$ is a quasi-isomorphism (as a tensor product of quasi-isomorphisms). 
For $S\subset S'\subset S''$, $ \tau_{S' S''}\tau_{S S'}=
\tau_{S S''}$.

If $K$ is the disjoint union of $K'$ and $K''$, 
$\vphi_K=\vphi_{K''}\vphi_{K'}$.  

If $K$ and $S'$ are disjoint the following commutes:
$$\begin{array}{ccc}
F(I\tbar S)&\mapr{\vphi_K} 
&F(I-K \tbar S) \\
\mapd{\tau_{SS'}}& &\mapdr{\tau_{SS'}} \\
F(I\tbar S')&\mapr{\vphi_K} &F(I-K \tbar S')\,\,.
\end{array}
$$
\smallskip 

(3)(Additivity of $F(X_1, \cdots, X_n\tbar S)$.)\quad
Assume that a variable $X_i$ is the direct sum of two objects,  
 $X_i=Y_i\bop Z_i$.
For  $i\not\in S$, let $S_1$, $S_2$ be the partition of $S$ by $i$, 
namely $S_1=S\cap (1, i)$ and $S_2=S\cap (i, n)$. 
Then we have the map 
$$\lsum_i(Y_i, Z_i): 
F(X_1,\cdots, X_{i-1}, Y_i, X_{i+1}, \cdots, X_n\tbar S)\to 
F(X_1,\cdots, X_{i-1}, Y_i\bop Z_i, X_{i+1},\cdots, X_n\tbar S)$$
defined as 
$$\begin{array}{rl}
1\ts \lsum_i(Y_i, Z_i)\ts 1:&
F(X_1,{\scriptstyle \cdots}, X_a\tbar S_1-\{a\})
\ts F(X_a,{\scriptstyle \cdots}, Y_i, {\scriptstyle \cdots}, X_b)\ts F(X_b, 
{\scriptstyle \cdots}, X_n\tbar S_2-\{b\})\\
\to &F(X_1,\cdots, X_{i-1}, Y_i\bop Z_i, X_{i+1},\cdots, X_n\tbar S)
\end{array}
$$
where $a$ is the largest element of $S_1$,  
$b$ is the smallest element of $S_2$, and $s_i(Y_i, Z_i)$ is the map in 
(\wqDG), (3). 
Similarly one has $\rsum_i(Y_i, Z_i): 
F(X_1,\cdots, Z_i, \cdots, X_n\tbar S)\to 
F(X_1,\cdots, Y_i\bop Z_i, \cdots, X_n\tbar S)$. 
If $1<i<n$,  we  define a map of complexes 
$$\pi_i(Y_i, Z_i): F(X_1,{\scriptstyle \cdots}, X_{i-1}, Y_i\tbar S_1)\ts F(Z_i,X_{i+1},
{\scriptstyle \cdots}
, X_n\tbar S_2) 
\to F(X_1,{\scriptstyle \cdots}, Y_i\bop Z_i,{\scriptstyle \cdots}, X_n\tbar S)\, $$
as the map 
$$\begin{array}{rl}
{1\ts \pi_i\ts 1}: & F(X_1,{\scriptstyle \cdots}, X_a\tbar S_1-\{a\})
\ts F(X_a,{\scriptstyle \cdots}, Y_i)\ts F(Z_i,{\scriptstyle \cdots}, X_b)\ts F(X_b, 
{\scriptstyle \cdots}, X_n\tbar S_2-\{b\})
\\
\to& F(X_1, \cdots, Y_i\bop Z_i, \cdots, X_n\tbar S)
\end{array}$$
where $a$, $b$ are as before,
and $\pi_i=\pi_i(Y_i, Z_i):F(X_a,{\scriptstyle \cdots}
, Y_i)\ts F(Z_i,{\scriptstyle \cdots}, X_b)\to F(X_a,{\scriptstyle \cdots}, Y_i\bop 
Z_i,, X_b)$ is the map given in (\wqDG), (3).
 
These maps  satisfy properties parallel to those for $\lsum_i$, $\rsum_i$, and 
$\pi_i(Y_i, Z_i)$ in 
(\wqDG), (3). 
In particular, if we define the map (in case $1<i<n$)
\begin{eqnarray*}
&\theta_i(Y_i, Z_i):& F(X_1, \cdots, Y_i, \cdots, X_n\tbar S) \bop
 F(X_1, \cdots, Z_i, \cdots, X_n\tbar S) \\
&& \bop F(X_1, \cdots, Y_i\tbar S_1)\ts F(Z_i, \cdots, X_n\tbar S_2) \\
&&\bop
  F(X_1, \cdots, Z_i\tbar S_1)\ts F(Y_i, \cdots, X_n\tbar S_2) \\
&\mapr{}&   
F(X_1, \cdots, Y_i\bop Z_i, \cdots, X_n\tbar S)\,,
\end{eqnarray*}
as the sum of the maps $\lsum_i$, $\rsum_i$, $\pi_i(Y_i, Z_i)$, and 
$\pi_i(Z_i, Y_i)$, then  it is a quasi-isomorphism. 
If $i=1$ or $n$, the map 
$$\theta_i(Y_i, Z_i): F(X_1,{\scriptstyle \cdots}, Y_i,{\scriptstyle \cdots}
, X_n\tbar S)\bop F(X_1,{\scriptstyle \cdots}, Z_i,{\scriptstyle \cdots}, X_n\tbar S)
\to F(X_1,{\scriptstyle \cdots}, Y_i\bop Z_i,{\scriptstyle \cdots}, X_n\tbar S)$$
defined as the sum of  $\lsum_i$ and $\rsum_s$, is a quasi-isomorphism. 
\bigskip

\sss{\HowqDG} {\it Homotopy category.}\quad 
A weak quasi DG category $\cC$ is not a category in the usual sense, since the 
composition is not defined. 
Nevertheless, one has composition in a weak sense.

For three objects $X$, $Y$ and $Z$, let 
$$\psi_Y: F(X, Y)\ts F(Y, Z)\to F(X, Z)$$
be the map in the derived category defined as the composition 
$\vphi_Y\scirc (\tau_Y)^{-1}$ where the maps are as in 
$$F(X, Y)\ts F(Y, Z)\mapl{\tau_Y} F(X, Y, Z)\mapr{\vphi_Y} F(X, Z)\,\,.$$
The map $\psi_Y$  is verified to be 
associative, namely the following commutes in the derived category:
$$\CD
F(X, Y)\ts F(Y, Z)\ts F(Z, W) &\mapr{\psi_Y\ts id} &F(X, Z)\ts F(Z, W) \\
\mapd{id\ts \psi_Z}  &    &\mapdr{\psi_Z} \\
F(X, Y)\ts F(Y, W)   &\mapr{\psi_Y}  &F(X, W)\,\,.
\endCD$$
To prove this identity, $\psi_Z(\psi_Y\ts id)=\psi_Y(id\ts \psi_Z)$, 
compose the quasi-isomorphism $(\tau_Y\ts id)\tau_Z=(id\ts \tau_Z)\tau_Y$ from right.
Using the commutation identities one has
\begin{eqnarray*}
\psi_Z(\psi_Y\ts id)(\tau_Y\ts id)\tau_Z&=&\psi_Z(\vphi_Y\ts id)\tau_Z \\
&=&\psi_Z\tau_Z\vphi_Y \\
&=&\vphi_Z\vphi_Y\,;
\end{eqnarray*}
similarly, $\psi_Y(id\ts \psi_Z)(id\ts \tau_Z)\tau_Y=\vphi_Y\vphi_Z$. 
Since $\vphi_Z\vphi_Y=\vphi_Y\vphi_Z$, the assertion follows. 

Let $H^0F(X, Y)$ be the $0$-th cohomology of $F(X,Y)$. 
$\psi_Y$ induces a map 
$$\psi_Y: 
H^0F(X, Y)\ts H^0F(Y, Z)\to H^0(\, F(X, Y)\ts F(Y, Z)\, )\mapr{H^0\psi_Y}
H^0F(X, Z)\,\,,$$
which is associative.
In particular $H^0F(X, X)$ is a ring. 
We often write $u\cdot v$ for $\psi_Y(u\ts v)$. 

The last condition (4) for a weak quasi DG category is:
\smallskip 

(4) For each $X$ there is an element $1_X\in H^0F(X, X)$ such that 
$1_X\cdot u =u$ for any $u\in H^0F(X, Y)$ 
and $u\cdot 1_X=u$ for $u\in H^0F(Y, X)$. 
(Call such $1_X$ the {\it identity}.) 
\smallskip  

Thus one can 
associate to $\cC$ an additive category, 
the associated {\it homotopy category}, denoted by 
$Ho(\cC)$. Objects of $Ho(\cC)$ are the same as the objects of $\cC$,
and $\Hom (X, Y):= H^0F(X, Y)$. Composition of arrows is the map 
induced from $\psi_Y$. The object $O$ is the zero object, and 
the direct sum $X\oplus Y$ is the direct sum in the  categorical 
sense. $1_X$ gives the identity $X\to X$. 
\bigskip 

{\it Remark.}\quad  More generally we  have maps 
$\psi_Y: H^mF(X, Y)\ts H^nF(Y, Z)\to 
H^{m+n}F(X, Z)$
for $m, n\in \ZZ$, defined in a similar manner. 
The groups $H^mF(X, Y)$ and the composition maps for them play roles 
when we consider quasi DG categories in the next subsection.
\bigskip 

\sss{\DefqDG} {\bf Definition.}\quad  A {\em   quasi $DG$ category $\cC$\/}
is a weak quasi DG category satisfying further conditions. 
Since it is a weak quasi DG category, it consists of data (i) - (iii) 
in (\wqDG), 
\smallskip 

(i) A symmetric monoidal groupoid $\cC_0$, 

(ii) Complexes $F(X_1, \cdots , X_n)$ for 
sequence of objects, 

(iii) Maps of complexes $\tau_k$ and $\vphi_\ell$, 
\smallskip 

\noindent satisfying the conditions (1)-(4) there:
\smallskip 

(1) Functoriality, 

(2) Commutation identities,

(3) Additivity, 

(4) Existence of the identity. 
\smallskip 

As discussed in (\HowqDG), we have the complexes 
$F(X_1, \cdots, X_n\tbar S)$ and maps $\tau_{S\, S'}$,  
$\vphi_K$ between them. 
For a quasi DG category, we impose
two more conditions (5) and (6) below. 
When necessary we will also impose additional data (iv)-(v), 
satisfying (7) and (8).
\smallskip 

(5)(Multiple complexes $F(X_1, \cdots , X_n|S)$.)\quad
For each sequence of objects $X_1, \cdots X_n$ ($n\ge 2$), 
and a subset $S$ of $(1, n)$, we are given a ($\ZZ$-free, bounded above) $a$-tuple 
complex of abelian groups $F(X_1, \cdots , X_n|S)$,
where $a=|S|+1$,  
a {\it surjective\/} quasi-isomorphism of complexes 
$$\sigma_S:F(X_1, \cdots, X_n)\to  \Tot F(X_1, \cdots , X_n|S)\,,$$
and an {\it injective\/} quasi-isomorphism of $a$-tuple complexes
$$\iota_S:  F(X_1, \cdots , X_n|S)\to F(X_1, \cdots , X_n\tbar S)$$
such that the map $\tau_S: F(X_1, \cdots, X_n)\to \Tot F(X_1, \cdots , X_n\tbar S)$
factors as 
$$F(X_1, \cdots, X_n)\mapr{\sigma_S}  \Tot F(X_1, \cdots , X_n|S)
\mapr{\iota_S}\Tot F(X_1, \cdots , X_n\tbar S)\,.\eqno{(\DefqDG.a)}$$
Note if $S=\emptyset$, one has $F(X_1, \cdots , X_n|\emptyset)
=F(X_1, \cdots, X_n)$. 

As in (\wqDG), for a finite ordered set $I$, a sequence of objects 
$(X_i)_{i\in I}$, and a subset $S\subset \ctop{I}$, 
one has a complex $F(X; I|S)$, also written 
$F(I| S)$ or $F(X|S)$.  
If $I_1, \cdots, I_c$ is the segmentation of $I$ by $S$, we also use 
the notation $F(I_1)\hts
\cdots\hts F(I_c)$  for $F(I|S)$. 
We will identify $F(I|S)$ with a subcomplex of $F(I\tbar S)$, and write
an element of $F(I|S)$ as a sum of $u_1\ts\cdots\ts u_c$ with $u_i\in F(I_i)$. 
\smallskip 

By means of the factorization (\DefqDG.a),  we have induced maps 
$\sigma_{S\, S'}$, $\vphi_K$ and $\iota_{S/T}$, with the following 
properties. 
For $S\subset S'$ one has a (unique) surjective quasi-isomorphism of 
multiple complexes
$$\sigma_{S\, S'}: F(X_1, \cdots , X_n| S)\to 
F(X_1, \cdots , X_n | S')\,\,,$$
or $F(I|S)\to F(I|S')$ for short,
such that 
the following diagram 
$$\begin{array}{ccc}
F(I| S)&\mapr{\iota_{S}}  &F(I\tbar S) \\
\mapd{\sigma_{S\,S'}}& &\mapdr{\tau_{S\,S'}} \\
F(I|S')&\mapr{\iota_{S'}} &F(I \tbar S')\,
\end{array}
$$
commutes (namely $\sigma$ and $\tau$ are compatible via the maps 
$\iota_S$). 
To be more precise, if $a=|S|+1$, the source
of $\sigma_{S\, S'}$ is an $a$-tuple complex while the 
target is an $(|S'|+1)$-tuple complex; the latter can be viewed as an $a$-tuple 
complex as in  (\FtbarS), and $\sigma_{S\, S'}$ is a map of $a$-tuple 
complexes. 
To show how the map $\sigma_{S\, S'}$ is obtained, let $c: [1, a]\to \{1\}$
and $c'=c\scirc f: [1, a']\to \{1\}$ be the maps to a one point set. 
The commutative diagram of complexes 
$$\begin{array}{ccc}
F(I)&\mapr{\tau_S}  &\Tot^c F(I\tbar S) \\
{\Vert}& &\mapdr{ \tau_{S\, S'} } \\
F(I)&\mapr{\tau_{S'}}  &\Tot^c\Tot^f F(I\tbar S') 
\end{array}
$$
induces another commutative diagram of complexes
$$\begin{array}{ccc}
\Tot^c F(I| S)&\mapr{\iota_S}  &\Tot^c F(I\tbar S) \\
\mapd{\sigma_{S\, S'}}& &\mapdr{ \tau_{S\, S'} } \\
\Tot^c\Tot^f F(I|S')&\mapr{\iota_{S'}}  &\Tot^c\Tot^f F(I\tbar S') 
\end{array}
$$
with a map of complexes $\sigma_{S\, S'}$. 
Since $\tau_{S\, S'}$ comes from a map of multiple complexes, and 
since $\iota_S$ and $\iota_{S'}$ are injections, one sees that 
$\sigma_{S\, S'}$ also comes from a map of multiple complexes
$\sigma_{S\, S'}: F(I| S)\to \Tot^f F(I|S')$. 
This type of reasoning will repeatedly appear in this subsection. 

Similarly one has a map of multiple complexes $\vphi_K: F(I|S)\to 
F(I-K|S)$ 
that is  compatible with the map $\vphi_K:F(I\tbar S)\to 
F(I-K\tbar S)$ via the maps $\iota_S$, namely $\vphi_K$ makes the following 
square commute:
$$\begin{array}{ccc}
F(I| S)&\mapr{\iota_S}  &F(I\tbar S) \\
\mapd{\vphi_K}& &\mapdr{ \vphi_{K} } \\
F(I-K| S)&\mapr{\iota_S}  &F(I-K \tbar S)\,\,. 
\end{array}
$$

For a subset $T\subset 
S$, if $I_1, \cdots, I_c$  is the segmentation corresponding to 
$T$, and $S_i=S\cap 
\ctop{I_i}$, there is an inclusion of multiple complexes 
$$\iota_{S/T}: F(I|S)\injto F(I_1|S_1)\ts \cdots\ts F(I_c|S_c)$$
such that the composition of $\iota_{S/T}$ and the tensor product of 
the inclusions $\iota_{S_i}:
F(I_i |S_i)\injto F(I_i\tbar S_i)$,
$$F(I_1|S_1)\ts \cdots\ts F(I_c|S_c)\injto 
F(I_1\tbar S_1)\ts \cdots\ts F(I_c\tbar S_c)=F(I\tbar S)$$
coincides with $\iota_S$. 
Notice that when $T=\emptyset$ or $T=S$, one has  $\iota_{S/\emptyset}=id$
or  $\iota_{S/S}=\iota_S$ respectively. 

The above compatibility of $\iota_S$ with $\sigma$ (resp. $\vphi$) 
generalizes to the compatibility of $\iota_{S/T}$ with $\sigma$ (resp. $\vphi$). 
The map $\sigma_{S\, S'}$ is compatible 
with the inclusion $\iota_{S/T}$:
If  $T\subset S\subset S'$, let $I_1, \cdots, I_c$ be the segmentation of 
$I$ by $T$, $S_i=S\cap \ctop{I_i}$ and $S'_i=S'\cap \ctop{I_i}$. 
Then the following commutes:
$$\begin{array}{ccc}
F(I| S)&\mapr{\iota_{S/T}}  &F(I_1| S_1)\ts\cdots\ts F(I_1| S_1) \\
\mapd{\sigma_{SS'}}& &\mapdr{\ts \sigma_{S_i\, S_i'}} \\
\phantom{\,\,.}F(I| S')&\mapr{\iota_{S'/T}} &F(I_1| S'_1)\ts\cdots\ts F(I_1| S'_1)\,\,.
\end{array}
$$
The map $\vphi_K$ is  compatible with  $\iota_{S/T}$:
With the same notation as above, and with $K_i=K\cap \ctop{I_i}$, 
the following commutes:
$$\begin{array}{ccc}
F(I| S)&\mapr{\iota_{S/T}}  &F(I_1| S_1)\ts\cdots\ts F(I_c| S_c) \\
\mapd{\vphi_K}& &\mapdr{\vphi_{K}} \\
F(I-K| S)&\mapr{\iota_{S/T}}  &F(I_1-K_1| S_1)\ts\cdots\ts F(I_c-K_c| S_c)\,\,.
\end{array}
$$
 
In (\RmkDefqDG) we discuss further properties of the complexes $F(I|S)$
that are consequences of (5) and (1)-(3). 
\bigskip 

(6)(Acyclicity of $\sigma$) For disjoint subsets $R, J$ of $\ctop{I}$ with $|J|\neq \emptyset$, 
the following sequence of 
complexes is exact, where the maps are alternating sums of 
$\sigma$, and $S$ varies over subsets of $J$:
$$F(I|R)\mapr{\sigma}
\moplus_{{|S|=1}\atop{S\subset J}} F(I| R\cup S)
\mapr{\sigma} \moplus_{{|S|=2}\atop{S\subset J}} F(I|R\cup S)
\mapr{\sigma} \cdots \to  F(I|R\cup J)\to 0\,\,.\eqno{(\DefqDG.b)}$$
\smallskip 

{\it Remarks to (6).}\quad  
$\bullet$\quad Since each $\sigma_{S\, S'}$ is a quasi-isomorphism, the total complex
of the double complex (\DefqDG.b) is acyclic. 

$\bullet$\quad From the remark above, the map 
$$\sigma: F(I|R)\to T:=\Tot[\moplus_{{|S|=1}\atop{S\subset J}} F(I| R\cup S)
\to \moplus_{{|S|=2}\atop{S\subset J}} F(I|R\cup S)
\to \cdots \to  F(I|R\cup J)\to 0]$$
is a quasi-isomorphism. 
This map factors as
$$F(I|R) \to \Ker[\moplus_{{|S|=1}\atop{S\subset J}} F(I| R\cup S)
\mapr{\sigma} \moplus_{{|S|=2}\atop{S\subset J}} F(I|R\cup S)]
\injto T\,;$$
the requirement implies that the first map is surjective, and the
first and the second maps are quasi-isomorphisms.

$\bullet$\quad  If $|J|=1$, the assumption says that  
$\sigma_{R\, R'}: F(I|R)\to F(I|R')$ is a surjective map. 
This was already required in (iii). 

$\bullet$\quad  One verifies by induction on $\sharp(R)$
 that the exactness of (\DefqDG.b) in case $R=\emptyset$ implies 
the exactness for any $R$. 
\bigskip

This concludes the definition of a quasi DG category. 
For the purpose of constructing a related quasi DG category $\cC^\Delta$
in subsequent sections,  
we need additional structure (iv) and (v) below. 
\bigskip

(iv) {\it Diagonal elements and diagonal extension}. \quad 
 We are given, for each object $X$ and a constant sequence of objects $i\mapsto X_i=X$
on a finite ordered set $I$ with $|I|\ge 2$, 
 a distinguished element, called the {\it diagonal element}, 
$$\bDelta_X(I)\in F(X; I)=F(X, \cdots, X)$$
of degree zero and coboundary zero.  
In particular for $|I|=2$ we write $\Delta_X=\bDelta_X(I)\in F(X, X)$.
In the sequel we will often  drop $X$ from $F( X; I)$. 

Let $(X_i)_{i\in I}$ be a sequence of objects on a finite ordered set $I$
with $|I|\ge 2$, and $\lambda: J\to I$ a surjective map of finite ordered sets. 
Then by $\lambda^* X$ we mean the sequence of objects $j\mapsto X_{\lambda(j)}$ on 
$J$. 
Assume given a map of complexes, called the {\it diagonal extension},
$$\lambda^*: F( X; I)\to F(\lambda^*X; J)$$
such that,  $\lambda^*=id $ if $\lambda=id$, and such that
if $\lambda': J'\to J$ is another surjective map, then 
${\lambda'}^*\lambda^*=(\lambda \lambda')^*$. 
When there is no confusion, we write $\diag(I, J)$ or $\diag$ for $\lambda^*$. 
One requires:
\smallskip 

(7)(Compatibility of the diagonal with the maps $\sigma$ and $\vphi$)\quad 
  The diagonal elements are assumed compatible with diagonal extension:
if $\lambda: J\to I$ is a surjective map of finite ordered sets, then 
$\lambda^*(\bDelta_X(I))=\bDelta_X(J)$. In particular, if 
$\lambda: I\to [1, 2]$ is any surjective map, one has 
$\bDelta_X(I)=\lambda^*(\Delta_X)$. 

The diagonal elements are compatible with 
the maps $\sigma$ and $\vphi$:
If $S\subset \ctop{I}$, and $I_1, \cdots, I_c$ is the corresponding 
segmentation, one has 
$$\tau_S(\bDelta_X(I)\, )=\bDelta_X(I_1)\ts\cdots\ts\bDelta_X(I_c)$$
in $F(I\tbar S)=F(I_1)\ts\cdots\ts F(I_c)$. For $K\subset\ctop{I}$, 
$$\vphi_K(\bDelta_X(I)\,)=\bDelta_X(I-K)$$
in $F(I-K)$. 
(These indeed follow from the compatibility of the diagonal extension and $\sigma$, 
$\vphi$, stated below.)
\smallskip 

The diagonal extension is compatible with the maps $\sigma$ and $\vphi$. 
For an element  $\ell\in J$ such that $\sharp \lambda^{-1}\{\lambda(\ell)\}=1$, 
the following square commutes:
$$\begin{array}{ccc}
F(X; I)   &\mapr{\lambda^*}  &F(\lambda^*X; J)  \\
\mapd{\vphi_{\lambda(\ell)}} & &\mapdr{\vphi_\ell} \\
F(X;  I-\{\lambda(\ell)\})   &\mapr{\lambda^*}  &F(\lambda^*X; J-\{\ell\} )  
 \end{array}
$$  
where the lower horizontal map is induced by the surjection 
$\lambda: J-\{\ell\}\to  I-\{\lambda(\ell)\}$. 
For $\ell\in J$ with $\sharp \lambda^{-1}\{\lambda(\ell)\}>1$, 
the diagram
$$\begin{array}{ccc}
F( X; I)   &\mapr{\lambda^*}  &F( \lambda^*X; J)  \\
 &\searrow &\mapdr{\vphi_\ell} \\
  &  &F( \lambda^*X; J-\{\ell\})  
 \end{array}
$$  
commutes, where the oblique arrow  $F( X;I)\to F(\lambda^*X; J-\{\ell\} )$ is induced by 
the surjection $\lambda: J-\{\ell\}\to  I$.

For $\ell\in \ctop{J}$ not over $\init(I)$ or $\term(I)$, 
let $J', J''$ be the segmentation of $J$ by $\ell$, 
and $I', I''$ be the segmentation of $I$ by $\lambda(\ell)$, and 
$\lambda': J'\to I'$, $\lambda'': J''\to I''$ be the surjections obtained 
by restricting $\lambda$. 
Then the following diagram commutes:
$$\begin{array}{ccc}
F(X;I)&\mapr{\lambda^*} &F( \lambda^*X; J) \\
\mapd{\tau_{\lambda(\ell)}}& &\mapdr{\tau_\ell} \\
F(X; I')\ts F(X; I'')&\mapr{{\lambda'}^*\ts {\lambda''}^*} &F(\lambda^*X; J')\ts F(\lambda^*X;
J'')\,\,.
\end{array}
$$
Here the lower horizontal map is the tensor product of the maps 
${\lambda'}^*$ and ${\lambda''}^*$. 
For $\ell\in \ctop{J}$ over $\init(I)$ or $\term(I)$, say $I=[1, n]$ and 
$\lambda(\ell)=n$, then 
with $J'$, $J''$ and $I'$, $I''$ as above, one has 
the following commutative diagram
$$\begin{array}{ccc}
F(  X; I)&\mapr{\lambda^*} &F(\lambda^*X; J) \\
\mapd{{\lambda'}^*}& &\mapdr{\tau_\ell} \\
F( \lambda^*X; J' )&\mapr{} &F( \lambda^*X; J')\ts F(X_n, \cdots, X_n)\,\,.
\end{array}
$$
Here the lower horizontal map is induced by $\lambda': J'\to I$, and 
the lower vertical map is $u\mapsto u\ts \bDelta_{X_n}(J'')$
with $\bDelta_{X_n}(J'')\in  F(X_n; J'')=F(X_n, \cdots, X_n)$. 
Similarly in case $\lambda(\ell)=\init(I)$. 
\bigskip 

{\it Remark to (iv).}\quad 
Given a surjection $\lambda: J\to I$ as above, the finite ordered set
$J$ is determined by $I$ and the function $m: I\ni i\mapsto m_i=\sharp(\lambda^{-1}(i))\in \ZZ_{>0}$. 
If $I=[1, n]$,  $J$ is isomorphic to the ordered set 
$$\{{1_1, \cdots, 1_{m_1}},
{2_1, \cdots, 2_{m_2}},  
\cdots,
{n_1, \cdots, n_{m_n}}\,\}\,\,$$
($i$ is repeated $m_i$ times). 
Then one has $F( X; I)=F(X_1, \cdots, X_n)$, and 
$$F(\lambda^*X; J )=
F(\overbrace{X_1, \cdots, X_1}^{\mbox{$m_1$ times}},
\overbrace{X_2, \cdots, X_2}^{\mbox{$m_2$ times}},
\cdots, 
\overbrace{X_n, \cdots, X_n}^{\mbox{$m_n$ times}})\,.$$
In the rest of this paper, instead of $J$ we will usually 
give $I$ and the multiplicities $\{m_i\}$. 
\bigskip

(v) {\it The set of generators, notion of proper intersection, and 
distinguished subcomplexes with respect to constraints}. 
\smallskip 

We impose that the complex $F( X; I)$ be free with a given basis set, 
and that there is given a notion of proper intersection for elements 
of $F( X; I)$.  

 For a sequence $X$ on $I$, the complex $F(I)=F(X; I)$ 
is assumed degree-wise $\ZZ$-free on a given set of 
generators $\cS_F(I)=\cS_F( X; I)$. 
More precisely $\cS_F(I)=\amalg_{p\in \ZZ} \cS_F(I)^p$, 
where $ \cS_F(I)^p$ generates $F(I)^p$. 
A non-zero element $u\in F(I)$ can be written $\sum c_\nu \al_\nu$, 
with $c_\nu$ non-zero integers and $\al_\nu\in \cS_F(I)$; then we say 
that $\al_\nu$ appears in the basis expansion of $u$. 

Let $I$ be a finite ordered set, and 
let $\{I_\al\}_{\al\in A}$ be a collection of subintervals (of cardinality $\ge 2$)
 of $I$,  indexed by a set $A$.  We say $\{I_\al\}$ is {\it almost disjoint} if for 
$\al$, $\be$ distinct, 
$I_\al\cap I_\be$ consists of at most one element. 
Then there is a total ordering $<$ on the set $A$ given by $\al < \be$ if and
only if $\init(I_\al)<\init(I_\be)$. 
When $\al < \be$, we have either $I_\al$ and $I_\be$ disjoint, or 
$\term(I_\al)=\init(I_\be)$.

 Let $I$ be a finite ordered set, and $I_1, \cdots, 
I_r$ be almost disjoint subintervals of $I$.  
Let $X$ be a sequence of objects on $I$.
We assume  given a subset $\BP(\{I_i\}_{1\le i\le r}) $ of the product set 
$$\prod_{1\le i\le r} \cS_F(I_i)\,,$$
satisfying the conditions (8) below. 
It will be useful to introduce a terminology.
For $\al_i\in \cS_F(I_i)$, a collection of elements indexed by  $\{1, \cdots, r\}$, 
we shall say  $\{\al_i\}$ is   {\it properly 
intersecting} when it is in $\BP$.
(The notation $\{\al_i\}$ is being used for the indexed set $\{\al_i\}_{i\in [1, r]}$,
 and it should not be confused with a set with elements $\al_i$.)
The terminology comes from examples related to  algebraic cycles, where the condition means 
algebraic cycles intersecting properly. 
\smallskip 
 
For the basis set $\cS_F(I)$ and the notion of proper intersection, the
following condition is to be satisfied.
\smallskip 

(8)(Notion of proper intersection and distinguished subcomplexes)\quad
For the condition of proper intersection we require three conditions:
\smallskip

(PI-1)\quad  If $\{\al_i \}_{ i\in [1, r]  }$
is properly intersecting, for any subset $J$ of $[1, r]$, 
$\{\al_i \}_{i\in  J}$
is properly intersecting. 

(PI-2)\quad   Assume that $[1, r]=J\amalg J'$, and that 
$\cup_{j\in J}I_j$ and $\cup_{j\in J'}I_j$ are disjoint. 
Then $\{\al_i \}_{ i\in [1, r]  }$
is properly intersecting if and only if $\{\al_i \}_{ i\in J }$
and $\{\al_i \}_{ i\in J' }$ are both properly intersecting. 

(PI-3)\quad  Assume $\{\al_1, \cdots, \al_r\}$ is properly intersecting, 
$i\in [1, r]$, and $\be_\nu\in \cS_F(I_i)$ appears in the basis expansion of $\partial \al_i$ 
(one has $\partial \al_i=\sum c_{i\nu}\be_\nu$ with 
$\be_\nu\in \cS_F(I_i)$ and $c_{i\nu}\neq 0$.)
Then the set 
$$\{\al_1, \cdots, \al_{i-1}, \be_\nu, \al_{i+1}, \cdots, \al_r\}$$
is properly intersecting.  
\smallskip 

The notion of proper intersection can be naturally extended to elements in the complexes $F(I)$:
for a collection of elements $u_i\in F(I_i)$ indexed by $i\in [1, r]$, 
we define $\{u_i\}$ to be properly intersecting if,  
for any choice of elements $\al_i\in \cS_F(I_i)$ appearing in the basis expansion of 
$u_i$, the indexed set $\{\al_i\}$ is properly intersecting. 
(We ignore those $i$ with $u_i=0$.)
\bigskip

Next we  explain conditions of constraint and distinguished subcomplexes.
Let $I$ be a finite ordered set, $I_1, \cdots, 
I_r$ be almost disjoint sub-intervals such that 
$\init(I_1)<\cdots <\init(I_r)$  and 
$\cup I_i=I$;
 equivalently, 
$\init(I_1)=\init (I)$, $\term(I_i)=\init(I_{i+1})$ or $\term(I_i)+1=
\init(I_{i+1})$ for $1\le i< r$,  and $\term(I_r)=\term(I)$. 
Assume given a sequence of objects $(X_i)_{i\in I}$ on $I$.
We shall consider a class of subcomplexes of 
$F(I_1)\ts \cdots \ts F(I_r)$ specified as follows. 

By a condition of {\it constraint} we mean a set of data
$$\cC=(I\injto \BI; X\mbox{ on } \BI;
  P; \{J_j\}_{j=1, \cdots, s}; \{f_j\in F(J_j)\}\,)\eqno{(\DefqDG.c)}$$
where 
\smallskip 

(a) $I\injto \bbI$ is an inclusion into another finite ordered set $\bbI$
such that the image of each $I_a$ is a sub-interval,

(b) $X$ is an extension of X to $\bbI$, still denoted by $X$, 

(c) $P\subset [1, r]$ is a (possibly empty) subset, 
 
(d)  $J_1, \cdots, J_s\subset \bbI$ is a (possibly empty)
set of sub-intervals of $\bbI$ such that 
the indexed family $\{I_i, J_j\}_{i,j}$  is almost disjoint
(namely, $I_i\cap J_j$ and $J_j\cap J_{j'}$ for $j\neq j'$ consist of at most one element),
and 

(e) $f_j\in F(J_j)$, $j=1, \cdots, s$, is a set of elements such that 
$\{f_j\}$ is properly intersecting. 
\smallskip 
  
Given such consider the subcomplex of $F(I_1)\ts \cdots \ts F(I_r)$ generated by 
$\al_1\ts\cdots \ts \al_r$, $\al_i\in \cS_F(I_i)$, such that the indexed set 
$$\{\{\al_i\}_{i\in P}, \{f_j\}_{1\le j\le s}\}$$
is properly intersecting
(it is a subcomplex by the last property of proper intersection).
If $P$ is empty, there is no condition.  

This subcomplex is denoted 
$$[F(I_1)\ts \cdots \ts F(I_r)]_\cC$$
specifying the condition of constraint, or more 
simply $[F(I_1)\ts \cdots \ts F(I_r)]_{\BI; f}$ or $[F(I_1)\ts \cdots \ts F(I_r)]_{f}$. 
A subcomplex of $F(I_1)\ts \cdots \ts F(I_r)$ obtained as a finite intersection of 
such subcomplexes is called a {\it distinguished subcomplex}. 
Thus a distinguished subcomplex is one of the form 
$$[F(I_1)\ts \cdots \ts F(I_r)]'=\bigcap_{i=1}^{m}[F(I_1)\ts \cdots \ts F(I_r)]_{\cC_i}$$
where $\cC_i$, $i=1, \cdots, m$ is a finite collection of conditions of constraint;
for distinct $i$'s, the inclusions $I\injto \BI_i$ are not assumed related. 

We require two conditions:
\smallskip 

(DS-1)\quad The inclusion of any 
 distinguished subcomplex 
$[F(I_1)\ts \cdots \ts F(I_r)]'\injto F(I_1)\ts \cdots \ts F(I_r)$
is a quasi-isomorphism.
\smallskip 

(DS-2)\quad 
As a special case, assume $I_1, \cdots, I_r$ is a segmentation of $I$, namely when 
 $\init(I_1)=\init (I)$, $\term(I_i)=\init(I_{i+1})$ for $1\le i<r$,  and $\term(I_r)=\term(I)$, the subcomplex of $F(I_1)\ts\cdots\ts F(I_r)$ generated by $\al_1\ts\cdots\ts\al_r$ with
$\{\al_i\}_{1\le i\le r}$  properly intersecting is a distinguished subcomplex;
it corresponds to the condition of constraint 
$\cC=(I=\BI, P=[1, r])$ with no $\{f_j\}$. 
We require that it coincides with $F(I|S)$ where   
 $S\subset \ctop{I}$ is the subset corresponding to 
the segmentation. 
\smallskip 

The second condition can be phrased in concrete terms as follows. 
Recall that $F(I_i)=\ZZ \cS_i$, the free abelian group on 
$\cS_i:=\cS_{F}(I_i)$. We have 
$$F(I_1)\ts\cdots\ts F(I_r)=\ZZ[\cS_1\times\cdots\times\cS_r]\,,$$
and $F(I|S)=\ZZ[\BP(\{I_i\}]$, with $\BP(\{I_i\})\subset \cS_1\times\cdots\times\cS_r$.
An element $u\in F(I_1)\ts\cdots\ts F(I_r)$ has a basis expansion 
$$u=\sum c_{\al_1, \cdots, \al_r}\al_1\ts\cdots\ts\al_r$$
with $(\al_1, \cdots, \al_r)\in \cS_1\times\cdots\times\cS_r$ and $c_{\al_1, \cdots, \al_r}$ non-zero integers. We then say that $\al_1\ts\cdots\ts\al_r$ appears 
in the basis expansion of $u$. 
So $u\in F(I|S)$ if and only if for each $\al_1\ts\cdots\ts\al_r$ appearing 
in its basis expansion, $\{\al_1, \cdots, \al_r\}$ is properly intersecting. 
In particular, for elements for $u_i\in F(I_i)$, the index set $\{u_1, \cdots, u_r\}$ is properly intersecting 
if and only if $u_1\ts\cdots\ts u_r\in F(I|S)$. 

In the definition of a constraint, there appears the set $P$. 
If we take as $P$ a smaller set, while keeping the other data, the corresponding
distinguished subcomplex becomes larger, by condition (PI-1). 

We also note that tensor product of distinguished subcomplexes is again 
a distinguished subcomplex. 
To explain it let  $[F(I_1)\ts \cdots \ts F(I_r)]'$
be a distinguished subcomplex as above. 
Let $I'$ be another finite ordered set, 
$I'_1, \cdots, I'_{r'}$ be almost disjoint set of sub-intervals such that 
$\init(I'_1)<\cdots <\init(I'_{r'})$  and 
$\cup I'_i=I'$; 
let $X'$ be a sequence of objects on $I'$.
Assume given a distinguished subcomplex  
$[F(I'_1)\ts \cdots \ts F(I'_{r'})]'$. 
Then the disjoint union $I\amalg I'$ is a finite totally ordered set 
(with $x< x'$ if $x\in I$ and $x'\in I'$), and 
$I_1, \cdots, I_r, I'_1, \cdots, I'_{r'}$ is an almost disjoint set of sub-intervals
with union $I\amalg I'$. 
Under these hypotheses, the tensor product 
$$[F(I_1)\ts \cdots \ts F(I_r)]'\ts [F(I'_1)\ts \cdots \ts F(I'_{r'})]'$$
is a distinguished subcomplex of $F(I_1)\ts \cdots \ts F(I_r)\ts F(I'_1)\ts \cdots \ts F(I'_{r'})$. (Proof uses (PI-2)\,). 
\smallskip 

{\it Remark to (v).}\quad   When we assume the structure (iv),  the condition (4) is redundant. 
From (7) and (8) it follows that $[\Delta_X]\in H^0F(X, X)$
is the identity in the sense of (\HowqDG). 
Indeed for $u\in H^0F(X, Y)$ take its representative $\underline{u}\in F(X, Y)^0$, then  
take its diagonal extension $\diag(\underline{u})\in F(X, X, Y)^0$. 
Then one has $\tau_2 (\diag(\underline{u}))=\Delta_X\ts \underline{u}$ and
$\vphi_2(\diag(\underline{u})\,) =\underline{u}$. 
The same argument shows the following, which is stronger than (4):
\smallskip 

(4)' For each $u\in H^nF(X, Y)$, $n\in \ZZ$, one has  
$1_X\cdot u =u$. 
Similarly for $u\in H^nF(Y, X)$,  $u\cdot 1_X=u$.
\bigskip 

\sss{\RmkDefqDG} {\it Remarks to (\DefqDG).}\quad 
The reader may skip this subsection until it is needed in (\GKL),(\HIGI), 
(\GISigma), and (\Sigmaprolong). 

From the condition (5), we easily deduce that 
the properties (1)-(3) for $F(X_1, \cdots, X_n)$ ``descend"
to the corresponding properties (1)'-(3)' for $F(X_1, \cdots, X_n|S)$, as we list 
below.
\smallskip 

(\RmkDefqDG.1) {\bf Properties of $F(I|S)$.}\quad 
(1)' (Functoriality.)\quad 
 The complex $F(X_1, \cdots, X_n|S)$ is covariantly 
functorial for morphisms:
for  a sequence of morphisms $f=(f_1, \cdots, f_n): (X_1, \cdots, X_n)$
\newline $\to (Y_1, \cdots, Y_n)$,
 there corresponds an isomorphism of multiple complexes
$$f_*: F(X_1, \cdots, X_n|S)\to F(Y_1, \cdots, Y_n|S)$$
that is covariantly functorial in $f$. 
If $X_i=O$ for some $i$, then $F(X_1, \cdots, X_n|S)$
equals zero.

The inclusion $\iota_{S/T}$ is  covariantly functorial:
For a sequence of morphisms $f$,  the square 
$$\begin{array}{ccc}
 F(X; I|S)   &\mapr{\iota_{S/T}}  &F(X; I_1|S_1)\ts \cdots\ts F(X; I_c|S_c)  \\
\mapd{f_*} &  &\mapdr{f_*} \\
  F(Y; I|S)   &\mapr{\iota_{S/T}}  &F(Y; I_1|S_1)\ts \cdots\ts F(Y; I_c|S_c)
 \end{array}
$$  
commutes. 

The map $\sigma_{S\, S'}$ is covariantly functorial:
For a morphism $f$,
 the following 
square commutes:
$$\begin{array}{ccc}
 F(X; I|S)   &\mapr{\sigma_{S\, S'}}  & F(X; I|S')   \\
\mapd{f_*} & &\mapdr{f_*} \\
F(Y, I|S)   &\mapr{\sigma_{S\, S'}}  &\phantom{\,,}F(Y, I|S')\,,
 \end{array}
$$ 
namely $f_*\sigma_{S\, S'}= \sigma_{S\, S'}f_*$.

The map $\vphi_K$ is covariantly functorial:
For a morphism $f$   
the square
$$\begin{array}{ccc}
 F(X; I|S)   &\mapr{\vphi_K }  & F(X; I-K|S)  \\
\mapd{f_*} & &\mapdr{ f_*} \\
F(Y; I|S)   &\mapr{\vphi_K }  & F(Y; I-K|S)
\end{array}
$$  
commutes, namely $f_*\vphi_K =\vphi_Kf_*$. 
\smallskip

(2)' (Commutation identities.)\quad 
The maps $\sigma_{S\, S'}$ and $\vphi_K$ are compatible
 with each other, as follows. 
 
We have that 
 $\sigma_{S\, S}=id$ and, for $S\subset S'\subset S''$, $\sigma_{S\,S''}=
\sigma_{S'\,S''}\sigma_{S\,S'}$.
If $K=K'\amalg K''$ then 
$\vphi_K=\vphi_{K''}\vphi_{K'}:
F(I\dbar S)\to  F(I-K\dbar S)$.

The maps $\sigma$ and $\vphi$ commute: For  $K$ and $S'$
 disjoint,  the following diagram commutes:
$$\begin{array}{ccc}
F(I| S)&\mapr{\vphi_K} 
&F(I-K | S) \\
\mapd{\sigma_{SS'}}& &\mapdr{\sigma_{SS'}} \\
F(I| S')&\mapr{\vphi_K} &F(I-K | S')\,\,.
\end{array}
$$
\smallskip 

(3)' (Additivity of $F(X_1, \cdots, X_n|S)$.)\quad 
For each $i$, $1\le i\le n$, and a sequence of objects 
$(X_1, \cdots, X_{i-1}, Y_i, Z_i, X_{i+1}, \cdots, X_n)$, 
one has maps of complexes 
$$\lsum_i(Y_i, Z_i): 
F(X_1,\cdots, X_{i-1}, Y_i, X_{i+1}, \cdots, X_n|S)\to 
F(X_1,\cdots, X_{i-1}, Y_i\bop Z_i, X_{i+1},\cdots, X_n|S)$$
and 
$$\rsum_i(Y_i, Z_i): 
F(X_1,\cdots, Z_i, \cdots, X_n|S)\to 
F(X_1,\cdots, Y_i\bop Z_i, \cdots, X_n|S)\,.$$
If  $1<i<n$ and $i\not\in S$, one also has a map of complexes
$$\pi_i(Y_i, Z_i): F(X_1, \cdots, X_{i-1}, Y_i|S_1)\ts
 F(Z_i,X_{i+1},  \cdots, X_n|S_2) 
\to F(X_1, \cdots, Y_i\bop Z_i, \cdots, X_n|S)\, $$
where $S_1=(1, i)\cap S$ and $S_2=(i, n)\cap S$. 
These maps satisfy conditions parallel to those in (3);
they are obtained in a obvious manner by replacing $F(X_1, \cdots, X_n)$
with $F(X_1, \cdots, X_n|S)$.  Among them we will mention below only less trivial  
ones. 
\smallskip 

$\bullet$\quad Compatibility with $\sigma$.\quad
The maps $s_i$ and $\sigma_{S\, S'}$ commute, namely the diagram 
$$\begin{array}{ccc}
F(X_1,{\scriptstyle\cdots}, Y_i,{\scriptstyle\cdots}, X_n|S)    &\mapr{s_i(Y_i, Z_i)}  &
F(X_1,{\scriptstyle\cdots}, Y_i\bop Z_i,{\scriptstyle\cdots}, X_n|S)  \\
\mapd{\sigma_{S\, S'}} && \mapdr{\sigma_{S\, S'}} \\
F(X_1,{\scriptstyle\cdots}, Y_i,{\scriptstyle\cdots}, X_n|S')    &\mapr{s_i(Y_i, Z_i)}  &
F(X_1,{\scriptstyle\cdots}, Y_i\bop Z_i,{\scriptstyle\cdots}, X_n|S')  \\
 \end{array}
$$  
commutes. 
The maps $\pi_i$ and $\sigma_{S\, S'}$ commute, namely if $S\subset S'$ and 
$i\not\in S'$, the following square commutes:
$$\begin{array}{ccc}
F(X_1,{\scriptstyle\cdots}, Y_i|S_1)\ts F(Z_i,{\scriptstyle\cdots}, 
X_n|S_2) &\mapr{\pi_i (Y_i, Z_i)} 
&F(X_1,{\scriptstyle\cdots}, Y_i\bop Z_i,{\scriptstyle\cdots}, X_n|S)  \\
\mapd{\sigma_{S_1\, S_1'}\ts \sigma_{S_2\, S_2'}} & & \mapdr{\sigma_{S\, S'}} \\
F(X_1,{\scriptstyle\cdots}, Y_i|S'_1)\ts F(Z_i,{\scriptstyle\cdots}, 
X_n|S'_2) &\mapr{\pi_i (Y_i, Z_i)} 
&F(X_1,{\scriptstyle\cdots}, Y_i\bop Z_i,{\scriptstyle\cdots}, X_n|S')\,.  
\end{array}
$$

If $i\in S$, the following diagram commutes:
$$\begin{array}{ccc}
F(X_1,{\scriptstyle\cdots}, Y_i|S_1)\ts F(Z_i,{\scriptstyle\cdots}, 
X_n|S_2) &\mapr{\pi_i (Y_i, Z_i)} 
&F(X_1,{\scriptstyle\cdots}, Y_i\bop Z_i,{\scriptstyle\cdots}, X_n|S)  \\
   \phantom{aaaaaaaaaaaaaaaaaaaaaa}  {\scriptstyle s\ts t}\!\!\!\!\!\!\!\!\!\!\! &\searrow      &\mapdr{\iota_{S/\{i\} }  } \\
    &&F(X_1, {\scriptstyle\cdots}, Y_i\bop Z_i|S_1)\ts F(Y_i\bop Z_i,{\scriptstyle\cdots}, X_n|S_2)\,.
\end{array}$$ 
\smallskip 

$\bullet$\quad Compatibility with $\vphi$.\quad
The maps $s_i$ and $\vphi_i$ commute, which means the commutativity of 
the following diagram:
$$\begin{array}{ccc}
F(X_1,{\scriptstyle\cdots}, Y_i, {\scriptstyle\cdots}, 
X_n|S) &\mapr{s_i (Y_i, Z_i)} 
&F(X_1,{\scriptstyle\cdots}, Y_i\bop Z_i,{\scriptstyle\cdots}, X_n|S)  \\
   \phantom{aaaaaaaaaaaaaaaaaaaaaa}  {{\scriptstyle\vphi_i}}\!\!\!\!\!\! &\searrow      &\mapdr{\vphi_i }   \\
    &&F(X_1, {\scriptstyle\cdots}, X_{i-1}, X_{i+1},{\scriptstyle\cdots}, X_n|S)\,.
\end{array}$$ 
\smallskip 

$\bullet$\quad  The  map  of additivity $\theta$.\quad 
The map  
\begin{eqnarray*}
&\theta_i(Y_i, Z_i):& F(X_1, \cdots, Y_i, \cdots, X_n|S) \bop
 F(X_1, \cdots, Z_i, \cdots, X_n|S) \\
&& \bop F(X_1, \cdots, Y_i|S_1)\ts F(Z_i, \cdots, X_n|S_2) \\
&&\bop
  F(X_1, \cdots, Z_i|S_1)\ts F(Y_i, \cdots, X_n|S_2) \\
&\mapr{}&   
F(X_1, \cdots, Y_i\bop Z_i, \cdots, X_n|S)\,,
\end{eqnarray*}
defined as the sum of the maps $s_i$, $t_i$, $\pi_i(Y_i, Z_i)$, and 
$\pi_i(Z_i, Y_i)$ (if $i=1$ or $n$, there are no ``cross terms"), 
is a quasi-isomorphism. 

We note that the map $\theta$ is compatible with $\sigma$ and $\vphi$.
\bigskip 

(\RmkDefqDG.2) {\bf Functorial properties of proper intersection.}\quad 
Assuming (v) and (8), we state functorial properties
of the notion of proper intersection. 
These are translations of functorial properties of 
the complex $F(I|S)$ discussed above. 
\smallskip 

(1)  Let $X$, $Y$ be sequences of objects indexed by $I=[1, n]$, and 
 $f: X\to Y$ be a sequence of isomorphisms of objects. 
 Let $I_1, \cdots, I_r$ be a segmentation of $I$.  
If $\{u_i\in F(X; I_i)\}_{i=1, \cdots, r}$ is properly intersecting, then 
$\{f_*(u_i)\in F(Y; I_i)\}_{i=1, \cdots, r}$ is also properly intersecting. 
(Indeed one has 
$u:=u_1\ts\cdots\ts u_r\in F(X|S)$ if $S$ corresponds to the segmentation.
Then $f_*(u)=(f_*u_1)\ts\cdots\ts (f_*u_r)\in F(Y|S)$, showing the assertion.)

(2) Assume $X=Y_i\oplus Z_i$ for each $i\in [1, n]$, and 
$s: F(Y; I) \to F(X; I)$ be the corresponding map for $I\subset [1, n]$. 
If $\{u_i\in F(Y; I_i)\}_{i=1, \cdots, r}$ is properly intersecting, then 
$\{s(u_i)\in F(X; I_i)\}_{i=1, \cdots, r}$ is also properly intersecting. 
(The verification is similar to that for (1), using 
$s(u_1\ts\cdots \ts u_r)=s(u_1)\ts\cdots\ts s(u_r)$. )

(3) Let $X$ be a sequence of objects on $I=[1, n]$, $I_1, \cdots, I_r$  a segmentation 
corresponding to $S$, and $K\subset (1, n)$ be disjoint from $S$. 
If $\{u_i\in F(X; I_i)\}_{i=1, \cdots, r}$ is properly intersecting, then letting 
$K_i=K\cap I_i$, the indexed set 
$\{\vphi_{K_i}(u_i)\in F(X: I_i-K)\}_{i=1, \cdots, r}$ is properly intersecting.

(4) To state that the notion of proper intersection is compatible with $\sigma$,  we need some preliminaries. 
With $X$ and $I_1, \cdots, I_r$ as above, 
a set of elements $v_j\in F(I_1)\ts\cdots\ts F(I_r)$ (for $j=1, \cdots, m$)
is said to be {\it basis-disjoint} if each basis $\al_1\ts\cdots\ts\al_r$ 
appears in at most one of $v_j$'s. 
Then one has that 
$\sum_{j=1, \cdots, m} v_j \in F(I|S)$ if and only if 
each $v_j\in F(I|S)$. 

Assume now given elements $u_i\in F(I_i)$ such that $\{u_1, \cdots, u_r\}$
is properly intersecting, an integer $j$ with $1\le j\le r$, and $k\in \ctop{I_j}$. 
Express 
$$\sigma_k(u_j)=\sum_\lambda u'_\lambda\ts u''_\lambda$$
where the set $\{u'_\lambda\ts u''_\lambda\}$ indexed by $\lambda$ is basis-disjoint.
Then for each $\lambda$, the set 
$$\{u_1, \cdots, u_{j-1}, u'_\lambda, u''_\lambda, u_{j+1}, \cdots, u_r\}$$
is properly intersecting. 
(For the proof, apply $\sigma_k$ to the element 
$u_1\ts\cdots\ts u_r\in F(I|S)$ to get 
$$\sum_\lambda u_1\ts\cdots\ts u_{j-1}\ts u'_\lambda\ts u''_\lambda\ts\cdots u_r
\in F(I|S\cup \{k\})\,.$$
Then the elements 
$\{u_1\ts\cdots\ts u'_\lambda\ts u''_\lambda\ts\cdots u_r\}_\lambda$ are basis-disjoint, 
so we have $u_1\ts\cdots\ts u'_\lambda\ts u''_\lambda\ts\cdots u_r\in F(I|S)$.)

(5) Let $X_1, \cdots, X_n$ be a sequence of objects,
with a given decomposition $X_e=Y_e\oplus Z_e$ for an element 
$e$ with $1<e<n$.  
For $r, s\ge 1$, assume given elements 
\begin{eqnarray*}
u_1\in& F(X_1, X_2, \cdots, X_{a-1}, X_a)\,, \\
u_2\in &F(X_a, X_{a+1}, \cdots, X_{b-1}, X_b)\,,\\
&\cdots  \\
u_r\in & F(X_d, X_{d+1}, \cdots, X_{e-1}, Y_e)\,,\\
u_{r+1}\in & F(Z_e, X_{e+1}, \cdots, X_{f-1}, X_f)\,,\\
u_{r+2}\in &F(X_f, X_{f+1}, \cdots, X_{g-1}, X_g)\,,\\
&\cdots  \\
u_{r+s}\in & F(X_i, X_{i+1}, \cdots, X_{n-1}, X_n)\,.
\end{eqnarray*}
(the variables are $X_1,\cdots, X_n$, except $Y_e$ or 
$Z_e$ in the $e$-th spot)
such that 
$\{u_1, \cdots, u_r\}$ and 
$\{u_{r+1}, \cdots , u_{r+s}\}$ are both properly intersecting. 
Then the indexed set 
$$\{u_1,  \cdots, u_{r-1}, s(u_r), t(u_{r+1}),
u_{r+2},  \cdots , u_{r+s}\}\,,$$
is properly intersecting, 
where $s(u_r)$ is the image by the map 
$s: F(X_d, X_{d+1}, \cdots, X_{e-1}, Y_e)
\to F(X_d, X_{d+1}, \cdots, X_{e-1}, X_e)$,
and similarly for $t(u_{r+1})\in F(X_e, X_{e+1}, \cdots, 
X_{f-1}, X_f)$. 
Also the indexed set 
$$\{u_1,  \cdots, u_{r-1}, \pi_e(u_r\ts u_{r+1}), 
u_{r+2},  \cdots , u_{r+s}\}$$
is properly intersecting, 
where $\pi_e(u_r\ts u_{r+1})$ is the image of $u_r\ts u_{r+1}$ by the map 
$$\pi_e:F(X_d, \cdots, Y_e)\ts F(Z_e,\cdots, X_f)\to 
F(X_d, X_{d+1}, \cdots, X_{e-1}, X_e, X_{e+1}, \cdots, X_{f-1}, X_f)\,.$$
(This is shown, if $r=s=2$, say, as follows. One has  
$u=u_1\ts u_2\in F(X_1,{\scriptstyle\cdots}, X_a,{\scriptstyle\cdots}, Y_e|\{a\}))$
and  
$u'=u_3\ts u_4 \in F(Z_e, {\scriptstyle\cdots}, X_f, {\scriptstyle\cdots}, X_n|\{f\})$. 
Let $v=\pi_e(u\ts u')\in F(X_1,{\scriptstyle\cdots},  X_e,
{\scriptstyle\cdots},  X_n|\{a, f\})$. 
Then we have 
$$ \pi_{\{a,e,  f\}}(v)=u_1\ts u_2\ts u_3\ts u_4$$
and 
$$ \pi_{\{a,  f\}}(v)=u_1\ts \pi_e(u_2\ts u_3)\ts u_4$$
showing that both $\{u_1, u_2, u_3,  u_4\}$ and 
$\{u_1, \pi_e(u_2\ts u_3),  u_4\}$ are properly intersecting.)
\bigskip

\sss{\ExqDG} {\bf Example.}\quad 
 A DG category can be viewed as a quasi DG category. 
It is a weak DG category, as already discussed in (\FtbarS).

The condition (5) is satisfied with 
$F(X_1, \cdots, X_n)=F(X_1, \cdots , X_n |S)=F(X_1, \cdots , X_n \tbar S)$
and $\sigma_S=\iota_S=id$. 
The condition (6) is also verified. 

Indeed we also have the additional structure (iv) and (v). 
We define 
$\Delta_X=id_X\in F(X, X)$ and 
$$\bDelta_X(I)=\Delta_X\ts\cdots\ts\Delta_X\in F(X; I)\,;$$
the diagonal extension $\lambda^*$ is defined in the evident manner. 
Then the condition (7) is satisfied. 
As for the condition of proper intersection, we declare {\it any\,}
indexed set $\{\al_i\}$ in $F(I_i)$  be properly intersecting. 
\bigskip

\sss{\SymbqDG} {\it The quasi DG category $\Symb(S)$.}\quad 
Let $S$ be a quasi-projective variety. 
The category of smooth varieties $X$ equipped with projective maps to $S$
will be denoted by $({\rm Smooth}/k\,, {\rm Proj}/S)$.  
A {\it symbol} over $S$ is an object the form 
$$\moplus_{\al\in A}(X_\al/S, r_\al)$$
where $X_\al$ is a collection of objects in $({\rm Smooth}/k\,, {\rm Proj}/S)$
indexed by a finite set $A$, and $r_\al\in \ZZ$. 
A morphism of symbols 
$$f: \moplus_{\al\in A}(X_\al/S, r_\al)\to \moplus_{\be\in  B}(Y_\be/S, s_\be)$$
is given by an isomorphism of sets $u: A\to B$ such that 
$r_\al =s_{u(\al)}$ for each $\al\in A$ and a collection of 
isomorphisms of varieties over $S$,  $f_\al: X_\al \to Y_{u(\al)}$ for 
$\al\in A$.  Composition of morphisms is given in the evident way.
The direct sum of two symbols is defined in an obvious manner, 
and the zero object corresponds to the symbol with $A=\emptyset$.

In \partI\,  we defined 
\begin{itemize}
\item{} the complexes $F(K_1, \cdots, K_n|S)$ for 
a sequence of symbols $K_i$ and $S\subset (1, n)$, 

\item{} the maps $\iota$, $\sigma$ and $\vphi$, 

\item{} the diagonal elements $\bDelta_K(I)$ and 
the diagonal extension, 

\item{} the notion of proper intersection on $F(K_1, \cdots, K_n)$, 
\end{itemize}
and  verified the conditions (1)-(8) for a quasi DG category. 
We refer to this as the quasi DG category $Symb(S)$. 

We have the relation to the cycle complex of S. Bloch, as follows. 
For $K=(X/S, r)$ and $(Y/S, s)$, there is a quasi-isomorphism
$$ \Z_{\dim Y-s+r}(X\times_S Y)\to F(K, L)\,. $$
The left hand side is the cycle complex of the fiber product $X\times_S Y$ 
in dimension $\dim Y-s+r$. 
Thus there is a canonical isomorphism
with the higher Chow group
$$H^{-n} F((X/S, r), (Y/S, s)\,)=
\CH_{\dim Y-s+r}(X\times_S Y, n)\,.$$
Via this isomorphism, the map $\psi$ in Remark to (\HowqDG) reads
$$\psi: 
\CH_{\dim Y-s+r}(X\times_S Y, n)\ts
\CH_{\dim Z-t+s}(Y\times_S Z, m)
\to \CH_{\dim Z-t+r}(X\times_S Z, n+m)\,.$$

If $f: X\to Y$ is a map over $S$, then its graph 
$\Gamma_f$ is an element of 
$\Z_{\dim X}(X\times_S Y) $ of degree zero, thus it gives a 
cocycle of degree zero in the complex $F((Y/S, 0), (X/S, 0)\,)$. 
Let $[\Gamma_f]\in H^0F((Y/S, 0), (X/S, 0)\,)$ be its cohomology class.
If $g: Y\to Z$ is another map, then one verifies that 
$\psi ([\Gamma_f]\ts [\Gamma_g])=[\Gamma_{g\scirc f}]$. 
\bigskip 

\sss{\Cdiag} {\it $C$-diagrams.}\quad 
Let $\cC$ be a quasi DG category. 
We will construct another quasi DG category $\cC^\Delta$ out of 
$\cC$. An object of  $\cC^\Delta$ is of the form 
$K= (K^m; f(m_1, \cdots  , m_\mu)\, )$, 
where $(K^m)$ is a sequence of objects of $\cC$
indexed by $m\in \ZZ$, almost all
of which are zero, and 
$$f(m_1, \cdots , m_\mu)\in F(K^{m_1}, \cdots , K^{m_\mu})^{-(m_\mu-m_1-\mu+1)}$$
is a collection of elements indexed by sequences $(m_1<m_2<\cdots<m_\mu)$
with $\mu\ge 2$. We require the following conditions:

(i) For each $j=2, \cdots, \mu-1$
$$\sigma_{K^{m_j}}(f(m_1, \cdots , m_\mu)\, )
=f(m_1, \cdots , m_j)\ts f(m_j, \cdots , m_\mu)\,\, $$
in $F(K^{m_1}, \cdots , K^{m_j})\ts F(K^{m_j}, \cdots , K^{m_\mu})$. 

(ii) For each $(m_1, \cdots , m_\mu)$, one has 
$$\bound f(m_1, \cdots , m_\mu)+\sum_{1\le t<\mu}\,
\sum_{m_t<k<m_{t+1}}
(-1)^{\,m_\mu+\mu+k+t}\vphi_{K^{m_k}}(f(m_1, \cdots ,m_t, k, m_{t+1}, \cdots,  m_\mu)\, )=0\,\,.$$
Here $\partial$ is the differential of the complex $F(K^{m_1}, \cdots, K^{m_\mu})$.
We  call an object of $\cC^\Delta$ a {\it $C$-diagram} in $\cC$.

In subsequent sections, 
for  $C$-diagrams $K_1, \cdots , K_n$ we will define complexes of abelian groups
$\BF(K_1, \cdots , K_n)$ together with maps $\sigma_{K_i}$  and $\vphi_{K_i}$. 
It will be shown that $\cC^\Delta$ forms a quasi $DG$ category, and that  
its homotopy category $Ho(\cC^\Delta)$
has the structure of a  triangulated category. 
\bigskip

\section{Function complexes $\BF (K_1, \cdots , K_n)$.}

In this section we keep the notation of \S 1. 
The operation $\Phi$ of (\Phicpx) is used in (\GKL).
A variant of the operation is discussed in (\htsPhicpx), which is needed in 
(\GISigma). 
In (\GISigma) we also refer to 
(\Tensorproductcpx) for tensor product of ``double" complexes.

Throughout this section, let $\cC$ be a quasi DG category having additional 
structure (iv) and (v) of (\DefqDG). 
We have defined the notion of $C$-diagrams in the 
category.  For a sequence of $C$-diagrams $K_1, \cdots, K_n$, 
we will define the complexes $\BF(K_1,\cdots, K_n)$ and  the maps 
$\vphi$ and $\sigma$ among them, and show that they satisfy the conditions for 
a quasi DG category (\DefqDG), except two of them that will be proven in \S 3. 
\bigskip

\sss{\secfcn.1} In this section 
a {\it sequence} is a pair $(M|M')$ consisting of a finite
 increasing sequence 
of integers $M=(m_1, \cdots, m_\mu )$ where
$m_1<\cdots <m_\mu$ with $\mu\ge 2$, and a subset $M'$
of $ M-\{m_1, m_\mu\}$. We allow $M'$ to be empty. 
For simplicity we also use the notation $\BM$ for $(M|M')$. 
When there is no confusion denote $(M|\emptyset)$ by $M$. 

Let $\init (M)=m_1$, $\term (M)=m_\mu$, $\ctop{M}=M-\{m_1, m_\mu\}$. 
Let $|M|=\mu$. For sequences $(M|M')$ and $(N|N')$
with  $\term (M)=\init(N)$, let 
$$\BM\scirc \BN= (M\cup N|M'\cup \{\term(M)\}\cup N')\,\,.$$

A {\it double sequence} a quadruple $(M_1|M'_1; M_2|M'_2)$. 
Here  $M_1$, $M_2$ are  finite sequences of integers, each of cardinality 
$\ge 1$, and 
$M'_1$ and $M'_2$ are subsets of $M_1-\{\init(M_1)\}$, $ M_2-\{\term (M_2)\}$, 
respectively. 
A double sequence may be 
viewed as a map defined on $[1, 2]$, which sends $i$ to  $\BM_i:=(M_i|M'_i)$.
To be specific, we will say it is a  double sequence on the set $[1, 2]$. 
(Note however that $\BM_i$ is not a sequence in the sense just defined, 
since $M_i$ may have cardinality one, and even if $|M_1|\ge 2$, $M'_1$ may 
contain $\term(M_1)$.)
We also use a single letter $A$ to denote a double sequence,  
$$A=(\BM_1; \BM_2)=(M_1|M'_1; M_2|M'_2)\,\,.$$

The following figure  on the left illustrates a sequence, where the line segment is $[m_1, m_\mu]$, the dots  indicate the set $M$ and the encircled dots the subset $M'$. 

\vspace*{0.5cm}
\hspace*{0.5cm}
\unitlength 0.1in
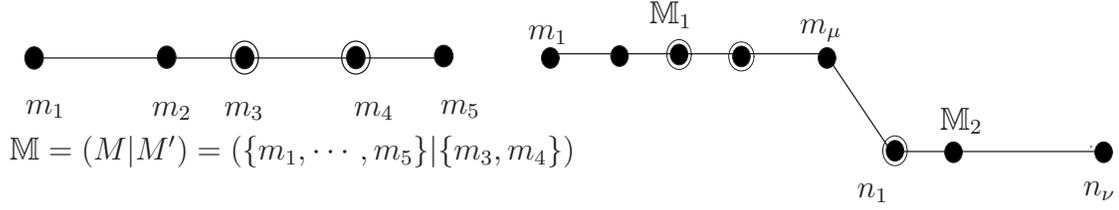
\begin{picture}( 57.7000,  9.1300)(  0.0000,-12.3300)
%
{\color[named]{Black}{%
\special{pn 8}%
\special{pa 136 608}%
\special{pa 2272 608}%
\special{fp}%
}}%
%
{\color[named]{Black}{%
\special{pn 0}%
\special{sh 1.000}%
\special{ia 822 604 48 58  0.0000000 6.2831853}%
}}%
{\color[named]{Black}{%
\special{pn 8}%
\special{ar 822 604 48 58  0.0000000 6.2831853}%
}}%
%
{\color[named]{Black}{%
\special{pn 0}%
\special{sh 1.000}%
\special{ia 128 608 48 58  0.0000000 6.2831853}%
}}%
{\color[named]{Black}{%
\special{pn 8}%
\special{ar 128 608 48 58  0.0000000 6.2831853}%
}}%
%
{\color[named]{Black}{%
\special{pn 0}%
\special{sh 1.000}%
\special{ia 2272 598 48 56  0.0000000 6.2831853}%
}}%
{\color[named]{Black}{%
\special{pn 8}%
\special{ar 2272 598 48 56  0.0000000 6.2831853}%
}}%
\put(0.8300,-9.2300){\makebox(0,0)[lb]{$m_1$}}%
\put(7.4500,-9.3400){\makebox(0,0)[lb]{$m_2$}}%
\put(11.2400,-9.3400){\makebox(0,0)[lb]{$m_3$}}%
\put(18.0400,-9.3400){\makebox(0,0)[lb]{$m_4$}}%
\put(22.5500,-9.2300){\makebox(0,0)[lb]{$m_5$}}%
\put(0.0000,-11.8000){\makebox(0,0)[lb]{$\BM=(M|M')=(\{m_1, \cdots, m_5\}|\{m_3, m_4\} )$}}%
%
{\color[named]{Black}{%
\special{pn 0}%
\special{sh 1.000}%
\special{ia 1230 606 48 58  0.0000000 6.2831853}%
}}%
{\color[named]{Black}{%
\special{pn 8}%
\special{ar 1230 606 48 58  0.0000000 6.2831853}%
}}%
%
{\color[named]{Black}{%
\special{pn 8}%
\special{ar 1230 606 72 86  0.0000000 6.2831853}%
}}%
%
{\color[named]{Black}{%
\special{pn 0}%
\special{sh 1.000}%
\special{ia 1812 606 48 58  0.0000000 6.2831853}%
}}%
{\color[named]{Black}{%
\special{pn 8}%
\special{ar 1812 606 48 58  0.0000000 6.2831853}%
}}%
%
{\color[named]{Black}{%
\special{pn 8}%
\special{ar 1812 606 74 86  0.0000000 6.2831853}%
}}%
%
{\color[named]{Black}{%
\special{pn 8}%
\special{pa 2820 586}%
\special{pa 4278 586}%
\special{fp}%
}}%
%
{\color[named]{Black}{%
\special{pn 8}%
\special{pa 4266 580}%
\special{pa 4624 1112}%
\special{fp}%
}}%
%
{\color[named]{Black}{%
\special{pn 8}%
\special{pa 5658 1072}%
\special{pa 5658 1072}%
\special{fp}%
}}%
%
{\color[named]{Black}{%
\special{pn 8}%
\special{pa 4640 1098}%
\special{pa 5736 1098}%
\special{fp}%
}}%
\put(33.5000,-4.5000){\makebox(0,0)[lb]{$\BM_1$}}%
\put(48.6900,-10.0600){\makebox(0,0)[lb]{$\BM_2$}}%
\put(41.3600,-5.2300){\makebox(0,0)[lb]{$m_\mu$}}%
\put(44.3600,-13.6300){\makebox(0,0)[lb]{$n_1$}}%
\put(27.1400,-5.3400){\makebox(0,0)[lb]{$m_1$}}%
\put(56.1900,-13.5300){\makebox(0,0)[lb]{$n_\nu$}}%
%
{\color[named]{Black}{%
\special{pn 0}%
\special{sh 1.000}%
\special{ia 2830 608 46 54  0.0000000 6.2831853}%
}}%
{\color[named]{Black}{%
\special{pn 4}%
\special{ar 2830 608 46 54  0.0000000 6.2831853}%
}}%
%
{\color[named]{Black}{%
\special{pn 0}%
\special{sh 1.000}%
\special{ia 4278 608 46 54  0.0000000 6.2831853}%
}}%
{\color[named]{Black}{%
\special{pn 4}%
\special{ar 4278 608 46 54  0.0000000 6.2831853}%
}}%
%
{\color[named]{Black}{%
\special{pn 0}%
\special{sh 1.000}%
\special{ia 3192 598 46 54  0.0000000 6.2831853}%
}}%
{\color[named]{Black}{%
\special{pn 4}%
\special{ar 3192 598 46 54  0.0000000 6.2831853}%
}}%
%
{\color[named]{Black}{%
\special{pn 0}%
\special{sh 1.000}%
\special{ia 5726 1102 46 54  0.0000000 6.2831853}%
}}%
{\color[named]{Black}{%
\special{pn 4}%
\special{ar 5726 1102 46 54  0.0000000 6.2831853}%
}}%
%
{\color[named]{Black}{%
\special{pn 0}%
\special{sh 1.000}%
\special{ia 4940 1102 46 54  0.0000000 6.2831853}%
}}%
{\color[named]{Black}{%
\special{pn 4}%
\special{ar 4940 1102 46 54  0.0000000 6.2831853}%
}}%
%
{\color[named]{Black}{%
\special{pn 0}%
\special{sh 1.000}%
\special{ia 3504 590 46 54  0.0000000 6.2831853}%
}}%
{\color[named]{Black}{%
\special{pn 4}%
\special{ar 3504 590 46 54  0.0000000 6.2831853}%
}}%
%
{\color[named]{Black}{%
\special{pn 4}%
\special{ar 3504 590 64 78  0.0000000 6.2831853}%
}}%
%
{\color[named]{Black}{%
\special{pn 0}%
\special{sh 1.000}%
\special{ia 3830 600 46 54  0.0000000 6.2831853}%
}}%
{\color[named]{Black}{%
\special{pn 4}%
\special{ar 3830 600 46 54  0.0000000 6.2831853}%
}}%
%
{\color[named]{Black}{%
\special{pn 4}%
\special{ar 3830 600 66 76  0.0000000 6.2831853}%
}}%
%
{\color[named]{Black}{%
\special{pn 0}%
\special{sh 1.000}%
\special{ia 4634 1094 46 54  0.0000000 6.2831853}%
}}%
{\color[named]{Black}{%
\special{pn 4}%
\special{ar 4634 1094 46 54  0.0000000 6.2831853}%
}}%
%
{\color[named]{Black}{%
\special{pn 4}%
\special{ar 4634 1094 64 78  0.0000000 6.2831853}%
}}%
\end{picture}%
\vspace{0.8cm}

\noindent The  figure on the right illustrates a double sequence. In the first line lies
$\BM_1$ which is a line segment with solid  and encircled dots, and in the second lies $\BM_2=(M_2=\{n_1, \cdots, n_\nu\}|M_2')$. (In the figure $m_\mu<n_1$, but there
are also cases $m_\mu=n_1$ and $m_\mu>n_1$.)
\bigskip

\sss{\secfcn.2} {\it The complex $F(M|M')$. }\quad 
Let $(K^m)$ be a sequence of objects in $\cC$
indexed by integers $m$, all but a finite
number of them being zero.  To a sequence $(M|M')$, one can 
associate the complex
$$F(\BM)=F(M|M'):= F(K^{m_1}, \cdots, K^{m_\mu}|M')\,\,.$$
If $\BM=(M|\emptyset)$, we simply write $F(M)=F(M|\emptyset)$;
in general, if $M_1, \cdots , M_r$ is the segmentation of $M$ given by $M'$, 
$F(M|M')=F(M_1)\hts\cdots\hts F(M_r)\subset F(M_1)\ts\cdots\ts F(M_r)$. 
In this section the differential of $F(\BM)$ is denoted $\partial$. 
For $k\in \ctop{M}-M'$ there is the corresponding map 
of complexes $\vphi_k: F(M|M')\to F(M-\{k\}|M')$. 
There is also the map $\sigma_k: F(M|M')\to F(M|M'\cup\{k\})$. 
The maps $\vphi_k$ commute with each other, $\sigma_k$ commute with 
each other, and $\vphi_k$ and $\sigma_\ell$ commute with other. 
\bigskip 

\sss{\FMM} {\it The complex $\bop F(M|M')$.}\quad
 We will define the structure of a complex
  on $\bop F(M|M')$, the direct sum over 
all sequences $(M|M')$. 

For $M=(m_1, \cdots, m_\mu)$, let 
$$ \quad \gamma(M)=m_\mu-m_1-\mu+1\,\,.$$
If $M$ and $N$ are sequences with $\term{M}=\init{N}$, then 
$\gamma(M\cup N)=\gamma(M)+\gamma(N)$. 

If $k\in \ctop{M}$, let 
$M_{\le k}:= \{m_i\in M \mid m\le k\}$. 
In this section,  for an integer $d$, write $\{d\}:=(-1)^{d}$;
this is useful when $d$ is a complicated expression. It will not 
lead to a confusion since braces are used exclusively for this. 
For $u\in F(M|M')$, let $|u|=\deg u$ (the degree in $F(M|M')$). 

We will define a map
$\bpartial: F(M|M')\to F(M|M')$ of degree 1, and a map 
$\bvphi: \bop F(M|M')\to \bop F(M|M')$ of degree 0. 
They are obtained from $\partial$ and $\vphi$ by putting appropriate 
signs. 

For $u\in F(M|\emptyset)$ define $\bpartial(u):=\partial u$.
For  $k\in \ctop{M}$, define
$\bvphi_k: F(M|\emptyset)\to F(M-\{k\}|\emptyset)$ by 
$$\bvphi_k(u):=\{|u|+\gamma(M_{\le k})\}\vphi_k(u)\,\,. $$
Note that the maps $u\mapsto \bpartial(u)$ and $u\mapsto \bvphi_k(u)$ increase the number
$|u|+\gamma(M)$ by one. 
In general let $M_1, \cdots, M_r$ be the 
segmentation of $M$ corresponding to $M'$, so that 
$F(M|M')=F(M_1)\hts\cdots\hts F(M_r)$. 
For $u=u_1\ts\cdots\ts u_r\in F(M|M')$, define 
$$\bpartial(u)=\sum_i\{\sum_{j>i}(|u_j|+\gamma(u_j)\,)\}u_1\ts\cdots\ts(
\bpartial u_i)\ts\cdots\ts u_r\,\,.$$
Here $\gamma(u_i):=\gamma(M_i)$ if $u_i\in F(M_i)$. 
For  $k\in \ctop{M}-M'$, let 
$i$ be such that $k\in M_i$, and 
define  $\bvphi_k: F(M|M')\to F(M-\{k\}|M')$ by 
$$\bvphi_k(u)=\{\sum_{j>i}(|u_j|+\gamma(u_j)\,)\}u_1\ts\cdots\ts\bvphi_k (
u_i)\ts\cdots\ts u_r\,\,.$$
Now let 
$$\bvphi(u):=\sum_{k}\bvphi_k(u)\,\,,$$
the sum over $k\in \ctop{M}-M'$. 

One verifies the following equalities:
$$\bpartial\bpartial(u)=0,\quad \bvphi\bvphi(u)=0,\quad \bpartial\bvphi(u)+
\bvphi\bpartial(u)=0\,\,. \eqno{(\FMM.a)}$$
For $u\ts v\in F(M|M')$, 
where $u\in F(M_1)\hts\cdots\hts F(M_s)$, $v\in F(M_{s+1})\hts\cdots\hts 
F(M_r)$, 
$$\bpartial(u\ts v)=\{|v|+\gamma(v)\}\bpartial u\ts v + u\ts \bpartial v\,\,,
\quad
\bvphi(u\ts v)=\{|v|+\gamma(v)\}\bvphi (u)\ts v + u\ts \bvphi (v)\,\,.\eqno{
(\FMM.b)}$$
Indeed the equalities (\FMM.b) are obvious. 
From these we obtain:
$$\bpartial\bpartial(u)=(\bpartial\bpartial u)\ts v + u\ts (\bpartial\bpartial v)\,\,,$$
$$\bvphi\bvphi(u)=(\bvphi\bvphi u)\ts v + u\ts (\bvphi\bvphi v)\,\,,$$
$$(\bpartial\bvphi +\bvphi\bpartial)(u\ts v)
=(\bpartial\bvphi +\bvphi\bpartial)u\ts v + u\ts (\bpartial\bvphi +\bvphi\bpartial)v\,\,.$$
 Using them, one may assume $u\in F(M|
\emptyset)$ to prove  the three equalities (\FMM.a). 
The verification for $\bpartial$ is obvious. 
For the second equality of (1), one must show, for  $k, l\in \ctop{M}$ with $k<l$ and $u\in F(M|\emptyset)$, 
$$(\bvphi_k\bvphi_l+\bvphi_l\bvphi_k)u=0\,\,.$$
But 
$$\bvphi_k\bvphi_l(u)=\{|u|+\gamma(M_{\le l})\}\{|u|+\gamma(M_{\le k})\}\vphi_k\vphi_l(u)\,\,,$$
$$\bvphi_l\bvphi_k(u)=\{|u|+\gamma(M_{\le k})\}\{|u|+\gamma((M-\{k\})_{\le l})\,\}\vphi_l\vphi_k(u)\,\,,$$
which add up to zero, since $\gamma((M-\{k\})_{\le l})=\gamma(M_{\le l})+1$. 
The verification of the second equality is similar. 
\bigskip 

Let $\bop F(M|M')$ be the direct sum over all sequences $(M|M')$, and  
$$\delta:=\bpartial +\bvphi:  \bop F(M|M')\to \bop F(M|M')\,\,.$$
Define the {\it 
first degree} of  $u\in F(M|M')$ by
$$\deg_1 (u) =|u|+\gamma(M)\,\,.$$
Then $\delta$  increases the first degree by 1. We have the following 
proposition, so  $\bop F(M|M')$ is a complex with degree $\deg_1$ and differential $\delta$, and the differential is compatible with  tensor product. 
The next proposition is a restatement of the identities 
(\FMM.a) and (\FMM.b). 
\bigskip  

(\FMM.1) {\bf Proposition.}\quad 
(1) $\delta\delta=0$. 

(2) For $u\ts v\in F(M|M')$, 
where $u\in F(M_1)\hts\cdots\hts F(M_s)$, $v\in F(M_{s+1})\hts\cdots\hts 
F(M_r)$, one has 
$\delta(u\ts v)=\{|v|+\gamma(v)\}\delta u\ts v + u\ts \delta v\,\,.$
\bigskip 

If one fixes $M'$ and takes the sum over $(M|M')$ where only $M$ varies with condition 
$M'\subset \ctop{M}$, 
one still obtains a complex; taking then the sum over $M'$ gives the complex
discussed above. 

In particular, 
$\bop F(M)=\bop F(M|\emptyset)$ is a complex, which appears in the following 
subsection.
\bigskip 

\sss{\DefCdiagram}
In the complex $\bop F(M)$, 
an element $f=(f(M))\in \bop F(M)$ is of first degree 0 if
$\deg f(M)+\gamma(M)=0$. 
It satisfies $\delta(f)=0$ if  for each $M=(m_1, \cdots, m_\mu)$, 
$$\partial f(M)+\sum_k \bvphi_k(f(M\cup\{k\})\, )=0$$
where $k$ varies over the set $[\init{M}, \term{M}]-M$.
Concretely 
$$\partial f(m_1, \cdots, m_\mu)+\sum_{t=1}^{\mu-1}
\sum_{m_t<k<m_{t+1}} (-1)^{m_\mu+\mu+k+t}\vphi_k(f(m_1, \cdots, m_t, k, m_{t+1}, \cdots, m_\mu)=0\,\,.$$
We now restate the definition of a $C$-diagram.
\bigskip

(\DefCdiagram.1) {\bf Definition.}\quad A {\it $C$-diagram} $K=(K^m; f(M))$ in 
the quasi DG category 
is a sequence of objects $K^m$ indexed by $m\in \ZZ$, all but a finite number
 of them being zero, and a set of elements $f(M)\in F(M)^{-\gamma(M)}$ indexed
 by $M=(m_1, \cdots, m_\mu)$, 
satisfying the following conditions:

(i) For each $k\in \ctop{M}$, $\sigma_k(f(M))=f(M_{\le k})\ts f(M_{\ge k})$ 
in $F(M_{\le k})\ts F(M_{\ge k})$.
(To be precise one should write $\tau_k$ for $\sigma_k$, but we may not make the 
distinction.)

(ii) $f=(f(M)\,)\in \bop F(M)$ satisfies $\delta(f)=0$. 
\bigskip 

For an object $X$  and $n\in \ZZ$,  there is a $C$-diagram $K$ with 
$K^n=X$, $K^m=0$ if $m\neq n$, and $f(M)=0$ for all $M$.
We write $X[-n]$ for this. 
\bigskip 

\sss{\Defbsigma} {\it The differential $\bsigma$.}\quad
Under the same assumption, we define, for each 
$k\in \ctop{M}-M'$, the map 
$\bsigma_k: F(M|M')\to F(M|M'\cup\{k\})$ as follows. 
For $u\in F(M|\emptyset)$, if
$\sigma_k(u)=\sum u'\ts u''$
where $u'\in F(M_{\le k})$ and $u''\in F(M_{\ge k})$,
let
$$\bsigma_k(u)=\sum \{\deg_1(u')\cdot \gamma (u'')\}u'\ts u''\,\,.$$
One has $\bsigma_k(u)\in F(M|\{k\})$. 
Indeed, since $F(M|\{k\})$ is a ``double" subcomplex of 
$F(M_{\le k})\ts F(M_{\ge k})$, 
writing $\sigma_k(u)=\sum u^{a, b}$, the sum of elements of bidegree 
$(a, b)$, each $u^{a, b}$ is in $F(M|\{k\})$; hence 
$\sum \{a\cdot \gamma(M_{\ge k})\} u^{a, b}$ is also in $F(M|\{k\})$. 

In general, let $M_1, \cdots, M_r$ be the segmentation of $M$ by $M'$. 
For $u=\sum u_1\ts u_2\ts\cdots\ts u_r\in F(M|M')$ with $u_i\in F(M_i)$, 
and for $k\in M_i-M'$, set 
$$\bsigma_k(u)=(-1)^{|{M'}_{>k}|} u_1\ts\cdots\ts\bsigma_k(u_i)\ts\cdots\ts u_r\,;$$
then we let  
$$\bsigma(u):= \sum_k \bsigma_k(u) \,\,$$ 
where $k$ varies over the set $\ctop{M}-M'$.  
(Here ${M'}_{>k}$ denotes the subset of $M'$ of elements $>k$.)
For $u\in F(M|M')$, define $\tau(u):=|M'|$.
Note $\bsigma$ increases $\tau(u)$ by one.
\bigskip

(\Defbsigma.1)
{\bf Proposition.}\quad (1) $\bsigma\bsigma(u)=0$. 

(2) $\delta\bsigma(u)=\bsigma\delta(u)$. 

(3) $\bsigma(u\ts v)=\{\tau(v)+1\}\bsigma(u)\ts v + u\ts \bsigma(v
)\,\,.$
\smallskip

{\it Proof.}\quad
(3) is obvious from the definitions. For (1), using (3) one may thus assume 
$u\in F(M|\emptyset)$. 
If $u\in F(M|\emptyset)$, for  $k, l\in \ctop{M}$ with $k<l$, write 
$$\sigma_k\sigma_l(u)=\sum u'\ts u''\ts u'''$$
with $u'\in F(M_{\le k})$, $u''\in F(M_{\ge k}\cap M_{\le l})$, and $u'''\in F(M_{\ge l})$. 
Then we have
$$\bsigma_k\bsigma_l (u)=\{(|u'|+|u''|+ \gamma(u')+\gamma(u'')\,)\cdot \gamma(u''')\}
\{1+(|u'|+\gamma(u'))\cdot\gamma(u'')\} \sigma_k\sigma_l(u)\,\,,$$
$$\bsigma_l\bsigma_k(u)=\{(|u'|+\gamma(u'))\cdot(\gamma(u'')+\gamma(u'''))\}
\{(|u''|+\gamma(u''))\cdot\gamma(u''')\} \sigma_l\sigma_k(u)\,\,.$$
The sum of them is zero. 

For (2), using (3) and (\FMM.1), (2), one may again assume $u\in F(M|\emptyset)$. 
We will then show 
$$\bpartial \bsigma_k (u)=\bsigma_k\bpartial (u)$$
and 
$$\bsigma_k\bvphi_l(u)= \bvphi_l \bsigma_k(u)$$
for $k\neq l$. 
For the first equality, writing  $\sigma_k(u)=\sum u'\ts u''$, 
one has 
$$\sigma_k(\partial u)=\partial \sigma_k (u)=\sum\{|u''|\} \partial u'\ts u''
+\sum u'\ts (\partial u'')$$
since $\sigma_k$ and $\partial$ commute. 
So 
$$\bsigma_k(\partial u)=\sum\{|u''|+(|\partial u'|+\gamma(u'))\cdot\gamma(u'')\}(\partial u')\ts u''
+\sum \{|u'|+\gamma(u'))\cdot \gamma(u'')\} u'\ts(\partial u'')\,\,.$$
On the other hand, 
$$\bpartial \bsigma_k (u)=\{ (|u'|+\gamma(u'))\cdot \gamma(u'')\} \cdot \bpartial (u'\ts u'')\,\,,$$
which coincides with the above. 
The proof of the second equality is similar: write 
$\sigma_k(u)=\sum u'\ts u''\ts u'''$ as in the proof of (\FMM.1), and 
compute $\bsigma_k\bvphi_l(u)$ and $\bvphi_l\bsigma_k(u)$. 
\bigskip 

\sss{\GKL} {\it The complexes $\BG(K, L)$ and $\BF(K, L)$.}\quad 
Let $K=(K^m; f_K(M))$ and $L$ $=(L^m; f_L(M))$ be $C$-diagrams. 
To a double sequence $A=(M|M'; N|N')$ one associates the complex 
$$F(A)=F(M|M'; N|N'):= F(K^{m_1}, \cdots, K^{m_\mu}; L^{n_1}, \cdots, L^{n_\nu}|
M'\cup N')\,\,.$$
To be precise, consider the finite ordered set $M\amalg N$
(where we give the ordering $m<n$ if $m\in M$ and $n\in N$),  the sequence of objects 
$(K^m)_{m\in M}$, $(L^n)_{n\in N}$ 
on it, and the subset $M'\amalg N'\subset M\amalg N$. 
 Then  the corresponding complex is the $F(A)$. 

If $M_1, \cdots, M_r$ is the segmentation of $M$ given by $M'$, and 
$N_1, \cdots, N_s$ that of $N$ given by $N'$, then 
$$F(M|M'; N|N')=F(M_1)\hts\cdots\hts F(M_{r}\cup N_1)\hts F(N_2)\hts\cdots
\hts F(N_s)\,\,.$$
We refer to $M_1, \cdots, M_{r-1}, M_r\cup N_1, N_2, \cdots, N_s$ as 
the segmentation of $M\cup N$ by $M'\cup N'$. 

Let us say the double sequence is {\it free} when $M'$ and $N'$ are empty; then 
the corresponding complex is free of tensor products. 

As for $F(M|M')$, one has maps $\partial, \vphi$ and $\sigma$ among 
the $F(A)$. 
Recall for $M=(m_1, \cdots, m_\mu)$, $\gamma(M):=m_\mu-m_1-\mu+1$. 
For $A=(\BM; \BN)$, set  
$|A|=|M|+|N|$, 
$$\gamma(A)=\gamma(M)+\gamma(N)+(n_1-m_\mu)$$
 and 
$\tau(A)=|M'|+|N'|$. 
Note $|A|$ and $\gamma(A)$ depend only on $(M; N)$, while 
$\tau(A)$ depends only on $(M'; N')$. 
The following obvious equalities will be repeatedly used. 
If $k\in A$ with  $k\not\in \{m_1, n_\nu\}\cup M'\cup N'$, then 
let $A-\{k\}$ be the double sequence obtained by removing $k$ from $M$ or $N$, 
and keeping $M'$ and $N'$.
Then 
$$\gamma(A-\{k\})=\gamma(A)+1\,\,.$$
Also, 
$$\gamma(A)=\gamma(A_{\le k})+ \gamma(A_{\ge k})$$
where $A_{\le k}$, for instance, is the (double or single) sequence consisting 
of elements $a\in A$  with $a\le k$. 

For $u\in F(A)$, let $\gamma(u)=\gamma(A)$ and  $\tau(u)=\tau(A)$. 
One can then define the maps $\bpartial$ and $\bvphi$ as before, 
as well as the sum 
$$\delta=\bpartial+\bvphi: \bop F(A)\to \bop F(A)\,\,.$$
Specifically if $u\in F(A)$ with $A=(M|\emptyset; N|\emptyset)$, 
then $\bpartial (u)=\partial (u)$. If 
$u=u_1\ts\cdots\ts u_r\ts u_{r+1}\ts\cdots\ts u_{r+s-1}
\in F(M_1)\hts\cdots\hts F(N_s)$ as above, then 
$$\bpartial (u)=\sum_i
\{\sum_{j>i} (|u_j|+\gamma(u_j)\,)\}u_1\ts\cdots (\bpartial u_i)
\ts\cdots \ts u_{r+s-1}\,\,.$$
Similarly for $\bvphi$. 

One also has the map 
$\bsigma: \bop F(A)\to \bop F(A)$. These maps satisfy the same identities as
before.  In addition, we will define  maps $\Bf_K$ and $\Bf_L$. 

For this purpose one needs to invoke (\DefqDG), (v)  
and take appropriate quasi-isomorphic subcomplexes. 
There is a distinguished subcomplex 
$[F(M|M'; N|N')]_f\injto F(M|M'; N|N')$ satisfying the following conditions:

\quad$\bullet$\quad If $P=(P|\emptyset)$ is a sequence with 
$\term(P)=\init(M)$, then the map 
$$f_K(P)\ts (-): [F(\BM; \BN)]_f\to [F(P\scirc \BM; \BN)]_f$$
is defined. 
Here $P\scirc \BM:=(P\cup M| \{\term(M)\}\cup M')$. 

\vspace*{0.5cm}
\hspace*{2cm}
\unitlength 0.1in
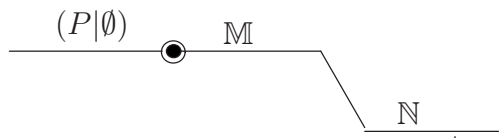
\begin{picture}( 25.6000,  6.3000)(  6.4000,-13.9000)
%
{\color[named]{Black}{%
\special{pn 8}%
\special{pa 1470 940}%
\special{pa 2270 940}%
\special{fp}%
}}%
%
{\color[named]{Black}{%
\special{pn 8}%
\special{pa 2270 940}%
\special{pa 2500 1340}%
\special{fp}%
}}%
%
{\color[named]{Black}{%
\special{pn 8}%
\special{pa 2970 1390}%
\special{pa 2970 1390}%
\special{fp}%
}}%
%
{\color[named]{Black}{%
\special{pn 8}%
\special{pa 2500 1360}%
\special{pa 3200 1360}%
\special{fp}%
}}%
\put(17.6000,-9.0000){\makebox(0,0)[lb]{$\BM$}}%
\put(26.7000,-13.2000){\makebox(0,0)[lb]{$\BN$}}%
%
{\color[named]{Black}{%
\special{pn 0}%
\special{sh 1.000}%
\special{ia 1500 940 38 34  0.0000000 6.2831853}%
}}%
{\color[named]{Black}{%
\special{pn 8}%
\special{ar 1500 940 38 34  0.0000000 6.2831853}%
}}%
%
{\color[named]{Black}{%
\special{pn 8}%
\special{pa 640 940}%
\special{pa 1440 940}%
\special{fp}%
}}%
\put(8.8000,-8.9000){\makebox(0,0)[lb]{$(P|\emptyset)$}}%
%
{\color[named]{Black}{%
\special{pn 8}%
\special{ar 1498 944 62 58  0.0000000 6.2831853}%
}}%
\end{picture}%
\vspace{0.8cm}

\quad$\bullet$\quad Similarly if $\term(N)=\init(Q)$, then
$$(-)\ts f_L(Q): [F(\BM; \BN)]_f\to [F(\BM; \BN\scirc Q)]_f$$
is defined. 

\quad$\bullet$\quad $[F(\BM; \BN)]_f$ is closed under the maps 
$\vphi$ and $\sigma$. 
\smallskip 

 For the existence of such a subcomplex, with reference to 
\Propintersection, we proceed as follows. 
Let $I=M\amalg N$, and 
$I_1, \cdots, I_c$ be its segmentation by $M'\amalg N'$
(same as $M_1, \cdots, M_{r-1}, M_r\cup N_1, N_2, \cdots, N_s$ at the beginning 
of this subsection). 
To specify a condition of 
constraint $\cC$ on this, let 
$$\BI:=[-w, \term(M)]\amalg [\init(N), w]$$
where $w$ is a positive integer large enough so that 
$K^m=O$ for $m<-w$ and $L^m=O$ for $m>w$. 
As a subset of $[1, c]$, we take  $P=[1, c]$. 
Choose a set of almost disjoint sub-intervals of $\BI$ such that
for each $j$, one has either $J_j\subset [-w, \init (M)]$ or
$J_j\subset [\term(N), w]$. 
Set 
$$f(J_j)=\left\{
\begin{array}{rl}
f_K(J_j) &\quad \mbox{if $J_j\subset [-w, \init (M)]$}\,,\\
f_L(J_j)&\quad \mbox{if $J_j\subset [\term(N), w]$}\,.
\end{array}
\right. 
$$
These give a constraint $\cC=(\BI; P; \{J_j\}; \{f(J_j)\})$. 
Thus the corresponding distinguished subcomplex is generated by 
$u_1\ts\cdots\ts u_c$ with $u_i\in F(I_i)$ such that 
$$\{u_1, \cdots, u_c,  \{f(J_j)\}\,\}\qquad \mbox{is properly intersecting.}
\eqno{(\GKL.a)}$$
We take all possible sets of sub-intervals $\{J_j\}$, and take the 
intersection of the corresponding distinguished subcomplexes;
this gives a subcomplex satisfying the required conditions, 
as can be shown using functorial properties of proper intersection,
(\RmkDefqDG.2).
In short, we could say that we
take $\{f_K(P)\}$ and $\{f_L(Q)\}$ as the set of constraints. 
In the rest of this paper we write $F(\BM; \BN)$ for $[F(\BM; \BN)]_f$
as long as no confusion is likely.  

For a sequence $P=(P|\emptyset)$ with $\term(P)=\init(M)$ and 
$u\in F(\BM; \BN)$, let 
$$\Bf_K(P)\ts u=\{|\tau(u)|+1\}f_K(P)\ts u\in F(P\scirc \BM; \BN)\,\,.$$
Let 
$\Bf_K\ts u:=\sum \Bf_K(P)\ts u$ where $P$ varies over sequences with 
$\tm(P)=\init (M)$. When there is no confusion, we also write 
$\Bf_K(u)=\Bf_K\ts u$ like an operator. 

If $\init(Q)=\tm (N)$ let 
$$u\ts \Bf_L(Q)=- u\ts f_L(Q)\in F(\BM; \BN\scirc Q)\,\,,$$
and $u\ts \Bf_L=\sum u\ts\Bf_L(Q)$, the sum over $Q$ with $\init (Q)=\tm (N)$. 
We also write $\Bf_L(u)=u\ts \Bf_L$.
These operations are subject to the following identities. 
\bigskip 

(\GKL.1)
\Prop{\it 
(1) $\delta\Bf_K=\Bf_K\delta$. $\delta\Bf_L=\Bf_L\delta$.

(2) One has 
$$\bsigma(\Bf_K\ts u)+\Bf_K\ts(\Bf_K\ts u )+ \Bf_K\ts \bsigma(u)=0\,\,.$$
Similarly 
$\bsigma(u\ts \Bf_L)+(u\ts \Bf_L)\ts\Bf_L +  \bsigma(u)\ts \Bf_L=0$.
(With the operator-like notation, $\bsigma\Bf_K+\Bf_K\bsigma+\Bf_K\Bf_K=0$, 
etc. )

(3) One has 
$\Bf_K\ts (u\ts v)=\{\tau(v)+1\}(\Bf_K\ts u)\ts v$. 
(Here  $u\in F(\BM)=F(K^{\BM})$, $v\in F(\BM'; \BN)$ with $\term(M)=\init(M')$; or either,
$u\in F(\BM; \BN)$, $v\in F(\BN')=F(L^{\BN'})$ with $\term(N)=\init(N')$. )
\newline
Similarly 
$(u\ts v)\ts\Bf_L=u\ts(v\ts\Bf_L)$. 

(4) $(\Bf_K\ts u)\ts \Bf_L+ \Bf_K\ts (u\ts \Bf_L)=0$.
(With the operator-like notation, $\Bf_K\Bf_L+\Bf_L\Bf_K=0$. )
}\bigskip 

{\it Proof.}\quad 
(1) We show $\delta\Bf_K=\Bf_K\delta$. 
Since $\delta(f)=0$, by (\FMM.1), (2), 
$$\delta(\Bf\ts u)=\{\tau(u)+1\}f\ts \delta(u)\,\,.$$
Noting $\tau(\delta(u))=\tau(u)$, one has 
$$\Bf\ts \delta(u)=\{\tau (u)+1\} f\ts \delta(u)$$
as well. 
The proof of $\delta\Bf_L=\Bf_L\delta$ is similar. 

(2) By $|f(P)|+\gamma(P)=0$, one has
$$\bsigma(f(P))=\sum f(P_1)\ts f(P_2)\,,$$
the sum over all segmentations of $P$ to $P_1$ and $P_2$. 
So, using (\Defbsigma.1), 
\begin{eqnarray*}
&&\bsigma(\Bf (P)\ts u)\\
&=& \{\tau(u)+1\} (\{\tau(u)+1\}\bsigma(f(P))\ts u
+f(P)\ts \bsigma(u) \,)\\
&=& \sum f(P_1)\ts f(P_2)\ts u
+\{\tau(u)+1\} f(P)\ts \bsigma(u)\,\,.
\end{eqnarray*}
Also, 
$$\Bf (P_1)\ts \Bf(P_2)\ts u= -f(P_1)\ts f(P_2)\ts u\,,$$
and 
$$\Bf(P)\ts \bsigma(u)=\{\tau(\sigma(u))+1\} f(P)\ts \bsigma(u)\,,$$
where $\tau(\sigma(u))= \tau (u)+1$. 
The equality hence follows. 

(3) Obvious from the  definitions. 

(4) One has 
$$(\Bf_K \ts u)\ts \Bf_L=-\sum \{\tau (u)+1\}(f_K(P)\ts u) \ts f_L(Q)\,,$$
and 
$$\Bf_K\ts (u\ts \Bf_L) =-\sum \{\tau(u)\} f_K(P)\ts (u\ts f_L(Q))\,,$$
which add to zero. 
\bigskip 

For $u\in F(A)$ define 
$$d'(u)=\bsigma(u)+\Bf_K(u)+\Bf_L(u) \,\,.$$
It increases $\tau(u)$ by one,  and by (\Defbsigma.1) and (\GKL.1), we have
$d'd'=0$  and $\delta d'=d' \delta$. 
In addition,  $d'(u\ts v)=\{\tau(v)+1\}(d'u)\ts v + u\ts (d'v)$. 
Thus the direct sum $\bop_{A}F(A)$ is a ``double" complex $H^{a, b}$
with respect to 
the two gradings 
$$a=\deg_1(u) +1, \quad b=\tau(u)\,, $$
with $\deg_1(u):=|u|+\gamma(A)$, 
and the commuting differentials $-\delta$, $d'$. 
Denote it by $H^{\bullet\, \bullet}(K, L)$:
$$H^{\bullet\, \bullet}(K, L)=\bop_A F(A)=\bop_A  [F(A)]_f\,.$$
Clearly $H^{a, b}=0$ unless $b\ge 0$. 

Let 
$$\BH(K, L)=\Tot( H^{\bullet\, \bullet}(K, L)   \, )$$
 be the associated total  complex. 
The total degree is  of $u\in F(A)$ is given by 
 $\deg_{\BH} (u)= \deg_1(u)+\tau(u)
+1$, and 
the total differential $d_{\BH}$ is given by
$$d_{\BH}(u)=-\delta(u)+ (-1)^{\deg_1(u)+1}d'(u)$$
 on $u\in F(A)$.  
\bigskip 

(\GKL.2)
\Prop{\it 
 For each $a$, the complex with respect to $d'$, 
$$H^{a, 0}\to H^{a, 1}\to \cdots$$
 is exact. }
 \smallskip 
 
{\it Proof.}\quad 
On the complex $(H^{a, \bullet}, d')$, consider the filtration given by the sum of terms $F(A)$
with $\init{M}\le p $ and $\term{N}\ge q$, for 
varying $p,q$; in the subquotients  the maps $\Bf_K$ and $\Bf_L$ are zero, 
so the differentials $d'$ are just $\bsigma$. So they are sums of the complexes 
$$F(\cI|\emptyset)\mapr{\sigma}
\moplus_{|S|=1} F(\cI|S)
\mapr{\sigma} \moplus_{|S|=2} F(\cI | S)
\mapr{\sigma} \cdots \to  F(\cI|\ctop{\cI})\to 0\,\,$$
where $\cI= M\amalg N$ with $\init{M}= p $ and $\term{N}= q$ given,  and 
$S$ varies over subsets of $\ctop{\cI}$.
The claim follows from the acyclicity axiom of $\sigma$, \Acyclicitysigma.
\bigskip 

Applying the operation $\Phi$ in (\Phicpx) 
we obtain a complex  
$$\BG^{\bullet}(K, L):=\Phi  H^{\bullet\, \bullet}(K, L)=\Ker (d': H^{\bullet\, 0}\to H^{\bullet\, 1})\,\,,$$
and the differential is the restriction of $d_{\BH}$ 
to $\Phi  H^{\bullet\, \bullet}(K, L)$. 
So the degree and differential are given by $\deg_\BG(u)= \deg_1(u)+1$ and 
$d_\BG=-\delta$ on $u\in F(A)$. 
Set finally $$\BF(K, L)=\BG(K, L)[1]\,\,.$$
The degree and the differential are given as follows:
$$\deg_{\BF}(u)=\deg_1 (u)\,\,,\quad d_{\BF}(u)=\delta(u). $$
The following facts are obvious from the definitions. 
\bigskip  

(\GKL.3)\quad 
An element $u\in H^{a, 0}$  consists of 
$u(M; N)\in F(K^M; L^N)$ with $\deg_1 u(M; N)=a-1$, where 
$(M; N)$ varies over free double sequences. 
It is  in  $\BG(K, L)^a$ 
if the following condition  is satisfied.

(i) For $k\in M$, $k\neq \init(M)$, 
$$\bsigma_k(u(M; N))=f_K(M_{\le k})\ts u(M_{\ge k}; N)\,\,, $$
equivalently,
$$\sigma_k(u(M; N))=f_K(M_{\le k})\ts u(M_{\ge k}; N)\,\,.$$

(ii) For
$k\in N$, $k\neq \term(N)$, 
$$\bsigma_k(u(M; N)\,)= 
u(M; N_{\le k})\ts f_L(N_{\ge k})\,, $$
equivalently, 
$$\sigma_k(u(M; N)\,)=\{(a-1)\cdot\gamma(N_{\ge k})\}\, 
u(M; N_{\le k})\ts f_L(N_{\ge k})\,.$$
\bigskip 

(\GKL.4) {\bf Example.}\quad 
If $X$, $Y$ are objects in $\cC$ and $K=X[0]$ and $L=Y[n]$ (see (2.4)),  we have 
$\BF(K, L)=F(X, Y)[n]$. 
\bigskip

\sss{\HIGI} {\it The complexes  $\BH(I)$ and $\BG(I)$.}\quad
Let $n\ge 2$.  Assume given a sequence of $C$-diagrams 
$K_i=(K_i^m; f_{K_i}(M))$
for $i=1, \cdots, n$. We will define complexes 
$\BH(K_1, \cdots, K_n)$ and $\BG(K_1, \cdots, K_n)$, 
generalizing $\BH(K, L)$ and $\BG(K, L)$ in case $n=2$. 
As in the case $|I|=2$, $\BG$ is a quasi-isomorphic subcomplex of $\BH$. 
The generalization of $\BF(K, L)$ to $\BF(K_1, \cdots, K_n)$
will be discussed later in this section.

A {\it multi-sequence} on $[1, n]$ is a $2n$-tuple of finite sequences 
$$A=(\BM_1; \cdots ; \BM_n)=(M_1|M'_1; M_2|M'_2; \cdots , M_{n-1}|M'_{n-1}; 
M_n|M'_n)$$ 
satisfying the following conditions: 
\smallskip 

\quad$\bullet$\quad  Each $M_i$ is non-empty;

\quad$\bullet$\quad
$M'_j\subset M_j$ for $j=1, \cdots, n$, and 
$\init(M_1)\not\in M'_1  $, 
and $ \term(M_n)\not\in  M'_n  $.  
\smallskip 

 The following illustrates a 
multi-sequence. The vertical direction is for  $[1, n]$, and 
the horizontal direction for the  $\BM_i$. 

\vspace*{0.5cm}
\hspace*{2cm}
\unitlength 0.1in
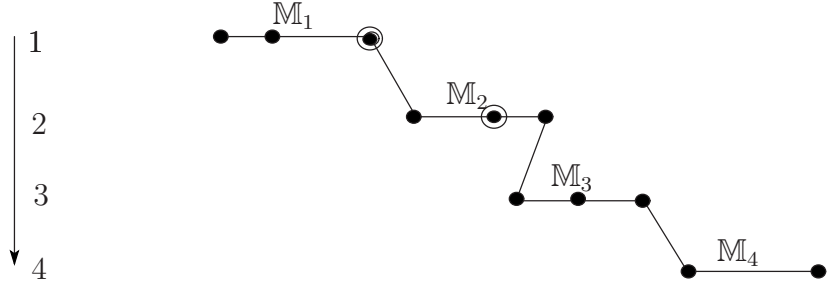
\begin{picture}( 42.7300, 14.3300)(  1.5300,-22.5300)
%
{\color[named]{Black}{%
\special{pn 8}%
\special{pa 1240 990}%
\special{pa 2040 990}%
\special{fp}%
}}%
%
{\color[named]{Black}{%
\special{pn 8}%
\special{pa 2040 990}%
\special{pa 2270 1390}%
\special{fp}%
}}%
%
{\color[named]{Black}{%
\special{pn 8}%
\special{pa 2970 1390}%
\special{pa 2970 1390}%
\special{fp}%
}}%
%
{\color[named]{Black}{%
\special{pn 8}%
\special{pa 2270 1410}%
\special{pa 2970 1410}%
\special{fp}%
}}%
%
{\color[named]{Black}{%
\special{pn 8}%
\special{pa 2970 1410}%
\special{pa 2810 1838}%
\special{fp}%
}}%
%
{\color[named]{Black}{%
\special{pn 8}%
\special{pa 2808 1850}%
\special{pa 3508 1850}%
\special{fp}%
}}%
%
{\color[named]{Black}{%
\special{pn 8}%
\special{pa 3468 1840}%
\special{pa 3708 2230}%
\special{fp}%
}}%
%
{\color[named]{Black}{%
\special{pn 8}%
\special{pa 3708 2220}%
\special{pa 4408 2220}%
\special{fp}%
}}%
\put(15.3000,-9.5000){\makebox(0,0)[lb]{$\BM_1$}}%
\put(24.4000,-13.7000){\makebox(0,0)[lb]{$\BM_2$}}%
\put(29.8800,-18.0000){\makebox(0,0)[lb]{$\BM_3$}}%
\put(38.5800,-21.9000){\makebox(0,0)[lb]{$\BM_4$}}%
%
{\color[named]{Black}{%
\special{pn 8}%
\special{pa 180 990}%
\special{pa 180 2180}%
\special{fp}%
\special{sh 1}%
\special{pa 180 2180}%
\special{pa 200 2114}%
\special{pa 180 2128}%
\special{pa 160 2114}%
\special{pa 180 2180}%
\special{fp}%
}}%
\put(2.5000,-10.7000){\makebox(0,0)[lb]{$1$}}%
\put(2.5000,-10.7000){\makebox(0,0)[lb]{$1$}}%
\put(2.7000,-15.0000){\makebox(0,0)[lb]{$2$}}%
\put(2.8000,-18.7000){\makebox(0,0)[lb]{$3$}}%
\put(2.7000,-22.6000){\makebox(0,0)[lb]{$4$}}%
%
{\color[named]{Black}{%
\special{pn 0}%
\special{sh 1.000}%
\special{ia 1260 990 38 34  0.0000000 6.2831853}%
}}%
{\color[named]{Black}{%
\special{pn 8}%
\special{ar 1260 990 38 34  0.0000000 6.2831853}%
}}%
%
{\color[named]{Black}{%
\special{pn 0}%
\special{sh 1.000}%
\special{ia 1530 990 38 34  0.0000000 6.2831853}%
}}%
{\color[named]{Black}{%
\special{pn 8}%
\special{ar 1530 990 38 34  0.0000000 6.2831853}%
}}%
%
{\color[named]{Black}{%
\special{pn 0}%
\special{sh 1.000}%
\special{ia 2270 1410 38 34  0.0000000 6.2831853}%
}}%
{\color[named]{Black}{%
\special{pn 8}%
\special{ar 2270 1410 38 34  0.0000000 6.2831853}%
}}%
%
{\color[named]{Black}{%
\special{pn 0}%
\special{sh 1.000}%
\special{ia 2960 1410 38 34  0.0000000 6.2831853}%
}}%
{\color[named]{Black}{%
\special{pn 8}%
\special{ar 2960 1410 38 34  0.0000000 6.2831853}%
}}%
%
{\color[named]{Black}{%
\special{pn 0}%
\special{sh 1.000}%
\special{ia 2808 1840 38 34  0.0000000 6.2831853}%
}}%
{\color[named]{Black}{%
\special{pn 8}%
\special{ar 2808 1840 38 34  0.0000000 6.2831853}%
}}%
%
{\color[named]{Black}{%
\special{pn 0}%
\special{sh 1.000}%
\special{ia 3468 1850 38 34  0.0000000 6.2831853}%
}}%
{\color[named]{Black}{%
\special{pn 8}%
\special{ar 3468 1850 38 34  0.0000000 6.2831853}%
}}%
%
{\color[named]{Black}{%
\special{pn 0}%
\special{sh 1.000}%
\special{ia 3708 2220 38 34  0.0000000 6.2831853}%
}}%
{\color[named]{Black}{%
\special{pn 8}%
\special{ar 3708 2220 38 34  0.0000000 6.2831853}%
}}%
%
{\color[named]{Black}{%
\special{pn 0}%
\special{sh 1.000}%
\special{ia 4388 2220 38 34  0.0000000 6.2831853}%
}}%
{\color[named]{Black}{%
\special{pn 8}%
\special{ar 4388 2220 38 34  0.0000000 6.2831853}%
}}%
%
{\color[named]{Black}{%
\special{pn 8}%
\special{ar 2050 1000 38 34  0.0000000 6.2831853}%
}}%
%
{\color[named]{Black}{%
\special{pn 0}%
\special{sh 1.000}%
\special{ia 3130 1840 38 34  0.0000000 6.2831853}%
}}%
{\color[named]{Black}{%
\special{pn 8}%
\special{ar 3130 1840 38 34  0.0000000 6.2831853}%
}}%
%
{\color[named]{Black}{%
\special{pn 0}%
\special{sh 1.000}%
\special{ia 2040 1004 36 30  0.0000000 6.2831853}%
}}%
{\color[named]{Black}{%
\special{pn 8}%
\special{ar 2040 1004 36 30  0.0000000 6.2831853}%
}}%
%
{\color[named]{Black}{%
\special{pn 8}%
\special{ar 2040 1000 66 60  0.0000000 6.2831853}%
}}%
%
{\color[named]{Black}{%
\special{pn 0}%
\special{sh 1.000}%
\special{ia 2690 1410 36 30  0.0000000 6.2831853}%
}}%
{\color[named]{Black}{%
\special{pn 8}%
\special{ar 2690 1410 36 30  0.0000000 6.2831853}%
}}%
%
{\color[named]{Black}{%
\special{pn 8}%
\special{ar 2690 1410 66 60  0.0000000 6.2831853}%
}}%
\end{picture}%
\vspace{0.8cm}

Associated to $A$ is a finite ordered set $M=M_1\cup\cdots\cup M_n$
(disjoint union) and $(K^m_i)$ defines a sequence of objects on it, 
$M_i\ni m\mapsto K^m_i$ on each $M_i$.
There 
corresponds the complex 
$$F(A)=F(\BM_1; \cdots ; \BM_n)=F(M_1\cup\cdots\cup M_n|M_1'\cup\cdots\cup M_n')\,\,.$$
An element $u\in F(A)$ can be written as a sum of 
$u_1\ts\cdots\ts u_r$ where $u_i\in F(A_i)$ and $(A_1, \cdots, A_r)$ is the segmentation
of $\amalg M_i$ by $\amalg M'_i$. 
Let $\partial$ be its differential. Set 
$|A|=\sum |M_i|$, 
$$\gamma(A)=\sum_{1\le i\le n}\gamma(M_i)+\sum_{1\le i<n} (\init(M_{i+1}) -\term(M_i)\,)\,\,,$$
and $\tau(A)=\sum |M'_i|$. 

Now consider the group
$$\bop_{A} F(A)$$
the direct sum over all multi-sequences $A$ on $[1, n]$.
We make it into a ``triple" complex, denoted $H\tbullet=H\tbullet(K_1, \cdots, K_n)$. 
It is analogous to the double complex $H\dbullet$ in the previous subsection.
Set
$$M'_{int}=\bigcup_{1<i<n} M'_i\,, \quad M'_{out}= M'_1\cup M'_n\,\,.$$
\smallskip 

\quad$\bullet$  The first degree (namely the number $a$ in $H^{a,b,c}$) 
is $\deg_1 (u)=|u|+\gamma(u)$, and the first differential is 
$d_1=-\delta$, where $\delta$ is to be defined as follows. 
 
As in the case $n=2$, one has the map $\bpartial$.
For each $k\in M_1-M_1'$ with $k\neq\init(M_1)$ (resp. 
$k\in M_n-M_n'$ with $k\neq\term(M_n)$), 
one has the map 
$$\bvphi_k: F( \BM_1; \cdots ; \BM_n)\to 
F( (M_1-\{k\}|M'_1); \BM_2; \cdots ; \BM_n)\,\,$$
(resp. $F( \BM_1; \cdots ; \BM_n)\to 
F(\BM_1; \BM_2; \cdots ; (M_n-\{k\}|M'_n))$\,). 
So $\bvphi=\sum \bvphi_k$ and  $\delta=\bpartial
+ \bvphi$ are endomorphisms of $\bop F(A)$. 

\quad$\bullet$   The second degree is $\deg_2(u)=|M'_{int}|+1$, and the second differential is $d_2=\bsigma_{int}$
defined as follows:
For  $k\in \cup_{1<i<n}  (M_i-M'_i)$, define
$$\bsigma_k(u)=(-1)^{| (M'_{int})_{>k}|} u_1\ts\cdots\ts\bsigma_k(u_j)\ts\cdots\ts u_r\,\,$$
where $\bsigma_k(u_j)$ is as in (\Defbsigma), 
and taking the sum of them, 
$$\bsigma_{int}(u):= \sum \bsigma_k(u)\,.$$

\quad$\bullet$  The third degree is $\deg_3(u)=|M'_{out}|$, and the third differential is $d_3=\bsigma_{out} +\Bf_{K_1}
+\Bf_{K_n}$, where the three operators are defined as follows. 

For $k\in (M_1-M'_1)\cup (M_n- M'_n)$, let 
$$\bsigma_k(u)=(-1)^{| (M'_{out})_{>k}|}u_1\ts\cdots\ts\bsigma_k(u_j)\ts\cdots\ts u_r\,\,$$
(so $j=1$ or $r$), and 
$$\bsigma_{out}(u):= \sum \bsigma_k(u)\,,$$
 the sum over $k\in (M_1-M'_1)\cup (M_n- M'_n)$.

For a sequence $P$ with $\term(P)=\init(M)$, set  
$$\Bf_{K_1}(P)\ts u =\{|M'_{out}|+1\}f_{K_1} (P)\ts u\quad\text{and}\quad 
\Bf_{K_1} (u)=\sum_{P}\Bf_{K_1}(P)\ts u\,\,.$$ 
Similarly set 
$$ u \ts \Bf_{K_n}(Q)= - u\ts f_{K_n}(Q)\quad\text{and}\quad
\Bf_{K_n} (u)=\sum_{Q} u \ts \Bf_{K_n}(Q)\,.$$
For these maps to be defined, one must take appropriate distinguished subcomplexes 
of $F(A)$ as for the case $n=2$. 
 \bigskip 

By the following result, which can be shown as in the previous subsection, we have a ``triple" complex.
\bigskip 

(\HIGI.1)
\Prop{\it (1) $\bsigma_{int}$ is a differential, and it commutes with $\delta$;
namely we have the following identities:
$$\bsigma_{int}\bsigma_{int}=0\,, \quad \bsigma_{int}\delta =\delta \bsigma_{int}\,.$$
Similarly for $\bsigma_{out}$:
$$\bsigma_{out}\bsigma_{out}=0\,, \quad \bsigma_{out}\delta =\delta \bsigma_{out}\,.$$
 The differentials $\bsigma_{int}$ and $\bsigma_{out}$ commute:
$$ \bsigma_{int}\bsigma_{out}=\bsigma_{out}\bsigma_{int}\,\,.$$
 
(2) The maps $d_1=-\delta$, $d_2=\bsigma_{int}$, and $d_3=\bsigma_{out} +\Bf_{K_1}
+\Bf_{K_n}$
are differentials, commuting with each other. 
 }\smallskip

{\it Proof.}\quad 
(1) Similar to the proof of (\Defbsigma.1). 

(2) Similar to the proof of (\GKL.1). 

(3) Follows from (1) and (2). 
\bigskip 
 
 Let $\BH(K_1, \cdots, K_n)$ be the total complex of $H\tbullet$.  The total differential is 
$$d_{\BH}=-\delta +(-1)^{\deg_1}d_2+ (-1)^{\deg_1+\deg_2}d_3\,\,.$$
As in the case $|I|=2$, we have the following claim; the proof is parallel to that
for (\GKL.2). 
\bigskip 

(\HIGI.2)
\Prop{\it
The complex 
$$H^{\bullet\,\bullet\,0}\mapr{d_3}H^{\bullet\,\bullet\,1}\mapr{d_3} \cdots$$
is exact. 
}\bigskip

Now set  $G\dbullet(K_1, \cdots, K_n)=\Phi (H\tbullet(K_1, \cdots, K_n)\,)$; it is a 
``double'' complex.  We also define $\BG(K_1, \cdots, K_n)$ to be its total complex.
Note the total degree and  differential of this complex are
$\deg_1+ \deg_2$ and 
$$d_\BG=-\delta +(-1)^{\deg_1}\bsigma_{int}$$
when acting on $u$.   
It is a quasi-isomorphic subcomplex of $\BH(K_1, \cdots, K_n)$.
When the sequence $K_1, \cdots, K_n$ is understood, write them as 
$G^{\bullet\,\bullet}([1, n])$ and 
 $\BG([1, n])$, respectively. 

The same construction applies to any finite totally ordered set
$I$ and a sequence of $C$-diagrams indexed by $I$, so we have the complex 
$\BG(I)$. 

If $|I|=2$, say $I=[1, 2]$, the construction in this subsection is 
related to that in (\GKL) as follows. 
Let $K=K_1$ and $L=K_2$. Then 
the second degree of $H\tbullet(K, L)$ is 
 concentrated in one, and $\bsigma_{int}=0$.
 Thus one sees that the partial totalization $\Tot_{12}(H\tbullet(K, L))$ coincides with 
 $H\dbullet(K, L)$ in (\GKL), and 
 that the totalization of $G\dbullet(K, L)$ coincides with 
 $\BG(K, L)$ in (\GKL).  

For the sake of reference we record the condition for an element of $H^{a, b, 0}(I)$
be in $G^{a, b}(I)$. 
\bigskip 

(\HIGI.3)($\sigma$-consistency)\quad 
Let $u\in H^{a, b, 0}([1, n])$ be an element with components 
$$u(M_1; \BM_2; \cdots; \BM_{n-1}; M_n)\in F(M_1; \BM_2; \cdots; \BM_{n-1}; M_n)$$
 (note the sequences $M_1$ and $M_n$ are free by the assumption $\deg_3(u)=0$). 
Then $u$ is in $G^{a, b}([1, n])$ if and only if the following equalities  are 
satisfied (we then say $u$ is {\it $\sigma$-consistent}.) Note the signs in the 
second equality. 
\smallskip 

(i) For each $k\in M_1-\{\init(M_1)\}$, one has 
$$-f_{K_1}((M_1)_{\le k})\ts u((M_1)_{\ge k}; \BM_2; \cdots; \BM_{n-1};M_n)
+\bsigma_k ( u(M_1; \BM_2; \cdots ; \BM_{n-1};M_n)\,)=0\,\,, $$
equivalently, 
$$-f_{K_1}((M_1)_{\le k})\ts u((M_1)_{\ge k}; \BM_2; \cdots; \BM_{n-1};M_n)
+\sigma_k ( u(M_1; \BM_2; \cdots ; \BM_{n-1};M_n)\,)=0\,\,.$$

(ii) For each $k\in M_n-\{\term(M_n)\}$, 
$$-u(M_1; \BM_2; \cdots; \BM_{n-1}; (M_n)_{\le k})\ts f_{K_n}((M_n)_{\ge k})
+ \bsigma_k( u(M_1; \BM_2; \cdots; \BM_{n-1};M_n)\,)     =0\,\,,$$
equivalently, writing  
$u(M_1; \BM_2; \cdots; \BM_{n-1}; M_n)=\sum u_1\ts\cdots\ts u_r\,,$ 
with $u_i\in F(A_i)$ as before, one has 
$$-u(M_1; \BM_2; \cdots; \BM_{n-1}; (M_n)_{\le k})\ts f_{K_n}((M_n)_{\ge k})
+ \sum\pm\sigma_k(u_1\ts\cdots\ts u_r)=0\,\,,$$
the sign for $\sigma_k$ being given by $(-1)^{\deg_1 (u_r)\cdot \gamma((M_n)_{\ge k})}$. 
\bigskip

The following proposition shows a difference between the cases $|I|=2$ and 
$|I|\ge 3$. 
\bigskip 

\sss{\acyclicGI} {\bf Proposition.}\quad{\it
$\BG(I)$ is acyclic if $|I|\ge 3$. }
\smallskip  

{\it Proof.}\quad  We show the acyclicity of $\BH(I)$ for $I=[1, n]$. 
For each pair of integers $(a, b)$, the sum 
$${\cal F}(a, b)=\moplus_{\init{M_1}\le a,\, \term{M_n}\ge b} F(\BM_1;\cdots; \BM_n)$$
is a subcomplex of $\BH(I)$, and gives its filtration (increasing in $a$, 
decreasing in $b$). 
In a successive quotient, which is the form 
$$\Gr_{\cal F}^{(a, b)}=\moplus_{\init{M_1}=a,\, \term{M_n}= b} F(\BM_1;\cdots; \BM_n)\,\,,$$
the maps $\Bf_{K}$ are zero, so the differential is a signed sum of 
$\partial$, $\vphi$, and $\sigma$. 
Consider an increasing filtration on it defined by
$$Fil_c=\moplus_{|M|\le c} F(\BM_1;\cdots; \BM_n)\,\,. $$
In a successive quotient $\Gr_{Fil}^c\Gr_{\cal F}^{(a, b)}$, one has 
$\bvphi=0$, so it is a sum of the total complexes of the 
following form:
$$0\to F(\cI|\emptyset)\mapr{\sigma}
\moplus_{|S|=1} F(\cI| S)
\mapr{\sigma} \moplus_{|S|=2} F(\cI | S)
\mapr{\sigma} \cdots \to  F(\cI|\ctop{\cI})\to 0\,\,$$
where $\cI=\amalg M_i$ and $S$ varies over subsets of $\ctop
{\cI}$.  
Since $|I|\ge 3$, one has $|\cI|\ge 3$; thus the total complex of the above is acyclic as shown by the spectral sequence argument using the following obvious lemma. 
\bigskip 

(\acyclicGI.1) {\bf Lemma.}\quad{\it
Let $A$ be an abelian group, $T$ be a non-empty finite set. 
For each subset $S$ of $T$ let $A_S=A$ be a copy of $A$, and for 
an inclusion $S\subset S'$, let $\al_{S\, S'}: A_S\to A_{S'}$ be the identity map. 
Then the sequence 
$$0\to A_\emptyset \to 
\moplus_{|S|=1}  A_S\to 
\moplus_{|S|=2}  A_S\to \cdots \to A_T\to 0 \,, $$
where the maps are alternating sums of $\al_{S\, S'}$, is exact. }
\bigskip

\sss{\htsPhicpx}{\it The subcomplexes $\Phi (A)\hts \Phi(B)$.}\quad
This subsection is a complement to (\Phicpx), from which we keep the notation. 
We explain a procedure to give a quasi-isomorphic subcomplex of the 
complexes $\Phi(A)\ts \Phi(B)$. (This will be needed in the next subsection.)
 
Assume  we have a quasi-isomorphic ``quadruple" subcomplex $A^{\dbullet}\hts B^{\dbullet} \injto A^{\dbullet}\ts B^{\dbullet}$,  see (\Multicpx.1). 
Assume further that  $A^{\dbullet}$ and $B^{\dbullet}$ are complexes of {\it free $\ZZ$-modules}. 
By the freeness assumption, the natural map 
$\Phi(A)^{\bullet}\ts \Phi(B)^{\bullet}\to A^{\bullet\, 0}\ts B^{\bullet\, 0}$ 
is an injection. 
We set 
$$\Phi(A)^{\bullet}\hts \Phi(B)^{\bullet}:= (\Phi(A)^{\bullet}\ts \Phi(B)^{\bullet}\,)\, \cap\,
(A^{\bullet\, 0}\hts B^{\bullet\, 0})$$
the intersection being taken in $A^{\bullet\, 0}\ts B^{\bullet\, 0}$.
Note one has the identity  
$$\Phi(A)^{\bullet}\hts \Phi(B)^{\bullet}= \Phi ( \Tot_{24}(A^{\dbullet}\hts B^{\dbullet})\,)\,.\eqno{(\htsPhicpx.a)}$$
Using this one sees that
the four inclusions in the following square are all quasi-isomorphisms. 
$$\begin{array}{ccc}
 \Phi(A)^{\bullet}\hts \Phi(B)^{\bullet}   &\injto  &A\dbullet\hts B\dbullet  \\
  \downarrow & &\downarrow  \\
 \Phi(A)^{\bullet}\ts \Phi(B)^{\bullet}  &\injto  &\phantom{\,\,.}A\dbullet\ts B\dbullet\,\,.
 \end{array}
$$
\smallskip 

 If the freeness assumption is not satisfied, one could 
take the identity (\htsPhicpx.a) as the 
definition of the left hand side; but in this paper the freeness assumption will be always satisfied.
\smallskip 

More generally, assume given a sequence of ``$n$-tuple" complexes 
 $A_1^{\bullet\cdots \bullet}, A_2^{\bullet\,\cdots \bullet}, 
\cdots, A_c^{\bullet\,\cdots \bullet}$, and a quasi-isomorphic 
multiple subcomplex $A_1^{\bullet\cdots \bullet}\hts A_2^{\bullet\,\cdots \bullet}\hts A_c^{\bullet\,\cdots \bullet}$.  One then has a ``$c(n-1)$"-tuple subcomplex 
$$\Phi(A_1)\hts\cdots\hts\Phi(A_c)\injto \Phi(A_1)\ts\cdots\ts\Phi(A_c)\,.$$

In the next subsection, we have ``triple" complexes $A\tbullet$, $B\tbullet$, and 
a quasi-isomorphic subcomplex $A\tbullet\hts B\tbullet$ is given.
One then has 
a ``quadruple" complex $\Phi(A)\dbullet\hts\Phi(B)\dbullet$. 
We have, as in (\Tensorproductcpx), a ``double" complex
$$\Phi(A)\dbullet\htimes\Phi(B)\dbullet:=\Tot_{13}\Tot_{24}(
\Phi(A)\dbullet\hts\Phi(B)\dbullet)\,\,,$$
and an isomorphism 
$$u: \Tot(\Phi(A)\dbullet)\hts \Tot(\Phi(B)\dbullet) \isoto 
\Tot (\Phi(A)\dbullet\htimes\Phi(B)\dbullet)\,.\eqno{(\htsPhicpx.b)}$$
\bigskip

\sss{\GISigma} {\it The complex $\BG(I|\Sigma)$}\quad
In the rest of this section we work under the following assumption:
 $I$ is a finite totally ordered set, and given 
a sequence of $C$-diagrams 
$K_i=(K_i^m; f_{K_i}(M))$
indexed by $I$. For simplicity of notation we often assume $I=[1, n]$. 

Let  $I_1, \cdots, I_c$ be a segmentation of $I$. 
We will define a quasi-isomorphic multiple subcomplex denoted 
$$\BG(I_1)\tildets\cdots\tildets \BG(I_c)\subset \BG(I_1)\ts\cdots\ts \BG(I_c)\,\,.$$
The usage of the symbol $\tildets$ is made in line with the convention 
(\Multicpx.1) for a quasi-isomorphic multiple subcomplex, but to avoid confusion 
with $\hts$. 

First consider the case $c=2$, and for  simplicity assume $I_1=[1,m]$ and $I_2=[m, n]$. 
Let $A_1=(\BM_1; \cdots ;\BM_m)$ and $A_2=({\BN}_m;
\cdots ;{\BN}_n)$ be multi-sequences on $[1, m]$ and $[m, n]$,
respectively.  

\vspace*{0.5cm}
\hspace*{0cm}
\input{new-pair-multiseq}
\vspace{0.8cm}

Let 
$$F(A_1)\tildets F(A_2)$$
be the distinguished subcomplex of $F(A_1)\ts F(A_2)$
defined as follows. 
\smallskip 

(i){\it Case   $\term(M_m)>\init(N_m)$\,.}\quad 
One has a distinguished subcomplex $[F(A_1)]_f\injto F(A_1)$
satisfying the same conditions with respect to $f_{K_1}(P)\ts (-)$
and $(-)\ts f_{K_n}(Q)$ 
as the complex 
$[F(M|M'; N|N')]_f$ in (\GKL).
Similarly one has $[F(A_2)]_f$. 
We set
 $$F(A_1)\tildets F(A_2)=[F(A_1)]_f \ts [F(A_2)]_f\,,
$$
which is a distinguished subcomplex as a tensor product of distinguished subcomplexes, 
see (\DefqDG), (8). 
\smallskip 

(ii){\it Case $\term(M_m)\le \init(N_m)$\,.}\quad 
Set $M=\amalg M_i$, $N=\amalg N_i$, $M'=\amalg M'_i$ and $N'=\amalg N'_i$ for 
convenience, and 
let $I=M\amalg N$, and $I_1, \cdots, I_r$ be the segmentation of $I$ obtained as the union of the segmentation of $M$ by $M'$ and that of $N$ by $N'$. 
The corresponding tensor product complex is
$F(I_1)\ts\cdots\ts F(I_r)$, and it contains 
$F(A_1)\ts F(A_2)$ as a distinguished subcomplex. 
To specify a condition of constraint on the tensor product complex,  let 
$$\BI=[-w, \term(M_1)]\amalg (\coprod_{2\le i\le m-1} M_i)\amalg  
[\term(M_m), \init(N_m)]\amalg(\coprod_{m+1\le i\le n-1} N_i)
 \amalg 
[\init(N_n), w]\,,$$
with $w$ large enough (see the figure on the right). 
Choose a set of almost disjoint sub-intervals $\{J_j\}$
such that each $J_i$ is contained in one of the intervals 
$[-w, \init(M_1)]$, $[\init(M_m), \term(N_m)]$, or 
$[\term(N_n), w]$. 
Correspondingly let 
$$f(J_j)=\left\{
\begin{array}{rl}
f_{K_1}(J_j)&\quad \mbox{if $J_j\subset [-w, \init(M_1)]$,}\\
f_{K_m}(J_j)&\quad \mbox{if $J_j\subset [\term(M_m), \init(N_m)]$,}\\
f_{K_n}(J_j)&\quad \mbox{if $J_j\subset [\term(N_n), w]$.}
\end{array}
\right. $$
The set of data $\cC=(\BI; P=[1, r]; \{J_j\}; \{f(J_j)\})$ gives a constraint;
the corresponding distinguished subcomplex is generated by 
$u_1\ts\cdots\ts u_r$ with $u_i\in F(I_i)$ such that 
$$\{u_1, \cdots, u_r,  \{f(J_j)\}\,\} \qquad\mbox{ is properly intersecting.}
\eqno{(\GISigma.a)}$$
To obtain 
our subcomplex $F(A_1)\tildets F(A_2)$, 
take the intersection of such distinguished subcomplexes 
for  all possible sets of sub-intervals $\{J_j\}$, then take the intersection 
with $[F(A_1)]_f \ts [F(A_2)]_f$ as in Case (i):
$$F(A_1)\tildets F(A_2)=\bigcap_{\cC} [F(A_1)\ts F(A_2)]_\cC\, \cap\, ([F(A_1)]_f \ts [F(A_2)]_f)\,.$$

The subcomplex $F(A_1)\tildets F(A_2)$ was so specified that 
\smallskip 

(a) $H\tbullet\tildets H\tbullet$ is a subcomplex of $H\tbullet\ts H\tbullet$ (see below for this), 
and 

(b) the map $\rho_m$ can be defined in the next subsection;
the following fact will be used there.  
If   $\term(M_m)=\init(N_m)=\{\ell\}$
(we will then say that $A_1$ and $A_2$ are {\it composable}) one has 
 $F(A_1)\tildets F(A_2)\subset [F(A_1\scirc A_2)]_f$, 
where we recall  $A_1\scirc A_2$ is the ``concatenation" of $A_1$ and $A_2$, namely
$$A_1\scirc A_2:=(\BM_1; \cdots ; \BM_{m-1}; (M_m\cup N_m|M_m'\cup\{\ell\}\cup N'_m);
\cdots ;\BN_n)\,. $$ 
Indeed $\cap_{\cC} [F(A_1)\ts F(A_2)]_\cC$ equals $[F(A_1\scirc A_2)]_f$
in this case. 
\smallskip

With these subcomplexes, let  
$$H\tbullet([1, m])\tildets H\tbullet([m, n])\subset H\tbullet([1, m])\ts H\tbullet([m, n])$$
be the ``6-tuple" subcomplex defined as the sum $\bop F(A_1)\tildets F(A_2)
\subset \oplus [F(A_1)]_f\ts [F(A_2)]_f$.
 It is a quasi-isomorphic 
subcomplex. To this subcomplex, 
we apply the construction of (\htsPhicpx) to obtain a quasi-isomorphic ``4-tuple"
subcomplex 
$$G\dbullet([1, m])\tildets G\dbullet([m, n])
=\Phi(H\tbullet([1, m])\,)\tildets \Phi(H\tbullet([m, n])\,)
\injto 
G\dbullet([1, m])\ts G\dbullet([m, n]) \,. $$
Then set 
$$\BG([1, m])\tildets \BG([m, n])=\Tot_{12}\Tot_{34}(G\dbullet([1, m])\tildets G\dbullet([m, n])\,)\,.$$ 
This is a ``double" complex, which can be turned into a double complex by changing signs as in (0.1).  

Also, as discussed in (\htsPhicpx), we have the ``double" complex 
$$G\dbullet([1, m])\tildetimes G\dbullet([m, n])=\Tot_{13}\Tot_{24}(G\dbullet([1, m])\tildets G\dbullet([m, n])\,)\,\,$$
and an isomorphism of complexes (\htsPhicpx.b)
$$u: \BG([1, m])\tildets \BG([m, n])
\to \Tot( G\dbullet([1, m])\tildetimes G\dbullet([m, n])\,)\,.$$

If $c>2$ the definition is similar.
For multi-sequences $A_1, \cdots, A_c$ on $I_1, \cdots, I_c$ respectively,
one can define a distinguished subcomplex 
$$F(A_1)\tildets F(A_2)\tildets\cdots\tildets F(A_c)$$
subject to similar conditions. 
Taking the sum of them we have the ``3$c$-tuple" subcomplex 
$$H\tbullet(I_1)\tildets\cdots\tildets H\tbullet(I_c)$$
(also denoted $H\tbullet(I|\Sigma)$ if $\Sigma$ corresponds to 
the segmentation $I_1, \cdots, I_c$)
and taking the $\Phi$-part of this  we get a ``2c-tuple" quasi-isomorphic subcomplex 
$$G\dbullet(I_1)\tildets\cdots\tildets G\dbullet(I_c)
\subset 
G\dbullet(I_1)\ts\cdots\ts G\dbullet(I_c)\,.$$
Then we get  a ``$c$-tuple" subcomplex
$$\BG(I|\Sigma)=\BG(I_1)\tildets\cdots\tildets\BG(I_c):= \Tot_{12}\Tot_{34}\cdots\Tot_{2c-1, 2c}
(G\dbullet(I_1)\tildets\cdots\tildets G\dbullet(I_c))\,\,
\eqno{(\GISigma.b)}
$$
of $\BG(I\tbar\Sigma)=\BG(I_1)\ts\cdots\ts\BG(I_c)$. 
We always view this as a $c$-tuple complex by changing signs as in (0.1). 
Let $\iota_\Sigma: \BG(I|\Sigma)\to \BG(I \tbar \Sigma)$ be the inclusion.
The complex $\BG(I|\Sigma)$ will play a major role in the rest of this section.
\bigskip 

\sss{\rhoPi} {\it The maps $\rho$ and $\Pi$.}\quad 
We define a map called the product map (for $m\in \Sigma$)
$$\rho_m: \BG(I|\Sigma)\to \BG(I|\Sigma-\{m\})\,\,.$$
This is to be a map of multiple complexes as follows. 
If $I_1, \cdots, I_a$ is the segmentation of $I$ by $\Sigma-\{m\}$, 
$I'_1, \cdots, I'_{a+1}$ the segmentation by $\Sigma$, and $f: [1, a+1]\to [1, a]$ the 
map such that $I'_{f(i)}\subset I_i$, then 
$\Tot^f \BG(I|\Sigma)$ and $\BG(I|\Sigma-\{m\})$ are
 $a$-tuple complexes,  
and $\rho_m$ is a map of $a$-tuple complexes. 

For simplicity of notation let us consider the case 
$\Sigma=\{m\}$, where the map is of the form
\newline
$\rho_m: \BG([1, m])\tildets \BG([m, n])\to \BG([1, n])$. 
Consider the map 
$$\rho_m: H\tbullet([1, m])\tildets H\tbullet([m, n])\to H\tbullet([1, n])$$
defined as follows. 
The source is the sum of $F(A_1)\tildets F(A_2)$, where $A_1$ and $A_2$ are multi-sequences on 
$[1, m]$ and $[m, n]$, respectively. 
On each direct summand $F(A_1)\tildets F(A_2)$ where $A_1$ and $A_2$ are composable, 
$\rho_m$ is the inclusion map 
$F(A_1)\tildets F(A_2)\to [F(A_1\scirc A_2)]_f$.
If $A_1$ and $A_2$ are not composable,   $\rho_m$ is set to be zero
on  $F(A_1)\tildets F(A_2)$. 

Observe that $\rho_m$ takes $H^{a,b, c}\tildets H^{a', b', c'}$ to $H^{a+a', b+b', c+c'}$. 
Indeed if $A_1=(\BM_1; \cdots ;\BM_m)$, $A_2=(\BN_m;\cdots ;\BN_n)$ are composable,  
$M_m=N_m=\{\ell\}$, and $u\in F(A_1)$, $v\in F(A_2)$, set
$$M'_{int}=\cup_{1<i\le m} M'_i\,, \quad N'_{int}=\cup_{m\le i<n} N'_i\,.$$
One has 
$$\deg_1(u\ts v)=|u|+|v|+\gamma(M_1)+\gamma(N_n)=\deg_1(u)+\deg_1(v)\,,$$
$$\deg_2(u\otimes v)=|M'_{int}\cup\{\ell\}\cup N'_{int}|+1= \deg_2 u+\deg_2 v\,,$$
and 
$$\deg_3(u\otimes v)=|M'_1\cup N'_m|= \deg_3(u)+\deg_3(v)\,.$$
By restriction one obtains a map 
$\rho_m: H^{\bullet\,\bullet\,0}([1, m])\tildets H^{\bullet\,\bullet\,0}([m, n])\to H^{\bullet\,\bullet\,0}([1, n])$.
\bigskip 

(\rhoPi.1) \Lem{\it 
(1) The  map $\rho_m$ above gives rise to a map of 
 ``double" complexes
$$\rho_m: G^{\bullet\, \bullet}([1, m])\tildetimes G^{\bullet\, \bullet}([m, n])
\to G^{\bullet\, \bullet}([1, n])\,\,.$$ 

(2) One has an induced map of complexes 
$$\rho_m: \BG([1, m])\tildets \BG([m, n])\to \BG([1, n])\,;$$
 it sends the sum of terms 
$u\ts v\in F(A_1)\tildets F(A_2)$, with $A_1$, $A_2$ composable, 
to the sum of $\{(\deg_1 u)\cdot (\deg_2 v)\}u\ts v$. 
}\bigskip 

{\it Proof.}\quad 
(1)  The meaning of the statement is that the degree-preserving map 
$$\rho_m: H^{\bullet\,\bullet\,0}([1, m])\tildets H^{\bullet\,\bullet\,0}([m, n])
\to H^{\bullet\,\bullet\,0}([1, n])$$
sends $G^{\bullet\, \bullet}\tildets G^{\bullet\, \bullet}$ into 
$ G^{\bullet\, \bullet}$, and  
the induced map 
$$G^{\bullet\, \bullet}\tildetimes G^{\bullet\, \bullet}:=\Tot_{13}\Tot_{24}(G^{\bullet\, \bullet}\tildets G^{\bullet\, \bullet}\,)
\to G^{\bullet\, \bullet}$$
is a map of ``double" complexes:   
One has, for a homogeneous element  $u\ts v\in G\dbullet\tildets G\dbullet$, 
equalities
$$\delta\rho(u\ts v)=(-1)^{\deg_1 v}\rho(\delta u\ts v)+\rho( u\ts \delta v)\,,$$
$$\bsigma_{int}\rho(u\ts v)=(-1)^{\deg_2 v}\rho(\bsigma_{int} u\ts v)+\rho( u\ts \bsigma_{int} v)\,\,.$$

For $u=(u(A_1)\,)\in G\dbullet([1, m])$ and $v=(v(A_2)\,)\in G\dbullet([m, n])$, 
one has $\rho_m(u\ts v)=\sum u(A_1)\ts v(A_2)$, the sum over composable pairs
 $A_1=(\BM_1; \cdots ; \BM_m)$
and $A_2=(\BN_m;\cdots ; \BN_n)$ both of third degree 0
($\BM_1=(M_1|\emptyset)$ and  $\BM_m, \BN_m, \BN_n$ are free sequences).
It is in $G\dbullet$
since the condition (\HIGI.3) holds; this shows the first assertion. 

The first equality holds, since for each composable $(A_1, A_2)$, 
one has 
$$\delta (u(A_1)\ts v(A_2))=(-1)^{\deg_1 v}(\delta u(A_1)\,)\ts v(A_2)+
 u(A_1)\ts (\delta v(A_2))\,.$$
by (\FMM.1). 

For the second equality, if $A_1$ and $A_2$ are composable with $\term(M_m)=
\init(N_m)=\{\ell\}$, 
one has from the definition
$$\bsigma_{int}(u(A_1)\ts v(A_2))=\bsigma'_{int}(u(A_1)\ts v(A_2))
+\bsigma''_{int}(u(A_1)\ts v(A_2))$$
where $\bsigma'_{int}:=\sum \bsigma_k$, the sum over $k\in \cup_{1<i<m}(M_i-M'_i)$
and $k\in \cup_{m<i<n}(N_i-N'_i)$, and 
$\bsigma''_{int}:=\sum \bsigma_k$, the sum over
 $k\in (M_m\cup N_m)-\{\ell\}$. 
We have 
$$\bsigma'_{int}(u(A_1)\ts v(A_2))=(-1)^{\deg_2 v}
\left(\bsigma_{int}( u(A_1))\ts v(A_2)+u(A_1)\ts \bsigma_{int} (v(A_2))\,\right)\,\,, $$
by the same proof as for (\Defbsigma.1), (3). 

As for the other term $\bsigma''_{int}(u(A_1)\ts v(A_2))$, we claim that the sum of 
them over composable pairs $(A_1, A_2)$ equals zero. 
If $(L|\{\ell, \ell'\})$ is a sequence with $\ell<\ell'$  in $L$, 
$L$ is segmented into three intervals $L_{\le \ell}$, $L_{[\ell, \ell']}$, and 
$L_{\ge \ell'}$. 
The elements $u(A_1)$, $v(A_2)$ associated to the composable pair 
$$A_1=(M_1; \BM_2; \cdots; \BM_{m-1}; L_{\le \ell}),\quad 
A_2=(L_{\ge \ell}; \BN_{m+1};\cdots;\BN_{n-1}; N_n)$$
will be denoted by $u(L_{\le \ell})$ and $v(L_{\ge \ell})$. 
Similarly we have elements $u(L_{\le \ell'})$ and $v(L_{\ge \ell'})$.
Now $\bsigma''_{int}(u\ts v)$ has in it terms 
$$\bsigma_{\ell'}(u(L_{\le \ell})\ts v(L_{\ge \ell})\,)
+\bsigma_{\ell}(u(L_{\le \ell'})\ts v(L_{\ge \ell'})\,)\,,
$$
that equals (with $\tau=|N'|$)
$$(-1)^\tau u(L_{\le \ell})\ts \bsigma_{\ell'}(v(L_{\ge \ell})\,)
+(-1)^{\tau+1} \bsigma_{\ell}(u(L_{\le \ell'})\,)\ts v(L_{\ge \ell'})\,.
$$
We find it is zero using the identities (see (\HIGI.3)\,)
$$\bsigma_{\ell'}(v(L_{\ge \ell})\,)=f_{K^m}(L_{[\ell, \ell']})\ts v(L_{\ge \ell'})\,,$$
$$\bsigma_{\ell}(u(L_{\ge \ell'})\,)=u(L_{\le \ell})\ts 
f_{K^m}(L_{[\ell, \ell']})\,. $$

(2) Apply the operation $\Tot$ to the map $\rho_m$ of (1) to  
get a map of complexes
$$\Tot(G^{\bullet\, \bullet}([1, m])\tildetimes G^{\bullet\, \bullet}([m, n])\,)
\to \BG([1, n])\,,$$
 and compose it with the isomorphism 
$u:\BG([1, m])\tildets \BG([m, n])\to \Tot(G\dbullet\tildetimes G\dbullet)$
recalled in (\htsPhicpx). 
This completes the proof of the lemma. 
\bigskip

The case $\Sigma$ contains more than one element is similar, so one has 
 $\rho_m: \BG(I|\Sigma)\to \BG(I|\Sigma-\{m\})$. 
 Explicitly, consider the map ($I_i\cup I_{i+1}=\{m\}$)
 $$\rho_m: H\tbullet(I_1)\tildets\cdots\tildets H\tbullet(I_c)
 \to H\tbullet(I_1)\ts\cdots\ts H\tbullet(I_i\cup I_{i+1})\tildets
 \cdots \tildets H\tbullet(I_c)$$
 which is the identity on $H\tbullet(I_j)$ with $j\neq i, i+1$, and which is, 
on the remaining factors, 
 the map 
 $$\rho_m: H\tbullet(I_i)\tildets H\tbullet(I_{i+1})\to 
 H\tbullet(I_i\cup I_{i+1})$$
 as in (\GISigma). 
 This induces by restriction  the map of ``$2(c-1)$-tuple" complexes
$$
\begin{array}{ll}
\rho_m: &G\dbullet(I_1)\tildets \cdots\tildets 
(G\dbullet(I_i)\tildetimes G\dbullet(I_{i+1})\,)\tildets\cdots\tildets
 G\dbullet(I_c)   \\
 \to &G\dbullet(I_1)\tildets \cdots\tildets G\dbullet(I_i\cup I_{i+1})
 \tildets \cdots\tildets G\dbullet(I_c)\,.
\end{array}$$
Taking the totalization $\Tot_{12}\cdots\Tot_{2c-3, 2c-2}$, 
obtains a map of ``$(c-1)$-tuple" complexes 
$$
\begin{array}{ll}
\rho_m:& \BG(I_1)\tildets \cdots\tildets 
\Tot(\BG(I_i)\tildetimes \BG(I_{i+1})\,)\tildets
\cdots\tildets
\BG(I_c)  \\
 \to & \BG(I_1)\tildets \cdots\tildets \BG(I_i\cup I_{i+1})
 \tildets \cdots\tildets \BG(I_c)\,.
\end{array}$$
This is also a map of  $(c-1)$-tuple complexes if 
we view the source and the target as ``$(c-1)$"-tuple complexes
as in (0.1). 
 Arguing as before on the factor $G\dbullet(I_i)\tildets G\dbullet(I_{i+1})$
 \newline
$ \to G\dbullet(I_i\cup I_{i+1})$, we obtain a map of 
$(c-1)$-tuple complexes
$$\rho_m: \BG(I_1)\tildets \cdots\tildets \BG(I_c)
 \to \BG(I_1)\tildets \cdots\tildets \BG(I_i\cup I_{i+1})
 \tildets \cdots\tildets \BG(I_c)\,.$$
 
 By definition it is obvious that, for $m$, $m'$ distinct, $\rho_m\rho_{m'}=\rho_{m'}\rho_m$.
So one can define, for $K\subset \ctop{I}$, a map 
$$\rho_K: \BG(I|\Sigma)\to \BG(I|\Sigma-K)\eqno{(\rhoPi.a)}$$
by composing $\rho_m$ for $m\in \Sigma$. 
 \bigskip 

For $k\in \ctop{I}-\Sigma$, define a map of complexes 
$$\Pi_k: \BG(I|\Sigma)\to \BG(I-\{k\}|\Sigma)\eqno{(\rhoPi.b)}$$
as follows. Let $I=[1, n]$ for simplicity. 

First assume $\Sigma=\emptyset$. 
On the direct sum $F(A)$ where $A=(\BM_1; \cdots ;\BM_n)$
is a multi-sequence on $[1, n]$, let $\Pi_k=0$ unless 
$\BM_k=(M_k|\emptyset)$ with $|M_k|=1$. 
For such $\BM_k$, letting $M_k=\{j\}$ and 
$A|_{[1, n]-\{k\}}=(\BM_1; \cdots; \widehat{\BM_k};
\cdots;\BM_n)$, 
we define 
$$\Pi_k: F(A)\to 
F(A|_{[1, n]-\{k\}})\,\,, $$
by $\Pi_k(u)= (-1)^j\vphi_j (u)\,\,$
({\it not} a typo for $(-1)^j\bvphi_j (u)$). 
The following figure is for the case $n=3$ and  $k=2$. 

\vspace*{0.5cm}
\hspace*{1.5cm}
\unitlength 0.1in
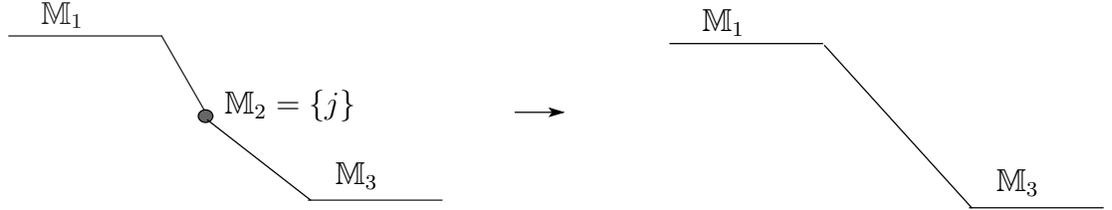
\begin{picture}(57.28,11.10)(1.70,-11.20)
%
\special{pn 8}%
\special{pa 170 220}%
\special{pa 970 220}%
\special{fp}%
%
\special{pn 8}%
\special{pa 970 220}%
\special{pa 1200 620}%
\special{fp}%
%
\special{pn 8}%
\special{pa 1738 1080}%
\special{pa 2438 1080}%
\special{fp}%
\put(3.4000,-1.8000){\makebox(0,0)[lb]{$\BM_1$}}%
%
\special{pn 8}%
\special{sh 0.600}%
\special{ar 1200 640 38 33  0.0000000 6.2831853}%
%
\special{pn 8}%
\special{pa 1210 660}%
\special{pa 1750 1080}%
\special{fp}%
\put(13.0000,-6.7000){\makebox(0,0)[lb]{$\BM_2=\{j\}$}}%
\put(18.8000,-10.2000){\makebox(0,0)[lb]{$\BM_3$}}%
%
\special{pn 8}%
\special{pa 3630 260}%
\special{pa 4430 260}%
\special{fp}%
%
\special{pn 8}%
\special{pa 5198 1120}%
\special{pa 5898 1120}%
\special{fp}%
\put(38.0000,-2.2000){\makebox(0,0)[lb]{$\BM_1$}}%
\put(53.4000,-10.6000){\makebox(0,0)[lb]{$\BM_3$}}%
%
\special{pn 8}%
\special{pa 4440 270}%
\special{pa 5210 1120}%
\special{fp}%
%
\special{pn 8}%
\special{pa 2820 620}%
\special{pa 3070 620}%
\special{fp}%
\special{sh 1}%
\special{pa 3070 620}%
\special{pa 3003 600}%
\special{pa 3017 620}%
\special{pa 3003 640}%
\special{pa 3070 620}%
\special{fp}%
\end{picture}%
\vspace{0.8cm}

\noindent
 Taking the sum over $A$'s, we obtain a map $\Pi_k: H\tbullet(I)\to H\tbullet(I-\{k\})$. 
Note that the map $u\mapsto \Pi_k(u)$ preserves the number $\gamma$, and thus also $\deg_1$. 
One verifies
\bigskip 

(\rhoPi.2) \Lem{\it (1) The map  $\Pi_k: H\tbullet(I)\to H\tbullet(I-\{k\})$
 is a map of ``triple" complexes. 

(2) If $k\neq k'$, $\Pi_k\Pi_{k'}=\Pi_{k'}\Pi_{k}$. 
}\bigskip 

{\it Proof.}\quad 
(1) For $A=(\BM_1; \cdots ;\BM_n)$ with $\BM_k=(M_k|\emptyset)$,
$|M_k|=1$ and $M_k=\{j\}$, if $u\in F(A)$, one has
$$\bpartial \vphi_j(u)= \vphi_j\bpartial(u)\,\,,$$
$$\bvphi_i \vphi_j(u)= \vphi_j\bvphi_i (u)\quad\text{for}\quad i\neq j.$$
To show the first equality, let $A_1, \cdots, A_r$ be the 
segmentation of $M=M_1\cup\cdots\cup M_n$ by $M'$, and assume 
$u=u_1\ts\cdots\ts u_r\in F(A)$ with $u_i\in F(A_i)$ as before.  
Then 
$$\bpartial(u)=\sum_i\{\sum_{a>i}\deg_1(u_a)\}u_1\ts\cdots\ts(
\partial u_i)\ts\cdots\ts u_r\,\,.$$
One also has, if $M_k$ is contained in $A_s$, 
$$\vphi_j(u_1\ts\cdots\ts u_r)=
u_1\ts\cdots\ts\vphi_j(u_s)\ts\cdots u_r$$
by the compatibility of $\vphi_j$ and the tensor product, (\DefqDG),(2).
Since $\deg_1(u_k)=\deg_1\vphi_j(u_k)$, the first equality
follows. The proof of the second equality is similar. 

We also have 
$$\bsigma_a\vphi_j(u)=\vphi_j\bsigma_a(u)\,\,.$$
Indeed if $a\in M-M'$, $a\neq j$, and if $\sigma_a (u)=\sum u'\ts u''$, then 
$\bsigma_a (u)=\sum \{\deg_1(u')\cdot \gamma(u'')\}u'\ts u''$. 
If $a>j$, $u'$ and $\vphi_j(u')$ have the same $\deg_1$; if
$a<j$, $u''$ and $\vphi_j(u'')$ have the same $\deg_1$. Hence the equality follows. 
Taking the sum over $a$, we obtain 
$$\bsigma_{int}\vphi_j(u)=\vphi_j\bsigma_{int}(u)\,\,,$$
$$\bsigma_{out}\vphi_j(u)=\vphi_j\bsigma_{out}(u)\,\,.$$
As for $\Bf$, one easily verifies 
$$\Bf_{K_1}\vphi_j(u)=\vphi_j\Bf_{K_1} (u)\,\,,
\quad 
\Bf_{K_n}\vphi_j(u)=\vphi_j\Bf_{K_n} (u)\,\,.$$

When $\BM_k=(M_k|\emptyset)$ with $M_k=\{j, j'\}$, 
 $j<j'$, then for 
$u\in F(A)$, 
$$((-1)^{j'}\vphi_{j'}\,)\bvphi_j+\bvphi_{j'}((-1)^{j}\vphi_{j}\,)=0\,\,.$$
(These terms appear in $\Pi_k\bvphi(u)$.)
To prove this assume $A_s$ contains $M_k$, so that 
$$\bvphi_j(u)=\sum_i\{\sum_{a>i}\deg_1(u_a) \}u_1\ts\cdots\ts
\bvphi_j (u_i)\ts\cdots\ts u_r\,\,, $$
and $\bvphi_j (u_i)= \{|u_i|+\gamma((A_s)_{\le j})\}\vphi_j(u_i)$. 
The identity follows from 
$\gamma((A_s)_{\le j'}) - \gamma((A_s)_{\le j})=j'-j-1$. 

Combining these, we conclude that $\Pi_k$ commutes with $\delta$, 
$\bsigma_{int}$, and $d_3=\bsigma_{out}+\Bf_{K_1}+ \Bf_{K_n}$. 

(2) is obvious from the definition. 
\bigskip 

Taking $\Phi$ and then $\Tot$ we get an induced map 
 $\Pi_k: \BG(I)\to \BG(I-\{k\})$.
 
We extend this to a map of complexes 
$\Pi_k: \BG(I|\Sigma)\to \BG(I-\{k\}|\Sigma)$ as follows. 
First consider the map (where $k\in I_i$)
$$H\tbullet(I_1)\tildets\cdots\tildets H\tbullet(I_c)\to 
H\tbullet(I_1)\tildets\cdots\tildets H\tbullet(I_i-\{k\})\tildets\cdots\tildets
H\tbullet(I_c)$$
which sends  $u=u_1\ts\cdots\ts u_c $ to 
$$u_1\ts\cdots\ts \Pi_k(u_i)\ts\cdots\ts u_c\,\,.$$
Applying $\Phi$, we obtain a map of ``$2c$-tuple" complexes
$$G\dbullet(I_1)\tildets\cdots\tildets G\dbullet(I_c)\to 
G\dbullet(I_1)\tildets\cdots\tildets G\dbullet(I_i-\{k\})\tildets\cdots\tildets
G\dbullet(I_c)\,,$$
and further applying $\Tot$, we get  a map of ``$c$-tuple" complexes 
$$\Pi_k: \BG(I|\Sigma)\to \BG(I-\{k\}|\Sigma
)\,.$$
The source and target can be viewed as $c$-tuple complexes, and then 
$\Pi_k$ is a map of $c$-tuple complexes.
It satisfies  the following properties. 
\bigskip

(\rhoPi.3) \Lem{\it (1) The 
$\Pi_k: \BG(I|\Sigma)\to \BG(I-\{k\}|\Sigma
)$ is a map of multiple complexes. 

(2) If $k\neq k'$, $\Pi_k\Pi_{k'}=\Pi_{k'}\Pi_{k}$. 

(3) If $m\neq k$,  $\Pi_k\rho_m=\rho_m\Pi_k$. 

(4) The map $\Pi_m\rho_m: \BG(I|\Sigma)\to \BG(I-\{m\}|\Sigma-\{m\})$ is zero. 
}\smallskip 

{\it Proof.}\quad (1) holds by definition. 
(2) holds  at the level of ``triple" complexes $H\tbullet(I|\Sigma)$,
the proof being parallel to that for (\rhoPi.2). 
(3) and (4) are obvious from the definitions. 
\bigskip

For a non-empty 
subset $K\subset \ctop{I}-\Sigma$, define 
$$\Pi_K: \BG(I|\Sigma)\to \BG(I-K|
\Sigma)$$
by composing $\Pi_k$ for $k\in K$. 
If $K$ is the disjoint union of $K'$ and $K''$, then 
$\Pi_K=\Pi_{K'}\Pi_{K''}$.
\bigskip 

\sss{\CompDoublecpx}{\it Complements on double complexes.}\quad 
A ``double" complex $(A\dbullet; d_1, d_2)$ is the same thing as a complex of complexes
$$\mapr{d_1} A^{p, \bullet} \mapr{d_1} A^{p+1, \bullet}\mapr{d_1}\cdots$$
where each $A^{p, \bullet}=(A^{p, \bullet}; d_2)$ is a complex. 
The total complex $\Tot(A)$, as defined in (\Multicpx), has a filtration by 
subcomplexes with graded quotients $A^{p, \bullet}[-p]$. 

Note that the total complex $A^\bullet =\Tot(A)$ satisfies the following properties. 
For each $i$, there is a direct sum decomposition $A^i=\moplus_p A^i_p$, where 
$A^i_p:=A^{p, i-p}$, and the differential $d$ is the sum of $d': A^i_p\to A^{i+1}_{p+1}$
and $d'': A^i_p\to A^{i+1}_{p}$. 

Conversely assume $(A^\bullet; d)$ is a complex of abelian groups and given a direct sum decomposition 
$$A^i=\moplus_{p\in \ZZ} A^i_p$$
for each $i$, such that $d$ decomposes as $d'+d''$ as above.  
Then if we set $A^{p, q}=A^{p+q}_p$, $d_1=d'$ and $d_2=(-1)^pd''$, then we obtain a
``double" complex $(A\dbullet; d_1, d_2)$. 
We recover $(A^\bullet, d)$ as the associated total complex. 

In the next we have a complex equipped with such a direct sum decomposition, 
which we will view as a ``double" complex. 
\bigskip

\sss{\GIT}{\it The complexes $\BG(I, T)$.}\quad 
For $I=[1, n]$, and a multi-sequence 
$A=(\BM_1; \cdots; \BM_n)$, the {\it type} of
$A$ is the subset of $\ctop{I}$ defined by 
$$T=\{i\in \ctop{I}|M'_i\neq \emptyset\}.$$
For $T\subset \ctop{I}$, consider the direct sum 
$$\moplus_{A\text{ of type } T} F(A)\,\,.$$
We form a ``triple" complex denoted $H\tbullet(I, T)$ as follows. 

The first degree and differential are the same as before in (\HIGI). 
The second degree is $\deg_2=|M'_{int}|$, and the second differential is 
$\bsigma_{int}=\sum \bsigma_k(u)$, 
the sum over $k \in M_i- M'_i$ with $i\in T$. 
The third degree and differential are the same  as before:
 $\deg_3(u)=|M'_{out}|$,  $d_3=\bsigma_{out} +\Bf_{K_1}
+\Bf_{K_n}$. 

The $H\tbullet(I, T)$ is a subquotient of $H\tbullet(I)$, as follows. 
The complex $H\tbullet(I)$ has a filtration by subcomplexes indexed by 
types. For a given type $T$ the corresponding subcomplex $Fil^T$ is given
as the sum $\bop F(A)$, where $A$ varies over multi-sequences $A$ with 
type containing $T$. The subquotient (at $T$) in the filtration is the 
complex $H\tbullet(I, T)$ introduced above. 

Temporarily let $\BG'(I, T)$ be the complex obtained from $H\tbullet(I, T)$
by applying $\Phi$
and taking $\Tot$:  $\BG'(I, T)=\Tot(\Phi H\tbullet(I, T)\,)$. 
Thus its differential is 
$-\delta+(-1)^{\deg_1}d_2$, 
with 
$$d_2=\sum_k\bsigma_k\,,\qquad\mbox{the sum over $k\in M_i-M'_i$ with $i\in T$.}$$
As a graded group $\BG(I)$ is the direct sum of $\BG'(I, T)$ for $T\subset \ctop{I}$. 
Recall that the differential of $\BG(I)$ is 
$-\delta+\sum(-1)^{\deg_1} \bsigma_k$, where the sum is over 
$k\in M_i-M'_i$, $i\neq 1, n$. 
Also note:
 \smallskip 

$\bullet$\quad 
The map $\delta$ takes $\BG'(I, T)$ to itself. 

$\bullet$\quad 
For $k\in M_i- M'_i$ with $i\in T$, the map  $\bsigma_k$ 
takes $\BG'(I, T)$ to itself.

$\bullet$\quad 
For $k\in M_i- M'_i$ with $i\not\in T$, 
$\bsigma_k$ takes $\BG'(I, T)$ to $\BG'(I, T\cup\{i\})$. 
\smallskip 

Now we can apply the discussion of (\CompDoublecpx) and present 
$\BG(I)$ as a ``double" complex.   
We define $\BG(I, T)$ to be the complex obtained from $\BG'(I, T)$ by
shifting the degree by 
$-|I|+|T|+2$:
$$\BG(I, T)=\Tot(\Phi H\tbullet(I, T)\,)
[\,-|I|+|T|+2]\,\,. $$
Namely 
$$\deg_{\BG(I, T)}(u)= \deg_{\BG(I)}(u)+|I| - |T|-2 $$
and the differential $d_{\BG(I, T)}$ is $(-1)^{-|I|+|T|+2}$ times the original differential 
of $\Tot(\Phi H\tbullet(I, T)\,)$. 

The complex $\BG(I)$ is the total complex of the 
``double" complex  
$$\BG(I, \emptyset)\mapr{\tilde{\bsigma} }\moplus_{|T|=1}
\BG(I, T)\mapr{\tilde{\bsigma} }
\cdots\mapr{\tilde{\bsigma}}\BG(I, \ctop{I})\,\,,
\eqno{(\GIT.a)}$$ where  $\BG(I, T)$ 
is  placed in first degree $-|I|+|T|+2$, and the first differential 
$\tilde{\bsigma}$ on $\BG(I, T)$ is given by 
$$\tilde{\bsigma}=\sum_k (-1)^{\deg_1} \bsigma_k\,,\qquad
\mbox{the sum over 
$k\in M_i- M'_i$ with $i\not\in T$.}\eqno{(\GIT.b)}$$

The following remark is obvious, but will be often used.
If $|I|=2$, then $T$ is an empty set, so $\BG(I, \emptyset)$ is placed in first degree 0, 
and 
$$\BG(I)=\BG(I, \emptyset)\,\,. $$  
If $|I|=3$, say $I=[1, 3]$, then the complex $\BG([1, 3])$ is  of the form 
$$\BG([1, 3], \emptyset)\mapr{\tilde{\bsigma}}\BG([1, 3], \{2\})$$
with the two terms placed in degrees  $-1$ and $0$, respectively. 
More generally, $\BG([1, n], T)$ is placed in degrees $-n+2, \ldots, 0$. 
\bigskip

We go one step further, and write $\tilde{\bsigma}$ as a signed sum of maps 
$\bsigma_{T, T'}$. 
For $T\subset T'$ with $|T'|=|T|+1$, let 
$$\bsigma_{T, T'}:=(-1)^{|T_{>i}|}\sum 
(-1)^{\deg_1}\bsigma_k: \BG(I, T)\to \BG(I, T')\,\,,
\eqno{(\GIT.c)}$$
the sum over $k\in M_i$, $T'=T\cup\{i\}$. 
Then
$\bsigma_{T, T'}$ is a map of complexes.
The map $\tilde{\bsigma}$ is the sum of $(-1)^{|T_{>i}|}\bsigma_{T, T'}$. 
\bigskip 

(\GIT.1) \Lem{\it 
(1) If $i, j$ are distinct elements not in $T$, letting $T_1=T\cup\{i\}$, 
$T_2=T\cup\{j\}$, and $T''=T\cup\{i, j\}$, 
one has 
$\bsigma_{T_1\, T''}\bsigma_{T\, T_1}= \bsigma_{T_2\, T''}\bsigma_{T\, T_2}$. 

 (2) $\bsigma_{T, T'}$ is a quasi-isomorphism.
}\bigskip   
  
{\it Proof.}\quad (1) Immediate from the definition, the property 
$\bsigma\bsigma=0$, (\Defbsigma.1), and 
the sign change by $(-1)^{|T_{>i}|}$. 

(2) For each $A=(\BM_1;\cdots; \BM_n)$ of type $T$, 
 $T'=T\cup\{i\}$, the total complex of the following double complex 
 is acyclic by \Acyclicitysigma:  
 
\begin{eqnarray*}
 F(A)
&\mapr{{\bsigma}}&\moplus_{|M'_i|=1} F(\BM_1; \cdots; (M_i|M'_i);
\cdots; \BM_n) \\
&\mapr{{\bsigma}}&\moplus_{|M'_i|=2} F(\BM_1; \cdots; (M_i|M'_i);
\cdots; \BM_n)
\mapr{{\bsigma}}\cdots\to 
 F(\BM_1; \cdots; (M_i|M_i);
\cdots; \BM_n)\,\,.
\end{eqnarray*}
Taking the sum over $A$ and applying $\Phi$, the first term $F(A)$ gives $\BG(I, T)$, 
and the terms from the second to the last give $\BG(I, T')$. 
Hence the assertion. 
\bigskip 

By (1) can define, for any pair $T, T'$ with $T\subset T'$, the map 
$\bsigma_{T\, T'}: \BG(I, T)\to \BG(I, T')$ by taking a chain 
$T=T_0\subset T_1\subset \cdots\subset T_a=T'$ with $|T_{i+1}|=|T_i|+1$ and 
setting 
$$\bsigma_{T\, T'}=\bsigma_{T_{a-1}, T_a}\cdots \bsigma_{T_1, T_2}\bsigma_{T_0, T_1}\,.$$
This is a quasi-isomorphism. 
\bigskip

\sss{\GITSigma}{\it The complex $\BG(I, T|\Sigma)$.}\quad
One can refine (\GISigma) taking the type into account. 
For  $\Sigma\subset \ctop{I}$ and $T\subset \ctop{I}-\Sigma$, 
letting $I_1, \cdots, I_c$ be the corresponding segmentation of $I$ and $T_i=T\cap 
\ctop{I_i}$, let 
\begin{eqnarray*}
H\tbullet(I, T|\Sigma)&=& H\tbullet(I_1, T_1)\tildets\cdots\tildets H\tbullet(I_c, T_c)  \\
&=&\moplus F(A_1)\tildets\cdots\tildets F(A_c)\,\,,
\end{eqnarray*}
the sum over $(A_1, \cdots, A_c)$ where $A_i$ is a multi-sequence 
on $I_i$ of type $T_i$. 

With  this subcomplex we proceed as follows.
\smallskip 

$\bullet$\quad Applying the construction of (\htsPhicpx), we get the
``$2c$-tuple" subcomplex
$$\Phi(H\tbullet(I_1, T_1)\,)\tildets\cdots\tildets
 \Phi(H\tbullet(I_c, T_c)\,)\injto H{\tbullet}(I, T|\Sigma)$$
 which we abbreviate to $\Phi(H\tbullet(I, T|\Sigma)\, )$. 
 
$\bullet$\quad 
Take the totalization $\Tot_{12}\cdots\Tot_{2c-1, 2c}$ to 
get a ``$c$-tuple" complex, then turn it into a $c$-tuple complex
by changing signs of the differentials as in (\Multicpx). 

$\bullet$\quad 
Finally  shift the 
degree by $[d_1, \cdots, d_c]$, with $d_i=-|I_i|+|T_i|+2$, to obtain 
a $c$-tuple complex denoted  $\BG(I, T|\Sigma)$. So 
$$\BG(I, T|\Sigma):=( \Tot_{12}\cdots\Tot_{2c-1, 2c}\Phi(H\tbullet(I, T|\Sigma)\, )\,\,)
[d_1, \cdots, d_c]\,\,.$$

The  differential of the complex $\BG(I|\Sigma)$ comes from the maps 
$\delta$, $\bsigma_k$ on each $\BG(I_i)$ by taking tensor product, namely
$$d=\sum \pm 1\ts\cdots\ts \delta\ts 1\ts\cdots\ts 1 +\sum \pm 1\ts\cdots\ts \bsigma_k\ts\cdots \ts 1\,\,.$$
The map $1\ts\cdots\ts \bsigma_k\ts\cdots \ts 1$
takes  $\BG(I, T|\Sigma)$ to itself if $k\in M_i-M'_i$ with $i\in T$, and 
takes $\BG(I, T|\Sigma)$  to $\BG(I, T\cup\{i\} |\Sigma)$ 
if $k\in M_i-M'_i$ with $i\not\in T$. So the situation is parallel to 
that in the (\GIT).
Thus the $c$-tuple complex $\BG(I|\Sigma)$ is identified with 
the total complex of 
the $(c+1)$-tuple complex 
$$\BG(I, \emptyset|\Sigma)\mapr{\tilde{\bsigma} }\moplus_{|T|=1}
\BG(I, T|\Sigma)\mapr{\tilde{\bsigma} }
\cdots\mapr{\tilde{\bsigma}}\BG(I, \ctop{I}-\Sigma|\Sigma)\,\,$$
where $\BG(I, T|\Sigma)$ is placed in degree $-|I|+|T|+|\Sigma|+2$ and the 
differential $\tilde{\bsigma}$ is a signed sum of the maps $\bsigma_k$. 

One can define the map of complexes 
$$\bsigma_{T, T'}: \BG(I, T|\Sigma)\to 
\BG(I, T'|\Sigma)$$
 for $T\subset T'$, $|T'|=|T|+1$ by changing the sign of 
$\tilde{\bsigma}$ by $(-1)^{|T_{>i}|}$. 
 It is, up to quasi-isomorphism, 
of the form  $\pm 1\ts\cdots\ts \bsigma_{T_l, T'_l}\ts\cdots 1$ where
$i\in I_l$. 
Since $\bsigma_{T_l, T'_l}$ is a quasi-isomorphism (\GIT.1), it follows 
$\bsigma_{T, T'}$ is also a quasi-isomorphism (recall all the complexes involved are
free).  We have thus obtained the following generalization of (\GIT.1).
\bigskip

(\GITSigma.1) \Lem{\it The map $\bsigma_{T, T'}: \BG(I, T|\Sigma)\to 
\BG(I, T'|\Sigma)$ satisfies the following properties.

(1) If $i, j$ are distinct elements not in $T$, letting $T_1=T\cup\{i\}$, 
$T_2=T\cup\{j\}$, and $T''=T\cup\{i, j\}$,
one has 
$\bsigma_{T_1\, T''}\bsigma_{T\, T_1}= \bsigma_{T_2\, T''}\bsigma_{T\, T_2}$. 

 (2) $\bsigma_{T, T'}$ is a quasi-isomorphism.
}\bigskip   

As before, one can define the map $\bsigma_{T\, T'}$ for any $T\subset T'$.

The maps $\iota_\Sigma$, $\rho$ and $\Pi$ for $\BG(I|\Sigma)$ defined in previous 
subsections decompose
according to the type as follows. 
\smallskip 

$\bullet$\quad  There is an injective quasi-isomorphism 
$$\iota_\Sigma: \BG(I, T|\Sigma)\to \BG(I, T\tbar \Sigma)\,\eqno{(\GITSigma.d)}$$  
The map $\iota_\Sigma:  \BG(I|\Sigma)\to \BG(I\tbar \Sigma)$
is the sum of them as a map of graded groups. 

$\bullet$\quad
The map $\rho_m: \BG(I|\Sigma)\to \BG(I|\Sigma-\{m\})$ defined in (\rhoPi) is, as a map of graded groups, 
the direct sum of the maps  
$\rho_m: \BG(I, T|\Sigma)\to \BG(I, T\cup\{m\}|\Sigma-\{m\})$.
Observe that the degree shift is the same for the source and the target
(indeed the shift was chosen for this to hold).
For a subset $K\subset \ctop{I}$ disjoint from $T$ and $\Sigma$, 
composing $\rho_m$'s for $m\in K$ one obtains the map  
$$\rho_K:  \BG(I, T|\Sigma)\to \BG(I, T\cup K|\Sigma -K)\,.\eqno{(\GITSigma.e)}
$$
The map $\rho_K: \BG(I|\Sigma)\to \BG(I|\Sigma-K)$ is the sum of them. 
 
$\bullet$\quad
The map $\Pi_k:\BG(I|\Sigma)\to \BG(I-\{k\}|\Sigma)$ defined in (\rhoPi) is the sum of 
$\Pi_k: \BG(I, T|\Sigma)\to \BG(I-\{k\}, T|\Sigma)[-1]$. 
If $K\subset \ctop{I}$ is disjoint from $\Sigma$, 
composing $\Pi_k$'s for $k\in K$ one obtains (with $c=|K|$)
$$\Pi_K: \BG(I, T|\Sigma)\to \BG(I-K, T|\Sigma)[-c]\,\,.\eqno{(\GITSigma.f)}$$
The map $\Pi_K:\BG(I|\Sigma)\to \BG(I-K|\Sigma)$ is the sum of them. 
\bigskip

\sss{\BFIS} {\it The complex $\BF(I|S)$.}\quad
For each subset $S\subset \ctop{I}$ we define a complex of 
abelian groups $\BF(I|S)$. The construction is a variant of the 
so-called {\it bar construction}. 
For simplicity assume $I=[1, n]$. 

First we consider the case $S=\emptyset$.
As an abelian group, 
$$
\BF(I)=\BF(I|\emptyset):=\moplus_{\Sigma} \BG(I|\Sigma)\,\,,$$
where $\Sigma$ varies over subsets of $\ctop{I}$. 
The degree of $u=u_1\ts\cdots\ts u_c\in \BG(I|\Sigma)=\BG(I_1)\tildets\cdots\tildets\BG(I_c)$
is 
$$\deg_\BF(u):=\deg_{\BG}(u)-c= \sum (\epsilon_j-1)\,\quad\, \epsilon_j=\deg_\BG u_j\,.$$
The  differential $d_{\BF}$ is the sum $\bar{d}_{\BG} + \bar{\rho}$
of the maps given as follows. 
If  $I_1, \cdots, I_c$ is the partition of $I$ corresponding to 
$\Sigma$, on an element  $u=u_1\ts\cdots\ts u_c\in 
\BG(I|\Sigma)$ with $\epsilon_j=\deg(u_j)-1$,  
$$\bar{d}_{\BG}(u_1\ts\cdots\ts u_c) =-\sum (-1)^{\sum_{j>i} \eps_j}\,u_1
\ts\cdots \ts u_{i-1}\ts 
d_{\BG}(u_i)\ts\cdots \ts u_c\in \BG(I|\Sigma)\,\,,\eqno{(\BFIS.a)}$$
$$\bar{\rho}(u_1\ts\cdots\ts u_c)= \sum_{1\le i\le c-1} 
(-1)^{\sum_{j\ge i} \eps_j} \rho_{k_{i}}(u)\in \moplus_{k\in \Sigma}
 \BG(I|\Sigma-\{k\}) \,\,\eqno{(\BFIS.b)}$$
with $k_{i}=\term(I_{i})$. 
One easily verifies $\bar{d}_{\BG}\bar{d}_{\BG}=0$, $\bar{\rho}\bar{\rho}=0$, 
and $\bar{d}_{\BG} \bar{\rho}+\bar{\rho}  \bar{d}_{\BG}=0$ so that 
$d_{\BF}$ is a differential that increases $\deg_\BF$ by one. 
Note that there is an increasing filtration of $\BF(I)$ by subcomplexes, 
in which $Fil_c$ is the sum of $\BG(I|\Sigma)$ with $|\Sigma|+1\le c$. 
The graded quotient $\Gr_c^{Fil}$ is, as a group, 
the sum of $\BG(I|\Sigma)$ with given $c=|\Sigma|+1$. 
As a complex, $\BG(I|\Sigma)$ has degree given by $\deg_{\BG}(u)-c$ and differential
given by $\bar{d}_{\BG}$; we denote it by $\BG(I|\Sigma)^{\shift}$. 
Then by (\BFIS.a), $\BG(I|\Sigma)^{\shift}$ is a subcomplex of 
the tensor product complex $\BG(I_1)[1]\ts\cdots\ts\BG(I_c)[1]$. 
There is a canonical isomorphism of complexes 
$$\BG(I_1)[1]\ts\cdots\ts\BG(I_c)[1]\isoto \BG(I\tbar \Sigma)[c]\,$$
(see \RefMacLtwo, \,\S\S 8 and 9, in particular the proof of Proposition 9.3, 
for the isomorphism $K[1]\ts L[1]\to (K\ts L)[2]$ for complexes $K$, $L$);
it induces by restriction an isomorphism of complexes
$$\BG(I|\Sigma)^{\shift}\isoto \BG(I| \Sigma)[c]\,.\eqno{(\BFIS.c)}$$

If $|I|=2$, $\BF(I)$  coincides with $\BF(K, L)$ defined before. 

In general, let $|I|=n$ and set 
$$G_{i,j}=\moplus_{i=|T|,\, j=|\Sigma|} \BG(I, T|\Sigma)\,$$
for $0\le i, j\le n-2$ and $0\le i+j\le n-2$. 
Then $\BF(I)$ may be displayed as follows: 
$$\begin{array}{cccccccc} 
  &&   &  &    &       &       &G_{0,n-2}  \\
  &&   &  &    &       &        &\mapdr{\rho}  \\
  &&   &  &   &G_{0,n-3}  &\mapr{\bsigma} &G_{1, n-3} \\
  &&   &  &    &\mapd{\rho}&        &\mapdr{\rho}  \\ 
  && &G_{0,n-4} &\mapr{\bsigma} &G_{1, n-4} &\mapr{\bsigma} &G_{2, n-4} \\  
        &&   &\mapd{}  &     &\mapd{}  &    &\mapdr{}  \\       &&G_{0,1}\to\ph{\cdots\to}&\vdots&&\vdots&&\vdots                                   \\
        &&\mapd{\rho}\phantom{\to\cdots\to}   &\mapd{\rho}  &     &\mapd{\rho}  &    &\mapdr{\rho}  \\ 
G_{0,0}&\mapr{\bsigma}&G_{1,0}\to\cdots\to&G_{n-4,0} &\mapr{\bsigma} &G_{n-3, 0} &\mapr{\bsigma} &\,\,G_{n-2, 0}\,\,.
  \end{array}$$

For $S\subset \ctop{I}$, the complex $\BF(I|S)$ is the quotient complex of $\BF(I)$ given by 
$$
\BF(I|S):=\moplus_{\Sigma\supset S} \BG(I|\Sigma)\,\,,$$
where $\Sigma$ varies over subsets containing $S$. 
Note  $\moplus_{\Sigma\not\supset S} \BG(I|\Sigma)$ is a subcomplex, and 
$\BF(I|S)$ is the quotient.  
Obviously $\BF(I|\ctop{I})=\BG(I|\ctop{I})[n-1]$ if $n=|I|$. 

One has the corresponding surjection 
$$\sigma_{S\, S'}: \BF(I|S)\to \BF(I|S')$$
for $S\subset S'$; it is easy to see $\sigma_{S\, S'}$ is a quasi-isomorphism
(since all summands $\BG(I|\Sigma)$ except $\BG(I|\ctop{I})$ are acyclic by 
(\acyclicGI)\, ).  

There is an injective quasi-isomorphism $\iota_S: \BF(I|S)\to \BF(I\tbar S)$, 
defined as follows. Let $I_1, \cdots, I_c$ be the segmentation of $I$ by $S$. 
For $\Sigma\supset S$, let $\Sigma_i=\Sigma\cap \ctop{I_i}$. Then 
\begin{eqnarray*}
\BF(I\tbar S)&=& \BF(I_1)\ts\cdots\ts \BF(I_c) \\
    &=&\moplus_{\Sigma\supset S} \BG(I_1|\Sigma_1)\ts\cdots \ts \BG(I_c|\Sigma_c)
\end{eqnarray*}
using the definition of $\BF(I_i)$. 
On the other hand $\BF(I|S)=\moplus_{\Sigma\supset S}\BG(I|\Sigma)$ by definition. 
The $\iota_S$ is defined to be the sum of inclusions
 $\BG(I|\Sigma)\injto \BG(I_1|\Sigma_1)\ts\cdots \ts \BG(I_c|\Sigma_c)$.
 One easily verifies $\iota_S$ is compatible 
with the differentials. 
\bigskip 

For $K\subset \ctop{I}$ disjoint from $S$, let 
$$\vphi_K:\BF(I|S)\to \BF(I-K|S)$$
  be the sum of the maps
$\Pi_K: \BG(I|\Sigma)\to \BG(I-K|\Sigma)$
for  $K$ is disjoint from $\Sigma$, 
and the zero maps on $\BG(I|\Sigma)$ with $K\cap \Sigma=\emptyset$. 
This is a map of complexes, as can be seen from the identities (see (\rhoPi.3)\,)
$$\Pi_K\rho_m=\rho_m\Pi_K: \BG(I|\Sigma)\to \BG(I-K|\Sigma-\{m\})$$
for $m\not\in K$ and 
$\Pi_K\rho_m=0: \BG(I|\Sigma)\to \BG(I-K|\Sigma-\{m\})$
for $m\in K$. 

From the definitions it is obvious that one has 
$\sigma_{S\, S}=id$ and, for $S\subset S'\subset S''$, $\sigma_{S\,S''}=
\sigma_{S'\,S''}\sigma_{S\,S'}$.
Also, if $K=K'\amalg K''$ then 
$\vphi_K=\vphi_{K''}\vphi_{K'}:
\BF(I| S)\to  \BF(I-K| S)$.
Further, the maps $\sigma$ and $\vphi$ commute. 
\bigskip 

\sss{\DefCdzero} {\bf Definition.}\quad 
Let $(\cC^\Delta)_0$ be the category defined as follows. 
The objects are the $C$-diagrams in $\cC$. 
For $C$-diagrams $(K^m; f_K(M))$ and $(L^m; f_L(M))$, a morphism
in $(\cC^\Delta)_0$ is a collection of morphisms 
$u^m: K^m\to L^m$ in $\cC_0$ such that, if $u: F(K^M)\to F(L^M)$
is the isomorphism of complexes induced by $u$, then 
$u(f_K(M))=f_L(M)$. Composition of morphisms is defined in the obvious 
manner. 
\bigskip 

Thus $(\cC^\Delta)_0$ is the category of $C$-diagrams and isomorphisms
between them. 
The complexs $\BF(I)=\BF(K_1, \cdots, K_n)$ and $\BF(I|S)$
defined  in the previous subsection is functorial on $(\cC^\Delta)_0$. 
Note that the maps $\sigma_{S, S'}$ and $\vphi_K$ are functorial. 

Let $\tau_S: =\iota_S\sigma_S: F(I)\to F(I\tbar S)$
(in particular we have $\tau_k$). 
We now have the category $(\cC^\Delta)_0$, the complexes $\BF(I)$, 
and the maps $\tau_k$ and $\vphi_\ell$, that satisfy the conditions 
 (1),(2),(5) of (\DefqDG) (the same as (1),(2),(5) of (\wqDG)\,).
 Indeed the factorization (5) is given in 
the construction. 
 
In the following proposition, we show that the condition (6) are 
also satisfied. 
The remaining condition (3) and (4) for a quasi DG category will be shown in \S 3.
\bigskip 

\sss{\propFI} \Prop{\it 
(1) $\BF(I)$ is a complex of free $\ZZ$-modules and 
$\BF(I|S)$ is a multiple complex of free $\ZZ$-modules. 
The category $(\cC^\Delta)_0$, the complexes $\BF(I)$, and the maps 
 the maps $\tau_k$ and $\vphi_\ell$ satisfy the conditions 
(1),(2),(5) of (\DefqDG).

(2) The condition (6)
of (\DefqDG) is satisfied: 
Let $R$, $J$ be disjoint subsets of $\ctop{I}$, with $J$ non-empty. 
Then the following  sequence of complexes is exact 
(the maps are alternating sums of 
the quotient maps $\sigma$) 
$$\BF(I|R)\mapr{\sigma}
\moplus_{\stackrel{S\subset J}{|S|=1} } \BF(I| R\cup S)
\mapr{\sigma} \moplus_{\stackrel{S\subset J}{|S|=2}} 
\BF(I|R\cup S)
\mapr{\sigma} \cdots \to  \BF(I|R\cup J)\to 0\,\,.$$
}\bigskip

\sss{\Proofacyclicity} {\bf Proof of (\propFI), (2).}\quad
We verify (\propFI), (3), that is the acyclicity axiom of the map 
$\bsigma$. 

We begin with noting that the complex $\BF(I|S)$ has a decreasing 
filtration by subcomplexes, denoted  $Fil^\Sigma$,  indexed by sets $\Sigma$
containing $S$, given by 
$$Fil^\Sigma=\moplus_{\Sigma'\supset\Sigma} \BG(I|\Sigma')$$
(as a group). 
The graded quotient at $\Sigma$ is 
$\Gr_{Fil}^\Sigma \BF(I|S)=\BG(I|\Sigma)[\, |\Sigma|+1\, ]$. 
The map $\sigma_{S\, S'}: \BF(I|S)\to \BF(I|S')$ preserves the filtration.

Thus the sequence 
in  (\propFI), (3) also has a filtration indexed by $\Sigma$ containing
$R$. For the term $\BF(I|R\cup S)$, the corresponding graded quotient 
$\Gr_{Fil}^\Sigma$ is non-zero if and only if 
$\Sigma\supset R\cup S$, equivalently, 
if $S\subset (\Sigma -R)\cap J$. 
So the graded quotient of the sequence in question is (setting  $T:=(\Sigma -R)\cap J$\,):
$$\BG(I|\Sigma) \mapr{\sigma}
\moplus_{\stackrel{S\subset T}{|S|=1} } \BG(I|\Sigma)
\mapr{\sigma} \moplus_{\stackrel{S\subset T}{|S|=2}} 
\BG(I|\Sigma)
\mapr{\sigma} \cdots \to  \BG(I|\Sigma)\to 0\,\,.$$
Here the maps $\sigma$ are the alternating sums of the identity maps. 
If $T$ is non-empty, 
the exactness of this sequence (even with ($0\to$) at left)
 is consequence of Lemma (\acyclicGI.1).
If $T=\emptyset$, the sequence is trivially exact. 
\bigskip

\section{The quasi DG category $\cC^\Delta$.}

\sss{\PlanCd} We keep the assumption from the previous section, so
 $\cC$ is a quasi DG category with additional structure (iv) and (v). 
It has been shown that the complex  
$\BF(K_1, \cdots, K_n)$ and maps $\tau$, $\vphi$ 
satisfy the conditions of (\DefqDG), (1),(2) and (5). 
In this section we will verify the the remaining conditions
in (\ZFKL)-(\AdditivityBF), 
establishing Theorem (\ThmCdqDG), 
which says that the $C$-diagrams in $\cC$ form 
a quasi DG category.  

We further give another theorem (\ThmCdqDGtriang), which states that 
the homotopy category of the quasi DG category of $C$-diagrams has 
the structure of a triangulated category. 
The proof takes the rest of this section (\Shiftfunctor)-(\Axiomprooftwo). 
\bigskip

\sss{\ZFKL} {\it The complex $\BF(K, L)$.}\quad
We will  examine the composition in the homotopy category. 
Let us first recall  $\BF(K, L)$ is given by
$$\BF(K, L)=\BG(K, L)[1]\,, \quad \BG(K, L)=\Phi H\dbullet \subset H^{\bullet\, 0}\, $$
where 
$$H^{\bullet\, 0}=\moplus F(M|\emptyset; N|\emptyset)$$
the sum over double sequences $(M|\emptyset;  N|\emptyset)$.
We will abbreviate $(M|\emptyset; N|\emptyset)$ to $(M; N)$.  
So an element of $H^{\bullet\, 0}$ is of the form 
$u= (u(M; N)\, )\in \bop F(M; N)$. 
The differential $d_{\BF}$ acting on $u$ is nothing but $\delta=\bpartial + \bvphi$ as defined in \S \secfcn. 
The $u$  has $\deg_\BF=0$ if
$u(M; N)\in F(M; N)^{-\gamma(M; N)}$. 
We will often make use of the following conditions.
\smallskip 

(i) An element $u\in H^{\bullet\, 0}$ is in  $\BF(K, L)^{\bullet-1}$ 
if the following condition ({\it $\sigma$-consistency})
 is satisfied.
For $k\in M$, $k\neq \init(M)$, 
$$\sigma_k(u(M; N))=f_K(M_{\le k})\ts u(M_{\ge k}; N)\,\,.$$
For
$k\in N$, $k\neq \term(N)$, 
$$\sigma_k(u(M; N)\,)=\{(\deg_1 u(M; N))\cdot\gamma(N_{\ge k})\}\, 
u(M; N_{\le k})\ts f_L(N_{\ge k})\,\,$$
(recall $\{i\}:=(-1)^i$ for $i\in \ZZ$). 

(ii) It is $\delta$-closed if
$$\partial (u(M; N)\,)+\sum \bvphi_k(u(M\cup\{k\}; N)\,)
     + \sum \bvphi_k(u(M; N\cup\{k\})\,)=0\,\,.$$
Here $k$ in the first sum varies over the set $[\init{M}, \infty)-M:=\{k\in \ZZ\,| \,k\ge \init(M)\,\,\text{and}\,\, k\not\in M\}$, 
and $k$ in the second sum over $(-\infty, \term{N}]- N$. 
\smallskip 
 
Recall from \S \secDG\, we have the homotopy category $Ho(\cC^\Delta)$, 
where the morphisms are given by 
$\Hom_{Ho}(K, L)=H^0\BF(K, L)$. 
A morphism $u: K\to L$ in the homotopy category 
is represented by a cocycle of degree 0, namely by an 
element $\underline{u} \in Z^0\BF(K, L)$. 
\bigskip

\sss{\GKLM}{\it The complex $\BF(K, L, M)$.}\quad 
Let $K$, $L$, $M$ be three $C$-diagrams. The complex 
$\BF(K, L, M)$ is of the form 
$$\begin{array}{ccc}
& &\BG([1, 3], \emptyset|\{2\}) \\
& &\mapdr{\bar{\rho} }\\
\BG([1, 3], \emptyset)&\mapr{-\tilde{\bsigma} } &\BG([1, 3], \{2\})\,\,
\end{array}
$$
which we also write 
$$\begin{array}{ccc}
& &\BG(K, L)\tildets \BG(L, M) \\
& &\mapdr{\bar{\rho} } \\
\BG(K, L, M, \emptyset)&\mapr{-\tilde{\bsigma}} &\BG(K, L, M, \{2\})\,\,.
\end{array}
\eqno{(\GKLM.a)}
$$

$\bullet$\quad The differential $d_\BF$ is equal to $\bar{d}_{\BG}+\bar{\rho}$, to be specified below. 

$\bullet$\quad Recall that the complex $\BG(I)$ has differential 
$u\mapsto d_{\BG}(u)=-\delta(u)+(-1)^{\deg_1 u}\bsigma_{int}(u)$. 
Also, $\BG(I)$ is the direct sum of $\BG(I, T)$, each 
$\BG(I, T)$ is a complex with differential $d_{\BG(I, T)}$, and  on $\BG(I, T)$
$$d_{\BG}=\sum (-1)^{|T|+|I|}d_{\BG(I, T)} +\tilde{\bsigma}\,\, $$
where 
$\tilde{\bsigma}=\sum (-1)^{\deg_1} \bsigma_k$,
 the sum over $k\in M_i$ with $i\not\in T$ (see (\GIT.b)\,). 
In particular, if $I=[1, 3]$,  $\BG([1, 3])$ is of the form 
$$\BG([1, 3], \emptyset)\mapr{\tilde{\bsigma}} \BG([1, 3], \{2\})\,\,.$$ 

$\bullet$\quad For $u\ts v\in \BG(K, L)^\bullet\tildets \BG(L, M)^\bullet$, 
$\deg_\BF(u\ts v)=\deg_\BG(u)+\deg_\BG(v)-2$, and 
$d_\BF(u\ts v)=\bar{d}_\BG(u\ts v)+\bar{\rho}(u\ts v)$, where 
$$\bar{d}_\BG(u\ts v) =-(-1)^{\deg_\BG(v)-1}d_\BG(u)\ts v- u\ts d_{\BG}(v)\,\,,$$
$$\bar{\rho}(u\ts v)=(-1)^{\deg_\BG(v)-1} \rho(u\ts v)\,\,.$$

Note if $u\ts v\in \BG(K, L)^1\tildets \BG(L, M)^1$ the signs simplify as follows:
$$\bar{\rho}(u\ts v)=\rho (u\ts v)\,,$$
and, since $\deg_1 u=\deg_1 v=0$, 
$$\rho (u\ts v)=u\ts v$$
if $\term(A)=\init(B)$, 
and zero otherwise. 
 
$\bullet$\quad 
For $W\in \BG(\emptyset):=\BG(K, L, M, \emptyset)$, $\deg_\BF(W)=\deg_\BG(W)-1=\deg_{\BG(\emptyset)}(W)-2$, and 
$$d_\BF(W)=-d_\BG(W)= d_{\BG(\emptyset)}(W)-\tilde{\bsigma}(W)\,\,.$$

$\bullet$\quad For $z\in\BG(\{2\}):=\BG(K, L, M, \{2\})$, $\deg_\BF(z)=\deg_\BG(z)-1=\deg_{\BG(\{2\})}(z)-1$, and 
$$d_\BF(z)=-d_\BG(z)= -d_{\BG(\{2\})}(z)\,\,.$$

In particular, for a degree 0 element of $\BF(K, L, M)$ of the form 
$(u\ts v, W, z)$ with 
$$ u\ts v\in \BG(K, L)^1\tildets \BG(L, M)^1\,,\,\,
W\in \BG(K, L, M, \emptyset)^1\,,\,\, z\in \BG(K, L, M,\{2\})^{2}\,\,, $$
 one has 
\begin{eqnarray*}
d_\BF(u\ts v, W, z)&=& -d_\BG(u)\ts v- u\ts d_{\BG}(v) + \rho(u\ts v)\\
&&+ d_{\BG(\emptyset)}(W)-\tilde{\bsigma}(W) \\
&&- d_{\BG(\{2\})}(z)\,\,.
\end{eqnarray*}
We also recall  there is a map of complexes  $\Pi_2= \Pi_L: \BG(K, L, M, \emptyset)[1]\to \BG(K, M)$. 
\bigskip 

\sss{\Compos} \Prop{\it Let $u:K\to L$ and $v: L\to M$ be morphisms. 
Assume $u$ is represented  by $u\in \BG(K, L)^1$,  $v$  is 
represented by $v\in \BG(L, M)^1$, and $u\ts v\in \BG(K, L)\tildets \BG(L, M)$. 

(1) There are an element $W\in \BG(K, L, M, \emptyset)^2$, $d$-closed in $\BG(K, L, M,\emptyset)$ such that 
$$\tilde{\bsigma}(W)=\rho(u\ts v)\,\,.$$
In this case $(u\ts v, W)\in \BF(K, L, M)^0$ is $d_\BF$-closed. 
 
(2) If (1) is satisfied,  the element
 $\Pi_L(W)\in \BG(K, M)^1$ is $d$-closed and represents
the morphism $u\cdot v :K\to M$. 
}\smallskip 

{\it Remark.}\quad If $W\in \BG(K, L, M, \emptyset)^2$, then 
$\deg_1 W=0$.
\smallskip 

{\it Proof.}\quad 
(1) For $i\ge 0$, let 
$$F_i =\moplus_{|M'_2|=i} F(K^{\BM_1}; L^{\BM_2}; M^{\BM_3})\,,$$
the sum over triple sequences $(\BM_1; \BM_2; \BM_3)$. 
It is a ``double"  complex with respect to the differentials 
$-\delta$ and $d_3$ of (\HIGI). 
One has a sequence of ``double" complexes 
(with $c$ sufficiently large)
$$F_0\mapr{\bsigma'} F_1\mapr{\bsigma'} \cdots \to  F_c\to 0$$
where the maps $\bsigma'$ are $\sum (-1)^{\deg_1} \bsigma_k$, the sum 
over $k\in M_2-M'_2$; note the first $\bsigma': F_0\to F_1$ is $\ti\bsigma$ of 
(\GIT). 
By (\DefqDG), (6) and the remark to it, 
this sequence is exact and induces a surjective 
quasi-isomorphism 
$$\ti{\bsigma}: F_0\to \Ker(\bsigma': F_1\to F_2)\,.$$

Applying the exact functor $\Phi$ of (\Phicpx), we obtain an exact sequence for
$\BG(K, L, M)_i := \Phi (F_i)$, 
$$\BG(K, L, M)_0\mapr{\bsigma'} \BG(K, L, M)_1\mapr{\bsigma'} \cdots \to  \BG(K, L, M)_c\to 0\,.$$
Thus we get a surjective quasi-isomorphism 
$$\ti{\bsigma}: \BG(K, L, M)_0\to \Ker[\bsigma': \BG(K, L, M)_1\to \BG(K, L, M)_2]\,.$$

Now $\rho(u\ts v)$ is a cocycle 
in $\Ker [\bsigma': \BG(K, L, M)_1\to \BG(K, L, M)_2)]$.
Thus we obtain the assertion by the following fact:
If $f: A\to B$ is a surjective quasi-isomorphism of complexes, then 
the induced map $f: Z^nA\to Z^nB $ on cocycles is surjective. 

(2) By definition the composition is defined by lifting the class
$[u\ts v]$ to a cohomology class in $H^0\BF(K, L, M)$, then mapping 
by $\vphi_L$ to $H^0\BF(K, M)$. So the assertion follows. 
\bigskip 

\sss{\Identity}{\it The identity map.}\quad 
Let $K=(K^m; f(M))$ be an object. 
For a non-decreasing 
sequence $M=(m_1, \cdots, m_\mu)$, $m_1\le m_2\le \cdots \le m_\mu$, 
with $\mu\ge 2$, 
let $M'$ be the increasing sequence 
obtained by eliminating repetitions.  
Thus if, say, $M=(1,2,2,4,5,5)$, then $M'=(1,2,4,5)$. 
Set $\gamma(M)=\gamma(M')$. 
Assume first that $M$ is a non-constant sequence. 
There is a natural surjection $\lambda: M\to M'$, and one has the corresponding
diagonal extension map (\DefqDG), (iv). 
$$\lambda^*: F(K; M' )\to F( \lambda^*K; M)\,;$$
this is also denoted $\diag: F(M')\to F(M)$. 
Let $f(M)\in F(M)^{-\gamma(M)}$ be the image of $f(M')\in F(M')^{-\gamma(M)}$. 
We say that $f(M)$ is obtained from $f(M')$ by means of diagonal extension.
If $M$ is a constant sequence $M=(m, \cdots, m)$, we take 
$$f(M)=\bDelta_M\in F(M)=F(K^m, \cdots, K^m)\,, $$
the diagonal element in (\DefqDG), (iv).

These elements
$$f(M)=f(m_1, \cdots, m_\mu)\in F(K^{m_1}, \cdots, K^{m_\mu})^{-\gamma(M)}$$
satisfy the following properties. 
They follow from the definition and the compatibility of $\lambda^*$ with 
$\sigma$ and $\vphi$. 
\smallskip 

(i) For each $M$, 
$$\partial f(M)+ \sum_a \bvphi_a(f(M\cup\{a\})\, )=0\,\,$$
where $a$ varies over $[\init(M), \term(M)]-M= [\init(M'), \term(M')]-M'$. 

(ii) For $k=m_i$ with $k\neq m_1, m_\mu$, one has 
$$\sigma_{m_i}(f(M)\, )= f(m_1, \cdots, m_i)\ts f(m_i,\cdots, m_\mu)\,.$$
which we also write 
$\sigma_k (f(M))= f(M_{\le k})\ts f(M_{\ge k})$ for short, 
when there is no danger of confusion. 

(iii) For $k=m_i$, $k\neq m_1, m_\mu$, 
 $$\vphi_{m_i}(f(m_1, \cdots, m_\mu)\, )= f(m_1, \cdots, \widehat{m_i}, \cdots, m_\mu)\,\,.$$

(iv) If $m_1=\cdots =m_\mu=m$, then $f(m, \cdots, m)=\bDelta(m, \cdots, m)\in F(K^{m_1}, \cdots, K^{m_\mu})$. 
\bigskip  

For a free double sequence 
$(M; N)=(m_1, \cdots, m_\mu; n_1, \cdots, n_\nu)$ with 
$m_1<\cdots<m_\mu$, $n_1<\cdots<n_\nu$, define 
the element $\ti f (M; N)$ in $F(M; N)^{-\gamma(M; N)}$ by 
\begin{eqnarray*}
\ti f (M; N)&=&\ti f(m_1, \cdots, m_\mu; n_1, \cdots, n_\nu) \\
&=&\begin{cases}
(-1)^{n_1}
f(m_1, \cdots, m_\mu,  n_1, \cdots, n_\nu) 
 &\text{if } m_\mu=n_1\\
0  &\text{if } m_\mu\neq n_1\,\,.
\end{cases}
\end{eqnarray*}
One may simply write $\ti f(m_1, \cdots, m_\mu; n_1, \cdots, n_\nu)
=(-1)^{n_1}\delta_{m_\mu\, n_1}f(m_1, \cdots, m_\mu, n_1, \cdots, n_\nu)$. 
Note the repetition of indices can occur only in the first case, and then
 for $m_\mu$ and $ n_1$.
The collection $(\ti f(M; N)\, )$, as $(M; N)$ varies,
gives an element in $H^{1, 0}(K, K)$. 
\bigskip 

(\Identity.1)\Prop{\it The element $(\ti f(M; N)\,)$ is contained  in $   \BF(K, K)^0$
and $\delta$-closed. 
}\bigskip 

{\it Proof.}\quad We verify the two conditions in (\ZFKL).
The first condition is obvious. 
To show the identity
$$\partial  (\ti f(M; N)\, )+\sum \bvphi_k (\ti f(M\cup\{k\}; N)\, )
+ \sum \bvphi_k (\ti f(M; N\cup\{k\})\, )=0\,\,,$$
there are cases   $m_\mu>n_1$, 
$m_\mu=n_1$ and  $m_\mu<n_1$. 
If $m_\mu>n_1$, all the three terms are zero, so the identity trivially holds. 
If $m_\mu=n_1$ it holds by property (i) for $f(M)$. 
If $m_\mu<n_1$, the first term is zero and the last two terms are
$$\bvphi_{n_1}(\{n_1\}f(m_1, \cdots, m_\mu, n_1, n_1, \cdots, n_\nu)\, )
+\bvphi_{m_\mu}(\{m_\mu \}f(m_1, \cdots, m_\mu, m_\mu, n_1, \cdots, n_\nu)\, )
$$
which is zero, since by properties (ii) and  (iii) above, we have
$$\bvphi_{n_1}(f(m_1, \cdots, m_\mu, n_1, n_1, \cdots, n_\nu)\, )
=\{-\gamma(M, N)+ \gamma(m_1, \cdots, m_\mu, n_1)\} f(M, N)$$
and 
$$\bvphi_{m_\mu}(f(m_1, \cdots, m_\mu, m_\mu, n_1, \cdots, n_\nu)\, )
=\{-\gamma(M, N)+ \gamma(m_1, \cdots, m_\mu)\} f(M, N)\,\,.$$
(Note $\gamma(m_1, \cdots, m_\mu, n_1)- \gamma(m_1, \cdots, m_\mu)
=n_1-m_\mu -1$.)
\bigskip 

Let $\iota_K: K\to K$ be the morphism represented by $(\ti f(M; N)\,)$. 
The following shows it is the identity map. 
\bigskip 

\sss{\IdentityProp} \Prop{\it For any morphism $u: K\to L$, one has 
$\iota_K\cdot u= u$. Similarly for any $u: K\to L$, $u\cdot \iota_L=u$.
}\bigskip 

{\it Proof.}\quad 
Let $u$ be represented by $\underline{u}=(\underline{u}(M; N))\in  \BF(K, L)^0$. 
Then by the definition of $\rho$ recalled in (\GKLM), 
$\rho(\ti f\ts \underline{u})$ equals 
$$\sum_{\term{N}=\init{\bar{N}} }\ti f(M; N)\ts \underline{u}(\bar{N}; R)
\in  \bop F(K^M; K^N)\tildets F(K^{\bar{N} }; L^R)
\subset H^{0,2, 0}(K, K, L, \{2\})\,\,.$$
Here $(M; N)$ and $(\bar{N}; R)$ are free double sequences. 
The restriction $\term{N}=\init{\bar{N}}$ occurs since $\rho$ acts as identity in that case, and 
as zero otherwise. 

\vspace*{0.5cm}
\hspace*{2cm}
\unitlength 0.1in
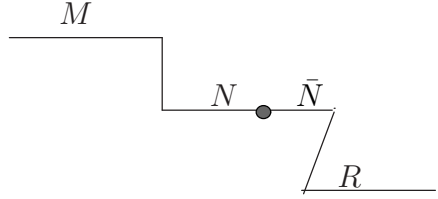
\begin{picture}( 22.3000, 10.5800)( 17.1000,-18.1800)
%
\special{pn 8}%
\special{pa 1710 1000}%
\special{pa 2510 1000}%
\special{fp}%
%
\special{pn 8}%
\special{pa 2510 1000}%
\special{pa 2510 1380}%
\special{fp}%
%
\special{pn 8}%
\special{pa 3410 1370}%
\special{pa 3410 1370}%
\special{fp}%
%
\special{pn 8}%
\special{pa 2510 1380}%
\special{pa 3400 1380}%
\special{fp}%
%
\special{pn 8}%
\special{pa 3410 1390}%
\special{pa 3250 1818}%
\special{fp}%
%
\special{pn 8}%
\special{pa 3240 1800}%
\special{pa 3940 1800}%
\special{fp}%
\put(19.7000,-9.3000){\makebox(0,0)[lb]{$M$}}%
\put(27.6000,-13.6000){\makebox(0,0)[lb]{$N$}}%
\put(34.2800,-17.8000){\makebox(0,0)[lb]{$R$}}%
%
\special{pn 8}%
\special{sh 0.600}%
\special{ar 3040 1390 38 34  0.0000000 6.2831853}%
\put(32.1000,-13.6000){\makebox(0,0)[lb]{$\bar{N}$}}%
\end{picture}%
\vspace{0.8cm}

\noindent  
The element $\rho(\ti f\ts \underline{u})$ is $\sigma$-consistent
(the condition (\ZFKL), (i) is obviously satisfied), namely contained in 
$\BG(K, K, L, \{2\})^1$, and  $d$-closed (as the image of the 
closed element $\ti f\ts \underline{u}$).

For each pair $(M; N)$ of non-decreasing sequences $M=(m_1\le  \cdots \le m_\mu)$
and $N=(n_1\le \cdots \le n_\nu)$ with $\mu\ge 1$, $\nu\ge 1$, 
 one has the complex $F(M; N)=F(K^M; L^N)$ just as we obtained 
$F(M)$ in (\Identity).  
To be specific, if $M'$ (resp. $N'$) is the sequence obtained from $M$ (resp. $N$)
 by eliminating repetitions, and $\lambda : M\to M'$ (resp. $\mu: N\to N'$) 
 is the natural surjection, we let 
 $$F(M; N)=F(K^M; L^N):=F((\lambda^*K)^M; (\mu^* L)^N)\,.$$
 There is the diagonal extension map 
$\diag: F(M'; N')\to F(M; N)$. 
Thus we have elements $\underline{u}(M; N)$ as the image of $\underline{u}(M'; N')$, 
and they satisfy properties similar to those for $f(M)$ in (\Identity). 

Let   
$$W=(W(M; N; R)\,)\in \moplus F(K^M; K^N; L^R)\subset H^{0,1,0}(K, K, L, \emptyset)
$$
where for each  free triple sequence $(M; N; R)$, 
\begin{eqnarray*}W(M; N; R)=\underline{u}(M_{\vartriangle} N; R)  
:=\begin{cases}
(-1)^\ell \underline{u}(M, N; R) &\text{ if } \term{M}=\ell=\init{N}\,\,,\\
0&\text{ if }  \term{M} \neq \init{N}\,\,.
\end{cases}
\end{eqnarray*}
The sign $_\vartriangle$ indicates applying diagonal extension 
at $K^\ell$ if
$\term{M}=\ell=\init{N}$ and putting sign $(-1)^\ell$. 
It is obvious that $W$ is $\sigma$-consistent in the sense of (\HIGI.3), namely 
\begin{eqnarray*}
\sigma_k (W(M; N; R)\, )         
&=&\begin{cases} 
  f(M_{\le k})\ts   W(M_{\ge k}; N; R)   &\text{ if }  k\in M-\{\init{M} \}\,,  \\
  W(M; N; R_{\le k})\ts f_L(R_{\ge k})    &\text{ if }  k\in R-\{\term{R}\}\,\,.
\end{cases}
\end{eqnarray*}
Also, $W$ is $\delta$-closed, as can be shown by the same argument 
as in the proof of the above proposition. 

Further one has 
$$\tilde{\bsigma}(W)=\rho(\ti f\ts \underline{u})\,\,.$$
Indeed, according to the definition of $\tilde{\bsigma}$ in (\GIT), 
$$\tilde{\bsigma}(W)=\sum_{(M; N; R)}\sum_k (-1)^{\deg_1 W} \bsigma_k(W(M; N; R)\, )$$
where $k$ varies over $k\in N$.
Since 
$$\sigma_k (W(M; N; R))
=\ti{f}(M; N_{\le k})\ts \underline{u}(N_{\ge k}; R)$$
for  $k\in N$, and since $\deg_1 \ti{f}(M; N_{\le k})=0$, one has 
$$\bsigma_k (W(M; N; R))=\sigma_k (W(M; N; R))\,.$$
Thus 
$$
\tilde{\bsigma}(W)=\sum_{(M; N; R)}\sum_k \ti{f}(M; N_{\le k})\ts \underline{u}(N_{\ge k}; R)\,.
$$
The right hand side coincides with $\rho(\ti f\ts \underline{u})$
in view of the description of $\rho(u\ts v)$ in (\GKLM). 
 
By Proposition (\Compos),  $\iota\cdot u$ is represented by $\Pi_L(W)(M; R)$, and they are:
\begin{eqnarray*}
\Pi_L(W)(M; R)&=&(-1)^\ell \vphi_\ell (W(M; \{\ell\}; R)\, )    \\
&=& \underline{u}(M; R)\,\,.
\end{eqnarray*}
Thus one has $\iota\cdot u=u$. 
The argument for the second statement is similar. 
\bigskip

\sss{\Sigmaprolong} {\it $\sigma$-consistent prolongation}. \quad 
For $a=1, 2$,  assume given a finite sequence of objects $k\mapsto X^k_a$. 
Let $X=X_1\oplus X_2$ be the direct sum sequence, namely the
sequence of objects given by $X^k=X^k_1\oplus X^k_2$. 
  For a finite set of integers $M$ with cardinality $\ge 2$ and a map $\al: M\to \{1, 2\}$, 
let ${X^M_\al}$ denote the sequence of objects $ X^k_{\al(k)}$ 
indexed by $k\in M$, and let 
 $F({X^M_\al}):=F(  X_{\al }; M \,)$ be the corresponding complex. 
One also has the complex $F(X^M)=F( X; M)$ associated to the sequence 
$X^k$ on $M$. 
Generalizing the map $\theta$ of (\DefqDG), we give a quasi-isomorphism
$$\theta:\quad \moplus F(X_{\al_1}^{M_1})\ts\cdots\ts F(X_{\al_c}^{M_c})
 \to F(X^M)
\eqno{(\Sigmaprolong.a)}$$
where the sum is over all segmentations $M_1, \cdots, M_c$ of $M$ and 
functions $\al_i:M_i\to \{1, 2\}$ satisfying the following 
condition:
\smallskip 

$$\text{at each}\quad k=\term{M_i}=
\init{M_{i+1}},\,\, 1\le i< c, \quad  
\al_i(k)\neq\al_{i+1}(k)\,.\eqno{(\Sigmaprolong.b)}$$ 
In other words, the maps $\al_i$ take distinct values at the overlaps of the segmentation. 
This map $\theta$ on the summand 
$F(X_{\al_1}^{M_1})\ts\cdots\ts F(X_{\al_c}^{M_c})$
is obtained as follows. 
Let $S\subset M$ correspond to the segmentation $M_1, \cdots, M_c$. 
On each $M_i$, consider a new sequence of objects $\ti X$ given by 
$$\ti X^k=\left\{
\begin{array}{ll}
X^k_{\al(k)}  &\mbox{ if } k\in S\,, \\
X^k\,&\mbox{ if } k\not\in S\,
\end{array}
\right.$$
(keep the $X^k_{\al(k)}$ for $k\in S$ and replace it with 
$X^k$ for $k\not\in S$). 
Recall from (\DefqDG), (3) that we have maps $s$ and $t$ 
for $X^k=X^k_1\oplus X^k_2$. 
Let 
$$F(X_{\al_i}^{M_i})\to F(\ti X^{M_i}), \quad u\mapsto u'\,,$$
be the map of complexes obtained by composing (in any order)
$s_k$ or $t_k$ (according as $\al(k)=1, 2$) 
for $k\in \ctop{M}_i$. 
One has the tensor product of these maps, 
$$F(X_{\al_1}^{M_1})\ts\cdots\ts F(X_{\al_c}^{M_c})
\to F(\ti X^{M_1})\ts\cdots\ts F(\ti X^{M_c}), \quad 
u_1\ts\cdots\ts u_c\mapsto u'_1\ts\cdots\ts u'_c\,.$$
Compose this with the map 
$$\pi: F(\ti X^{M_1})\ts\cdots\ts F(\ti X^{M_c})
\to F(X^M)$$
which is the composition (in any order) of $\pi_k$ 
in (3) of (\DefqDG) for $k\in S$. 
We thus have the map 
$$\theta: F(X_{\al_1}^{M_1})\ts\cdots\ts F(X_{\al_c}^{M_c})
\to F(X^M)\,,$$
$\theta(u_1\ts\cdots\ts u_c)=\pi (u'_1\ts\cdots\ts u'_c)$. 
Taking the sum of these, we have the map $\theta$ of (\Sigmaprolong.a).

Note that $\theta$ is compatible with the maps $\sigma$ and $\vphi$.
The compatibility with $\sigma$ means, 
for $u_1\ts\cdots\ts u_c\in F(X_{\al_1}^{M_1})\ts\cdots\ts F(X_{\al_c}^{M_c})$, 
one has 
$$\sigma_k(\theta(u_1\ts\cdots\ts u_c)\,)=
\theta(u_1\ts\cdots\ts u_i)\ts \theta(u_{i+1}\ts\cdots\ts u_c)
$$
if $k=\term(M_i)$, and 
$$\sigma_k(\theta(u_1\ts\cdots\ts u_c)\,)=
\sum \theta(u_1\ts\cdots\ts u_{i-1}\ts u')\ts \theta(u''\ts\cdots\ts u_c)$$
if $k\in \ctop{M_i}$ and 
$ \sigma_k(u_i)=\sum u'\ts u''$ with $u'\in F((M_i)_{\le k})$, 
$u''\in F((M_i)_{\ge k})$.

In particular, an element $u\in F(X^M_\al)$ for $\al: M\to \{1, 2\}$
(namely an element in a summand with $c=1$) -- called a
{\it primary component} -- gives rise to the element $u'\in F(X^M)$. 

It is convenient to introduce a variant  of the complex 
$\bop F(X^M)$ in (\FMM). 
We consider the group,  called the {\it primary part},
defined by 
$$\moplus_{(\al, M)} F(X^M_\al)\,,$$
the sum over sequences 
$M$ and $\al: M\to \{1, 2\}$; this is a double complex where the 
two differentials are 
$$\partial: F(X^M_\al)\to F(X^M_\al)$$
and the sum of 
$$\bvphi_k: F(X^M_\al)\to F(X^{M-\{k\}}_{\al'})$$
for $k\in \ctop{M}$ 
where $\al'$ is the restriction of $\al$ to $M-\{k\}$. 
From this one obtains a complex with degree given by 
$\deg_1(u)=|u|+\gamma(M)$ and 
differential $\delta=\partial +\bvphi$. 

An element  $(f(\al, M))$ of in $\bop F(X^M_\al)$ is said 
{\it $\sigma$-consistent} if 
for each pair $(\al, M)$ and 
 $k\in \ctop{M}$, one has 
$$\sigma_k(f(\al, M)\,)=f(\al_{\le k}, M_{\le k})\ts f(\al_{\ge k}, M_{\ge k})\,\,,
$$
where $\al_{\le k}$ (resp.  $\al_{\ge k}$) is the restriction of $\al$ to $M_{\le k}$ (resp. $M_{\ge k})$). 
We say $(f(\al, M))$ is $\delta$-closed if it is closed in the complex 
$\bop F(X^M_\al)$, namely if for each $(\al, M)$, 
$$\partial f(\al, M)+\sum_{(k, \tilde{\alpha})}\bvphi_k\left( f(\tilde{\alpha}, M\cup\{k\})
\right)=0$$
where the sum is over elements 
$k\in [\init M, \term M]-M$ and maps $\tilde{\alpha}: M\cup\{k\}
\to \{1, 2\}$ extending $\al$. 
One has:
\bigskip 

(\Sigmaprolong.1)\Prop{\it (1) Assume given a $\sigma$-consistent set of elements 
of first degree zero, 
\newline $(f(\al, M))$ in $\bop F(X^M_\al)$. Define  
 elements $f(M)\in F(X^M)$ by 
$$f(M)= \theta\left(\sum f(\al_1, M_1)\ts\cdots\ts f(\al_c, M_c)\,\right)\,\,,$$
the sum over  segmentations $M_1, \cdots, M_c$ and functions 
$\al_i: M_i\to \{1, 2\}$ satisfying (\Sigmaprolong.b). 
Then $f(M)\in F(X^M)$ is $\sigma$-consistent, namely it satisfies for 
$k\in \ctop{M}$,
$\sigma_k ( f(M)\, )= f(M_{\le k})\ts f(M_{\ge k})$ in $F(M_{\le k})\ts F(M_{\ge k})$. 
 We call $(f(M)\,)$ the {\rm $\sigma$-consistent prolongation} of $(f(\al, M)\,)$.

(2) Assume in addition the element  $(f(\al, M)\,)\in 
\bop F({X^M_\al})$ is $\delta$-closed.  Then $(f(M))\in \bop F(X^M)$ 
is $\delta$-closed.}
\smallskip 

{\it Proof.}\quad 
(1) This follows from the compatibility of $\sigma$ with the map $\theta$. 

(2) Let $M_1, \cdots, M_c$ be a segmentation of $M$, and  
$\al_i: M_i\to \{1, 2\}$ be functions 
satisfying the condition (\Sigmaprolong.b). 
Assume given elements 
$u_i\in F(X_{\al_i}^{M_i})$. Then one has 
$u_1\ts\cdots\ts u_c\in F(X_{\al_1}^{M_1})\ts\cdots\ts F(X_{\al_c}^{M_c})$. 
One verifies the following identities using the definitions and 
the compatibility of the map $\vphi$ with the decomposition (\Sigmaprolong.a):
$$\partial(\, \theta(u_1\ts\cdots\ts u_c)\,)=
\theta\left(\sum_i (-1)^{\nu_i} u_1\ts\cdots\ts u_{i-1}\ts(\partial u_i)
\ts u_{i+1}\ts \cdots \ts u_c\right)$$
with $\nu_i=\sum_{j>i}|u_j|$, 
and 
$$\bvphi\left(\theta(u_1\ts\cdots\ts u_c)\,\right)=\theta\left(\sum(-1)^{\mu_i} u_1\ts\cdots\ts u_{i-1}\ts \bvphi( u_i)
\ts u_{i+1}\ts \cdots \ts u_c\,\right)$$
with $\mu_i=\sum_{j<i}(\deg_1 u_j)+\sum_{j>i} |u_j|$.
The latter is shown by a repeated application of the equality
(which follows from the definition)
$$\bvphi(\theta(u\ts v))=(-1)^{|v|} \theta(\bvphi(u)\ts v)
+(-1)^{\deg_1 u} (u\ts \bvphi(v)\,)\,.$$
Thus if $\deg_1 (u_i)=0$ for $i<c$, then 
$\nu_i=\mu_i$ for all $i$ and 
$$\delta\left(\,\theta(u_1\ts\cdots\ts u_c)\right)=\theta\left(\sum_i (-1)^{\nu_i} u_1\ts\cdots\ts \delta ( u_i)
\ts \cdots \ts u_c\right)\,\,.$$
\noindent 
Applying these to the elements $u_i= f(\al_i, M_i)$, the assertion follows. 
\bigskip 

(\Sigmaprolong.2) {\bf Definition.}\quad 
Let $(K_1^m)$ and $(K_2^m)$ be finite sequences of objects indexed by 
$m\in \ZZ$.

If $(f(\al, M))\in \bop F(K^M_\al)$
is  of degree zero, $\sigma$-consistent and $\delta$-closed, 
so that the element $(f(M))\in \bop F(K^M)$ given by the 
above procedure is $\sigma$-consistent and 
$\delta$-closed,   
$(K^m; f(M))$ is called the $C$-diagram obtained by means of 
$\sigma$-consistent prolongation from $(f(\al, M)\,)$. 

If $(K_1, f_{K_1}(M))$ and $(K_2, f_{K_2}(M))$ are $C$-diagrams, 
and if one defines elements $f(\al, M)\in F(K_\al^M)$ by setting 
$f(\underline{1}, M)=f_{K_1}(M)$, $f(\underline{2}, M)=f_{K_2}(M)$
for constant functions $\al=\underline{1}$ and $\underline{2}$ and 
$f(\al, M)=0$ otherwise, they give a set of
$\sigma$-consistent and $\delta$-closed elements
(we then say that $(f(\al, M)\,)$ is {\it split}). 
The resulting $C$-diagram is called the {\it direct sum} of 
the $C$-diagrams $K_1$ and 
$K_2$. 
\bigskip

{\it Remark.}\quad (1) In all the  examples we encounter, 
the set of elements $(f(\al, M)\,)$ is either split, or 
the following weaker condition will be satisfied:
 $f(\al, M)=0$ unless $\al$ is non-decreasing, i.e., of the form
$\al=(1, \cdots, 1, 2, \cdots, 2)$. 

(2) It is clear that one can generalize the argument to the case $K$ is the direct sum of 
more than two objects.
\bigskip 

\sss{\Propintersectiontheta}\quad  With regard to the map $\theta$ in the previous subsection, 
we state a proposition which will be used in the rest of this section. 
Recall that we have $X^m=X^m_1\oplus X^m_2$, and for a function 
$\al: M\to \{1, 2\}$, we have the complexes 
$F(X^M_\al)=F(\al, M)$ and $F(X^M)= F(M)$. 

If $I$ is a finite set of integers, $\al: I\to \{1, 2\}$ is a function, and 
 $I_1, \cdots, I_r$ is a segmentation of $I$, then for elements 
 $u_i\in F(\al, I_i)$, $i=1, \cdots, r$, we have the condition 
of proper intersection for $\{u_i\}$, with respect to the sequence 
$m\mapsto X^m_{\al(m)}$.

Let $I$ be a finite set of integers, 
$T\subset S\subset \ctop{I}$, and 
$I_1, \cdots, I_r$ (resp. $J_1, \cdots, J_s$)
be the segmentation of $I$ by $S$ (resp. $T$). 
Thus there are integers $1=i_1<\cdots< i_s<i_{s+1}=r+1$
such that 
$$J_j=I_{i_j}\cup \cdots\cup I_{i_{j+1}-1}\,,
\qquad\mbox{for $j=1, \cdots, s$.}$$
Let $\al_i: I_i\to \{1, 2\}$, $i=1, \cdots, r$ be functions such that
if $ k=\term(I_i)\in T$ (resp. $\in S-T$) then 
$\al_i(k) = \al_{i+1}(k)$ (resp. $\al_i(k)\neq\al_{i+1}(k)$ )
(in other words, $\al_i$ are consistent on $T$, and inconsistent on $S-T$). 
Assume given elements $u_i\in F(\al_i, I_i)$ for $i=1, \cdots, r$. 

Before giving the statement, we give a few remarks. 
Note we have the maps
$$\theta: F(\al_{i_j}, I_{i_j})\ts\cdots\ts F(\al_{i_{j+1}-1}, I_{i_{j+1}-1})
\to F(J_j)$$
for $j=1, \cdots, s$. 

Assume that $k<\ell$ and $i_k, \cdots, i_\ell$ are consecutive integers so that 
$J_j=I_{i_j}$ for $j=k, \cdots, \ell$.
Then $\BI:=I_{i_k}\cup\cdots\cup I_{i_\ell}=
J_k\cup\cdots\cup J_\ell$, and the functions $\al_i$ on $I_i$ glue to define a 
function $\al$ on $\BI$.  
By what we said before, it then makes sense to ask whether the set 
$\{u_{i_k}, \cdots, u_{i_\ell}\}$ is properly intersecting. 

By the property (\propPIfive) of proper intersection in (\RmkDefqDG.2), we have:
\bigskip 

{\bf Proposition.}\quad{\it Under the above hypothesis, assume that 
the elements  $u_i\in F(\al_i, I_i)$ satisfy the following
condition: 

If $k<\ell$ and $i_k, \cdots, i_\ell$ are consecutive integers, then 
$\{u_{i_k}, \cdots, u_{i_\ell}\}$ is properly intersecting. 
\smallskip 

\noindent Then if we set $v_j=\theta(u_{i_j}\ts\cdots\ts u_{i_{j+1}-1})\in F(I_i)$ for 
$j=1, \cdots, s$, the set $\{v_1, \cdots, v_s\}$ is properly intersecting. }
\bigskip 

\sss{\AdditivityBGKL} {\it The additivity of $\BG(K, L)$.}\quad 
Assume given two sequences of objects $K_a=(K_a^{m} )$ for $a=1, 2$, 
and a set of elements  $(f(\al, M))\in \bop F(K^M_\al)$
 of first degree zero which is $\sigma$-consistent and $\delta$-closed,
 see (\Sigmaprolong).
Let $(K^m=K_1^m\oplus K_2^m; f(M))$ with $f(M)\in  
F(K^M)$ be the $C$-diagram
obtained by $\sigma$-consistent prolongation (\Sigmaprolong). 

 Let $(L; g(N))$ be another $C$-diagram. 
As a variant of the complex $H\dbullet(K, L)$ in (\GKL), consider the 
``double" complex (the $P$ in the notation indicates ``primary part")
$$\PH\dbullet=\PH\dbullet(K_1\bop K_2; L):=\bop [F(A)]_f\,,$$
the sum over $A=(\alpha, \BM; \BN)$, where $(\BM; \BN)$ is a double sequence
and $\al: M\to \{1, 2\}$ is a function, and 
$F(A)=F(\al, \BM; \BN):= F(K^{\BM}_\al; L^{\BN})$;
we denote by $[F(A)]_f$ a distinguished subcomplex.

The constraint for the distinguished subcomplex $[F(A)]_f$ is given by modifying 
that in (\GKL) as follows. 
Let $I=M\amalg N$, the segmentation $I_1, \cdots, I_c$, and $\BI$ be the same
as just before (\GKL.1). 
Take a function $\al: [-w, \term(M)]\to \{1, 2\}$ that takes
the same values as $\al: M\to \{1, 2\}$ on $M$.  
Let $\{J_j\}$ be a set of almost disjoint sub-intervals 
$\BI$ such that
for each $j$, one has either $J_j\subset [-w, \init(M)]$ or
$J_j\subset [\term(N), w]$, and set  
$$f(J_j)=\left\{
\begin{array}{cl}
f(\al, J_j) &\quad \mbox{if $J_j\subset [-w, \init (M)]$}\,,\\
g(J_j)&\quad \mbox{if $J_j\subset [\term(N), w]$}\,.
\end{array}
\right. 
$$
Then one has a condition of constraint
 $\cC=(\BI; [1, c] ; \{J_j\}; \{f(J_j)\})$, 
and the corresponding distinguished subcomplex 
$[F(A)]_\cC$.  
Consider all possible functions $\al$ extending $M$ and  
all possible sets of sub-intervals $\{J_j\}$, and take the 
intersection of the corresponding distinguished subcomplexes, 
and denote it by $[F(A)]_f$.
The constraint has been  so chosen that the
operations $\bpartial$, $\bvphi$, $\bsigma$, and $\Bf_K$, 
$\Bf_L$ below are defined on $[F(A)]_f$.

The bigrading $\PH^{a, b}$ is given by 
$a=\deg_1(u) +1=|u|+\gamma(M; N)+1$,  $b=\tau(u)$ as in (\GKL). 
The first differential is  $\bpartial+\bvphi$, where 
$\bvphi$ is the sum of 
$$\bvphi_k: [F(\al, \BM; \BN)]_f\to [F(\al, \BM-\{k\}; \BN)]_f$$
for $k\in \ctop{M}-(M'\cup\{\init(M)\})$
and 
$$\bvphi_k: [F(\al, \BM; \BN)]_f\to [F(\al, \BM; \BN-\{k\})]_f$$
for $k\in \ctop{N}-(N'\cup\{\term(N)\})$. 
The second differential is $d'=\bsigma+\Bf_K +\Bf_L$, 
where $\bsigma$ and $\Bf_L$ are the same as in  (\GKL), and 
$\Bf_K$ is the sum of 
$$[F(\al', \BM; \BN)]_f\ni u\mapsto \{|\tau (u)|+1\} f(\al, P)\ts u
\in [F(\al'\scirc \al, P\scirc \BM; \BN)]_f$$
for all pairs $(\al', P)$ with $\term(P)=\init(M)=\ell$ and $\al'(\ell)=\al(\ell)$, 
where  
$P\scirc \BM:= (P\cup M| \{\ell\}\cup M')$ as before 
and $\al'\scirc \al: P\cup M \to \{1, 2\}$ is the function that extends both 
$\al'$ and $\al$. 
From now we will write $F(A)$ for $[F(A)]_f$,
except when we need to exercise caution. 

For the ``double" complex $\PH\dbullet(K_1\oplus K_2, L)$, we have an analogue of the claim (\GKL.2), so that we can form the complex 
$$\PBG(K_1\bop K_2, L):= \Phi \PH\dbullet(K_1\bop K_2; L)\,.$$

We shall study the relationship to the complex $\BG(K, L):
=\Phi\PH\dbullet(K, L)$. 
Recall that an element $u$ of $\BG(K, L)$ is a $\sigma$-consistent collection 
$u=(u(M; N))$ with $u(M; N)\in F(K^M; L^N)$ for free 
double sequences $(M; N)$. 
As in (\Sigmaprolong), one has a quasi-isomorphism
$$
\theta: \moplus F(\al_1, K^{ M_1})\ts F(\al_2, K^{ M_2} )\ts\cdots\ts F(\al_r, K^{ M_r}; L^N)\to F(K^M; L^N)$$
where the sum is over segmentations $M_1, \cdots, M_r$ of $M$ and functions $\al_i: M_i\to \{1, 2\}$
taking distinct values at the overlaps  of the segmentation, 
$M_i\cap M_{i+1}$, $1\le i<r$. 
We call the summands  
$\bop F(\al, K^{M}; L^N)$  the {\it primary part}. 
\bigskip 

(\AdditivityBGKL.1){\bf Proposition.}\quad{\it
There is a quasi-isomorphism  of complexes 
$$\cP: \PBG(K_1\bop K_2, L)\to \BG(K, L)\,.$$
}\smallskip 

{\it Proof.}\quad 
To an element $u$ of $\PBG(K_1\bop  K_2, L)$, namely for a collection
$u=(u(\al, M; N))$, $u(\al, M; N)\in F(\al, K^M; L^N)$ for free 
double sequences $(M; N)$ and $\al: M\to \{1, 2\}$ satisfying the 
$\sigma$-consistency, associate an element $v=(v(M; N))$ in 
$\bop F(K^M; L^N)$ given by 
$$v(M; N)=\theta\left(
\sum f(\al_1, M_1)\ts \cdots \ts f(\al_c, M_c)\ts u(\al', M'; N)\right)\,$$
the sum over segmentations $(M_1, \cdots, M_c, M')$ of $M$ and functions 
$\al_i: M_i\to \{1, 2\}$, $\al': M'\to \{1, 2\}$ that take distinct values on the 
overlaps of the segmentation. 
By (\Propintersectiontheta), one has $v(M; N) \in [F(K^M; L^N)]_f$. 
Further $v$ is $\sigma$-consistent (the condition in (\ZFKL) is clearly satisfied), 
so it is an element of $\BG(K, L)$. Set $\cP(u)=v$. 

The assignment $u\mapsto \cP(u)$ is a map of complexes, since 
by the same computation 
as in the proof of (\Sigmaprolong.1), one has 
\begin{eqnarray*}
&
\delta(\theta\left(
\sum f(\al_1, M_1)\ts {\cdots} \ts f(\al_c, M_c)\ts u(\al', M'; N)
\right)\,) \\
=&\theta\left( 
\sum f(\al_1, M_1)\ts {\cdots} \ts f(\al_c, M_c)\ts \delta(u(\al', M'; N)\right)\,)\,, 
\end{eqnarray*}
in other words, $\delta(\cP(u))=\cP(\delta u)$.

We next show that $\cP$ is a quasi-isomorphism. 
On the ``double" complex $H\dbullet (K, L)$ consider a filtration by subcomplexes
indexed by pairs of integers $(m, n)$
$$F^{m, n}H\dbullet (K, L)=\moplus_{{\term(M)\le m},\, {\init(N)\ge n}}
F(K^\BM; L^\BN)\,.$$
One has an analogue of (\GKL.2) for $F^{m, n}H\dbullet (K, L)$ with the same proof, so one can apply the 
operation $\Phi$ to obtain a subcomplex 
$$F^{m, n}\BG(K, L):=
\Phi\left(F^{m, n}H\dbullet (K, L)\right) \subset \BG(K, L)\,,$$
giving a filtration of $\BG(K, L)$ by subcomplexes. 
Since the operation $\Phi$ is exact, 
a subquotient in the filtration is 
$$\Gr_F^{m, n}\BG(K, L)= \Phi\left(\Gr_F^{m, n}H\dbullet (K, L)\right)\,,$$
where
$$\Gr_F^{m, n}H\dbullet (K, L)
\moplus_{{\term(M)=m},\,{\init(N)= n}}
F(K^\BM; L^\BN)$$
is a ``double" complex with differentials $-\delta$ and $d'$, 
where $\delta=\partial+\sum \bvphi_k$, the sum over
 $k\in M-(M'\cup \{m\})$ and $k\in N-(N'\cup\{n\})$. 

On $\Gr_F^{m, n}H\dbullet (K, L)$ we consider the filtration 
analogous to the one in the proof of (\acyclicGI), 
indexed by $(a, b)$ with $a\le m$ and $b\ge n$, 
$$\cF(a, b)\Gr_F^{m, n}H\dbullet (K, L)
=\moplus_{{\init{M}\le a}, {\term{N}\ge b}} F(K^\BM; L^\BN)\,.$$
We show that the graded quotient 
$$\Gr_\cF^{a,b}\Gr_F^{m, n} H\dbullet (K, L)=
\moplus F(K^\BM; L^\BN)\,,$$
 the sum over $(\BM; \BN)$ satisfying $\init{M}= a\le \term{M}=m$, 
$\init{N}=n\le \term{N}=b$,
is acyclic unless $(a, b)=(m, n)$. 
Indeed consider the filtration $Fil_c$ on it given by
$Fil_c=\oplus F(K^\BM; L^\BN)$, the sum over $(\BM; \BN)$
satisfying $|M|+|N|\le c$. 
Then the graded quotient 
$\Gr_{Fil}^c \Gr_\cF^{a,b}\Gr_F^{m, n}$ is a sum of the total complexes
of the complexes of the form
$$0\to F(\cI|\emptyset)\mapr{\sigma}
\moplus_{|S|=1} F(\cI| S)
\mapr{\sigma} \moplus_{|S|=2} F(\cI | S)
\mapr{\sigma} \cdots \to  F(\cI|\ctop{\cI})\to 0\,\,$$
where $\cI=M\amalg N$ and $S$ varies over subsets of $\ctop
{\cI}$.  When $a<m$ or $b>n$, we have $|\cJ|\ge 3$, so the total complex is
acyclic by Lemma (\acyclicGI.1). 

Thus the natural surjection 
$$\Gr_F^{m, n}H\dbullet (K, L)\to \Gr_\cF^{m, n}
\Gr_F^{m, n}H\dbullet (K, L)=F(K^m; L^n)$$
is a quasi-isomorphism; consequently the induced surjection 
$$\Gr_F^{m, n}\BG(K, L)\to F(K^m; L^n)$$
is a quasi-isomorphism. 

We have a similarly defined filtration on the ``double"
complex $\PH\dbullet(K_1\oplus K_2, L)$:
$$F^{m, n}\PH\dbullet (K_1\oplus K_2,L)=\moplus_{{\term(M)\le m},{\init(N)\ge n}}
F(\al, K^\BM; L^\BN)\,, $$
and a proof parallel to the one above shows that there is a quasi-isomorphism
$$\Gr_F^{m, n}\PBG(K_1\oplus K_2, L)\to F(K_1^m; L^n)\oplus F(K_2^m; L^n). $$
The map 
$\cP: \PBG(K_1\oplus K_2, L)\to \BG(K, L)$ respects the filtrations $F^{m,n}$, 
so there is an induced map on $\Gr_F^{m, n}$, 
and we have a commutative square
$$\begin{array}{ccc}
\Gr_F^{m, n}\PBG(K_1\oplus K_2, L)    &\mapr{\cP}  & \Gr_F^{m, n}\BG(K, L) \\
\mapd{} &     &\mapdr{} \\
 F(K_1^m; L^n)\oplus F(K_2^m; L^n)   &\mapr{\theta}  &F(K^m; L^n)
 \end{array}
$$  
where the vertical maps are the above quasi-isomorphisms.
Since the lower horizontal map $\theta$ is a quasi-isomorphism, it follows that 
the upper horizontal arrow $\cP$ is a quasi-isomorphism. 
Hence the map  $\cP: \PBG(K_1\bop K_2, L)\to \BG(K, L)$
is a quasi-isomorphism. 
\bigskip 

Assume now given two sequences of objects $L_a=(L_a^{m} )$ for $a=1, 2$, 
and a set of elements  $(g(\al, N))\in \bop F(L^N_\al)$
of first degree zero which is $\sigma$-consistent and $\delta$-closed;
let  $(L^m=L_1^m\oplus L_2^m; g(N))$ with $g(N)\in  F(L^N)$ be the 
$C$-diagram given by $\sigma$-consistent prolongation. 
Let $K=(K^m; f(M))$ be another $C$-diagram. 
We define the ``double" complex 
$$\PH\dbullet =\PH\dbullet(K; L_1\bop L_2)$$
as the sum $\bop [F(A)]_f$; the sum is over $A=(\BM; \al, \BN)$, 
where $(\BM; \BN)$ is a double sequence,  $\al: N\to \{1, 2\}$ is a function, 
and $F(A)=F(K^{\BM}; L_\al^\BN)$;
by $[F(A)]_f$ we denote an appropriately defined distinguished subcomplex.
One makes it into a ``double" complex 
as one did with $\PH\dbullet(K_1\bop K_2; L)$, and 
then form the complex 
$\PBG(K; L_1\bop L_2):=\Phi \PH\dbullet(K; L_1\bop L_2)$. 
With this we have: 
\bigskip 

(\AdditivityBGKL.2){\bf Proposition.}\quad{\it
There is a quasi-isomorphism of complexes 
$$\cP: \PBG(K,  L_1\bop L_2)\to \BG(K, L)\,.$$
}\smallskip 

{\it Proof.}\quad Take an element  $u$ of $\PBG(K, L_1\bop L_2)$, 
namely  a $\sigma$-consistent collection of elements
$u=(u(M; \al, N))$, $u(\al, M; N)\in F( K^M;\al,  L^N)$ for free 
double sequences $(M; N)$ and $\al: N\to \{1, 2\}$, 
 and define  an element $v=(v(M; N))$ in 
$\bop F(K^M; L^N)$ given by 
$$v(M; N)= \theta\left(\sum(-1)^k u(M; \al', N')\ts g(\al_1, N_1)\ts\cdots \ts g(\al_c, N_c)\,\right)\,,
$$
$$\mbox{with}\quad k=(\deg_1 u(M; \al', N')\,)\cdot (\gamma(N_1)+\cdots +\gamma(N_c)\,)\,,$$
the sum over segmentations $(N', N_1, \cdots, N_c)$ of $N$ and functions 
$\al_i: N_i\to \{1, 2\}$, 
$\al': N'\to \{1, 2\}$ that take distinct values on the 
overlaps of the segmentation.
Here
$$\theta: \moplus F(K^M; \al', L^{N'})\ts F(\al_1, L^{N_1})\ts\cdots\ts 
F(\al_c, N_c)\to F(K^M; L^N)$$
is the quasi-isomorphism  as in (\Sigmaprolong).
Notice the sign assigned for each term.
We have $v(M; N)\in [F(K^M; L^N)]_f$ by (\Propintersectiontheta), 
thus $v\in H\dbullet(K; L)$;
also $v$ is $\sigma$-consistent, namely $v\in \BG(K, L)$.  
Set  $\cP(u)=v$. 
One verifies that $\cP$ is a map of complexes and is a quasi-isomorphism
as in the proof for the previous proposition.
\bigskip

We next specialize to the split case.
\bigskip 

(\AdditivityBGKL.3){\bf Proposition.}\quad{\it
Let  $(K_1, f_{K_1}(M))$ and $(K_2, f_{K_2}(M))$ be $C$-diagrams, 
and $(K; f(M))$ be the direct sum of them, see (\Sigmaprolong.2).
Then there is a canonical quasi-isomorphism
$$\cP: \BG(K_1, L)\oplus \BG(K_2, L)\to \BG(K, L)\,.$$
Similarly, one has a canonical quasi-isomorphism
$\BG(K, L_1)\oplus \BG(K, L_2)\to \BG(K, L)$. 
}
\smallskip 

{\it Proof.}\quad 
In the first case, one has an inclusion 
$H\dbullet(K_1, L)\oplus H\dbullet(K_2, L)\to \PH\dbullet(K_1\oplus K_2, L)$.
Indeed each summand $[F(K_1^{\BM}; L^\BN)]_f$ maps into 
$[F(K^{\BM}; L^\BN)]_f$ by (\Propintersectiontheta). 
Composing it with the map $\cP: \PH\dbullet(K_1\oplus K_2, L)\to H\dbullet(K, L)$ in (\AdditivityBGKL.1), we obtain a map 
$$H\dbullet(K_1, L)\oplus H\dbullet(K_2, L)\to H\dbullet(K, L)\,.$$
It induces a map $\BG(K_1, L)\oplus \BG(K_2, L)\to \BG(K, L)$. 
This is shown to be a quasi-isomorphism by an argument
parallel to that for 
(\AdditivityBGKL.1). 
The proof for the latter case is similar.
\bigskip

{\bf Remarks.}\quad 
$\bullet$ We refer to the maps $\cP$ in (\AdditivityBGKL.1), (\AdditivityBGKL.2)
and (\AdditivityBGKL.3) as the $\sigma$-consistent prolongation. 

$\bullet$ It is obvious how to generalize the above to the 
case where the object $K$ (or $L
$) is the direct sum of a finite number of 
 objects, 
$K_1\oplus\cdots \oplus K_r$.

$\bullet$ Combining (\AdditivityBGKL.1) and (\AdditivityBGKL.2), 
one can consider $\sigma$-consistent prolongation when 
$K=K_1\oplus K_2$ as in (\AdditivityBGKL.1), and 
$L=L_1\oplus L_2$ as in (\AdditivityBGKL.2). 
Then we should consider a $\sigma$-consistent set of elements 
$$u(\al, M; \be, N)\in F(K_\al^M; L_\be^N)$$
for free double sequences $(M; N)$ and functions $\al$ on $M$, 
$\be$ on $N$.   
\bigskip 

\sss{\AdditivityBG} {\it The additivity of $\BG(K_1, \cdots, K_n|\Sigma)$.}\quad 
Let $n\ge 2$, $I=[1, n]$, 
and $K_1, \cdots, K_n$ be a sequences of $C$-diagrams. 
For $\Sigma\subset (1,  n)$, we have the complex 
$\BG(I|\Sigma)=\BG(K_1, \cdots, K_n|\Sigma)$, as defined in (\GISigma).
For this, we give statements of 
additivity, generalizing (\AdditivityBGKL). 
Assume  $i$  is an integer 
with $1\le i\le n$, and we are given a direct sum decomposition 
$K_i= K'_i\oplus K''_i$. 
\bigskip 

{\bf Proposition.}\quad{\it
(1) For any $i$ with $1\le i\le n$, there is a map of complexes 
$$s_i: \BG(K_1, {\scriptstyle\cdots},  K'_i,  {\scriptstyle\cdots}, K_n|\Sigma)
\to \BG(K_1, {\scriptstyle\cdots},  K_i,  {\scriptstyle\cdots}, K_n|\Sigma)
\eqno{(\AdditivityBG.a)}
$$
(similarly $t_i: \BG(K_1, {\scriptstyle\cdots},  K''_i,  {\scriptstyle\cdots}, K_n|\Sigma)
\to \BG(K_1, {\scriptstyle\cdots},  K_i,  {\scriptstyle\cdots}, K_n|\Sigma)
$). 
The map $s_i$  is compatible with $\rho_m$ (for $m\in\Sigma$, allowing 
$m=i$), namely the square 
$$\begin{array} {ccc}
\BG(K_1, {\scriptstyle\cdots},  K'_i,  {\scriptstyle\cdots}, K_n|\Sigma)
&\mapr{s_i}
&\BG(K_1, {\scriptstyle\cdots},  K_i,  {\scriptstyle\cdots}, K_n|\Sigma)\\
\mapd{\rho_m} & &\mapdr{\rho_m} \\
\BG(K_1, {\scriptstyle\cdots},  K'_i,  {\scriptstyle\cdots}, K_n|\Sigma-\{m\})
&\mapr{s_i}
&\BG(K_1, {\scriptstyle\cdots},  K_i,  {\scriptstyle\cdots}, K_n|\Sigma-\{m\})
\end{array}
$$
commutes.
It is also compatible with $\Pi_k$ (for $k\not\in \Sigma$, allowing $k=i$), namely 
the following square commutes:
$$\begin{array}{ccc}
 \BG(K_1,\cdots, K_n|\Sigma)   &\mapr{s_i}  &\BG(K_1,  \cdots, K_n|\Sigma)  \\
\mapd{\Pi_k} &     &\mapdr{\Pi_k} \\
 \BG(K_1, {\scriptstyle\cdots},\widehat{K_k},{\scriptstyle\cdots},  K_n|\Sigma)  &\mapr{s_i}  
&\ph{\, .}\BG(K_1, {\scriptstyle\cdots},\widehat{K_k},{\scriptstyle\cdots},  K_n|\Sigma)\,.
 \end{array}
$$ 
Also, in the sense as in (\DefqDG),(3), 
$s_i$ and $s_j$ commute for $i\neq j$,
and $s_i$ is compatible with the constraint maps.

(2) In case $i\in \Sigma$, letting $\Sigma_1=\Sigma\cap (1, i)$ and 
$\Sigma_2=\Sigma\cap (i, n)$, one has a map of complexes
$$\pi_i(K'_i, K''_i): \BG(K_1, {\scriptstyle\cdots},  K'_i|\Sigma_1)
\ts \BG(K''_i,  {\scriptstyle\cdots}, K_n|\Sigma_2)
\to \BG(K_1, {\scriptstyle\cdots},  K_i,  {\scriptstyle\cdots}, K_n|\Sigma)
\,. \eqno{(\AdditivityBG.b)}$$
It is compatible with $\rho_m$ (for $m\in \Sigma$ with $m\neq i$), namely 
the following diagram commutes:
$$\begin{array}{ccc}
\BG(K_1, {\scriptstyle\cdots},  K'_i|\Sigma_1) \ts \BG(K''_i,  {\scriptstyle\cdots}, K_n|\Sigma_2)   &\mapr{\pi_i}  & \BG(K_1, {\scriptstyle\cdots},  K_i,  {\scriptstyle\cdots}, K_n|\Sigma) \\
\mapd{1\ts\rho_m} &     &\mapdr{\rho_m} \\
\BG(K_1, {\scriptstyle\cdots},  K'_i|\Sigma_1) \ts \BG(K''_i,  {\scriptstyle\cdots}, K_n|\Sigma_2-\{m\})   &\mapr{\pi_i}  & \BG(K_1, {\scriptstyle\cdots},  K_i,  {\scriptstyle\cdots}, K_n|\Sigma-\{m\}) 
 \end{array}
$$  
(the diagram is for $m>i$). 
It is compatible with $\Pi_k$ for $k\not\in \Sigma$. 
Also, $\pi_i$ and $\pi_j$ commute for $i\neq j$, 
and $s_i$ and $\pi_j$ commute (cf. (\DefqDG),(3)\,). 

(3) Assume $i\not\in \Sigma$. Then the map 
$$\theta_i(K'_i, K''_i): \BG(K_1, {\scriptstyle\cdots},  K'_i,  {\scriptstyle\cdots}, K_n|\Sigma_1) \bop
\BG(K_1, {\scriptstyle\cdots},  K''_i,  {\scriptstyle\cdots}, K_n|\Sigma_2)
\to 
\BG(K_1, {\scriptstyle\cdots},  K_i,  {\scriptstyle\cdots}, K_n|\Sigma)\,,$$
defined as the sum of the maps $s_i$ and $t_i$, is a quasi-isomorphism. 

(4) Assume $i\in\Sigma$. Then the map 
\begin{eqnarray*}
&\theta_i(K'_i, K''_i):& \BG(K_1, {\scriptstyle\cdots},  K'_i,  {\scriptstyle\cdots}, K_n|\Sigma) \bop
\BG(K_1, {\scriptstyle\cdots},  K''_i,  {\scriptstyle\cdots}, K_n|\Sigma)\\
&& \bop \BG(K_1, {\scriptstyle\cdots},  K'_i|\Sigma_1)\ts 
\BG(K''_i,  {\scriptstyle\cdots}, K_n|\Sigma_2)\\
&&\bop \BG(K_1, {\scriptstyle\cdots},  K''_i|\Sigma_1)\ts 
\BG(K'_i,  {\scriptstyle\cdots}, K_n|\Sigma_2) \\
&\mapr{}&   
\BG(K_1, {\scriptstyle\cdots},  K_i,  {\scriptstyle\cdots}, K_n|\Sigma)\,,
\end{eqnarray*}
defined as the sum of the maps $s_i$, $t_i$, $\pi_i(K'_i, K''_i)$, and 
$\pi_i(K''_i, K'_i)$, is a quasi-isomorphism. }
\smallskip 

{\it Proof.}\quad 
(1) First assume $\Sigma=\emptyset$ and $i=1$ or $n$. 
If $n=2$ the map $s_i$ was constructed in the previous subsection
(\AdditivityBGKL.3), and the case $n\ge 3$ is similar. 

Next assume $\Sigma=\emptyset$ and $1<i<n$. 
For simplicity, we consider the case $n=3$, $\Sigma=\emptyset$ and $i=2$;
thus we have $K_2=K_2'\oplus K''_2$. 
Recall that the ``triple" complex $H\tbullet({K_1}; {K_2}; {K_3})$
is the sum 
$$\moplus [F({K_1}^{\BM_1}; {K_2}^{\BM_2}; {K_3}^{\BM_3})\,]_f$$
over triple sequences $(\BM_1; \BM_2; \BM_3)$.
Here $[-]_f$ denotes the distinguished complex with respect to 
$\{f_{K_1}(P)\}$ and $\{f_{K_3}(Q)\}$, see (\GKL) and (\HIGI). 
Using the notation before (\GKL.a), it is generated by $u_1\ts\cdots\ts u_c$ with 
$u_i\in F(K_1, K_2, K_3; I_i)$ such that 
$$\{\{f_{K_1}(J_j)\}, u_1, \cdots, u_c, \{f_{K_3}(J_j)\}\}
\quad\mbox{is properly intersecting.}\quad
\eqno{(\AdditivityBG.c)}$$
The subcomplex $[F({K_1}^{\BM_1}; {K'_2}^{\BM_2}; {K_3}^{\BM_3})\,]_f$
is also given by the same condition, except then 
$u_i\in F(K_1, K'_2, K_3; I_i)$. 
For each $(\BM_1; \BM_2; \BM_3)$, one has the map  
obtained by composing $s_i$ for $i\in M_2$,
$$s: F({K_1}^{\BM_1}; {K'_2}^{\BM_2}; {K_3}^{\BM_3})
\to F({K_1}^{\BM_1}; {K_2}^{\BM_2}; {K_3}^{\BM_3})\,;$$
it induces a map between the distinguished subcomplexes 
$$s: [F({K_1}^{\BM_1}; {K'_2}^{\BM_2}; {K_3}^{\BM_3})]_f
\to [F({K_1}^{\BM_1}; {K_2}^{\BM_2}; {K_3}^{\BM_3})]_f\,.$$
because the condition $(\AdditivityBG.c)$ is respected by the map $s$ 
by (\propPItwo) in (\RmkDefqDG.2).  
Taking the sum over $(\BM_1; \BM_2; \BM_3)$, we obtain a degree
preserving map 
$$s: H\tbullet ({K_1}; {K'_2}; {K_3})
\to H\tbullet({K_1}; {K_2}; {K_3})\,.$$
This is a map of ``triple" complexes because 
(a) $s$ is compatible with $\partial$ since $s$ is a map of complexes, 
(b) $s$ is compatible with $\vphi$ and $\sigma$ since $s_i$ compatible with $\vphi$ and 
$\sigma$ by (\RmkDefqDG.2), 
and 
(c) $s$ is compatible with with $\Bf_{K_1}$ 
and $\Bf_{K_3}$, since $s_i
$ is compatible with tensor product.
Taking the operation $\Phi$ and then $\Tot$ as in (\HIGI), 
one obtains a map of complexes 
$$s: \BG({K_1}, {K'_2}, {K_3})
\to \BG({K_1}, {K_2}, {K_3})\,.$$
If $n$ is arbitrary and $\Sigma=\emptyset$, the argument is the same. 

When $\Sigma$ is non-empty and $i\not\in \Sigma$, 
we slightly extend the argument as follows. 
Assume, say, $n=4$, $\Sigma=\{3\}$ and $i=2$ (the general case is the same). 
Then one already has given the map
 $s: H\tbullet(K_1, K'_2, K_3)\to H\tbullet(K_1, K_2, K_3)$. 
We claim that there is a map 
$$s\ts 1: H\tbullet(K_1, K'_2, K_3)\tildets H\dbullet(K_3, K_4) \to 
H\tbullet(K_1, K_2, K_3)\tildets H\dbullet(K_3, K_4)\,.$$
This holds because the condition of constraint defining $H\dbullet\tildets H\dbullet$
in (\GISigma), for either Case (i) or (ii), 
 is respected by the map $s$ by (\RmkDefqDG.2). 
This induces the map as desired
$s\ts 1:  \BG(K_1, K'_2, K_3)\tildets \BG(K_3, K_4) \to 
\BG(K_1, K_2, K_3)\tildets \BG(K_3, K_4)$. 

When $\Sigma$ is non-empty and $i\in \Sigma$,  say $n=3$, $\Sigma=\{2\}$ and $i=2$, we argue as follows. 
We have the map of (\AdditivityBGKL.3),
$s: H\dbullet({K_1}; {{K'}_2})\to H\dbullet({K_1}; {{K}_2})$, and 
similarly 
$s: H\dbullet({K'_2}; {{K}_3})\to H\dbullet({K_2}; {{K}_3})$.
For simplicity we will drop the double or triple dots for $H$. 
They induce the map 
$$s\ts s: H({K_1}; {{K'}_2})\tildets H({K'_2}; {{K}_3})
\to H({K_1}; {{K}_2})\tildets H({K_2}; {{K}_3})
$$
since the condition defining $H\tildets H$
 is respected by the map $s\ts s$ by (\Propintersectiontheta).
It induces the map 
$$\BG({K_1}, {{K'}_2})\tildets \BG({K'_2}, {{K}_3})
\to \BG({K_1}, {{K}_2})\tildets \BG({K_2}, {{K}_3})\,.$$
The compatibility with $\rho_m$ and $\Pi_k$ is  obvious from the definitions. 

(2) Assume $n=3$, $\Sigma=\{2\}$ and $i=2$ (the general case being similar).
We have the maps
$s: H({K_1}; {{K'}_2})\to H({K_1}; {{K}_2})$ and 
$t: H({K''_2}; {{K}_3})\to H({K_2}; {{K}_3})$.
They induce the map 
$$s\ts t: H({K_1}; {{K'}_2})\tildets H({K''_2}; {{K}_3})
\to H({K_1}; {{K}_2})\tildets H({K_2}; {{K}_3})
$$
since the condition defining $H\tildets H$
 is respected by the map $s\ts s$ by (\Propintersectiontheta).
It induces the map 
$$\pi_2: \BG({K_1}, {{K'}_2})\tildets \BG({K''_2}, {{K}_3})
\to \BG({K_1}, {{K}_2})\tildets \BG({K_2}, {{K}_3})\,.$$

(3) If $\Sigma\neq (1, n)$, the assertion is obvious,
because then the  complex $\BG([1, n]|\Sigma)$ is acyclic.
If $\Sigma=(1, n)$, and $i=1$, say, then the assertion follows from the 
additivity of $\BG(K_1, K_2)$ in $K_1$, proven in (\AdditivityBGKL.3).

(4) If $\Sigma\neq (1, n)$, the assertion is obvious, because the  
complexes involved are all acyclic.  
If $\Sigma=(1, n)$ and $1<i<n$, then 
the assertion follows from the additivity of  $\BG(K_{i-1}, K_i)$
in $K_i$, and of  $\BG(K_{i}, K_{i+1})$ in $K_i$. 
\bigskip 

{\it Remark.}\quad Note in (3) there are no 
cross terms involved.
It thus differs from the additivity for 
the complex $F(I|S)$ as in (\RmkDefqDG.1). 
\bigskip

\sss{\AdditivityBF} {\it Additivity of $\BF(I)$.}\quad 
We keep the assumption of the previous subsection. 
Recall in (\BFIS) we have defined the complex $\BF(I)$
as well as the maps $\tau$ and $\vphi$.
\bigskip

{\bf Proposition.}\quad{\it
(1) For any $i$ with $1\le i\le n$, there is a map of complexes 
$$s_i: \BF(K_1, {\scriptstyle\cdots},  K'_i,  {\scriptstyle\cdots}, K_n)
\to \BF(K_1, {\scriptstyle\cdots},  K_i,  {\scriptstyle\cdots}, K_n)
\eqno{(\AdditivityBF.a)}
$$
(similarly $t_i: \BF(K_1, {\scriptstyle\cdots},  K''_i,  {\scriptstyle\cdots}, K_n)
\to \BF(K_1, {\scriptstyle\cdots},  K_i,  {\scriptstyle\cdots}, K_n)
$). 
The map $s_i$ is  compatible with $\vphi_k$ (for $1<k<n$), namely 
the following square commutes:
$$\begin{array}{ccc}
 \BF(K_1,{\scriptstyle\cdots}, K'_i,  {\scriptstyle\cdots} K_n)   &\mapr{s_i}  &\BF(K_1,{\scriptstyle\cdots}, K_i,  {\scriptstyle\cdots}, K_n)  \\
\mapd{\vphi_k} &     &\mapdr{\vphi_k} \\
 \BF(K_1, {\scriptstyle\cdots},\widehat{K_k},{\scriptstyle\cdots},  K_n)  &\mapr{s_i}  
&\ph{\, .}\BF(K_1, {\scriptstyle\cdots},\widehat{K_k},{\scriptstyle\cdots},  K_n)\,.
 \end{array}
$$ 
In the sense as in (\DefqDG.3), 
$s_i$ commutes with $\tau_j$ for $i\neq j$, 
$s_i$ and $s_j$ commute for $i\neq j$,
and $s_i$ is compatible with the constraint maps.

(2) For $i$ with $1<i<n$, there is a map of complexes 
$$\pi_i(K'_i, K''_i): \BF(K_1, {\scriptstyle\cdots},  K'_i)
\ts \BF(K''_i,  {\scriptstyle\cdots}, K_n)
\to \BF(K_1, {\scriptstyle\cdots},  K_i,  {\scriptstyle\cdots}, K_n)
\,. \eqno{(\AdditivityBF.b)}$$
It is compatible with $\vphi_k$ for $k\neq i$, namely the following diagram commutes:
$$\begin{array}{ccc}
\BF(K_1, {\scriptstyle\cdots},  K'_i) \ts \BF(K''_i,  {\scriptstyle\cdots}, K_n)   &\mapr{\pi_i}  & \BF(K_1, {\scriptstyle\cdots},  K_i,  {\scriptstyle\cdots}, K_n) \\
\mapd{\vphi_k\ts 1} &     &\mapdr{\vphi_k} \\
\BF(K_1, {\scriptstyle\cdots}, \widehat{K_k}, {\scriptstyle\cdots}, K'_i) \ts \BF(K''_i,  {\scriptstyle\cdots}, K_n)   &\mapr{\pi_i\ts 1}  & \BF(K_1, {\scriptstyle\cdots},  \widehat{K_k},  {\scriptstyle\cdots}, K_n) 
 \end{array}
$$  
(the diagram is for $k<i$). 
Also, in the sense of (\DefqDG.3), 
$\pi_i$ is compatible with $\tau_j$, 
$\pi_i$ and $\pi_j$ commute for $i\neq j$, 
and $s_i$ and $\pi_j$ commute.

(3) The map 
\begin{eqnarray*}
&\theta_i(K'_i, K''_i):& \BF(K_1, {\scriptstyle\cdots},  K'_i,  {\scriptstyle\cdots}, K_n) \bop
\BF(K_1, {\scriptstyle\cdots},  K''_i,  {\scriptstyle\cdots}, K_n)\\
&& \bop \BF(K_1, {\scriptstyle\cdots},  K'_i)\ts 
\BF(K''_i,  {\scriptstyle\cdots}, K_n)\\
&&\bop \BF(K_1, {\scriptstyle\cdots},  K''_i)\ts 
\BF(K'_i,  {\scriptstyle\cdots}, K_n) \\
&\mapr{}&   
\BF(K_1, {\scriptstyle\cdots},  K_i,  {\scriptstyle\cdots}, K_n)\,,
\end{eqnarray*}
defined as the sum of $s_i$, $t_i$, $\pi_i(K'_i, K''_i)$, and 
$\pi_i(K''_i, K'_i)$, is a quasi-isomorphism.}
\smallskip 

{\it Proof.}\quad 
(1) Recall $\BF(I)=\bop\BG(I|\Sigma)$ by definition. 
We take the sum of the maps $s_i$ of (\AdditivityBG.a) to define  $s_i$  of (\AdditivityBF.a). 
Since the former is compatible with $\rho_m$, it follows that the the latter is a map of 
complexes. The compatibility of the former with $\Pi_k$ implies the compatibility of the 
latter with $\vphi_k$. 
The other compatibilities are obvious from the definitions
(recall that $\tau_j$ is the sum of the identity maps 
 on some summands $\BG(I|\Sigma)$).

(2) Take the sum of the maps $\pi_i$ of (\AdditivityBG.b) to define  $\pi_i$  of (\AdditivityBF.b).

(3) Consider the filtration $Fil_c$ on the complexes defined as in (\BFIS). 
Taking $\Gr_{Fil}^c$ of the map $\theta_i$, one obtains the sum, over $\Sigma$ with $|\Sigma|=c$, 
of the maps $\theta_i$ of (\AdditivityBG), (3) if $i\not\in\Sigma$, and 
$\theta_i$ of (\AdditivityBG), (4) if $i\in\Sigma$. Since those $\theta_i$ are quasi-isomorphisms,
the assertion follows. 
\bigskip

The groupoid of $C$-diagrams defined in (\DefCdzero) is a symmetric monoidal 
category with respect to the direct sum  of (\Sigmaprolong.2). 
The complexes $\BF(K_1, \cdots, K_n)$ and the maps $\tau$, $\vphi$
satisfy conditions (\DefqDG), (3) and (4) by (\AdditivityBF)
and (\IdentityProp), respectively. 
We have thus proven the following theorem. 
(Statement (2) is in (\GKL.4), and statement (3) easily follows from 
(\Compos).)
\bigskip

\sss{\ThmCdqDG} {\bf Theorem.}\quad{\it
(1) Let $\cC$ be a quasi DG category possessing  the additional structure 
(iv), (v) of (\DefqDG). 
Then the symmetric groupoid $(\cC^\Delta)_0$ of $C$-diagrams in $\cC$, 
 the complexes $\BF(K_1, \cdots, K_n)$ and the maps $\tau$, $\vphi$ as in 
\S 2 satisfy the conditions of a quasi DG category. 
(This quasi DG category will be denoted by $\cC^\Delta$.)

(2) If $X$, $Y$ are objects in $\cC$ and $m, n\in \ZZ$, then 
we have $\BF(X[m], Y[n])=F(X, Y)[n-m]$. 

(3) In the homotopy category $Ho(\cC^\Delta)$, the composition
of morphisms
between the objects $X[m]$, $Y[n]$, $Z[\ell]$ (where 
$X, Y, Z$ are objects of $\cC$ and $m, n, \ell \in \ZZ$), 
$$H^0\BF(X[m], Y[n])\ts H^0\BF(Y[n], Z[\ell])\to H^0\BF(X[m], Z[\ell])$$
is identified via the isomorphisms (2) with 
the map 
$$H^{n-m}F(X, Y)\ts H^{\ell-n}F(Y, Z)\to H^{\ell -m}F(X, Z)$$
induced by $\tau$ and $\vphi$. (See Remark to (\HowqDG) for the latter map). 
}\bigskip 

The proof of the next theorem will take the rest of this section.
\bigskip 

\sss{\ThmCdqDGtriang} {\bf Theorem.}\quad{\it
Under the same assumption and notation, the homotopy category $Ho(\cC^\Delta)$
of the quasi DG category $\cC^\Delta$ has the structure of a triangulated 
category.
}\bigskip

\sss{\Shiftfunctor}
{\it Shifting functor.}\quad 
For an increasing sequence $M=(m_1, \cdots, m_\mu)$, 
let $M[1]=(m_1+1, \cdots, m_\mu+1)$. 
For a sequence $\BM=(M|M')$ as in \S \secfcn, 
let $$\BM[1]=(M[1]|M'[1])\,\,.$$
Likewise for a double sequence $(\BM; \BN)$, another double sequence
$(\BM[1]; \BN[1])$ is defined.

Let $K=(K^m; f(M))$ be a $C$-diagram.
Define another $C$-diagram 
$$K[1]=(K[1]^m; (f[1])(M)\,)$$
 by 
$(K[1])^m= K^{m+1}$ and $(f[1])(M)=f(M[1])\in F(M[1])$
(namely $(f[1])(m_1, \cdots, m_\mu)=f(m_1+1, \cdots, m_\mu+1)$\,). 

Let $K$ and $L$ be $C$-diagrams. Recall 
$$H\dbullet(K; L)=\moplus_{(\BM; \BN)} F(K^{\BM}; L^{\BN})=\bop F(\BM; \BN)\,\,.$$
So an element $u$ has $(\BM; \BN)$-component $u(\BM; \BN)\in F(\BM; \BN)$. 
Similarly 
$$H\dbullet(K[1], L[1])=\moplus_{(\BM; \BN)} F(\BM[1]; \BN[1])\,\,.$$ 
Define a map 
$$\shift: H\dbullet(K; L)\to H\dbullet(K[1]; L[1])$$
by sending $u=(u(\BM; \BN)\, )$ to $-u[1]$, where
$(u[1])(\BM; \BN)=u(\BM[1]; \BN[1])\in F(\BM[1], \BN[1])$. 
This map preserves $|u|$, $\gamma (u)$ and $\tau(u)$, and 
 compatible with the maps $\delta$ and 
$d'$. Thus it is an isomorphism of ``double" complexes.
Taking $\Phi$ we get an isomorphism of complexes
$\shift : \BG(K; L)\to \BG(K[1], L[1])$. 
 
Similarly for $n$ $C$-diagrams $K_1, \cdots, K_n$, one has an isomorphism
of ``triple" complexes
$$\shift: H\tbullet (K_1;\cdots ; K_n) \to H\tbullet(K_1[1]; \cdots ; K_n[1])
$$
by sending $u$ to $\shift (u)=(-1)^{n+1} u[1]$  where  
$(u[1])(\BM_1; \cdots; \BM_n)=u(\BM_1[1]; \cdots; \BM_n[1])$. 
Taking $\Phi$ and then $\Tot$, we get an isomorphism of complexes 
$$\shift: \BG(K_1, \cdots, K_n)\to \BG(K_1[1], \cdots, K_n[1])\,.$$
For $\Sigma$ a subset of $(1, n)$ with cardinality $r-1$, let $I_1, \cdots, 
I_r$ be the corresponding segmentation of $[1, n]$, and we let  
$$\shift: \BG(K_1, \cdots, K_n|\Sigma)\to \BG(K_1[1], \cdots, K_n[1]|\Sigma)\eqno{(\Shiftfunctor.a)}$$
be given by 
$$u=u_1\ts\cdots\ts u_r\mapsto (-1)^{n+1} u_1[1]\ts\cdots u_r[1]
=\shift(u_1)\ts \cdots\ts \shift(u_r)\,.$$
\bigskip 

(\Shiftfunctor.1) {\bf Proposition.}\quad {\it The map $\shift$ commutes with $\rho_m$ and $\Pi_k$, and it is compatible with additivity of (\AdditivityBG).
}
\smallskip 

{\it Proof.}\quad  The compatibility with $\rho_m$ is obvious by the definitions. 
The compatibility with the additivity, which means the 
commutativity with the maps $s_i$ and $\pi_i$ in (\AdditivityBG), 
is also obvious. 

The compatibility with $\Pi_k$ means the identity 
$\Pi_k\scirc \shift =\shift\scirc \Pi_k$. 
To show this, we reduce to the case $\Sigma$ is empty, 
where the verification is immediate recalling by definition
\begin{eqnarray*}
&&(\Pi_k u) (\BM_1; \cdots, \BM_{k-1}; \BM_{k+1};
\cdots; \BM_n)\\
&=&\sum_j (-1)^j \vphi_{K_k^j}\left(
u(\BM_1; \cdots, \BM_{k-1}; (\{j\}|\emptyset);\BM_{k+1}; \cdots; \BM_n)
\right)\,.
\end{eqnarray*}
\bigskip

Taking the sum of $\shift$ of (\Shiftfunctor.a) for $\Sigma$ containing $S$, we obtain 
the maps as in the following proposition. 
\bigskip 

(\Shiftfunctor.2) 
{\bf Proposition.}\quad {\it 
There is an isomorphism of complexes 
$$\shift: \BF(K_1, \cdots, K_n|S)\to \BF(K_1[1], \cdots, K_n[1]|S)
\eqno{(\Shiftfunctor.b)}$$
that is compatible with the maps $\sigma_{S\, S'}$ and $\vphi_K$. 
It is also compatible with  additivity of (\AdditivityBF). 
}\bigskip 

In particular, the isomorphism $\shift: \BF(K_1, \cdots, K_n)\to \BF(K_1[1], \cdots, K_n[1])$ is compatible with $\tau_k$ and $\vphi_\ell$. One thus has
an isomorphism of abelian groups 
$$\shift: H^0\BF(K, L)\to H^0\BF(K[1], L[1])\eqno{(\Shiftfunctor.c)}$$
that is compatible with composition 
$$H^0\BF(K_1, K_2)\ts H^0\BF(K_2, K_3)\to H^0\BF(K_1, K_3)$$
as defined in (\HowqDG) using $\tau$ and $\vphi$. 
We thus have a self-equivalence $\shift$ 
of the additive category $Ho(\cC^\Delta)$ given on 
objects by $K\mapsto K[1]$ and on morphisms by (\Shiftfunctor.c).
If $u: K\to L$ is a morphism represented by $\underline{u}\in 
Z^0\BF(K, L)$, 
then $\shift (u):K[1]\to L[1]$ is the morphism represented by $\shift(\underline{u})
=-\underline{u}[1]\in  Z^0\BF(K[1], L[1])$. 

In the sequel of this paper, we shall write $u[1]$ for $\shift(u)$, 
bearing in mind that it is indeed represented by $-\underline{u}[1]$. 
\bigskip

\sss{\Conemorphism} 
{\it The cone of a morphism.}\quad 
Let $u: (K^m; f(M)\,)\to (L^m; g(M)\,)$
be a morphism. Take a representative 
$\underline{u}=(\underline{u}(M; N)\,)$ in $Z^0\BF(K, L)$. 
We will define the  {\it cone } of $u$ as the 
$C$-diagram $C_u=(C_u^m; h(M)\,)$, with 
$$C_u^m=K^{m+1}\bop L^m\,\,, $$
and with elements $h(M)\in F(C_u^M)^{-\gamma(M)}$ to be given
below. 

According to (\Sigmaprolong), we shall specify primary elements 
$$h(\al, M)\in F(\al, C_u^M)$$
and take its $\sigma$-consistent prolongation.  

For $\al=\underline{1}$ the corresponding group is 
$F(\underline{1}, C_u^M)=F(K[1]^M)$. 
One has an identification 
$F(K^{M[1]})=F(K[1]^M)$, where the former is a summand of $\bop F(K^M)$, 
and the latter is a summand of $\bop F(\al, C_u^M)$;
the identification preserves the number $\gamma$ in respective groups, namely 
$\gamma(M[1])=\gamma(M)$, so it also preserves $\deg_1(u)=|u|+\gamma(M)$. 
Consider the element 
$f(M[1])\in F(K^{M[1]})$.
According to the identification we view it as an element of $F(K[1]^M)$, where  $\deg_1 f(M[1])=0$
still holds.   We take $h(\underline{1}, M)$ to be 
this element: 
$$h(\underline{1}, M)=f(M[1])\in F(K[1]^M)\,.$$
For $\al=\underline{2}$, let 
$$h(\underline{2}, M)=g(M)\in F(L^M)\,.$$
Next let $M', M''$ be a segmentation 
of $M$ (in the sense of (\Finiteorderedset)\,)
with $M', M''$ both non-empty.
By $\al=(\underline{1}^{M'}, \underline{2}^{M''})$ we denote the function 
on $M$ such that $\al(i)=1$ on $M'$ and $\al(i)=2$ on $M''$; it 
corresponds to the group 
$F(\al, C_u^M)=F(K[1]^{M'}, L^{M''})$. 
We have an identification 
$F(K^{M'[1]}; L^{M''})= F(K[1]^{M'}, L^{M''})$
where the former is a summand of $H\dbullet (K, L)$, and 
the latter is a summand of $\oplus F(\al, C_u^M)$. 
It preserves the number $\gamma$ attached to the respective summands, namely
$\gamma(M'[1]; M'')=\gamma (M)$. 
We take the element $\underline{u}(M'[1]; M'')\in 
F(K^{M'[1]}; L^{M''})$,  via the above identification 
view it as an element of the latter group, and take
$h(\al, M)$ to be this element. Thus
$$h(\al, M)= \underline{u}(M'[1]; M'')\in 
F(K^{M'[1]}; L^{M''})= F(K[1]^{M'}, L^{M''})\,\,(=F(\al, C_u^M)\,)\,.$$
Note that $\deg_1 h(\al, M)=0$. 
For functions $\al$ other than the above, we take $h(\al, M)=0$. 
To be concise,  following matrix expression is useful:
$$\renewcommand{\arraystretch}{1.3}
h(M)=
\begin{array}{cc}
    &K[1] \phantom{aaaaaaaa}L \\
\begin{array}{c}
K[1]\!\!\!\!\\
 L \!\!\!\!
\end{array}
&\left[
\begin{array}{cc}
f(M[1])&0  \\
\underline{u}(M'[1]; M'') & g(M)
\end{array}
\right]
\end{array} 
\eqno{(\Conemorphism.a)}
$$
where $M', M''$ are obtained from $M$ by partition. 

Since $f(M)$, $\underline{u}(M; N)$, and $g(M)$ are $\sigma$-consistent and 
$\delta$-closed,  the set of elements $(h(M))\in \bop F(C_u^M)$ obtained by 
means of $\sigma$-consistent prolongation is also $\sigma$-consistent and 
$\delta$-closed, defining the $C$-digram $C_u$. 
\bigskip 

There are canonical morphisms $\al(u): L\to C_u$ and $\beta(u): C_u\to
K[1]$ as we shall define. 

We begin with a complement to (\AdditivityBGKL.2). 
With notation there, assume that $g(\al, N)\in F(K^M; \al, L^N)$ 
satisfies the following condition: 
$g(\al, N)=0$ unless $\al$ is non-increasing. 
Then there is a canonical map of complexes
$$\BG(K, L_2)\to \BG(K, L)\,.\eqno{(\Conemorphism.b)}$$
Indeed one has an inclusion 
$H\dbullet(K, L_2)\to \PH\dbullet(K, L_1\oplus L_2)$, and it is a map 
of ``double" complexes under the assumption (see proof of (\AdditivityBGKL.3)
for an analogous argument). It induces a map of complexes 
$\BG(K, L_2)\to \PBG(K, L_1\oplus L_2)$, and composing with the map 
$\cP$ of (\AdditivityBGKL.2), we obtain the map (\Conemorphism.b). 

 The above assumption is satisfied for $C_u=K[1]\oplus L$, thus 
we have a map of complexes 
$\BF(L, L)\to \BF(L, C_u)$. 
Let $\underline{\al}=\underline{\al}(u)\in Z^0\BF(L, C_u)$ be 
the image of 
$(\tilde{g}(M; N))\in Z^0\BF(L, L)$.  
Concretely it is obtained by $\sigma$-consistent prolongation from 
the set of primary elements
$\underline{\al}(M; \gamma, N)\in F(L^M; \gamma, C_u^N)$ given by 
$\underline{\al}(M; \underline{2}, N)=\ti{g}(M; N)$ and 
$\underline{\al}(M; \gamma, N)=0$ for $\gamma\neq \underline{2}$; in matrix 
form, 
$$\underline{\al}(M; N) =
\begin{array}{cc}
   &\!\!\!\!L \\
\begin{array}{c}
K[1]\!\!\!\!\\
L\!\!\!\!\!\!\!\!
\end{array}
&\left[
\begin{array}{c}
 0\\
\ti{g}(M; N)
\end{array}
\right].
\end{array}
$$
Explicitly, 
$$\underline{\al}(M; N)=
\theta\left(\sum \ti{g}(M; N')\ts h(\gamma_1, N_1)\ts\cdots\ts h(\gamma_c, N_c)
\right),
\eqno{(\Conemorphism.c)}$$
the sum over segmentations $(N', N_1, \cdots, N_c)$ of $N$ and 
functions $(\gamma_1, \cdots, \gamma_c)$ taking distinct values at
the overlaps of the segmentation and $\gamma_1(\init(N_1))=1$. 
Let $\al(u): L\to C_u$ be the morphism represented by 
$\underline{\al}$. 

Similarly we have a map of complexes 
$$\BF(K[1], K[1])\to \BF(C_u, K[1])\,.$$
We have $(f[1])\,\tilde{}\in Z^0\BF(K[1], K[1])$, the element 
representing the identity of $K[1]$, and we take 
$\underline{\be}\in Z^0\BF(C_u, K[1])$ to be the image 
by the above map. 
In other words, it is obtained by $\sigma$-consistent 
prolongation from the set of primary elements
$$\renewcommand{\arraystretch}{1.3}
\underline{\be}(M; N)=\begin{array}{cc}
   &\phantom{aaa}K[1]\phantom{aaaaaa} L \\
 K[1] &\left[\,\,
(f[1])\,\tilde{}\,(M; N) \phantom{aa} 0
\right]. 
\end{array}$$
The explicit form is 
$$\underline{\be}(M; N)=\theta\left(
\sum h(\gamma_1, M_1)\ts\cdots\ts h(\gamma_c, M_c)
\ts (f[1])\,\tilde{}\,(M'; N)\,\right)\eqno{(\Conemorphism.d)}$$
the sum over segmentations $(M_1, \cdots, M_c, M')$ of $M$ and 
functions $(\gamma_1, \cdots, \gamma_c)$ taking distinct values at the 
overlaps of the segmentation and $\gamma_c(\term(M_c))=2$. 
We let $\be(u): C_u\to K[1]$ be the morphism represented by 
$\underline{\be}$. 

One verifies that the composition of $\al$ and $\be$ is zero. 
The argument follows the line of the proof for (\IdentityProp). 
Define the element    
$$W=(W(M; N; R)\,)\in \moplus F(L^M;C_u^N; K[1]^R)\subset H^{0,1,0}(L, C_u, K[1], \emptyset)
$$
by 
$$W(M; N; R)=\theta\left(
\sum\underline{\al}(M; N')\ts h(\al_1, N_1)\ts\cdots\ts h(\al_r, N_r)\ts 
\underline{\be}(N''; R)\right)\,.$$
The sum is over segmentations $(N', N_1,\cdots, N_r, N'')$ of $N$ and
functions $(\al', \al_1, \cdots, \al_r, \al'')$ on respective subsets 
which take distinct values at the overlaps of $N_1, \cdots, N_r$ and 
satisfy $\al_1(\init(N_1))=1$ and $\al_r(\term(N_r))=2$. 
We then have 
$\tilde{\bsigma}(W)=\rho(\underline{\al}\ts \underline{\be})$. 
If $N$ consists of one element $\ell$, then 
$$W(M; N; R)=\theta (\underline{\al}(M; N)\ts \underline{\be}(N; R)\,)$$
where $\theta$ is the map $\pi_\ell$. Thus 
we have $\Pi_{C_u}(W)(M; R)=0$ since the composition of $\pi_\ell$ with 
$\vphi_\ell$ is zero by (\DefqDG), (3).
This shows the assertion.

One thus has a triangle
$$K\overset{u}\to L\to C_u \overset{[1]}\to \,\,.$$
Such a triangle is called a standard distinguished triangle. 
We declare the {\it distinguished triangles} 
to be  the ones isomorphic to the standard ones.

For the axioms of a triangulated category, 
see for example \RefVe\, \RefWe\, or \RefKaSh, \S1. 
The verification of them is parallel to the case of the homotopy category of 
complexes; see e.g.  \RefKaSh, \S 1.4 for a detailed exposition.
The arguments are similar to the DG case, done in \RefHaonetwo \, in detail. 
The above definitions of  $h(M)$, $\alpha$, $\beta$ are motivated by the DG case. 
We will show two of the axioms in the following propositions. 
The other axioms are easily verified as in [KS], \S 1.4. 
\bigskip

\sss{\Axiomproof} \Prop{\it There exists an isomorphism 
$\phi: K[1]\to C_{\alpha (u)}$ such that the following diagram commutes:
$$\CD
 L @>\alpha (u) >> C_u @>{\beta (u)} >>
K[1] @>-u[1] >> L[1] \\
@V id VV  \ph{id} @V id VV \ph{\phi}@V \phi  VV
\ph{id} @V id VV
\\
L@>{\alpha (u)}>> C_u @>{\alpha (\alpha (u))}>>
C_{\alpha (u)} @>\beta(\alpha(u))>> \ph{\,\,.} L[1]\,\,.
\endCD
$$
}\bigskip 

{\it Proof.}\quad  To avoid conflict of notation, write 
$D=C_{\al(u)}$. 
Note $D^m=L[1]^m\oplus K[1]^m\oplus L^m$ as a sequence of objects. 
The structural elements $k(M)\in F(D^M)$
are given as the $\sigma$-consistent prolongation (with respect to 
the decomposition $L[1]\oplus C_u$)
of the set of primary elements (see (\Conemorphism.a) )
$$\renewcommand{\arraystretch}{1.3}
\begin{array}{cc}
    &L[1] \phantom{aaaaaaaa}C_u \\
\begin{array}{c}
L[1]\!\!\!\!\\
 C_u \!\!\!\!
\end{array}
&\left[
\begin{array}{cc}
g(M[1])&0  \\
\underline{\al(u)}(M; N) & h(M)
\end{array}
\right]
\end{array} 
$$
which is also given as the $\sigma$-consistent prolongation, with respect to 
the three summands $L[1]\oplus K[1]\oplus L$,
from the set of primary elements $k(\al, M)\in F(\al, D^M)$ for 
functions $\al: M\to \{1, 2, 3\}$, 
given as follows.
\smallskip 

(a) Type (1).\quad 
If $\al=\underline{1}$, then 
$k(\underline{1}, M)=g(M[1])\in F(L[1]^M)$.

(b) Type (13).\quad 
If $\al=({1}^{M'}, {3}^{M''})$ for a partition 
$M', M''$ of $M$ (both $M'$, $M''$ are non-empty)
such that $\term(M')+1=\init(M'')$, then
$$k(({1}^{M'}, {2}^{M''}), M)=g(M'[1]_\vartriangle M'')
\in F(L[1]^{M'}, L[1]^{M''})\,.$$
Here we let  $_\vartriangle$ mean, if $\term(M')+1=\init(M'')=\ell$, taking the diagonal extension at $L^{\ell}$ and 
putting the sign $(-1)^\ell$.
For convenience we set $g(M'[1]_\vartriangle M'')=0$ if
$\term(M')+1\neq \init(M'')$. 

(c) Type (2).\quad 
If  $\al=\underline{2}$, then 
$k(\underline{1}, M)=f(M[1])\in F(K[1]^M)$.

(d) Type (23).\quad If $\al=({2}^{M'}, {3}^{M''})$ for a partition 
$M', M''$ of $M$ (both $M'$ and $M''$ are non-empty)
 then 
$$k(({2}^{M'}, {3}^{M''}), M)=\underline{u}(M'[1];
 M'')\in F(K[1]^{M'}, L^{M''})\,.$$
 
(e) Type (3).\quad 
If $\al=\underline{3}$, then 
$k(\underline{1}, M)=g(M)\in F(L^M)$.

\noindent Take $k(\al, M)=0$ in all other cases.
One may express this as 
$$\renewcommand{\arraystretch}{1.3}
k(M)=
\begin{array}{cc}
    &\phantom{a}L[1]\phantom{aaaaaaaaa} K[1] \phantom{aaaaaaa}L \\
\begin{array}{c}
L[1]\!\!\!\!\\
K[1]\!\!\!\!\\
 L \!\!\!\!
\end{array}
&\left[
\begin{array}{ccc}
g(M[1]) &0 &0 \\
0& f(M[1])&0  \\
g(M'[1]_\vartriangle M'') &\underline{u}(M'[1]; M'') & g(M)
\end{array}
\right]
\end{array}
$$
These primary elements are $\sigma$-consistent and $\delta$-closed, and one 
takes the $\sigma$-consistent prolongation to obtain $k(M)\in F(D^M)$

Define morphisms 
$\psi: D=C_{\al(u)}\to K[1]$ and $\phi: K[1]\to D$  as follows. 

We shall give a set of elements  representing $\psi$,  
$(\underline{\psi}(M; N))\in Z^0\BF(D, K[1])$, 
according to (\AdditivityBGKL.1), as the 
$\sigma$-consistent prolongation of the set of elements 
$\underline{\psi}(\al, M; N)\in F(\al, D^M; K[1]^N)$ for 
functions $\al: M\to \{1, 2, 3\}$ given as follows. 
If $\al=\underline{2}$, and $\term(M)=\init(N)=\ell$, then 
$$\underline{\psi}(\underline{2}, M; N)=(f[1])(M_\vartriangle N)\in F(K[1]^M; K[1]^N)$$
where $_\vartriangle$ means taking diagonal extension 
at $K^{\ell+1}$ and putting 
$(-1)^\ell$\, ({\it not\/} $(-1)^{\ell+1}$). 
For other functions $\al$, we take $\underline{\psi}(\al, M; N)=0$. 
In matrix form, one may express it as 
$$\renewcommand{\arraystretch}{1.3}
\underline{\psi}(M; N)= \begin{array}{cc}
   &L[1] \phantom{aaa}K[1]\phantom{aaaaaa} L \\
 K[1] &\left[\,\, \,0\,\,\,
 \ph{\,\,\,}f[1](M_\vartriangle N)\ \phantom{aaa} 0
\right]  \\
 & 
\end{array}  
\,\, \in F(D^M;  K[1]^N)\,\,.$$
It is obvious that these form a  $\sigma$-consistent and 
$\delta$-closed set of elements. By means of (\AdditivityBGKL.1) we obtain 
$(\underline{\psi}(M; N)\,)\in \BF(D, K[1])$ which is also $\delta$-closed. 

We next give elements $\underline{\phi}(M; N)\in F(K[1]^M; D^N)$
representing $\phi$ as the $\sigma$-consistent prolongation of the 
set of elements $\underline{\phi}(M;\al,  N)\in F(K[1]^M; \al, D^N)$
for functions $\al:N\to \{1, 2, 3\}$, given as follows. 
\smallskip 

(a) Type (1).\quad When $\al=\underline{1}$, we take the element 
$\underline{u}(M[1]; N[1]) \in F(K^{M[1]}; L^{N[1]})$, and via 
the identification $F(K^{M[1]}; L^{N[1]})=F(K[1]^{M}; L[1]^{N})$
view it as an element of the latter group. 

We take $\underline{\phi}(M; \underline{1}, N)$ to be 
$\underline{u}(M[1]; N[1])$. 
\smallskip 
 
(b) Type (13).\quad  When $\al=(1^{N'}, 3^{N''})$ for a partition $N=N'\amalg N''$ such that $N'$, $N''$ both non-empty, 
 and $\term{N'}+1=\init{N''}=\ell$, 
we take the element $\underline{u}(M[1]; N'[1]\cup N'')\in F(K^{M[1]};
L^{N'[1]\cup N''})$ (here $N'[1]\cup N''$ is the set theoretic union, i.e.,  $\ell$ is counted 
once), take its image by the diagonal extension  at $L^\ell$
$$\diag:F(K^{M[1]}; L^{N'[1]\cup N''} )\to F(K^{M[1]}; L^{N'[1]}, L^{N''})\,,$$
put the sign $(-1)^\ell$, and via the identification 
$$F(K^{M[1]}; L^{N'[1]}, L^{N''}) =F(K[1]^{M}; L[1]^{N'}, L^{N''})\,\,(=
F(K[1]^M; \al, D^N)\,)$$
view it as an element of the latter group and 
denote it by $\underline{u}(M[1]; N'[1]_\vartriangle N'')$. 
Since $\underline{u}(M[1]; N'[1]$
\newline $\cup N'')$ is of first degree 0, i.e., 
of degree $-\gamma(M[1]; N'[1]\cup N'')$, and since 
$\gamma(M[1]; N'[1]\cup N'')=\gamma(M; N)$ (use $N'$ non-empty to show this), 
the element $\underline{u}(M[1]; N'[1]_\vartriangle N'')$ is also of first 
degree 0. 
We take $$\underline{\phi}(M;(1^{N'}, 3^{N''}),   N)=\underline{u}(M[1]; N'[1]_\vartriangle N'')\,.$$  
When $\term{N'}+1\neq \init{N''}$, we set 
$\underline{u}(M[1]; N'[1]_\vartriangle N'')=0$ for convenience.
\smallskip 

(c) Type (2).\quad When $\al=\underline{2}$ and $\term{M}=\init{N}$, 
we take the element 
$f[1](M\cup N)\in F(K^{M[1]\cup N[1]})$, 
take its image under the diagonal 
extension map 
$$F(K^{M[1]\cup N[1]})\to F(K^{M[1]}, K^{N[1]}) $$
put the sign $(-1)^\ell$, and view it as 
an element of $F(K[1]^{M}; K[1]^{N})$;
denote it by $f[1](M_\vartriangle N)$. 
Let
$\underline{\phi}(M; \underline{2}, N)=f[1](M_\vartriangle N)$. 
\smallskip 

(d) Type (23).\quad  When $\al=(2^{N'}, 3^{N''})$, where
$\term{M}=\init{N}$, and $N', N''$ is a partition of $N$ such that 
$N'\neq \emptyset$ and $N''\neq \emptyset$, 
we take the element 
$\underline{u}(M[1]\cup N'[1]; N'')\in F(K^{M[1]\cup N'[1]}; L^{N''})$, 
take its image by the diagonal extension at $K^{\ell+1}$, 
$$F(K^{M[1]\cup N'[1]}; L^{N''})\to F(K^{M[1]}, K^{N'[1]}; L^{N''})\,,$$
put the sign $(-1)^\ell$ ({\it not\/} $(-1)^{\ell+1}$), and by 
the identification (notice the place of semi-colon changes)
$$F(K^{M[1]}, K^{N'[1]}; L^{N''})=F(K[1]^{M}; K[1]^{N'},L^{N''})
\,\,(=F(K[1]^M; \al, D^N)\,)
$$
view it as an element of the latter group.
Denote it by $\underline{u}(M[1]_\vartriangle N'[1]; N'')$.
Since $\underline{u}(M[1]\cup N'[1]; N'')$ is of first degree 0, i.e., of 
degree $-\gamma(M[1]\cup N'[1]; N'')$, 
and $\gamma(M[1]\cup N'[1]; N'')=\gamma(M; N)$ by the assumption $N', N''$ both non-empty, 
the element $\underline{u}(M[1]_\vartriangle N'[1]; N'')$ is also of first 
degree 0. 
Take 
$$\underline{\phi}(M;(2^{N'}, 3^{N''}),   N)=\underline{u}(M[1]_\vartriangle N'[1]; N'')\,.$$ 
When $\term{M}\neq\init{N}$ we set $\underline{u}(M[1]_\vartriangle N'[1]; N'')=0$.
\smallskip 

(e) For other functions $\al$, we take $\underline{\phi}(M; \al, N)=0$. 
\smallskip 

\noindent By a direct case-by-case verification, one shows 
these primary elements are 
$\sigma$-consistent and $\delta$-closed. 
The proof of the latter proceeds as follows. 

The $\delta$-closedness means that for each $(M; \al, N)$ we have:
$$\partial \underline{\phi}(M; \al, N)
+\sum_k \bvphi_k ( \underline{\phi}(M\cup\{k\}; \al, N)\,)
+\sum_{(k, \tilde{\al})} \bvphi_k ( \underline{\phi}(M; \tilde{\al}, N\cup\{k\})\,)
=0$$
where  $k\in (-\infty, \term{M}] -M$ in the first sum, and where
$k\in [\init(N), +\infty)-N$
  and $\tilde\al$ extension of $\al$ in the second 
 sum. 
To show the identity holds, we argue according to the  type of $\al$. 

(a) If $\al=\underline{1}$ or 
$\al=(1^{N'}, 3^{N''})$ with $N'\neq \emptyset$, 
$N''\neq \emptyset$ and $\term{N'}+1=\init{N''}$, 
it holds by the $\delta$-closedness of $\underline{u}(M; N)$. 

(b) If $\al=(1^{N'}, 3^{N''})$ with $N'\neq \emptyset$, 
$N''\neq \emptyset$ and $\term{N'}+1\neq \init{N''}$, 
the first term of  the identity is zero, and 
it reduces to the easily verifiable identity
(with $n'=\term{N'}$, $n''=\init{N''}$)
$$\bvphi_{n''-1}\left(
\underline{u}(M[1]; (\,(N'\amalg \{n''-1\})[1]_\vartriangle N'')
\right)
+
\bvphi_{n'+1}\left(
\underline{u}(M[1]; (\,(N'[1])_\vartriangle (\{n'+1\}\cup N'')\,)
\right)=0\,.$$

(c) If $\al=\underline{2}$, it reduces to the $\delta$-closedness of $f(M_\vartriangle N)$, shown in (\IdentityProp).

(d) If $\al=(2^{N'}, 3^{N''})$ with $N'\neq \emptyset$, 
$N''\neq \emptyset$, the identity follows from the $\delta$-closedness 
of $\underline{u}(M; N)$.

(e) If $\al=\underline{3}$, the first term of the identity is zero, and 
letting $m_\mu=\term{M}$ and $n_1=\init{N}$, the identity reduces to:
$$\bvphi_{m_\mu} 
\left(
\underline{u}(M[1]_\vartriangle \{m_\mu\}[1]; N)
\right)
+
\bvphi_{n_1-1} 
\left(
\underline{u}(M[1];\{n_1-1\}_\vartriangle N)
\right)=0\,.$$
For other functions $\al$, the three terms are all zero. 

By the $\sigma$-consistent prolongation (\AdditivityBGKL.2) we obtain 
a set of elements
$(\underline{\phi}(M; N)\,)\in \BF(K[1], D)$ of degree 0
which is $\delta$-closed.
\bigskip 

One can prove 

(i) $\al(\al (u))\cdot \psi= \be(u)$. 

(ii) $\phi\cdot \be(\al(u))=-u[1]$. 

(iii) $\phi\cdot \psi=id$.

(iv) $\psi\cdot \phi=id$. 

\noindent We write out the proof of (i) and (ii) only.  
First note that according to the definition $\al(\al(u))$ is represented 
by $\vecta(M; N)\in F(C_u^M; C_{\al(u)}^N)$, which is the $\sigma$-consistent 
prolongation of the set of elements 
$\vecta(\al, M; \be, N)$ $\in F(\al, C_u^M; \be, C_{\al(u)}^N)$
for functions $\al: M\to \{1, 2\}$ and $\be: N\to \{1,2,3\}$ given as follows. 
\smallskip

(a) If $\al=\underline{1}$ and $\be=\underline{2}$, then 
$\vecta (\underline{1}, M; \underline{2}, N)= (f[1])(M_\vartriangle N)
\in F(K[1]^M; K[1]^N)$\,
 (the same element as the type (2) element for $\underline{\phi}$). 
 
(b) If $\al=\underline{2}$ and $\be=\underline{3}$, then 
$\vecta (\underline{2}, M; \underline{3}, N)= g(M_\vartriangle N)
\in F(L^M; L^N)$. 

(c) In case  $\term{M}=\init{N}=\ell$ and $\al=(1^{M'}, 2^{M''})$ with $M', M''$ segmentation of $M$ with 
$M'$ non-empty (but $M''=\emptyset$ allowed)
and $\be=\underline{3}$, then 
$$\vecta((1^{M'}, 2^{M''}), M; \underline{3}, N) \in F(K[1]^{M'}, L^{M''}; L^N)$$
is 
$\underline{u}(M'[1]; M''_\vartriangle N)$, the element  obtained from 
$\underline{u}(M'[1]; L^{M''\cup N})$ by diagonal extension  at $K^{\ell+1}$
and putting 
sign $(-1)^\ell$.

In case $\term{M}=\init{N}=\ell$ and $\al=\{1\}$ and 
$\be=(2^{N'}, 3^{N''})$ with $N', N''$ segmentation of $N$ with $N''$ non-empty
(but $N'=\emptyset$ is allowed)
$$\vecta(\underline{1}, M; (2^{N'}, 3^{N''}), N\,)=
\underline{u}((M_\vartriangle N')[1]; N'')\in F(K[1]^M; K[1]^{N'}, L^{N''})\,,$$
the element 
obtained from $\underline{u}(M[1]\cup N'[1]; N'')\in F(K[1]^{M\cup N}; L^{N''})$
by diagonal extension at $K^{\ell+1}$ and putting $(-1)^\ell$. 

(d) For other pairs of functions $\al, \be$, one has $\vecta(\al, M;\be, N)=0$. 
\smallskip 

This  set of elements is $\sigma$-consistent and $\delta$-closed, and 
we have 
$$\vecta(M; N)=\theta\left(
\sum h(\gamma_1, M_1)\ts\cdots \ts h(\gamma_t, M_t)\ts 
\vecta(\gamma', M'; \al', N')
\ts  k(\al_1, N_1)\ts\cdots k(\al_r, N_r)
\right)$$
where the sum is over $(M_1, \cdots, M_t, M')$, $t\ge 0$,  
a segmentation of $M$, and $(\gamma_1, \cdots, \gamma_t, \gamma')$
is a set of functions taking distinct values at the overlaps of the segmentation, 
as well as over $(N', N_1, \cdots, N_r)$, $r\ge 0$ and $(\al', \al_1, \cdots, \al_r)$.
Similarly, 
$$\underline{\psi}(\bar{N}; R)=
\theta\left(\sum k(\bar{\al}_1, \bar{N}_1)\ts\cdots\ts k(\bar{\al}_s, \bar{N}_s)
\ts\underline{\psi}(\bar{\al}', \bar{N}'; R)\,
\right)$$
the sum over segmentations $(\bar{N}_1, \cdots, \bar{N}_s, \bar{N}')$ of 
$\bar{N}$ and functions $(\bar{\al}_1, \cdots, \bar{\al}_s, \bar{N}')$ 
taking distinct values at the overlaps. 
Thus for $\vecta=(\vecta(M; N))$ and 
$\underline{\psi}=(\underline{\psi}(\bar{N}; R))$,
their tensor product 
$\vecta\ts \underline{\psi}\in \BG(C_u, D)\tildets \BG(D, K[1])$
is equal to 
$$\sum \vecta(M; N)\ts \underline{\psi}(\bar{N}; R)\in 
\moplus F(C_u^M; D^N)\tildets F(D^{\bar{N}}; K[1]^R)\,,$$
the sum over all free double sequences $(M; N)$ and $(\bar{N}; R)$. 
Consequently 
$$\rho(\vecta\ts \underline{\psi} )=\sum \vecta(M; N)\ts \underline{\psi}(\bar{N}; R)
\in  \BG(C_u, D, K[1], \{2\})
\,,$$
the sum over those pairs of double sequences with the condition
$\term{N}=\init{\bar{N}}$. 

We shall give an element $W\in \BG(C_u, C_{\al(u)} , K[1], \, \emptyset)$
such that 
$\tilde{\bsigma}(W)= \rho(\vecta\ts \underline{\psi} )$. 
For a free triple sequence $(M; N; R)$ and functions 
$\gamma$ on $M$, $\al$ on $N$,  we define elements  
$$W(\gamma, M; \al, N; R)\in F(\gamma, C_u^M; \al,D^N; K[1]^R)$$
as follows. 
When $\gamma=\underline{2}$ and $\al=\underline{2}$, let  
$$W(\underline{2}, M;\underline{2}, N; R)=
f[1](M_\vartriangle N_\vartriangle R)\,;$$
the right hand side is defined below. 
In all other cases, we take $W(\gamma, M; \al, N; R)=0$. 

To define $f[1](M_\vartriangle N_\vartriangle R)$, 
assume  $\term{M}=\init{N}=\ell$ and
$\term{N}=\init{R}=\ell'$. Consider
the element 
$f[1](M\cup N\cup R)\in F(K^{M[1]\cup N[1]\cup R[1]})$,
take its image under the diagonal extension 
at $K^{\ell+1}$ and $K^{\ell'+1}$, 
$$\diag: F(K^{M[1]\cup N[1]\cup R[1]})\to F(K^{M[1]}; K^{N[1]}; K^{R[1]})\,,$$
and put the sign $(-1)^{\ell+\ell'}$; the result is the 
$f[1](M_\vartriangle N_\vartriangle R)$. 
Note that if $\ell=\ell'$, then $K^{\ell+1}$ is counted three times in the diagonal extension, and the sign is $(-1)^{\ell+\ell}=+1$. 
If $\term{M}\neq \init{N}$ or $\term{N}\neq \init{R}$, we let
$f[1](M_\vartriangle N_\vartriangle R)$  be zero.

One then obviously has the identities 
$$
\sigma_k(W(\gamma, M; \al, N; R)\,)
=\left\{
\begin{array}{ll} 
h(\gamma_{\le k}, M_{\le k}) \ts  W(\gamma_{\ge k}, M_{\ge k}; \al, N; R) &\mbox{if $k\in M-\{\init{M}\}$\,,} \\
\vecta(\gamma, M; \al_{\le k}, N_{\le k})\ts \underline{\psi}(\al_{\ge k}, N_{\ge k}; R)      &\mbox{if $k\in N$\,,}  \\
W(\gamma, M; \al, N; R_{\le k})\ts f[1](R_{\ge k}) 
&\mbox{if $k\in R-\{\term{R}\}$\,.}  
\end{array}\right.
$$
Then we define 
the element $W(M; N; R)\in F(C_u^M; D^N; K[1]^R)$ by 
the formula
$$\begin{array}{cl}
W(M; N; R) 
&= \theta\left(
\sum h(\gamma_1, M_1)\ts {\scriptstyle\cdots} \ts 
h(\gamma_t, M_t)\ts \vecta(\gamma', M'; \al', N') \right.\\
&\left.\phantom{aaaaaaaaaaaaaa}
\ts k(\al_1, N_1)\ts {\scriptstyle\cdots} \ts k(\al_r, N_r)
\ts \underline{\psi}(\bar{\al}, \bar{N}; R)\,
\right)   \\
&+ \theta\left(\sum h(\gamma_1, M_1)\ts \cdots h(\gamma_t, M_t)\ts W(\gamma', M'; \al, N; R)
\right)\,. 
\end{array}\eqno{(\Axiomproof.a)}$$
The sum is over segmentations of $M$ and $N$ and functions taking distinct 
values at the overlaps.  In the first sum, $N$ is segmented into at least 
two sub-intervals $N', N_1, \cdots, N_r, \bar{N}$, while in the second sum 
$N$ is segmented into just itself. 
These elements satisfy
$$
\sigma_k(W( M; N; R)\,)
=\left\{
\begin{array}{cl} 
h(M_{\le k}) \ts  W(M_{\ge k}; N; R) &\mbox{if $k\in M-\{\init{M}\,\}$\,,} \\
\vecta( M; N_{\le k})\ts \underline{\psi}(N_{\ge k}; R)    &\mbox{if $k\in N$\,,}  \\
W(M; N; R_{\le k})\ts f[1](R_{\ge k}) 
&\mbox{if $k\in R-\{\term{R}\,\}$\,.}  
\end{array}\right.
$$
The first and third identities show that 
$W=(W( M; N; R)\,)\in H^{010}(C_u; D; K[1])$ is contained in 
$\BG(C_u; D; K[1]\,, \emptyset)$, and the second identity shows
$\tilde{\bsigma}(W)= \rho(\vecta\ts \underline{\psi} )$. 
By (\Compos), $\al(\al(u))\cdot \psi$ is represented by 
$\Pi_D(W)$.
Recall by definition 
$$(\Pi_D(W))(M; R)=\sum_j (-1)^j \vphi_{D^j} (W(M; \{j\}; R)\, )\,.$$
One has by (\Axiomproof.a)
\begin{eqnarray*}
W(M; \{j\}; R)
&=&\theta\left( \sum h(\gamma_1, M_1)\ts {\scriptstyle\cdots} \ts 
h(\gamma_t, M_t)\ts \vecta(\gamma', M'; \al', \{j\})
\ts \underline{\psi}(\bar{\al}, \{j\}; R) \right)  \\
&+&\theta\left(\sum h(\gamma_1, M_1)\ts \cdots \ts h(\gamma_t, M_t)\ts
 W(\gamma', M'; \al, \{j\}; R)
\right)\,.
\end{eqnarray*}
In the first sum, $\al'(j)\neq \bar{\al}(j)$; since the composition of
 $\pi_j$ and $\vphi_j$ is zero, 
the first sum yields zero by $\vphi_{D^j}$. 
For the second sum, using the identity 
$$(-1)^j \vphi_{K^{j+1}} (f[1](M'_\vartriangle\{j\}_\vartriangle R))
=
\left\{
\begin{array}{cl}
(f[1])(M'_\vartriangle R)\, &\mbox{if $\term{M'}=\init{R}=j$}, \\
0 &\mbox{otherwise},
\end{array}
\right.$$
and the compatibility of $\theta$ and $\vphi$, one has 
$$(-1)^j \vphi_{D^j} (W(M; \{j\}; R)\, )
=
\theta\left(\sum h(\gamma_1, M_1)\ts \cdots \ts h(\gamma_t, M_t)\ts
(f[1])(M'_\vartriangle R)\right)
$$
if $\term{M}=\init{R}=j$, and zero otherwise. 
Therefore
$$(\Pi_D(W))(M; R)
=\theta\left(\sum h(\gamma_1, M_1)\ts \cdots \ts h(\gamma_t, M_t)\ts
 (f[1])(M'_\vartriangle R)
\right)\,;
$$
the sum is over segmentations of $M$ and functions $(\gamma_1, \cdots, \gamma_t)$ taking distinct values at the overlaps and  $\gamma_t(\term(M_t))\neq 2$. 
These are the representatives of $\be(u)$, so we have shown (i). 

To prove (ii),  the map $\be(\al)$ is represented by elements 
$\vectb(M; N)\in F(C_u^M; K[1]^N)$ given by 
$$\vectb(M; N)
=\theta\left(
\sum k(\gamma_1, M_1)\ts\cdots\ts k(\gamma_c, M_c)\ts (g[1])(M'_\vartriangle N)\,
\right)\,, $$
the sum over segmentations of $M$ and functions such that $\gamma_c(\term(M_c))\neq 1$. 
On the other hand, 
$$\underline{\phi}(M; N)=
\theta\left(
\sum \underline{\phi}(M; \al', N')\ts 
k(\al_1, N_1)\ts\cdots\ts k(\al_r, N_r)
\right)\,.$$
Hence $\rho(\underline{\phi}\ts \underline{\vectb})\in 
\BG(K[1], D, L[1], \{2\})$ equals 
$$
\begin{array}{l}
\theta\left(
\sum \underline{\phi}(M; \al', N')\ts 
k(\al_1, N_1)\ts{\scriptstyle\cdots}\ts k(\al_r, N_r)  
\right. \\
\phantom{ \underline{\phi}(M; \al', N')\ts 
k(\al_1, N_1)\ts{\scriptstyle\cdots} } 
\left.
\ts k(\gamma_1, \bar{N}_1)\ts{\scriptstyle\cdots}
\ts k(\gamma_c, \bar{N}_c)\ts (g[1])(N''_\vartriangle R)\,
\right)\,, 
\end{array}
$$
the sum over free double sequences $(M; N)$, $(\bar{N}, R)$ 
such that $\term{N}=\init{\bar{N}}$, and
segmentations $(N', N_1, \cdots, N_c, N'')$ of $N$ and functions 
on them. 
Let $W\in \BG(K[1], D, L[1], \emptyset)$ be the element with components 
$$\begin{array}{cl}
W(M; N; R)&=
\theta\left(\sum \underline{\phi}(M; \al', N')\ts k(\al_1, N_1)\ts
\cdots\ts k(\al_c, N_c)\ts (g[1])(N''_\vartriangle R)\,
\right)   \\
&\phantom{aa}+\theta\left(\underline{u}(M[1]; N[1]_\vartriangle R[1])\,\right)
\end{array}$$
in $F(K[1]^M, D^N, L[1]^R)$. 
The first sum is over segmentations $(N', N_1, \cdots, N_c, N'')$ of $N$
and functions taking distinct values at the overlaps and 
$\al_c(\term{N_c})\neq 1$. 
One has 
$\tilde{\bsigma}(W)=\rho(\underline{\phi}\ts \underline{\vectb})$. 

We now show $(\Pi_D (W))(M; R)=\underline{u}(M[1]; R[1])$. 
Indeed when $N$ consists of one element $N=\{j\}$, 
$$\begin{array}{cl}
W(M; \{j\}; R)&=
\theta\left(\underline{\phi}(M; \al', \{j\})\ts 
 (g[1])(\{j\}_\vartriangle R)\,
\right)   \\
&\ph{aa}+\theta\left(\underline{u}(M[1]; \{j\}[1]_\vartriangle R[1])\,\right)\,.
\end{array}$$
Therefore, 
$$(-1)^j \vphi_{D^j} (W(M; \{j\}; R)\,)
=(-1)^j\theta\left(\vphi_{L^{j+1}}(\underline{u}(M[1]; \{j\}[1]_\vartriangle R[1])  )\,  \right)$$
which equals 
$\underline{u}(M[1]; R[1])$ if $j=\term(R)$ and zero otherwise. 
Taking the sum over $j$, the assertion follows. 

Recall from (\Conemorphism) that 
$u[1]$ is represented by$- \underline{u}(M[1]; N[1])$, thus $\phi\cdot \be(\al)=-u[1]$. 
\bigskip 

\sss{\Axiomprooftwo} \Prop{\it 
A commutative square 
$$\begin{array}{ccc}
 K   &\mapr{u}  &L  \\
\mapd{v} &     &\mapdr{w} \\
 K'  &\mapr{u'}  &L'
 \end{array}
$$  
in $Ho(\cC^\Delta)$ can be extended to a morphism of triangles 
$$\begin{array}{cccccccc}
          K &\mapr{u} &L &\mapr{\al(u)} &C_u &\mapr{\be(u)} &K[1] \\
     \mapd{v}  & &\mapd{w} & &\mapd{\zeta}  & &\mapdr{v[1]} \\
 K' &\mapr{u'}&L' &\mapr{\al(u')} &C_{u'} &\mapr{\be(u')} &\phantom{\,.}K'[1] \,.
 \end{array}
 $$
}\smallskip 
{\it Proof.}\quad 
Let $\underline{u}=(\underline{u}(M; N)\,)$ be a representative of $u$, and 
similarly for $u'$, $v$, and $w$. 
By (\Compos), there is a cocycle $W=(W(M; N; R)\,)$ in $\BG(K, L, L', \emptyset)^2$
such that 
$$\ti\bsigma(W)=\rho(\underline{u}\ts \underline{w})\,, $$
and  there is a cocycle  $W'=(W'(M; N; R)\,)$ in  $\BG(K, K', L', \emptyset)^2$
such that 
$$\ti\bsigma(W')=\rho(\underline{v}\ts \underline{u'})\,. $$
Let $\zeta: C_u\to C_{u'}$ be the morphism represented by 
$(\underline{\zeta}(M; N)\,)$ with
$\underline{\zeta}(M; N)\in F(C_u^M; C_{u'}^N)$,  obtained by $\sigma$-consistent
prolongation from the $\sigma$-consistent, $\delta$-closed set of elements 
$\zeta(\al, M; \be, N)$ in $F(\al, C_u^M;\be, C_{u'}^N)$
given as follows. 
\smallskip 

(a) If  $\al=(1^{M'}, 2^{M''})$ with $M', M''$ both non-empty and 
$\be=\underline{2}$, let 
$$\zeta((1^{M'}, 2^{M''}),  M; \underline{2}, N)=
W(M'[1]; M''; N)\,.$$

(b) If $\al=\underline{1}$ and $\be=(1^{N'}, 2^{N''})$
with $N'$, $N''$ both non-empty, let 
$$\zeta(\underline{1},  (1^{N'}, 2^{N''}), N)=
W'( M; N'[1]; N'')\,.$$

(c)  If $\al=\underline{1}$ and $\be=\underline{1}$, then 
$\zeta(\underline{1}, M; \underline{1}, N)=v(M[1]; N[1])$. 
\smallskip 

(d)  If $\al=\underline{2}$ and $\be=\underline{2}$, then 
$\zeta(\underline{2}, M; \underline{2}, N)=w(M; N)$. 
\smallskip 

 It can be can be shown that both compositions $\al(u)\cdot \zeta$ and 
$w\cdot \al(u')$ coincide with the morphism $x: L\to C_{u'}$
represented by $(\underline{x}(M; N)\,)$ with 
$\underline{x}(M; N)\in F(L^M; C_{u'}^N)$, that is  obtained 
by the $\sigma$-consistent prolongation from the $\sigma$-consistent, 
$\delta$-closed set of elements given by 
$$\underline{x} (M; \al, N)=
\left\{
\begin{array}{cl}
\underline{w}(M; N) &\mbox{if $\al=\underline{2}$, } \\
0 &\mbox{otherwise\,.}
\end{array}
\right.$$
The verification is analogous to the ones in the proof of 
the previous proposition. 

Similarly, both $\be(u)\cdot v[1]$ and $\zeta\cdot \be(u')$  
coincide with the morphism $y: C_u\to K'[1]$ represented by 
$\underline{y}(M; N)\in F(C_u^M; K'[1]^N)$ obtained by 
the $\sigma$-consistent prolongation from the $\sigma$-consistent, 
$\delta$-closed set of elements given by 
$$\underline{y} (\al, M;N)=
\left\{
\begin{array}{cl}
\underline{v}(M[1]; N[1]) &\mbox{if $\al=\underline{1}$, } \\
0 &\mbox{otherwise\,.}
\end{array}
\right.$$
This concludes the proof. 
\bigskip

\section{The triangulated category of motives over a variety}
\bigskip

\sss{\DefDS} 
Let $S$ be a quasi-projective variety over a field $k$. 
 We take the quasi DG category 
of symbols $Symb(S)$, recalled in (\SymbqDG), and 
apply the constructions of the previous sections. 
We obtain the quasi DG category $Symb(S)^\Delta$, and then its 
homotopy category $Ho(Symb(S)^\Delta)$. 
\bigskip 

\Def We set $\cD(S)=Ho(Symb(S)^\Delta)$. 
This is a triangulated category by (\ThmCdqDGtriang). 
We call this the {\it triangulated category of mixed motives over $S$}.
\bigskip 

Recall that $({\rm Smooth}/k\,, {\rm Proj}/S)$ denotes the category of 
smooth varieties over $k$ equipped with projective maps to $S$. 
For $X$ in $({\rm Smooth}/k\,, {\rm Proj}/S)$ and $r\in \ZZ$, there
corresponds an object 
$$h(X/S)(r):=(X/S,  r)[-2r]$$
 of $\cD(S)$. 
We write $h(X/S)$ for $h(X/S)(0)$. 

If $f: X\to Y$ is a map over $S$, one has
the class of the graph of $f$, $[\Gamma_f]\in \CH_{\dim X}(Y\times_S X)$. 

The next theorem follows from (\SymbqDG) and (\ThmCdqDG). 
\bigskip 

\sss{\ThmDS} \Thm{\it 
The triangulated category $\cD(S)$ has the following properties.

(1) For two  objects $h(X/S)(r)$ and $h(Y/S)(s)$ where 
$X$, $Y$ are in $({\rm Smooth}/k\,, {\rm Proj}/S)$ and $r, s \in \ZZ$, 
 we have a canonical isomorphism 
$$\Hom_{\cD(S)}(h(X/S)(r)[2r], h(Y/S)(s)[2s-n])=\CH_{\dim Y-s+r}(X\times_S Y, n)$$
the right hand side being the higher Chow group of the fiber product 
$X\times_S Y$. 

(2) The composition of morphisms between three such objects 
(with shifts), 
$$
\begin{array}{c}
\Hom(h(X/S)(r)[2r], h(Y/S)(s)[2s-n])
\ts \Hom(h(Y/S)(r)[2s-n], h(Z/S)(t)[2t-n-m])\\
\to \Hom(h(X/S)(r)[2r], h(Z/S)(t)[2t-n-m])
\end{array}$$
is identified via the isomorphism in (1) with 
the map 
$$\psi: 
\CH_{\dim Y-s+r}(X\times_S Y, n)\ts
\CH_{\dim Z-t+s}(Y\times_S Z, m)
\to \CH_{\dim Z-t+r}(X\times_S Z, n+m)\,$$
in (\SymbqDG). 

(3) There is a functor 
$$h: ({\rm Smooth}/k\,, {\rm Proj}/S)^{opp}\to 
\cD(S)$$
that sends $X$ to $h(X/S)$, and a map  $f: X\to Y$ to 
the map $f^* : h(Y/S)\to h(X/S)$ corresponding to 
$[\Gamma_f]\in \CH_{\dim X}(Y\times_S X)$.
}\bigskip

{\bf References.}
\smallskip




\RefBlone Bloch, S. : Algebraic cycles and higher $K$-theory, 
Adv. in Math. 61 (1986), 267 - 304. 

\RefBltwo  ---  : The moving lemma for higher Chow groups, J. Alg. Geom. 
3 (1994), 537--568. 

\RefBlthree --- :  
Some notes on elementary properties of higher chow groups, including functoriality properties and cubical chow groups, preprint on Bloch's home page.

\RefCH Corti, A. and Hanamura, M. : Motivic decomposition and intersection Chow
groups I, Duke Math. J. 103 (2000), 459-522.


\RefFu Fulton, W. : Intersection Theory, Springer-Verlag, Berlin, New York, 1984.



\RefHaoneone  Hanamura, M.  :  
Mixed motives and algebraic cycles I, Math. Res.
Letters 2(1995), 811-821. 

\RefHaonetwo  Hanamura, M.  :
Mixed motives and algebraic cycles II, 
 Invent. Math. 158 (2004), 105-179. 

\RefHaonethree  Hanamura, M.  :
Mixed motives and algebraic cycles III,
Math. Res. Letters  6(1999), 61-82.

\RefHaDescent   --- :  Homological and cohomological motives of algebraic varieties,
Invent. Math. 142(2000), 319-349. 

\RefHaRelcorresp  --- : Cycle theory of relative correspondences, preprint. 

\RefKa  Kapranov, M. M. : On the derived categories of coherent sheaves on some homogeneous spaces,  Invent. Math.  92  (1988), 479--508.

\RefLev  Levine, M.,  : Mixed Motives, Mathematical Surveys and Monographs, vol. 57,
1998, American Math. Soc.

\RefKaSh  Kashiwara, M., Schapira, P.: Sheaves on Manifolds, Springer, Berlin, New York,
1990. 

\RefMacL Mac Lane, S.:  Categories for the Working Mathematician, Second edition. Graduate Texts in Mathematics, 5. Springer-Verlag, New York, 1998. 

\RefMacLtwo Mac Lane, S.:   Homology. Reprint of the 1975 edition. Classics in Mathematics. Springer-Verlag, Berlin, 1995.

\RefTera Terasoma, T. : DG-categories and simplicial bar complexes, 
Moscow Math. J. vol. 10 (2010), 231--267.

\RefVe Verdier, J.-L. :
Des cat\'egories d\'eriv\'ees des cat\'egories ab\'eliennes, 
Ast\'erisque No. 239 (1996). 

\RefVo Voevodsky, V. : Triangulated categories of motives over a field,
in: Cycles, Transfers, and Motivic Cohomology Theories (by V. Voevodsky, A. Suslin, and 
E. M. Friedlander), Princeton Univ. Press (2000). 

\RefWe  Weibel, C. A. : An introduction to homological algebra. Cambridge Studies in Advanced Mathematics, 38. Cambridge University Press, Cambridge (1994).
\bigskip 
{}\bigskip 

{\it Department of Mathematics,}

{\it Tohoku University}

{\it  Aramaki Aoba-ku, 980-8587, Sendai, Japan} 

\end{document}